\documentclass{siamart171218}

\usepackage{amsmath}
\usepackage{float}
\usepackage{hyperref}
 \usepackage{array,multirow}
\usepackage{sidecap}
\usepackage{bbold}
\usepackage{graphicx}
\usepackage{wrapfig}

\newtheorem{rmk}{Remark}

\newcommand{\ep}{\epsilon} 
\newcommand{\Aep}{A^{\epsilon} }
\newcommand{\bep}{b^{\epsilon} }
\newcommand{\cep}{c^{\epsilon}}
\newcommand{\Fep}{F^{\epsilon} }
\newcommand{\dt}{\Delta t }
\newcommand{\At}{\tilde{A} }
\newcommand{\bt}{\tilde{b} }
\newcommand{\ct}{\tilde{c} }

\begin{document}

\title{Stability Analysis and Performance Evaluation of Mixed-Precision Runge-Kutta Methods\thanks{Submitted to the editors on
December 15, 2022.
\funding{Burnett and Gottlieb's work was partially supported by 
ONR UMass Dartmouth Marine and UnderSea Technology (MUST) grant N00014-20-1-2849 
under the project S31320000049160, by DOE grant DE-SC0023164
sub-award RC114586 - UMD, and by AFOSR grants FA9550-18-1-0383 and FA9550-23-1-0037.
Grant's work was partially supported by Michigan State University, by
AFOSR grants FA9550-19-1-0281 and FA9550-18-1-0383 and  by 
DOE grant DE-SC0023164.
}
}}

\author{Ben Burnett\thanks{Center for Scientific Computing and Data Science Research,
UMass Dartmouth. \email{bburnett@umassd.edu}}
\and
Sigal Gottlieb*\thanks{Corresponding author. Center for Scientific Computing and Data Science Research,
UMass Dartmouth. \email{sgottlieb@umassd.edu}}
\and Zachary J. Grant \thanks{CSME, Michigan State University.
  \email{zack.j.grant@gmail.com}}
}

\maketitle

\begin{abstract}
{Additive Runge-Kutta methods designed for preserving highly accurate solutions in 
mixed-precision computation were proposed and analyzed in \cite{Grant2022}. 
These specially designed methods use reduced precision for the implicit  computations
and full precision for the explicit computations.
In this work we analyze the stability properties of these methods and their sensitivity to the low precision
rounding  errors, and demonstrate their performance in terms of accuracy and efficiency.
We develop codes in FORTRAN and Julia to solve nonlinear systems of ODEs and PDEs using 
the mixed precision additive Runge-Kutta (MP-ARK) methods.
The convergence, accuracy, runtime, and energy consumption of these methods is explored.
We show that for a given level of accuracy, suitably chosen
MP-ARK methods may provide significant reductions in runtime.}
\end{abstract}

\section{Introduction}
 
High precision numerical simulations may be desirable when highly accurate solutions are required, or
for reducing rounding errors (e.g. in long-time  simulations). However, the use of high precision is often 
inefficient. Mixed precision simulations aim to combine the {\em efficiency} of low precision with the {\em accuracy} 
of high precision computations. However, care must be taken  to prevent the low precision computations 
from degrading the overall accuracy of the method, or the high precision computations from adversely 
impacting the efficiency.

The use of mixed precision to enhance performance of numerical algorithms
has become increasingly common in numerical linear algebra operations \cite{b1}.
In \cite{perf1}, the authors explored several mixed precision algorithms for solving 
sparse linear systems and identified algorithms that benefit from a mixed precision 
implementation. However, while the performance of mixed precision algorithms for linear algebra applications
has been well-studied \cite{perf1} \cite{perf2} \cite{perf3} \cite{perf4} \cite{perf5} \cite{perf6},
mixed precision algorithms for solving ordinary differential equations (ODES) have only recently been proposed
in \cite{Grant2022} and \cite{Croci2022}. 
In \cite{Grant2022} Z. Grant proposed a numerical analysis framework for development of 
order conditions for mixed-precision  Runge--Kutta methods by using an additive Runge--Kutta method 
approach and framing the  use of multiple  precisions as a perturbation. 
Using this rigorous approach, he developed  novel  mixed-precision  Runge--Kutta additive (MP-ARK)
methods that reduce the cost of the   computationally expensive implicit stages in the 
 Runge--Kutta methods  by employing a low precision computation. 
 Meanwhile, the structure of the MP-ARK is designed to suppress the low precision errors either 
 by introducing inexpensive explicit high-precision correction terms,
or by designing novel methods that internally suppress the low precision perturbations.

The work in \cite{Grant2022} focused on the development of a framework for order analysis of MP-ARK methods, 
and the development and convergence verification of such MP-ARK methods, using a simplified "chopping"
routine to simulate low precision in MATLAB.
In this work, we revisit the novel methods in \cite{Grant2022} to rigorously study their stability properties, 
numerically validate their convergence, and evaluate their efficiency on several clusters, 
CPU  architectures, and programming languages. We explore half,  single, double, and quadruple precision,
and the related  mixed precision codes. We review the development of  perturbed additive Runge--Kutta
for mixed precision implementations in Section \ref{sec:review} and discuss the order conditions that must be satisfied
so that the convergence of the methods has perturbation error of orders $ O(\epsilon^m \dt^n)$ for integers 
$m= 1, 2$ and $n=1, 2, 3$.  In Section \ref{sec:convergence} we  confirm that these methods
converge as expected from the design order and investigate the impact of the stiffness of the problem
on the mixed precision convergence rate. 
In Section \ref{sec:stability} we analyze  the stability properties of the  MP-ARK
and in Section \ref{sec:sensitivity}  we explore the sensitivity of the mixed precision methods to roundoff.
Finally, in Section \ref{sec:performance} we numerically investigate the performance of the 
mixed precision methods and we demonstrate that for a given level of accuracy, suitably chosen
MP-ARK methods may provide significant reductions in runtime (and by extension energy consumption).

\section{Review of  perturbed additive Runge--Kutta methods for mixed precision implementations}
\label{sec:review}
\subsection{Motivating Example}
A simple mixed precision approach is illustrated by the implicit midpoint rule  written as a
two stage Runge--Kutta method:
  \begin{subequations} \label{IMPButcherMP}
  \begin{align}
{y_\ep}^{(1)}&= u^n + \frac{\Delta t}{2} \Fep \left(y_\ep^{(1)} \right)  \label{IMPButcherMP1}\\
{u}^{n+1}&= u^n + \Delta t F\left( y_\ep^{(1)}\right) \label{IMPButcherMP2}
\end{align}
\end{subequations}
  where we use a lower-precision $\Fep$ in the first stage and an accurate high precision 
  evaluation of $F$ in the second stage.
  The motivation here is that the first stage requires a potentially costly implicit solve, and
  the use of the lower-precision $\Fep$ (possibly combined with a lower tolerance in the implicit solve)
  can significantly mitigate this cost. To recover the accuracy lost by this lower precision implicit stage, 
 we compute   the cheaper explicit second stage in higher precision.
 
 The lower precision operator $\Fep$ can be seen as a perturbation of the high fidelity operator $F$, where
 the perturbation is given by 
  \begin{eqnarray} \label{tau}
F(u) - \Fep(u) = \ep  \tau(u)  ,
   \end{eqnarray}
 here $\tau$ is a (possibly a non-smooth) function of order one.
The  $O(\ep \dt)$ error that results from the use of $\Fep$ in the first stage is further dampened by $\dt$ in the second stage.
The error between the full precision and mixed precision implementation over one time-step  is thus of order   $O(\ep \dt^2 )$,
which over the full time evolution will result in a  global error contribution of $O(\ep \dt )$.

We observe that it is possible to add in correction stages that will damp out the error further:
\begin{subequations} \label{IMPButcher_MPcorr} 
\begin{align} 
y^{(1)}_{[0]} &= u^n+ \frac{1}{2} \dt  \Fep(y^{(1)}_{[0]}) \label{IMPButcher_MPcorr1} \\ 
y^{(1)}_{[k]} & = u^n+ \frac{1}{2}  \dt  F(y^{(1)}_{[k-1]}) ,  \; \; \;  k=1:m-1 \label{IMPButcher_MPcorr2}\\
u^{n+1}&= u^n + \Delta t F \left(y^{(1)}_{[m-1]})  \right) . \label{IMPButcher_MPcorr3}
\end{align}
\end{subequations}
Each correction damps the initial stage perturbation error by a factor of $\dt$, so that after $m$ high precision evaluations 
(the $m-1$ corrections and the final reconstruction stage) we
obtain an error over one step that is of order $O(\ep \dt^{m+1} )$.
Over the full time evolution will result in a  global error contribution of $O(\ep \dt^{m} )$.

\subsection{Order conditions for mixed precision Runge--Kutta methods} \label{sec:OC}
A more systematic approach to analyzing mixed precision implementations of Runge--Kutta methods
was developed \cite{Grant2022} by considering these as additive Runge Kutta methods of the form:
\begin{subequations} \label{RK-Butcher}
\begin{align}
y^{(i)}= u^n + \dt \sum_{j=1}^{s}A_{ij}F(y^{(j)}) +  \dt \sum_{j=1}^{s}\Aep_{ij}\Fep (y^{(j)})\\
u^{n+1}= u^n + \dt \sum_{j=1}^{s}b_{j}F(y^{(j)}) +  \dt \sum_{j=1}^{s}\bep_{j}\Fep (y^{(j)}).
\end{align}
\end{subequations}
Using the fact that $\Fep(y)$  is an approximation to  $F(y)$ as in \eqref{tau}, we 
rewrite the scheme \eqref{RK-Butcher} in terms of  $F$ and its perturbation $\tau$ (defined in \eqref{tau}):
\begin{subequations} 
\begin{align}
\label{PerturbedMethoda}
y^{(i)}&= u^n + \dt \sum_{j=1}^{s}\tilde{A}_{ij}F(y^{(j)}) +  \ep \dt \sum_{j=1}^{s}\Aep_{ij} \tau(y^{(j)})\\
u^{n+1}&= u^n + \dt \sum_{j=1}^{s}\tilde{b}_{j}F(y^{(j)}) + \ep  \dt \sum_{j=1}^{s} \bep_{j} \tau (y^{(j)}) \label{Perturbed Methodb}
\end{align}
\end{subequations}
where $\tilde{A}_{ij}= A_{ij}+\Aep_{ij}$ and $\tilde{b} = b_{j}+\bep_{j}$.

This form allows us to use the additive B-series representation to track the evolution of $ F$ 
and its interaction with  the perturbation  $\tau$. For example, a second order expansion is:
\begin{align} \label{expansion}
u^{n+1} &= \underbrace{ u^n +  \dt\bt eF(u^n) + \dt^2\bt \ct F_y(u^n)F(u^n)}_{scheme} \\
& +   \underbrace{\ep\dt \left( \bep e\tau(u^n) +    \dt \left( \bep \ct \tau_y(u^n)F(u^n)  +\bt \cep F_y(u^n) \tau(u^n)
 +  \ep  \bep \cep \tau_y(u^n)\tau(u^n) \right) \right) }_{perturbation}
\nonumber  \\
& + O(\dt^3) \nonumber 
\end{align} 
where $\ct = \tilde{A} e$.
As we expect, there are two sources of error: the approximation errors of the scheme and the perturbation errors from the low precision
operator $\tau$. The order conditions arising from the scheme are obtained by setting the above expansion equal to the Taylor expansion:
\[ \bt e = 1, \; \; \; \bt \ct = \frac{1}{2}.\]
The leading order perturbation errors can be zeroed out by setting the coefficients of $\ep\dt $ equal to zero. However, when the perturbation is 
{\bf not smooth}, as is expected to be the case when it results from rounding errors,
 we cannot rely on cancellation errors to obtain these zeroes, we must instead require that each term be zero
\[ \left| \bep \right| e = 0 \; \; \; \Rightarrow \; \; \; \bep_j = 0. \]
This  cancels out many of the $\ep\dt^2 $ terms! To zero out the rest, we require
\[  \left| \bt \right|   \left| \cep \right| =0 \; \; \; \Rightarrow \; \; \; \bt_j  \cep_j =0 \; \; \forall j.\]

The order conditions for higher order smooth and non-smooth perturbation errors 
were derived in \cite{Grant2022}. These order conditions allow us to ensure that 
methods will provide the correct order for both the numerical methods
and the perturbation. In the table below we repeat the  conditions on the perturbation 
coefficients that will ensure the correct order 
for smooth and non-smooth perturbations.

\[ \begin{array}{|lcc|} \hline 
& \multicolumn{2}{c|}{\mbox{Order conditions for:} } \\
\mbox{Expansion } & \mbox{non-smooth} & \mbox{smooth } \\ 
\mbox{Terms}  & \mbox{perturbations} & \mbox{perturbations} \\ \hline
& & \\
\ep \dt \; \; \; \; & \left| \bep \right| e      & \bep e \\ 
\ep \dt^2 &  \left| \bep \right| \left|  \ct \right|   & \bep \ct \\
\ep \dt^2 & \left| \bt \right| \left| \cep \right| & \bt \cep \\ 
\ep^2\dt^2  &  \left| \bep  \right|  \left| \cep \right| & \bep \cep \\
\ep \dt^3 & \left| \bep  \right|   \left| \At  \right|  \left| \ct   \right|     & \bep \At\ct \\ 
\ep \dt^3 &	\left|  \bt  \right|  \left| \Aep  \right|   \left| \ct   \right|    & \bt \Aep \ct   \\ 
\ep \dt^3 & 	\left|  \bt  \right| \left| \At  \right| \left| \cep \right|   &  \bt \At \cep  \\
\ep \dt^3 &  \left| \bep  \right| (\ct \cdot \ct)    &  \bep(\ct \cdot \ct) \\
\ep \dt^3  & \left| \bt \right|( \left|\ct \right| \cdot \left| \cep \right|)   &  \bt(\ct \cdot \cep)  \\
\ep^2 \dt^3  & \left|  \bep \right| \left| \Aep \right| \left| \ct \right| & \bep \Aep \ct  \\
 \ep^2 \dt^3  &   \left| \bep \right|  \left| \At \right| \left| \cep  \right|   & \bep \At \cep \\
  \ep^2 \dt^3  & \left| \bt \right|  \left| \Aep \right| \left| \cep \right|     &  \bt \Aep \cep   \\
    \ep^2 \dt^3  & \left|  \bep \right|( \left| \cep \right| \cdot \left| \ct \right|)    &    \bep(\cep \cdot \ct)  \\
     \ep^3\dt^3 &   \left| \bep \right|  \left| \Aep \right|  \left|\cep  \right|   & \bep \Aep \cep \\
     \ep^3 \dt^3  &   \left| \bep \right|( \cep \cdot \cep)    & \bep(\cep \cdot \cep)  \\ \hline 
\end{array} \]

\subsection{Mixed precision methods used in this work} \label{sec:methods}
In the current work, we consider three methods: the implicit midpoint rule, a third order singly diagonally implicit 
Runge--Kutta method, and a novel four-stage third order additive Runge--Kutta method designed in \cite{Grant2022}.

\smallskip 

\noindent{\bf [1] IMR: Implicit midpoint rule with corrections}
This is the method given above in \eqref{IMPButcher_MPcorr}.

\smallskip 

\noindent{\bf [2] SDIRK: Third order singly diagonally implicit method with corrections}
We take the two stage  third order singly diagonally implicit method \cite{Alexander1977} 
\begin{subequations} \label{SDIRK-MP}
\begin{eqnarray}
     y^{(1)} &=& u^n + \gamma \dt \Fep(y^{(1)})  \\
     y^{(2)} &=& u^n + \left( 1- 2 \gamma \right)  \dt F(y^{(1)})  + \gamma \dt \Fep(y^{(2)})   \\
      u^{n+1}&=&u^n+ \frac{1}{2} \dt F(y^{(1)})  + \frac{1}{2} \dt F(y^{(2)}) 
\end{eqnarray}
\end{subequations}
with $\gamma =  \frac{\sqrt{3}+3}{6}$
The consistency conditions are satisfied to order three. The coefficient of the 
highest order non-zero perturbation term is
$ \bt \cep = \gamma $
 so we have a perturbation error $E_{per} = O(\dt^2 \ep)$ at each time step.
 If the method is stable to the roundoff error, we would expect a global error
 \[ Error = O(\dt^3) + O(\ep \dt  ) .\]
 
 Adding explicit high precision correction stages  results in the method:
  \begin{subequations} \label{SDIRK-MPfix}
\begin{eqnarray}
y_{[0]}^{(1)} &=& u^n+\gamma \dt \Fep(y^{(1)}_{[0]}) \\
y_{[k]}^{(1)}  &=& u^n+\gamma \dt  F(y^{(1)}_{[k-1]}) \; \; \; \mbox{for} \; \; k = 1,...,  m-1 \\
y_{[0]}^{(2)}  &=& u^n+(1-2\gamma) \dt F(y^{(1)}_{[m-1]})+\gamma \dt F(y^{(2)}_{[0]}) \\
y_{[k]}^{(2)}  &=& u^n+(1-2\gamma) \dt F(y^{(1)}_{[m-1]})+\gamma \dt F(y^{(2)}_{[k-1]}) \; \; \; \mbox{for} \; \; k = 1,...,  m-1 \\
  u^{n+1}&=&u^n+ \frac{1}{2} \dt F(y^{(1)}_{[m-1]}))  + \frac{1}{2} \dt F(y^{(2)}_{[m-1]})) .
  \end{eqnarray}
\end{subequations}
If $m=1$ we have the uncorrected method, which should have global $ O(\dt^3) + O(\ep \dt  )$. Each correction will raise the perturbation 
order by one, so that we obtain a global error of the form $O(\dt^3) + O(\ep \dt^m)$.

 For the uncorrected method, if  $\dt^2 > \ep$, we expect to see the usual third order convergence. However, if  $\dt^2 < \ep$, 
 the perturbation error will dominate. For this reason, we suggest including a correction term if $ \dt < \sqrt{\ep}$.
 The method with one correction ($m=2$) will have  global error of order $O(\dt^3) + O(\ep \dt^2)$, 
 so that if $\dt > \ep$ we should see third order convergence. If $ \dt < \ep$ we will see the perturbation error
 dominating. To fix this, we suggest that a second correction be added ($m=3$), at which point we should observe
 a global error of order $ O(\dt^3)$.

\smallskip

\noindent{\bf [3]  NovelA: Novel Four stage third order ARK method.}
This method, given in \cite{Grant2022} is given by 
\begin{eqnarray}\label{MP-4s3pA}
y^{(1)} & = & u^n + \Delta t A^{\epsilon}_{1,1} F^{\epsilon} \big(y^{(1)}\big)  \nonumber  \\
y^{(2)} & = & u^n + \Delta t A_{2,1} F \big(y^{(1)}\big)  \nonumber \\
y^{(3)} & = & u^n + \Delta t \big[ A_{3,1} F \big(y^{(1)}\big) + A_{3,2} F \big(y^{(2)}\big)  + A^{\epsilon}_{3,1} F^{\epsilon} \big(y^{(1)}\big) + A^{\epsilon}_{3,3} F^{\epsilon} \big(y^{(3)}\big) \big] \nonumber  \\
y^{(4)} & = & u^n + \Delta t \big[ A_{4,1} F \big(y_{(1)}\big) + A_{4,2} F \big(y^{(2)}\big) 
+ A_{4,3} F \big(y^{(3)}\big) \big] \nonumber \\
u^{n+1} & = & u^{n} + \frac{1}{2}\Delta t \big[ \big(y^{(2)}\big) + \big(y^{(4)}\big) \big]
\end{eqnarray}
The matrix $A$ is given by
\[ A_{2,1}=0.211324865405187, \; 
A_{3,1}=0.709495523817170, \; 
A_{3,2}=-0.865314250619423,\]
\[ A_{4,1}= 0.705123240545107, \; 
A_{4,2} =0.943370088535775, \; 
A_{4,3}= -0.859818194486069,\]
and the matrix$\Aep$ is given by
\[\Aep_{1,1}= 0.788675134594813,  \; \;
\Aep_{3,1}= 0.051944240459852, \; \;
\Aep_{3,3} = 0.788675134594813,\]
where all unspecified values   $A_{i,j} = \Aep_{i,j} = 0$.
The final row vectors are given by
\[ b= ( 0, \frac{1}{2}, 0 , \frac{1}{2}) , \; \; \; \; 
\bep = ( 0, 0, 0 , 0) .\]
This method performs as a third order method with a second order perturbation error $O(\dt^3) + O(\ep \dt^2)$. 
This means that when $\dt > \ep$ we should see third order convergence. (It is, of course, possible to add explicit corrections to
this method as well).
It will be interesting to compare this method to the  SDIRK method \eqref{SDIRK-MPfix} with $m=2$, which has the same order, and the 
SDIRK method with $m=3$ which has $O(\dt^3) + O(\ep \dt^3)$.
We note that if we have a smooth perturbation, the order of this NovelA method
is higher: $O(\dt^3) + O(\ep \dt^3)$, because the order conditions for this case may rely on cancellations, as in \cite{Grant2022}.

\section{Order verification for the MP-ARK methods} \label{sec:convergence}
To verify the order of convergence of the numerical methods, we consider the three methods in Section
\ref{sec:methods}, with corrections where relevant.
We observe that the corrections act as expected, to raise the order of the mixed-precision method up to 
that of the full precision. Similarly, we observe that as the theory predicts,
the NovelA method behaves much like the SDIRK method with one explicit correction.
Another aspect of this convergence study is the effect of the stiffness of the problem
on the convergence of the mixed precision method. We observe that when the problem is stiffer,
the errors may be slightly larger,  and require generally smaller time-steps to settle down to the asymptotic regime.

\subsection{Convergence of the mixed precision method for the van der Pol system}
\label{sec:vdpCONV}
Our first test problem is the van der Pol system
\begin{subequations} \label{vdp}
\begin{eqnarray}
y_1' & = & y_2 \\
y_2' & = &  \alpha y_2 (1 - y_1^2) - y_1
\end{eqnarray}
\end{subequations}
 with  initial conditions $y_1(0) = 2.0$, $y_2(0) = 0.0$, for time $t = [0,1]$. A reference solution using an explicit fourth order Runge--Kutta at a small time step was used for calculating the errors. The reference solution was computed entirely in quadruple precision and all the numerical results, regardless of computational precision, were cast to quadruple precision in order to compute the error between a method and this reference solution.

\begin{figure}[htb]
\begin{center} 
{\includegraphics[width=0.24\textwidth]{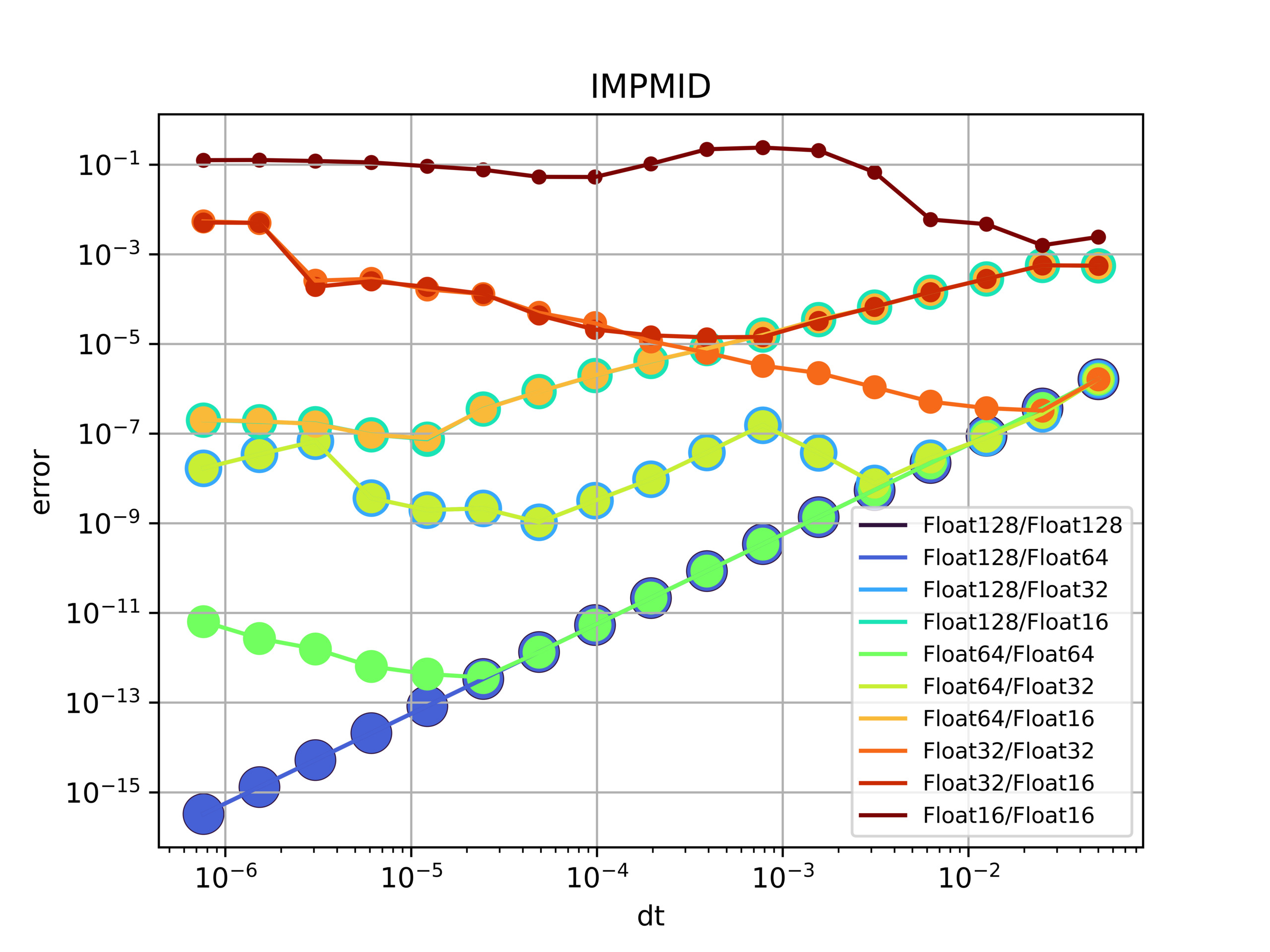}}
{\includegraphics[width=0.24\textwidth]{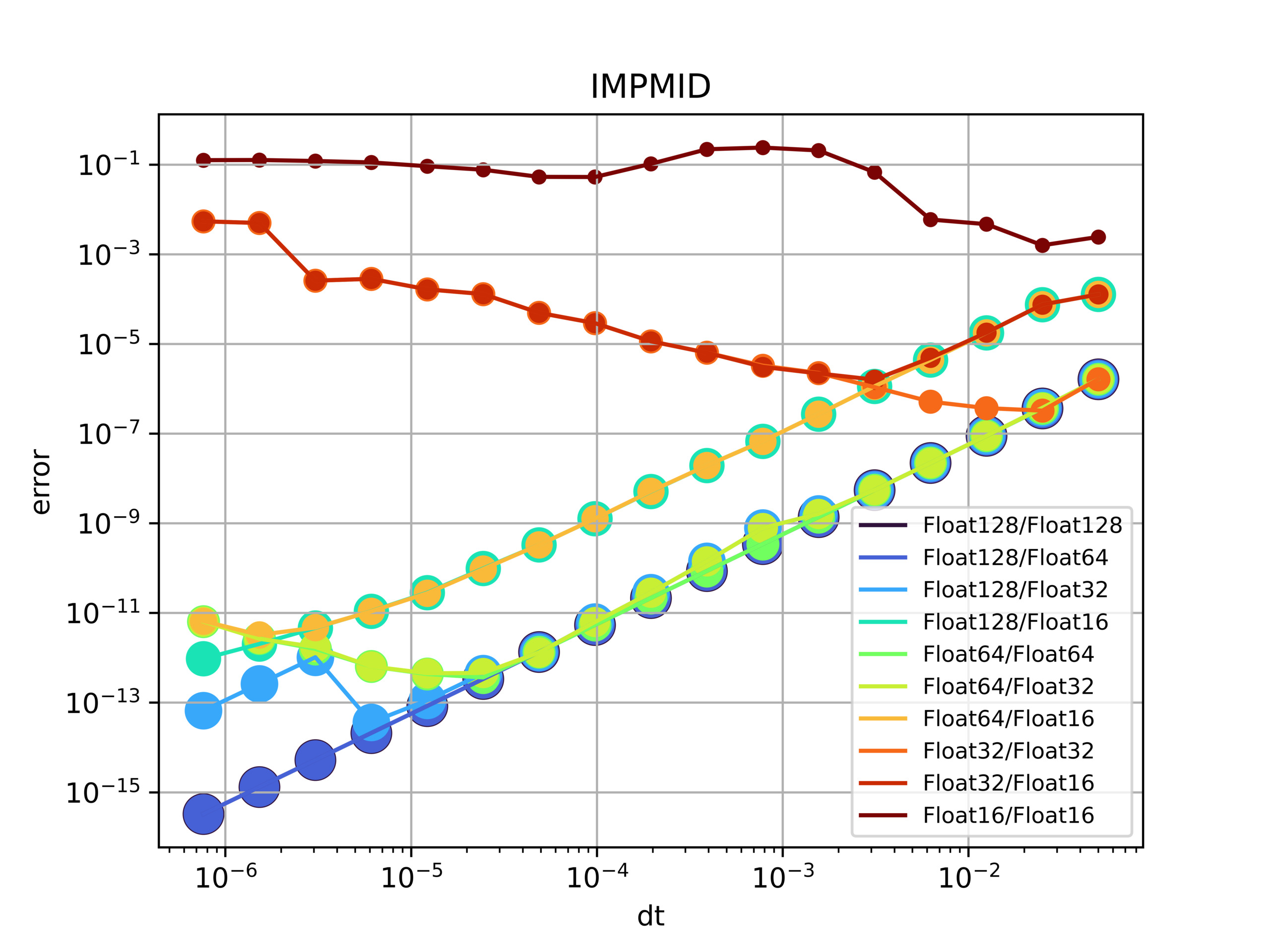}} 
{\includegraphics[width=0.24\textwidth]{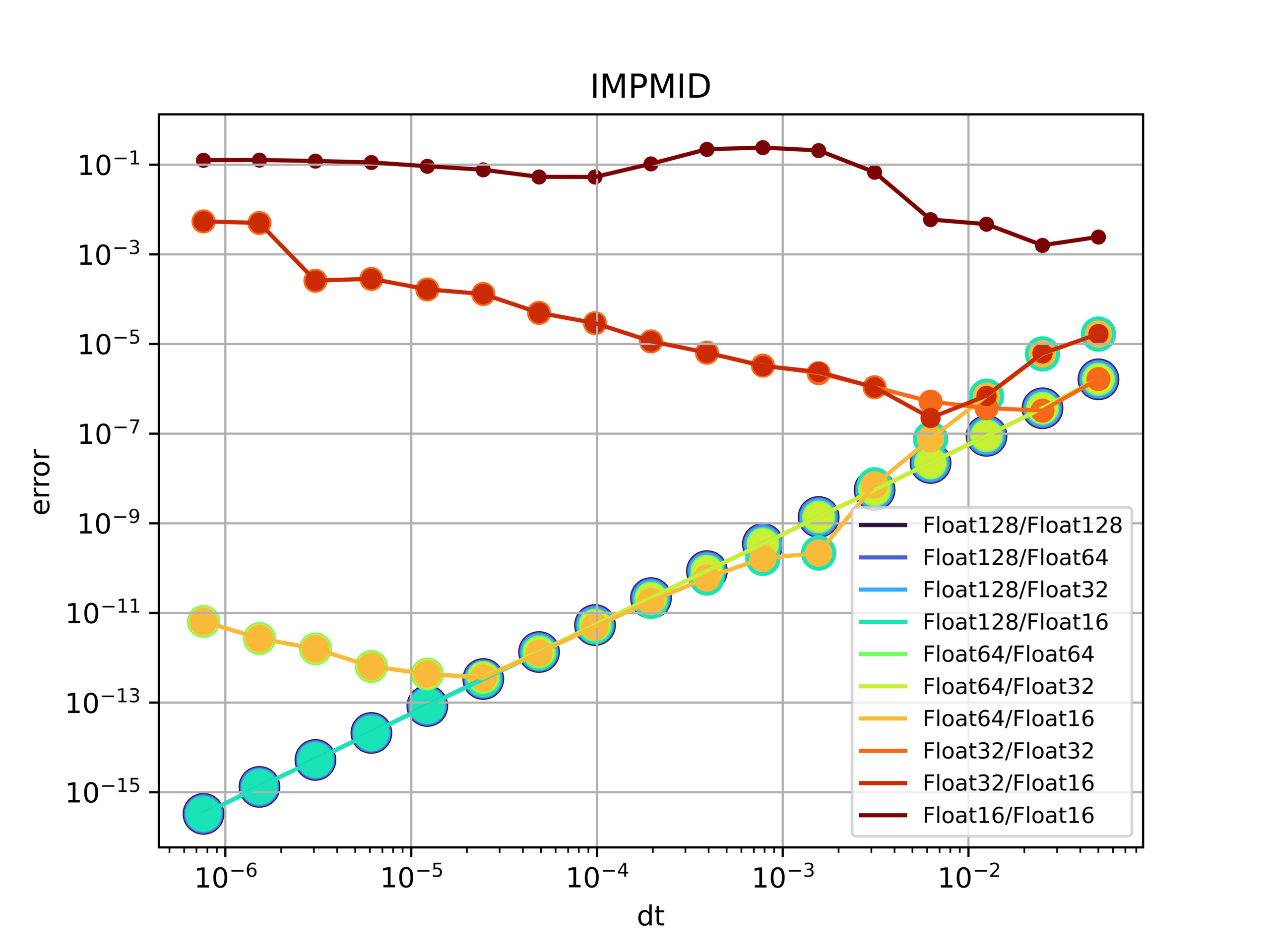}}
{\includegraphics[width=0.24\textwidth]{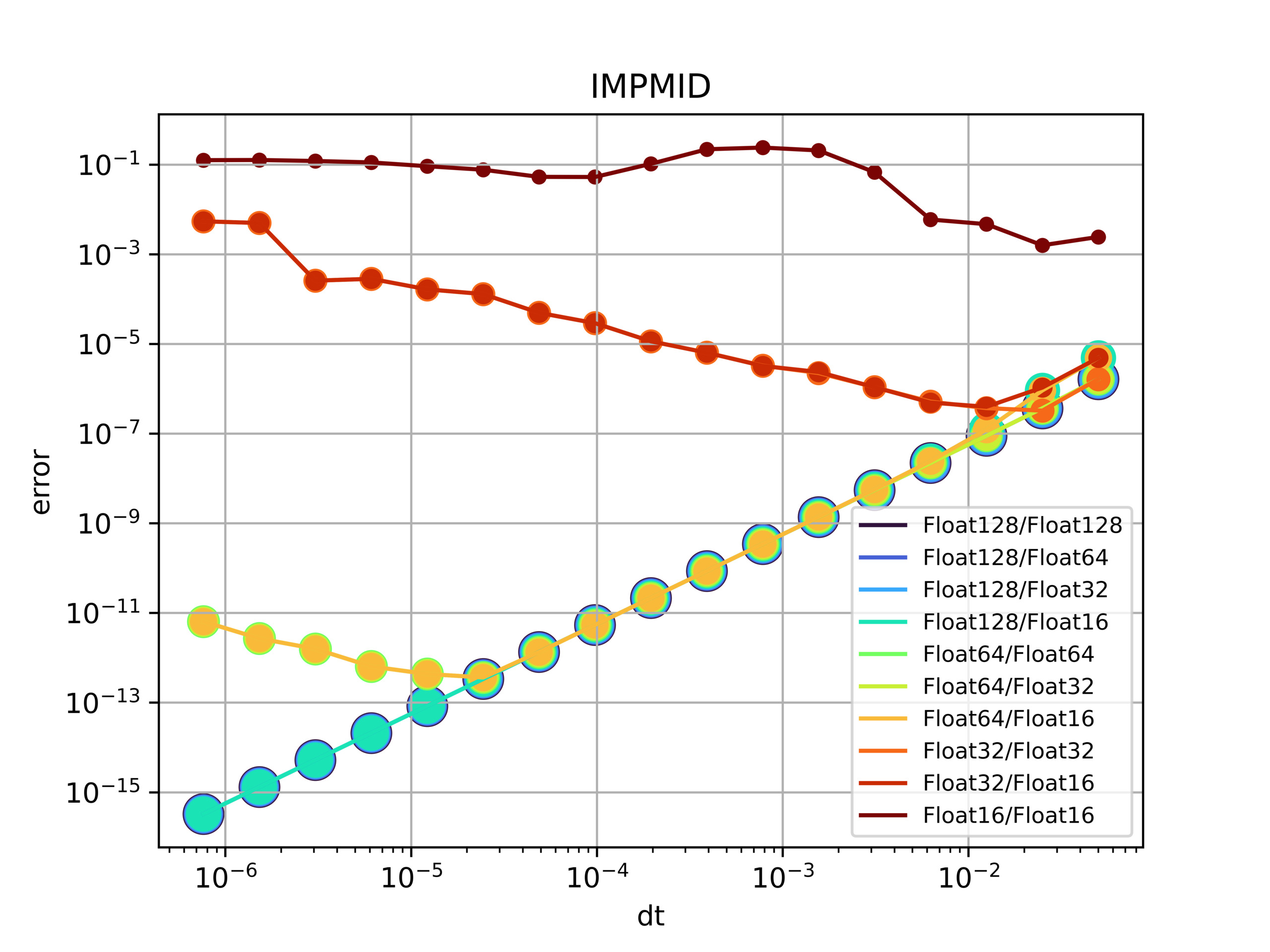}} 
\caption{Converges plots for the implicit midpoint rule for the van der Pol equation with $\alpha =3.0$.
No corrections (left), one correction (middle-left), two corrections (middle-right), three corrections (right).
\label{vdpIMRconv} 
}
\end{center}
\end{figure}

Figure \ref{vdpIMRconv} shows the convergence of the mixed precision and full precision methods 
for the implicit midpoint rule for the van der Pol equation with $\alpha =3.0$.
On the left, we observe that with no corrections, the half precision (Float16)  and single precision (Float32)  codes
do not converge as $ \Delta t$ is refined.
The mixed precision quad/half (Float128/Float16) and double/half codes  (Float64/Float16) 
both converge at first order for $10^{-5} \leq \Delta t \leq 10^{-2}$, after which convergence levels off.
The mixed precision single/half (Float32/Float16) behaves like the quad/half and double/half codes
for $\Delta t$ large enough. However, convergence starts to level off when
$10^{-4} \leq \Delta t \leq 10^{-3}$, and the errors start to behave like the single precision (Float32)  code.
This is not surprising, as we would expect the mixed precision single/half  code to behave no better than the 
single precision code.
The mixed precision quad/single (Float128/Float32) and double/single codes  (Float64/Float32) 
both converge at second order, following the same line as the double precision and quad precision 
codes for $\frac{1}{320} \leq \Delta t$,  after which convergence levels off.
The double precision (Float64) code continues to converge at second order while  $ \Delta t \geq  \frac{1}{40,960}$,
and the mixed precision quad/double (Float128/Float64) and quad precision code continue to converge.

The middle-left image in Figure \ref{vdpIMRconv} shows the behavior of the errors when when explicit correction is added
to the mixed-precision codes (but not to the full precision codes).
We see two major improvements over the non-corrected codes: the first is that the 
mixed precision quad/single (Float128/Float32) and double/single codes  (Float64/Float32) 
both converge at second order, following the same line as the double precision, the mixed quad/double, 
and the quad precision  codes until the double precision bifurcation point at $ \Delta t \approx  \frac{1}{40,960}$.
The second major improvement is that the mixed precision quad/half (Float128/Float16) and double/half codes  (Float64/Float16) 
both converge at second order while $ \Delta t \geq \frac{1}{163,840}$,  where the error level is below $10^{-11}$.

The middle-right image in Figure \ref{vdpIMRconv} shows the continued improvement in the behavior of the errors
when {\em two} explicit correction are added to the mixed-precision codes (but not to the full precision codes).
With two corrections, the full precision quad and  double codes, and the mixed precision quad/double, 
quad/single, quad/half, double/single and double/half codes all converge at second order with the same level of errors
until the double precision bifurcation point at $ \Delta t \approx  \frac{1}{40,960}$ (error level $\approx 10^{-12}$).
After this point, all the errors from the mixed precision quad/double,  quad/single, and quad/half codes look 
like the errors from the quad precision codes,
and all the errors from the  mixed precision double/single and double/half codes look like the 
errors from the double precision codes. An additional explicit correction, as shown in the right image in
Figure \ref{vdpIMRconv}  does not improve the convergence rates or overall behavior any more, but the 
errors from the mixed precision quad/half and double/half errors track more exactly those from the quad 
precision codes, with no oscillations.

\begin{figure}[htb] 
\begin{center}
{\includegraphics[width=0.24\textwidth]{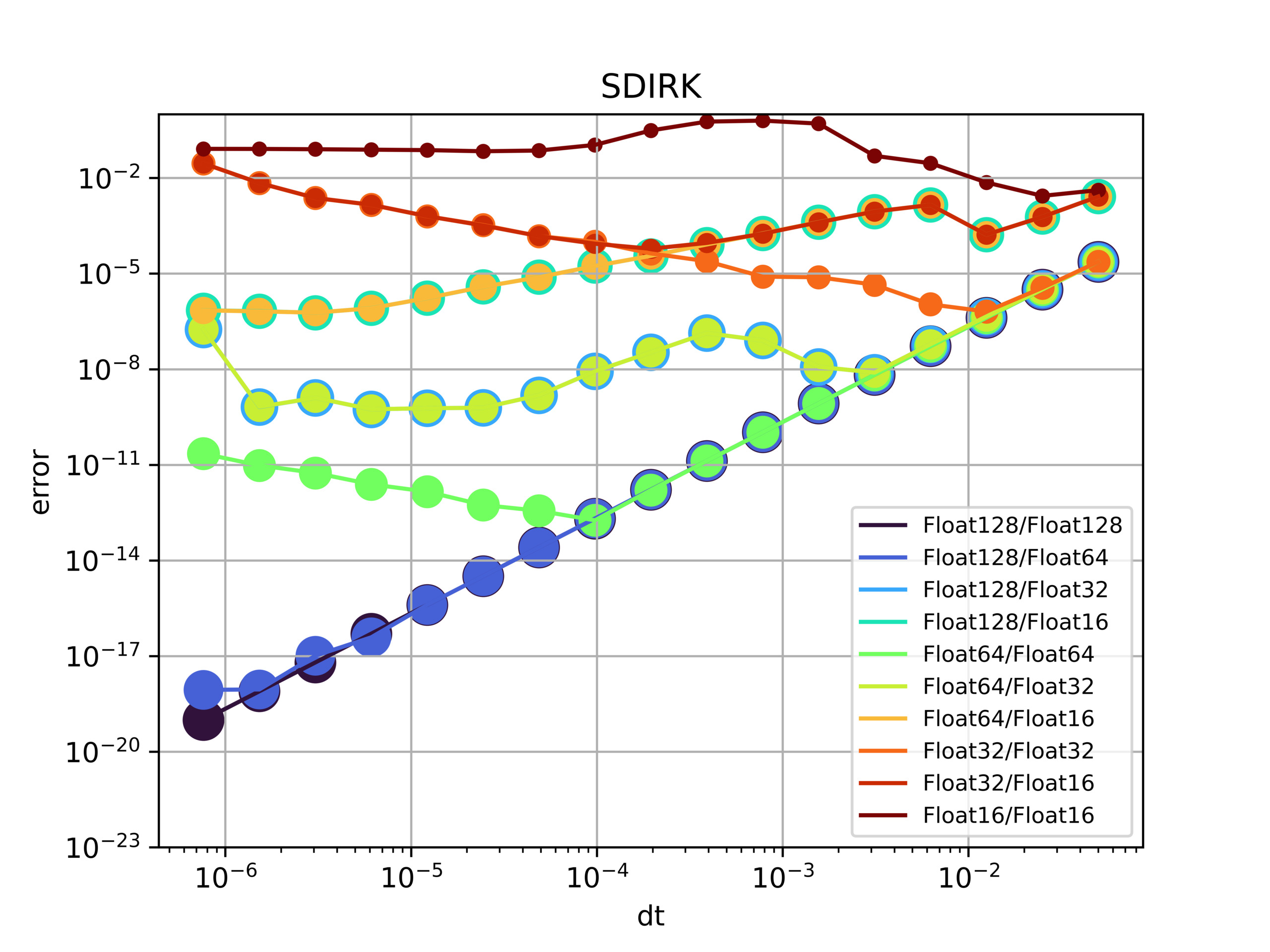}}
{\includegraphics[width=0.24\textwidth]{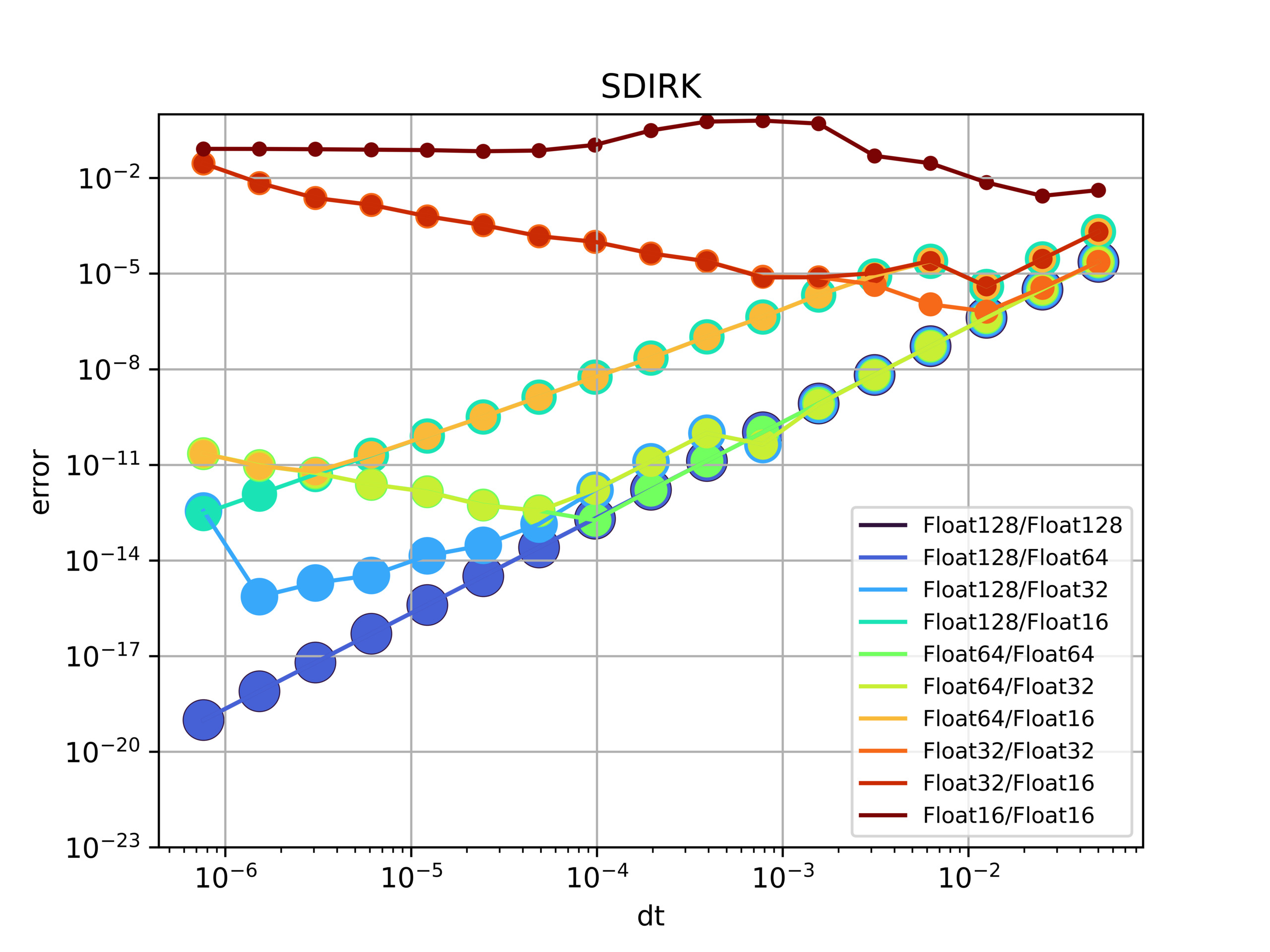}} 
{\includegraphics[width=0.24\textwidth]{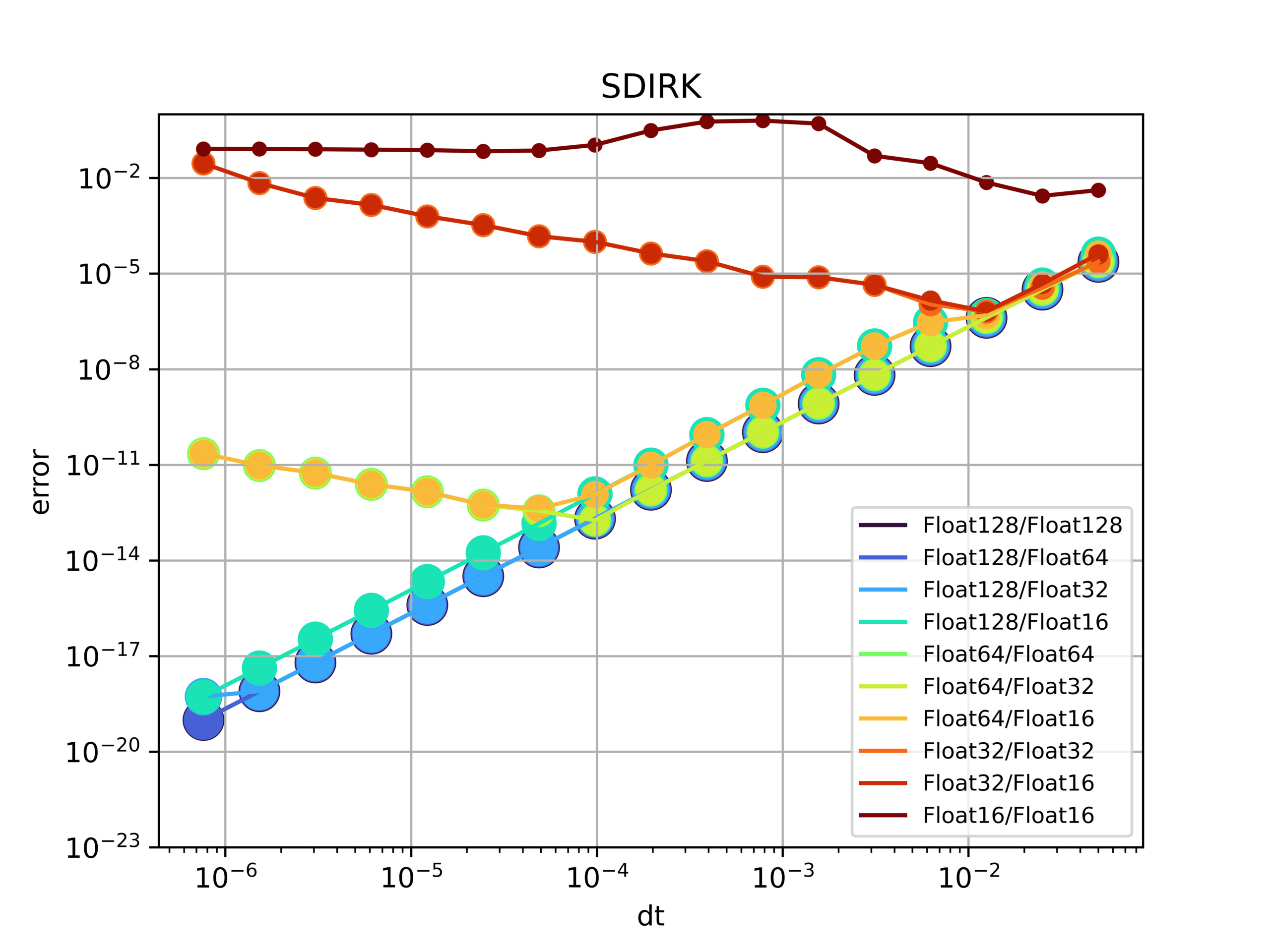}}
{\includegraphics[width=0.24\textwidth]{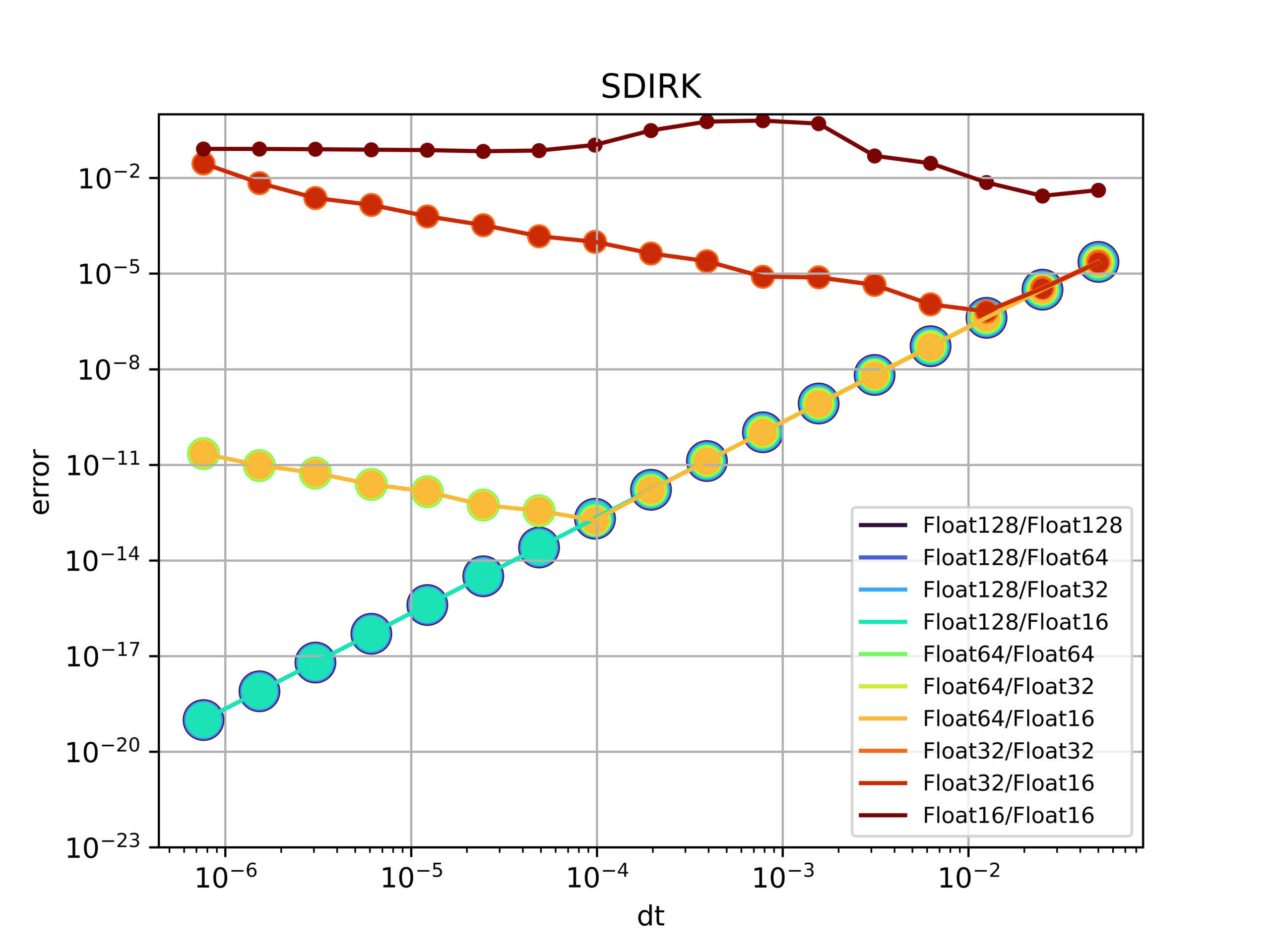}} \\
{\includegraphics[width=0.24\textwidth]{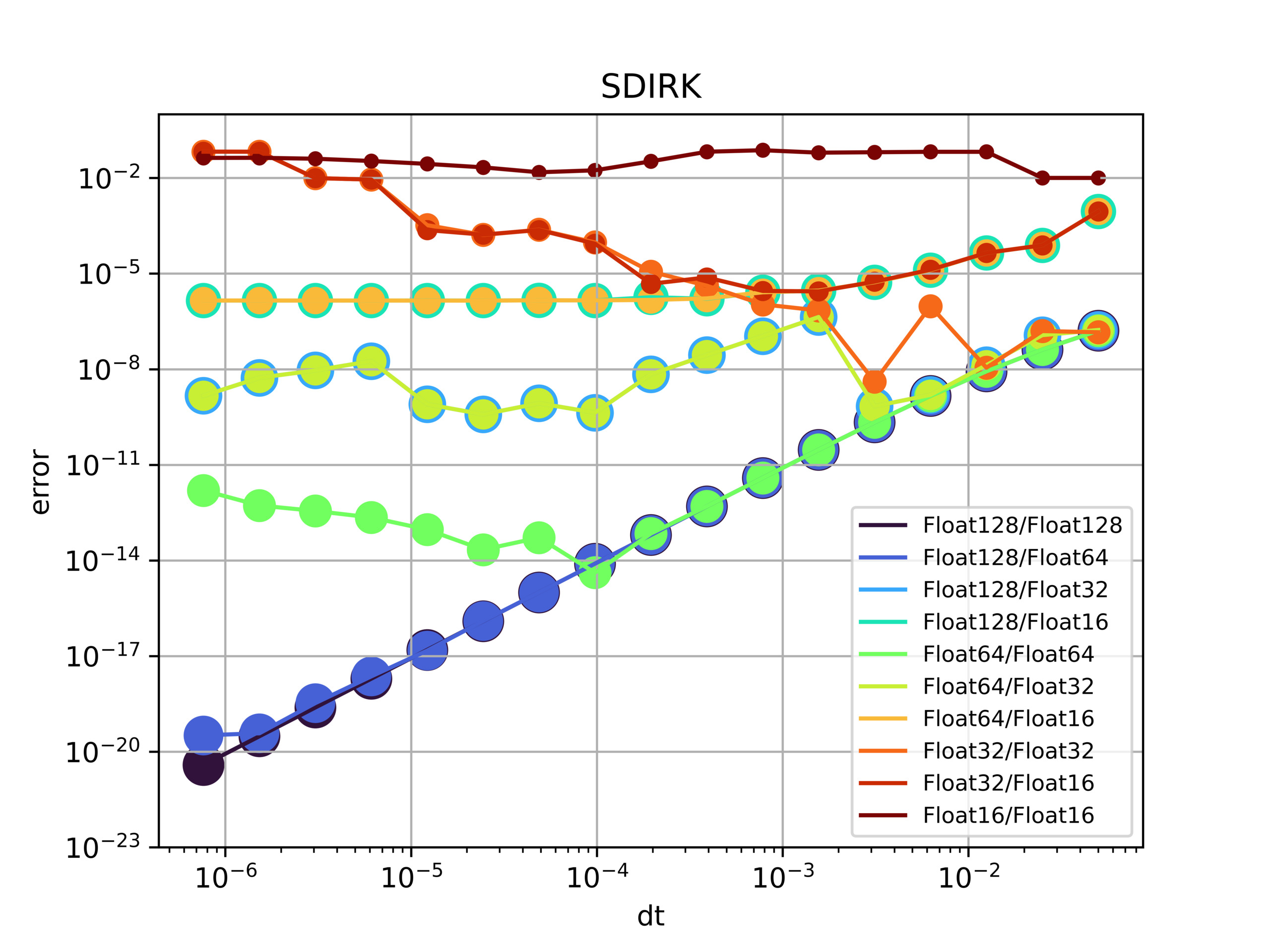}}
{\includegraphics[width=0.24\textwidth]{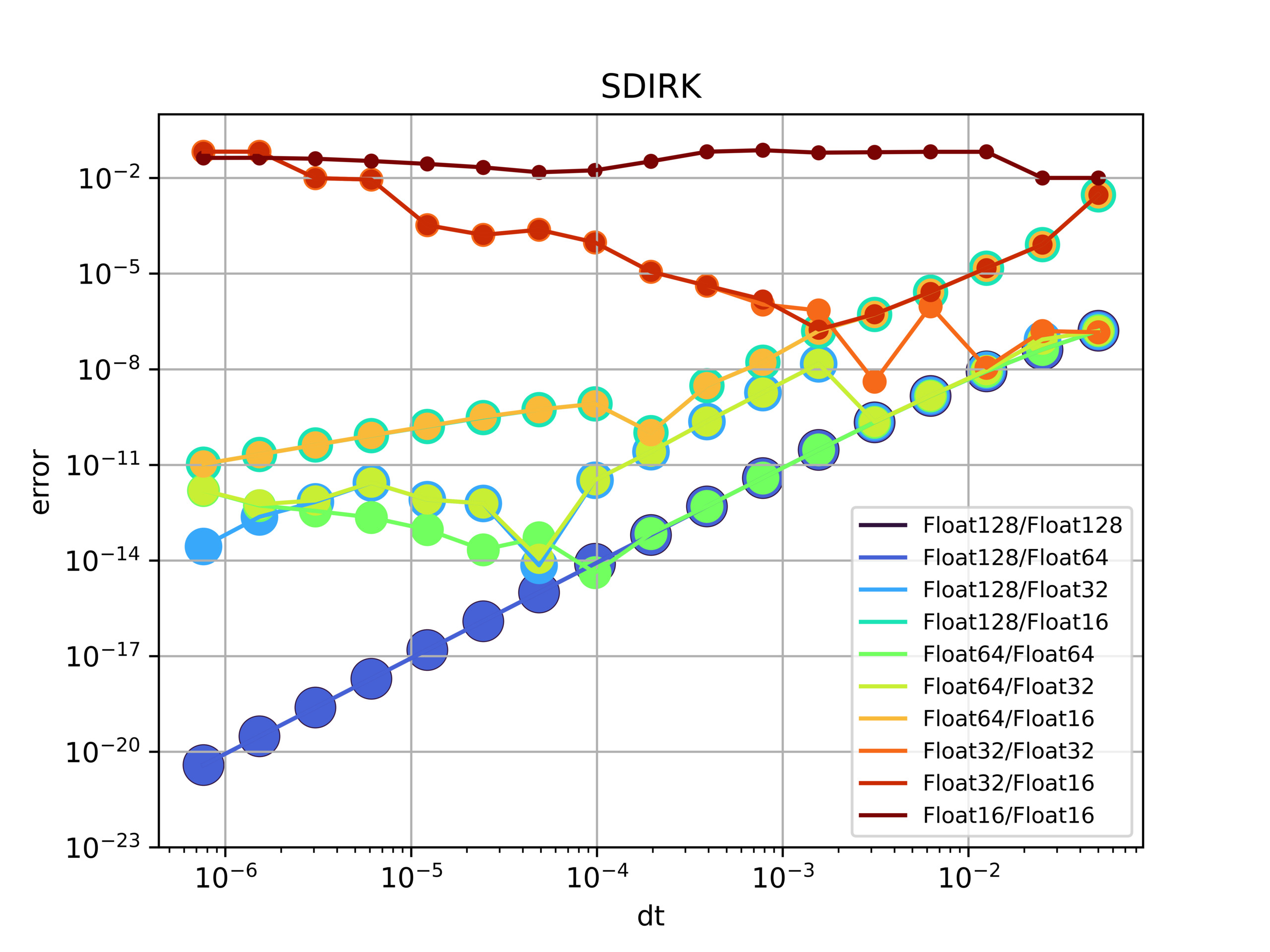}} 
{\includegraphics[width=0.24\textwidth]{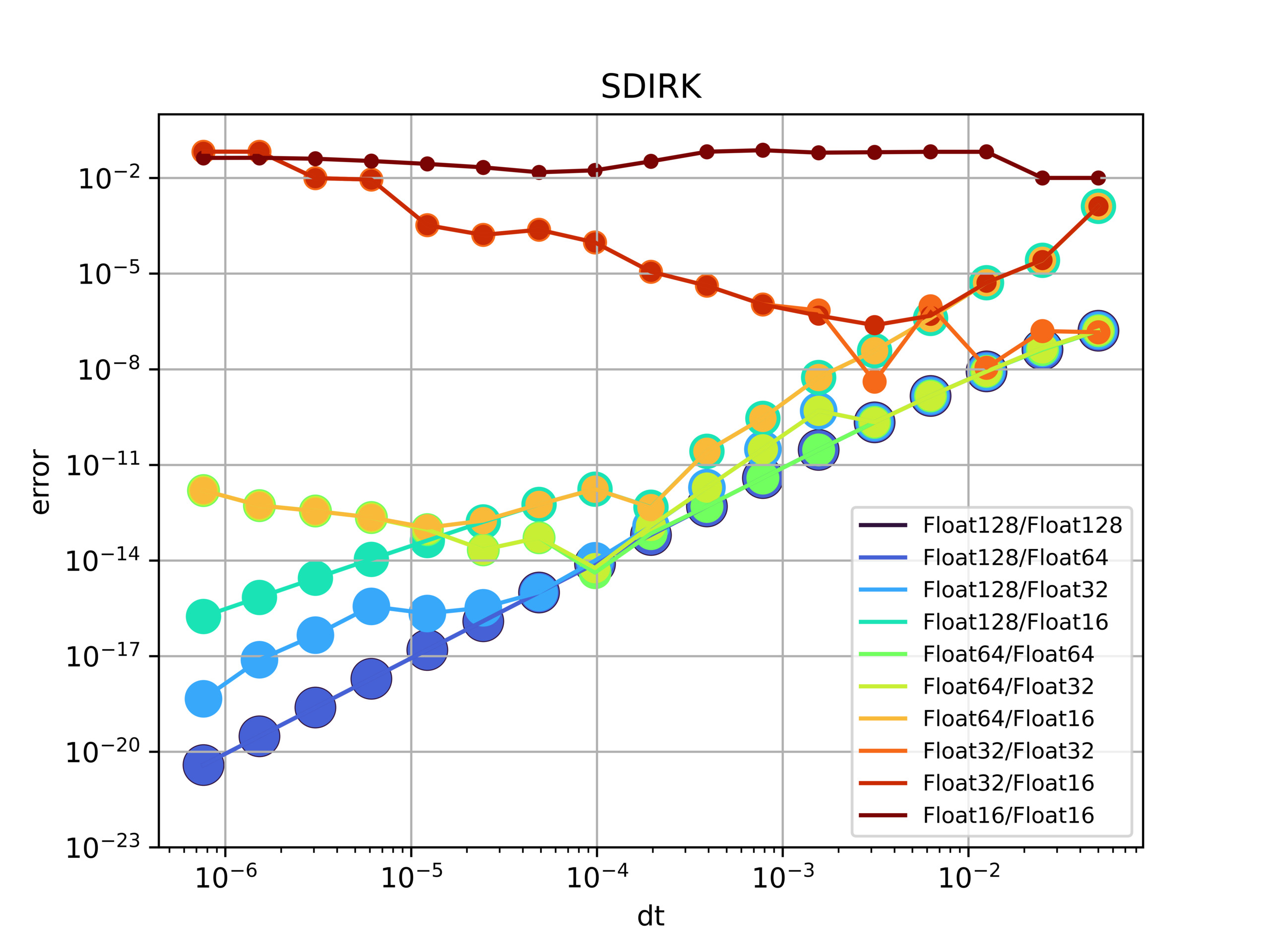}}
{\includegraphics[width=0.24\textwidth]{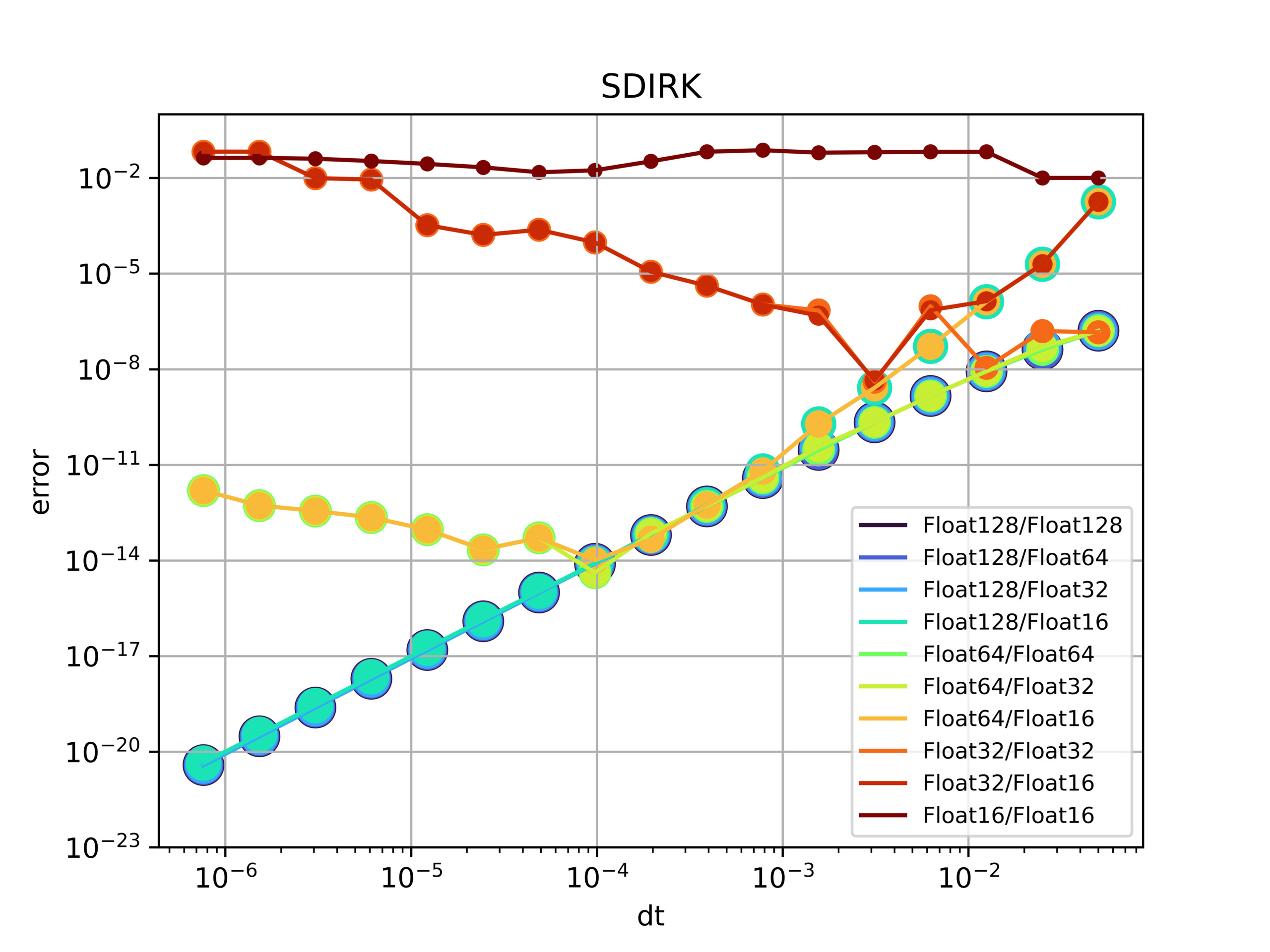}} \\
{\includegraphics[width=0.24\textwidth]{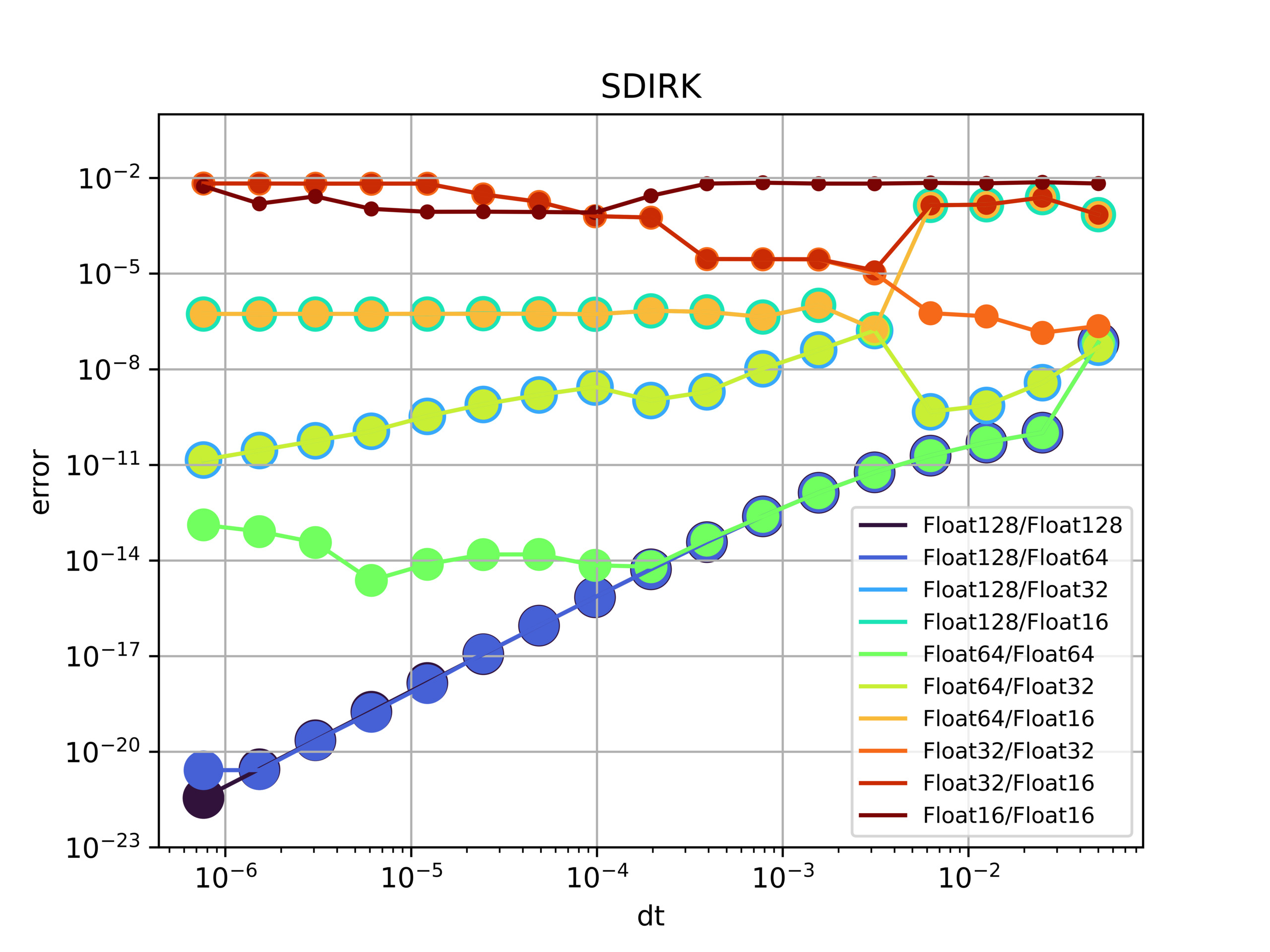}}
{\includegraphics[width=0.24\textwidth]{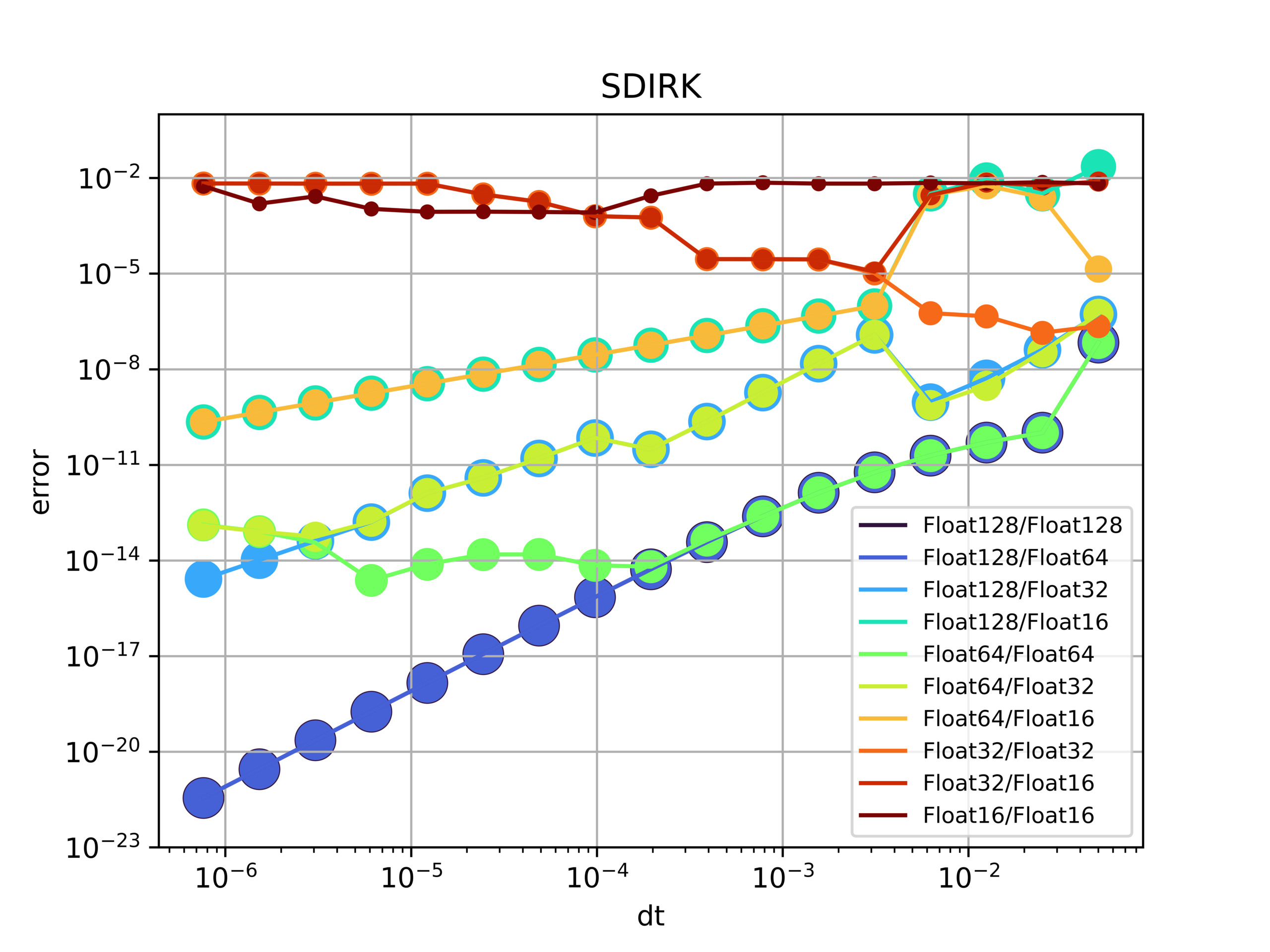}} 
{\includegraphics[width=0.24\textwidth]{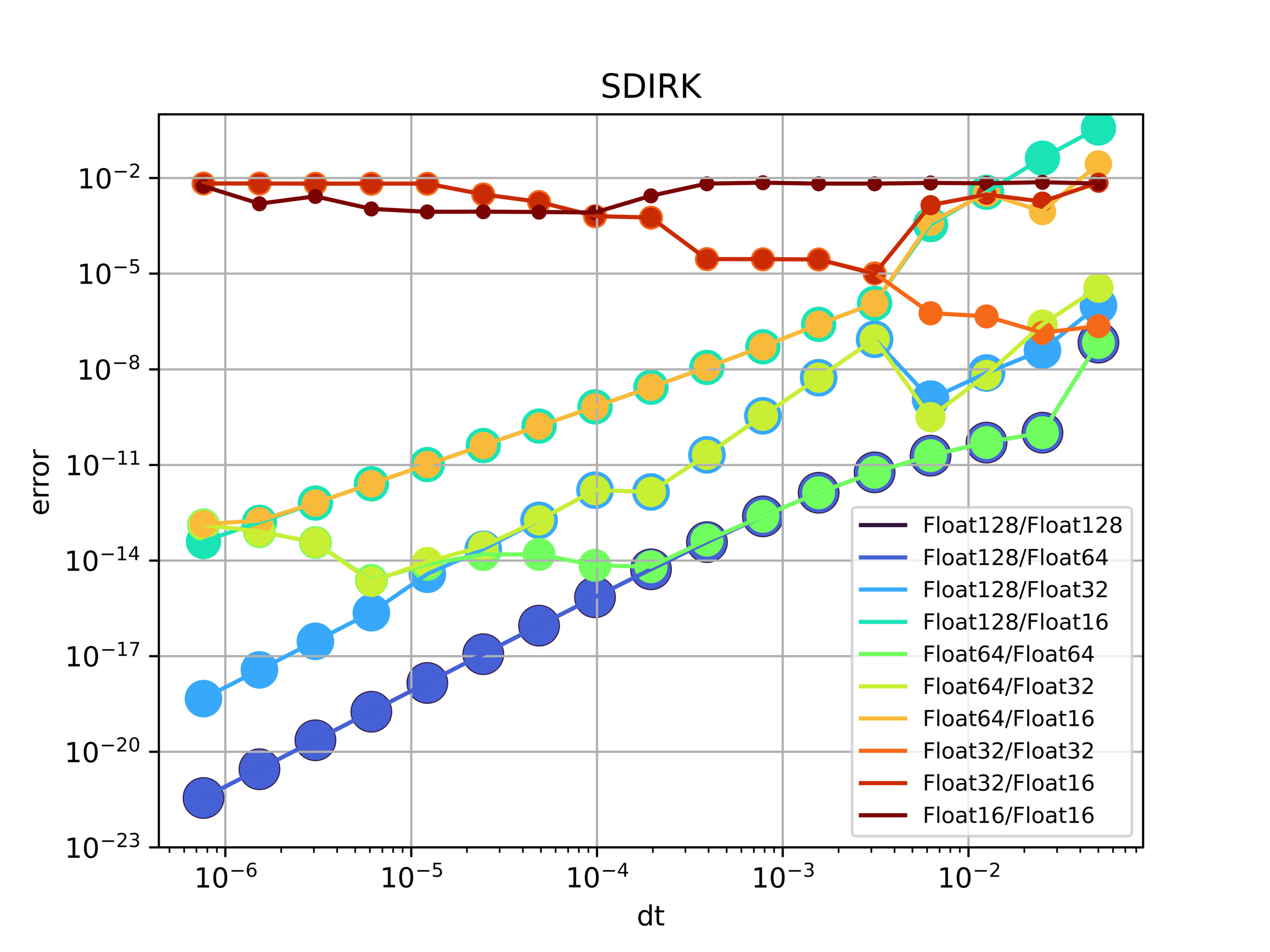}}
{\includegraphics[width=0.24\textwidth]{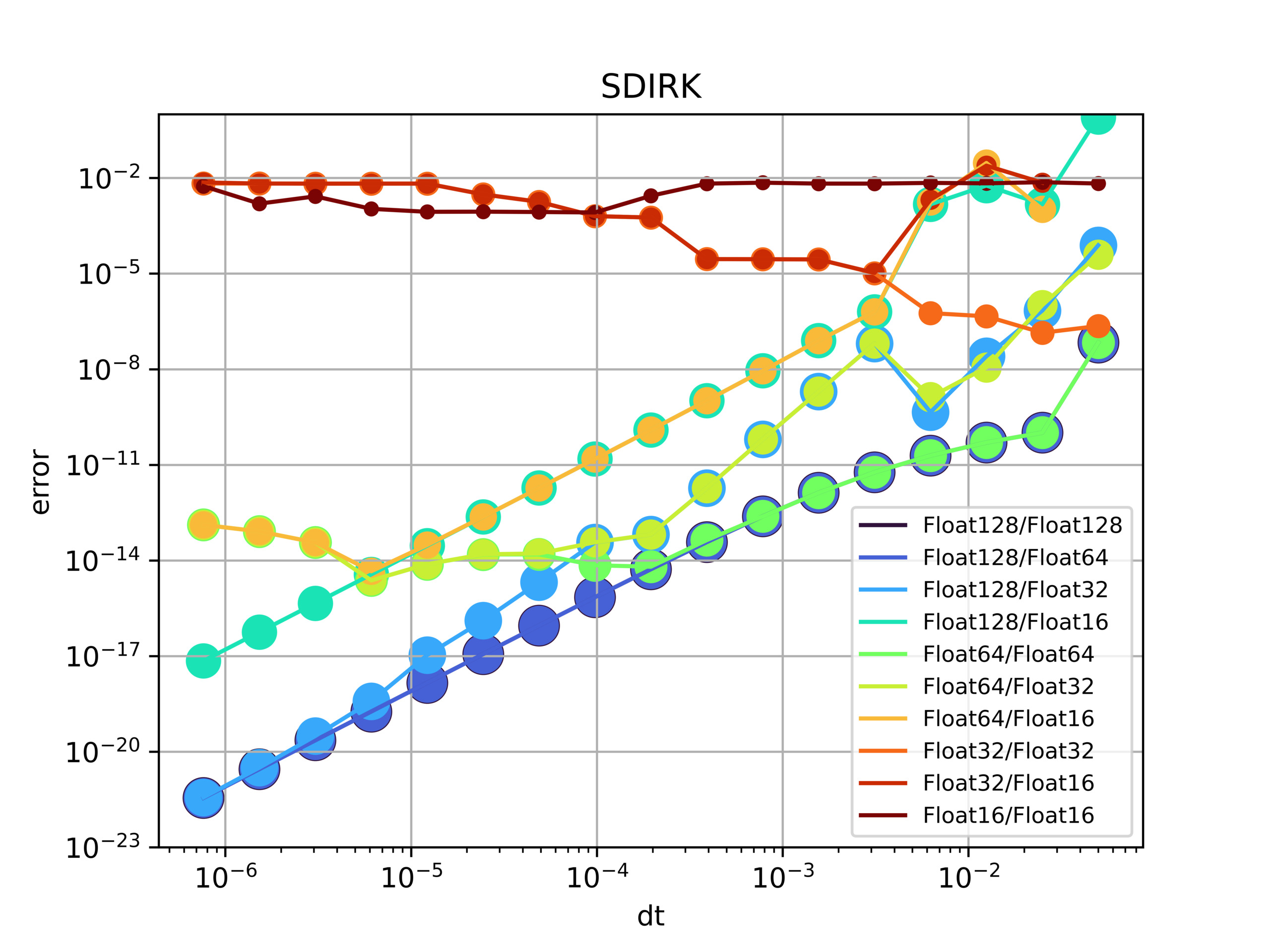}} \\
\caption{Convergence history of the  SDIRK method for the van der Pol equation with with no corrections (left),
one correction (left-middle), two corrections (right-middle), and three corrections (right),
for $\alpha =1,10,100$ (top to bottom).
\label{vdpSDIRKconv} }
\end{center}
\end{figure}

A more detailed set of images is shown in Figure \ref{vdpSDIRKconv}, which 
shows the convergence of the mixed precision and full precision methods 
for the SDIRK methods with up to three corrections for the van der Pol equation with a selection of 
increasing stiffness values  (top to bottom: $\alpha=1, 10, 100$).
For $\alpha =1$, the behavior of the methods with up to three corrections is qualitatively similar to that of 
the implicit midpoint rule in Figure \ref{vdpIMRconv} above. In this case, we see that the convergence of 
the corrections improve the order of convergence of the mixed precision methods to third order.
The effect of the stiffness parameter is notable in that the errors behave worse and require generally smaller
time-steps to settle down -- this is why the bottom right image looks "messy" compared to the ones above it.
An interesting feature is seen in the bottom row images: the errors for the mixed quad/half and double/half codes
are a horizontal line with no corrections, first order with one correction, second order with two corrections, and 
third order with two corrections, albeit with higher errors than the double or quad codes produce.
The mixed precision  quad/single and double/single start at first order with no correction, but look like fourth order with three corrections, which enable the  mixed precision quad/single code to give errors as small as the quad precision code, for a small
enough time-step.

\begin{figure}[htb]
\begin{center}
{\includegraphics[width=0.45\textwidth]{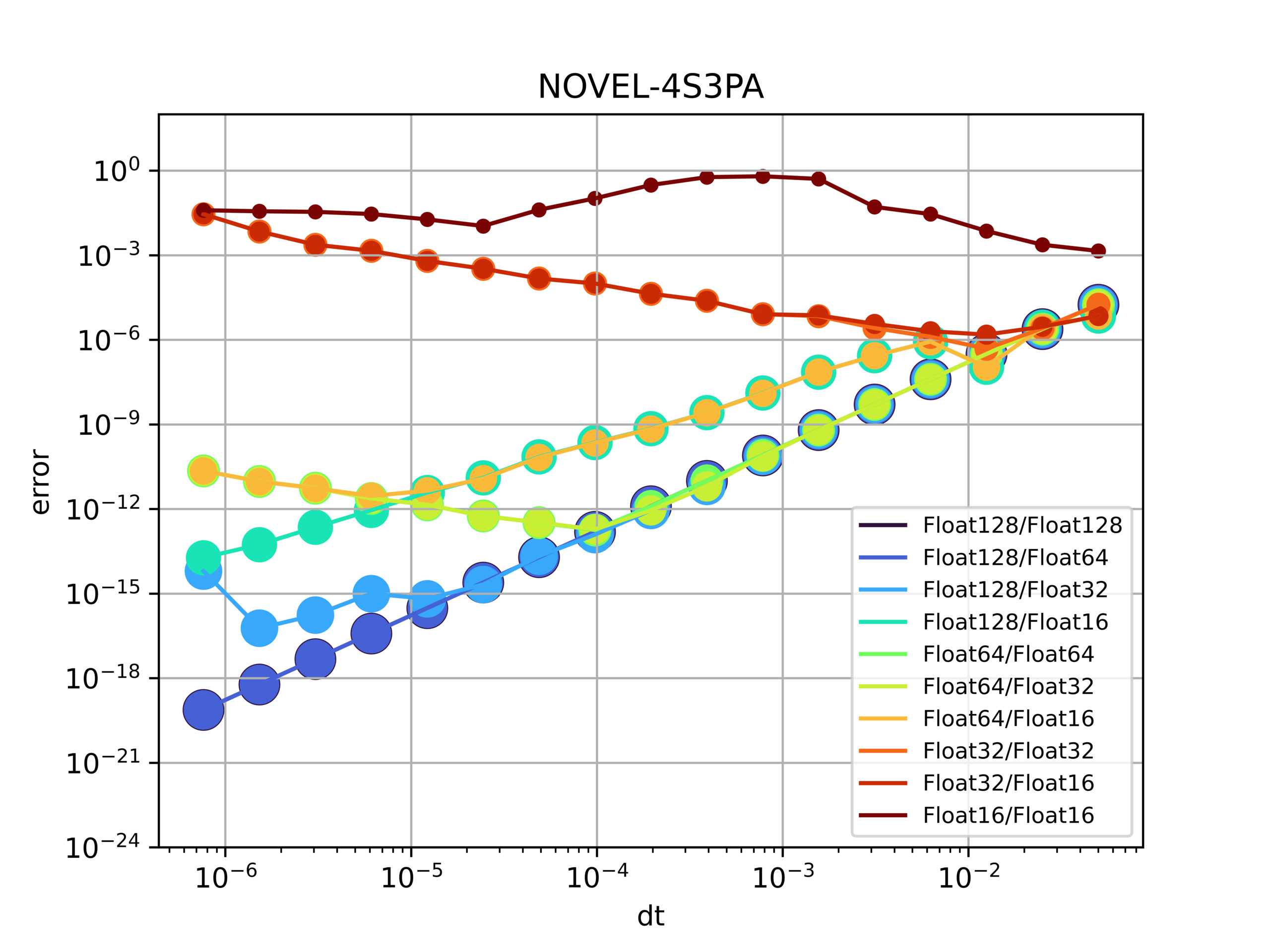}}
{\includegraphics[width=0.45\textwidth]{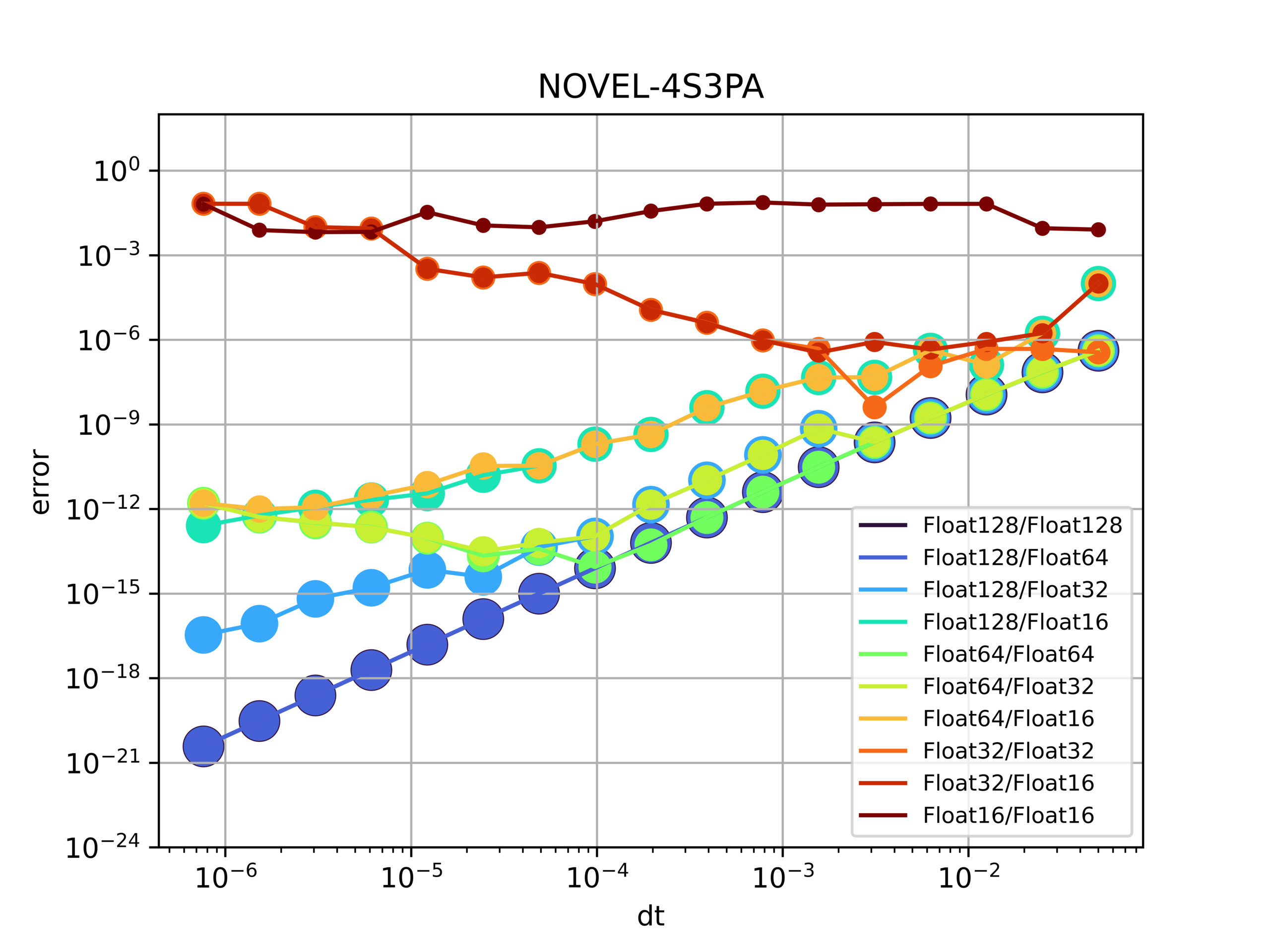}}\\
{\includegraphics[width=0.45\textwidth]{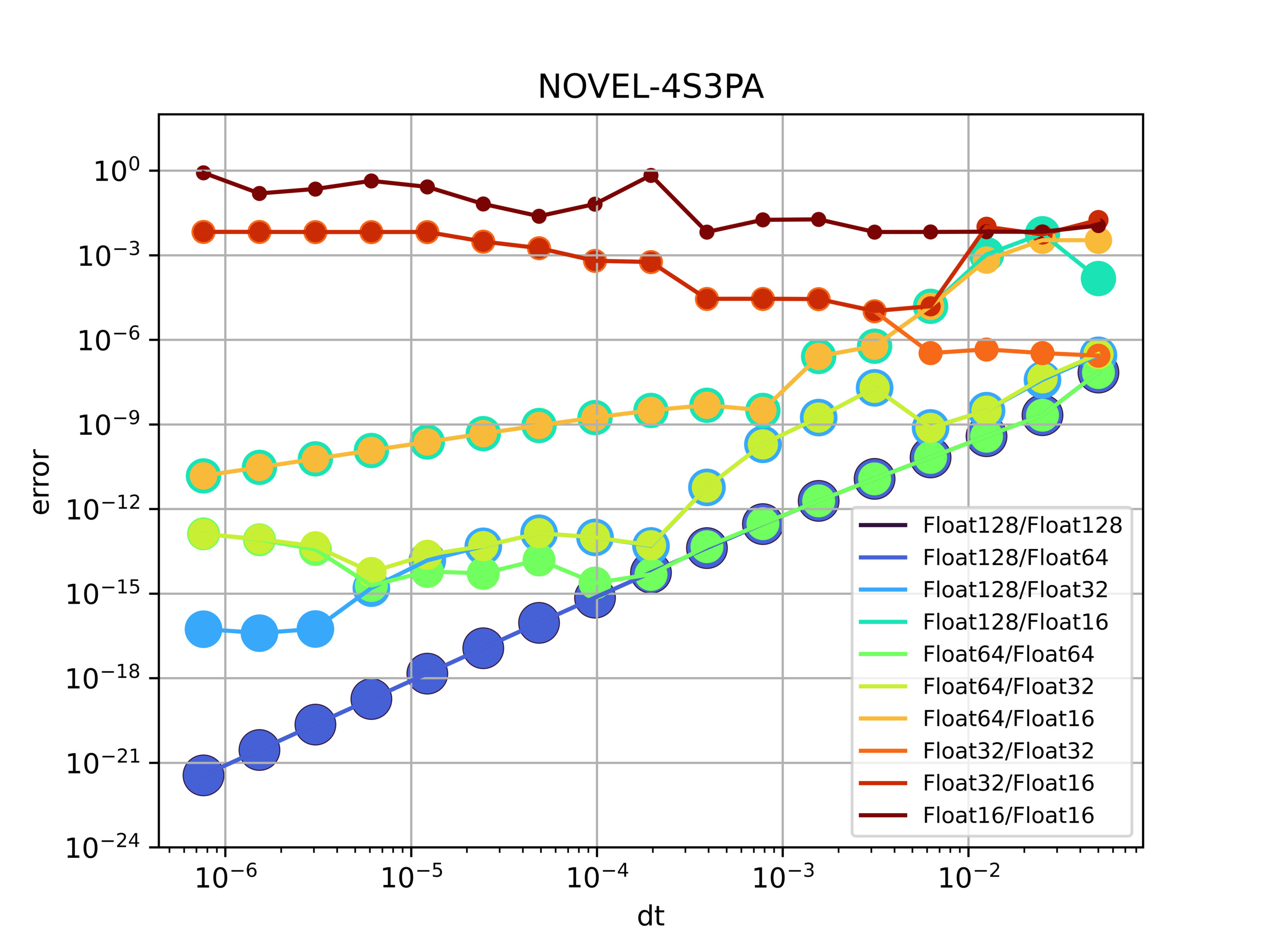}}
{\includegraphics[width=0.45\textwidth]{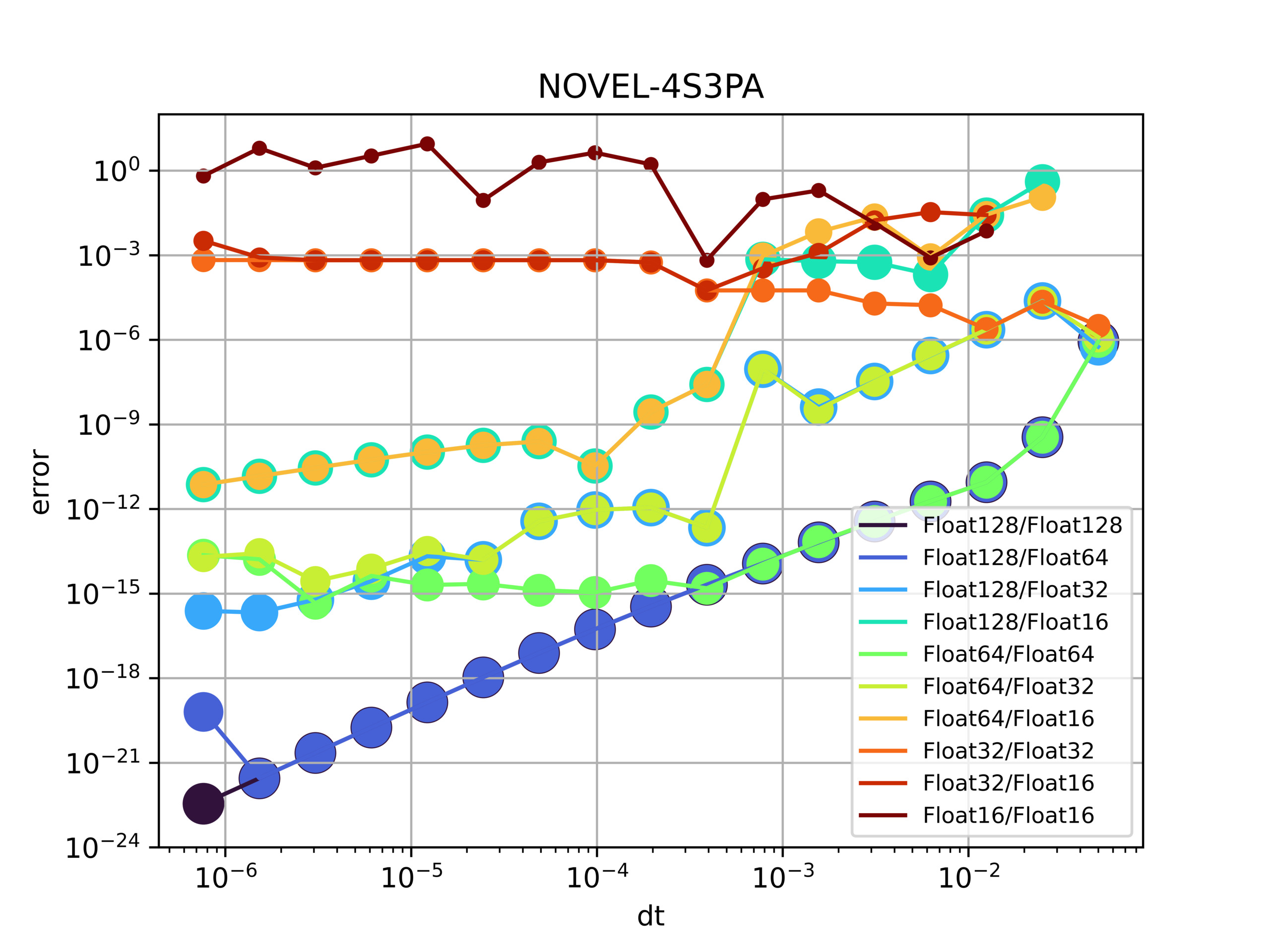}} 
\caption{Convergence history of the  NovelA method  \eqref{MP-4s3pA} for the van der Pol equation  for $\alpha =1,10,100, 1000$ (left to right).}
\end{center}
\label{vdpNOVELconv}
\end{figure}

Similar results can be seen with the NovelA method (Eqn. \eqref{MP-4s3pA})
for the van der Pol equation with a selection of 
increasing stiffness values  (top to bottom: $\alpha=1, 10, 100,1000$).
In Figure \ref{vdpNOVELconv} we observe that as expected, this method performs similarly 
to the SDIRK method with one explicit correction.
The increasing stiffness does not result in dramatically different results, but the 
errors are slightly worse and require generally smaller time-steps to settle down, and the 
error graphs generally look more messy as the stiffness is increased.

\subsection{Convergence of the mixed precision method for a viscous Burgers' equation system}
\label{sec:vBCONV}
For the second convergence test we consider the 1D viscous Burgers' equation
\begin{eqnarray}
u_t + \left(\frac{1}{2} u^2 \right)_x = \frac{1}{100} u_{xx}
\end{eqnarray}
on $x \in (0,1)$ with initial condition $u(x,0) = \sin(2 \pi x)$, and  boundary conditions $u(0,t) =u(1,t)  = 0$.
We begin by semi-discretizing the problem with $N_x$ points in space, where $N_x=50, 100, 200$.
The nonlinear term is approximated using 
a forward difference  and the diffusion term is approximated using  a centered difference. 
The resulting system of ordinary differential equations is evolved to final time $T_f =1.0$
using one of the methods described above. 
The reference solution is computed for $N_x =50, 100, 200$ points in space in quad precision 
using the classical fourth order Runge--Kutta method with $\Delta t = 10^{-7}$. Note that we focus only on the 
errors from the time evolution, so we treat the semi-discretized ODE, not the PDE.

\begin{figure}[htb]
\begin{center}
{\includegraphics[width=0.45\textwidth]{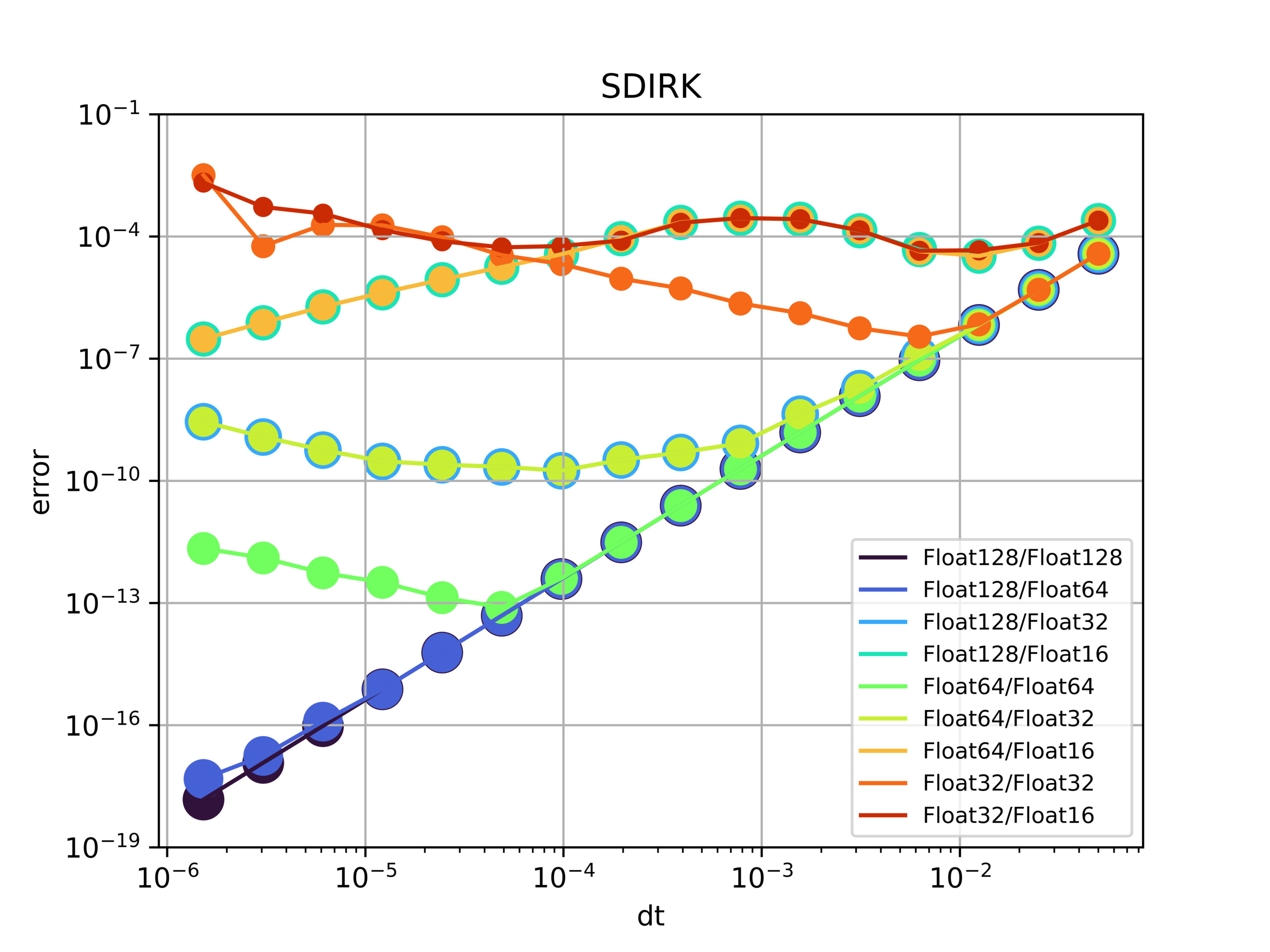}}
{\includegraphics[width=0.45\textwidth]{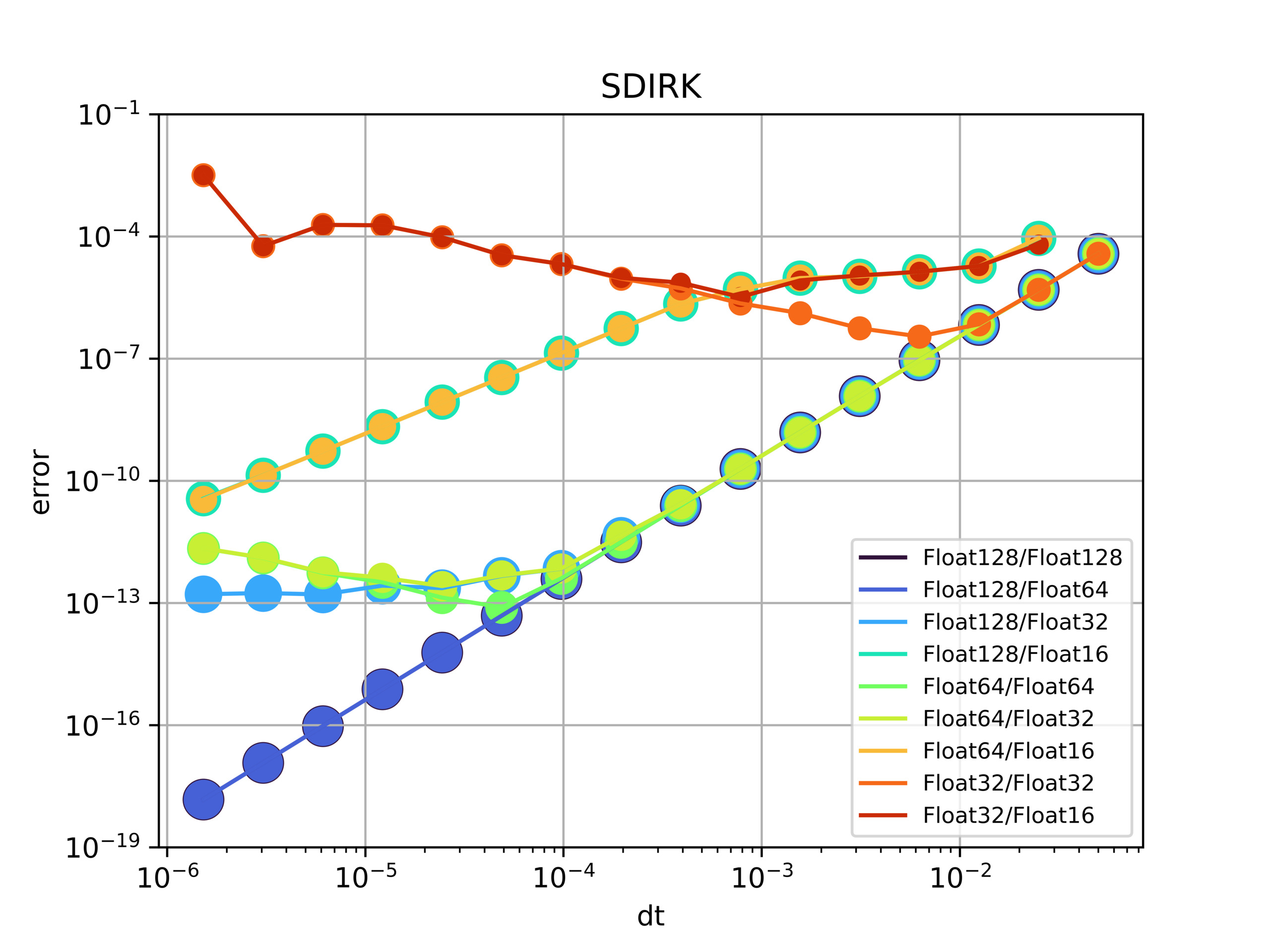}} \\
{\includegraphics[width=0.45\textwidth]{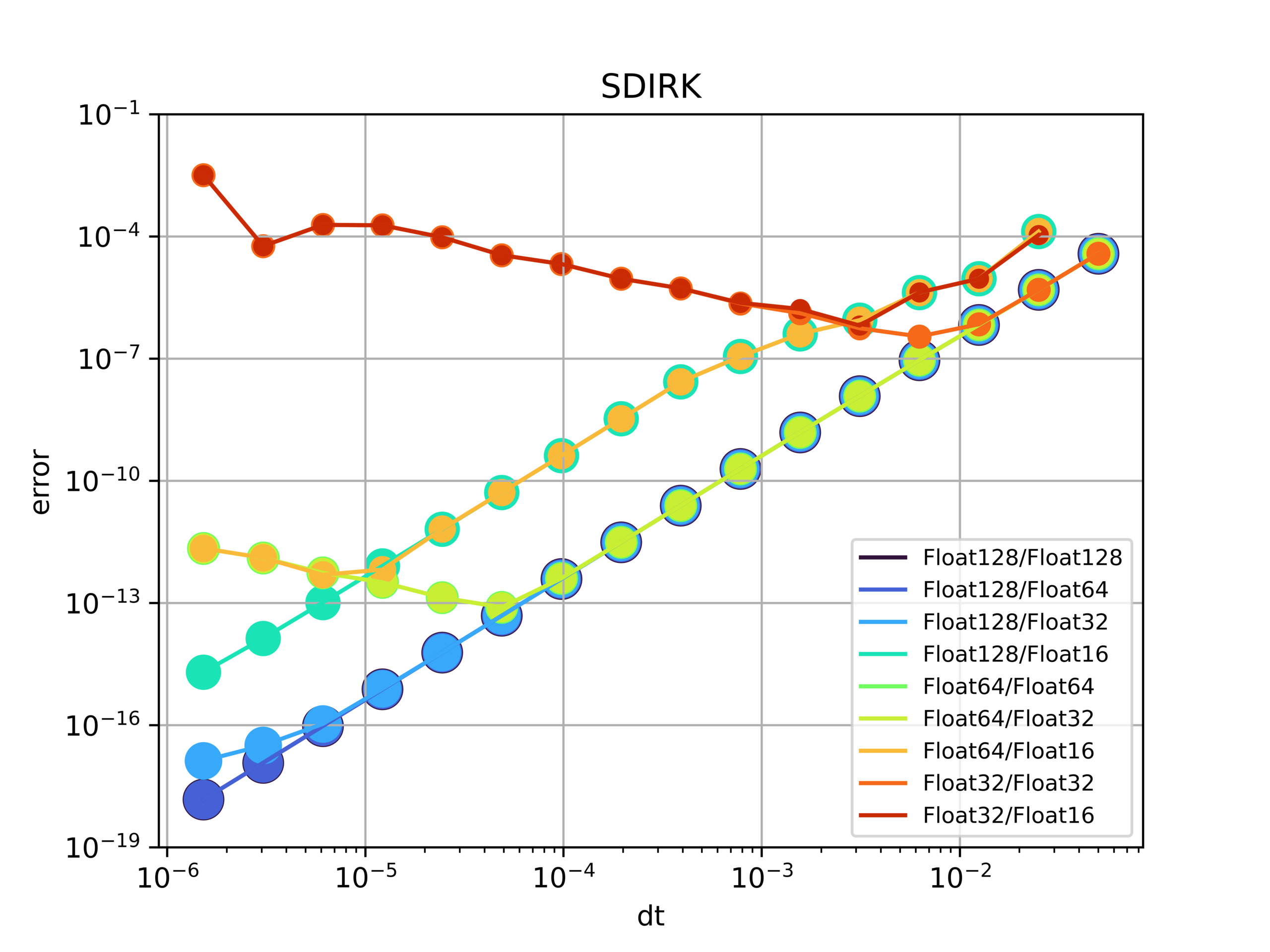}}
{\includegraphics[width=0.45\textwidth]{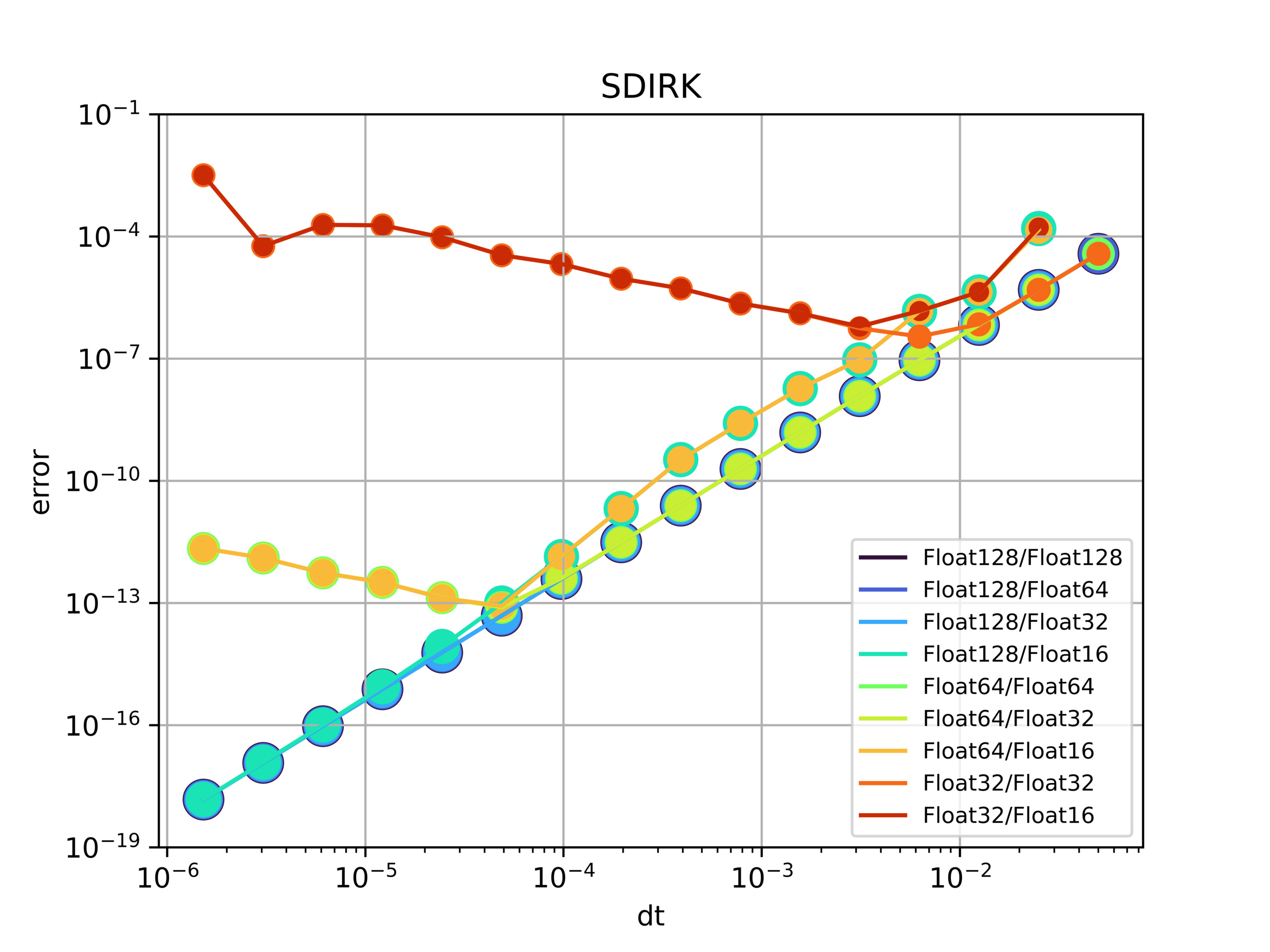}} \\
\caption{Convergence history of the   SDIRK method for the viscous Burgers' equation with $N_x=50$ with no corrections (top left),
one correction (top right), two corrections (bottom left), three corrections (bottom right).
\label{BurgersSDIRKconv}
}
\end{center}
\end{figure}

Figure \ref{BurgersSDIRKconv} shows the convergence history for SDIRK with up to three corrections for $N_x =50$.
The results show the effect of the explicit corrections much like those seen in Subsection \ref{sec:vdpCONV}.
On the top left, we observe that with no corrections, the half precision (Float16)  does not converge,
and the  single precision (Float32)  and mixed precision single/half (Float32/Float16) codes
do not converge as the time-step is refined below $ \Delta t \leq 10^{-2}$ .
The mixed precision quad/half (Float128/Float16) and double/half codes  (Float64/Float16) 
both converge at first order when the time-step is sufficiently refined,  $\Delta t \leq 10^{-3.5}$.
The mixed precision quad/single (Float128/Float32) and double/single   (Float64/Float32) codes
follow the same third order convergence line as the double precision and quad precision  
codes initially (for $ 3\times 10^{-3}  \leq \Delta t$),
but the errors do not continue to decay, and in fact start to grow, as $\Delta t$ is refined below $10^{-3}$.
The double precision (Float64) code continues to converge at third order until $ \Delta t \approx 5 \times10^{-5}$,
and the mixed precision quad/double (Float128/Float64) and quad precision code continue to converge at third order
until much further down. When the errors fall below the level of $10^{-16}$ the errors from the mixed quad/double code 
start to separate from the full quad precision code. However, one explicit correction (top right) changes this: once a
correction is applied, the errors from the quad precision and mixed quad/double precision codes look identical.
One correction also causes the The mixed precision quad/single (Float128/Float32) and double/single   (Float64/Float32) codes
follow the same third order convergence line as the double precision and quad precision  
codes all the way down  to   $  \Delta t \approx 10^{-4} $, which is approximately the where the errors from the 
double precision code bifurcates from those of the quad precision codes. Two explicit corrections 
reduce the errors for the quad/single (Float128/Float32) code to those of the quad codes, and three explicit 
corrections reduce the errors for the quad/single (Float128/Float32) code to those of the quad codes as well.

\begin{figure}[htb]
\begin{center}
{\includegraphics[width=0.32\textwidth]{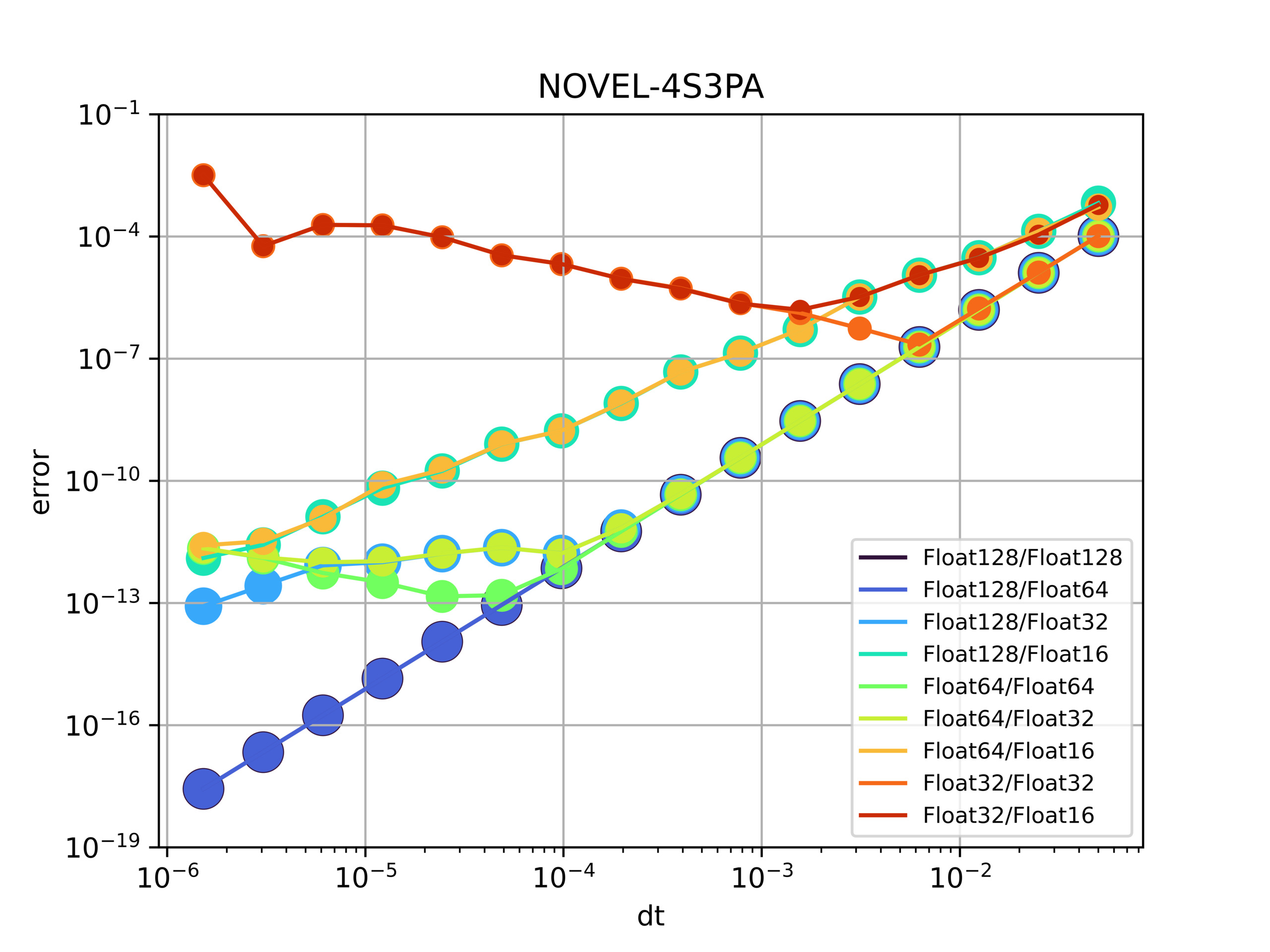}}
{\includegraphics[width=0.32\textwidth]{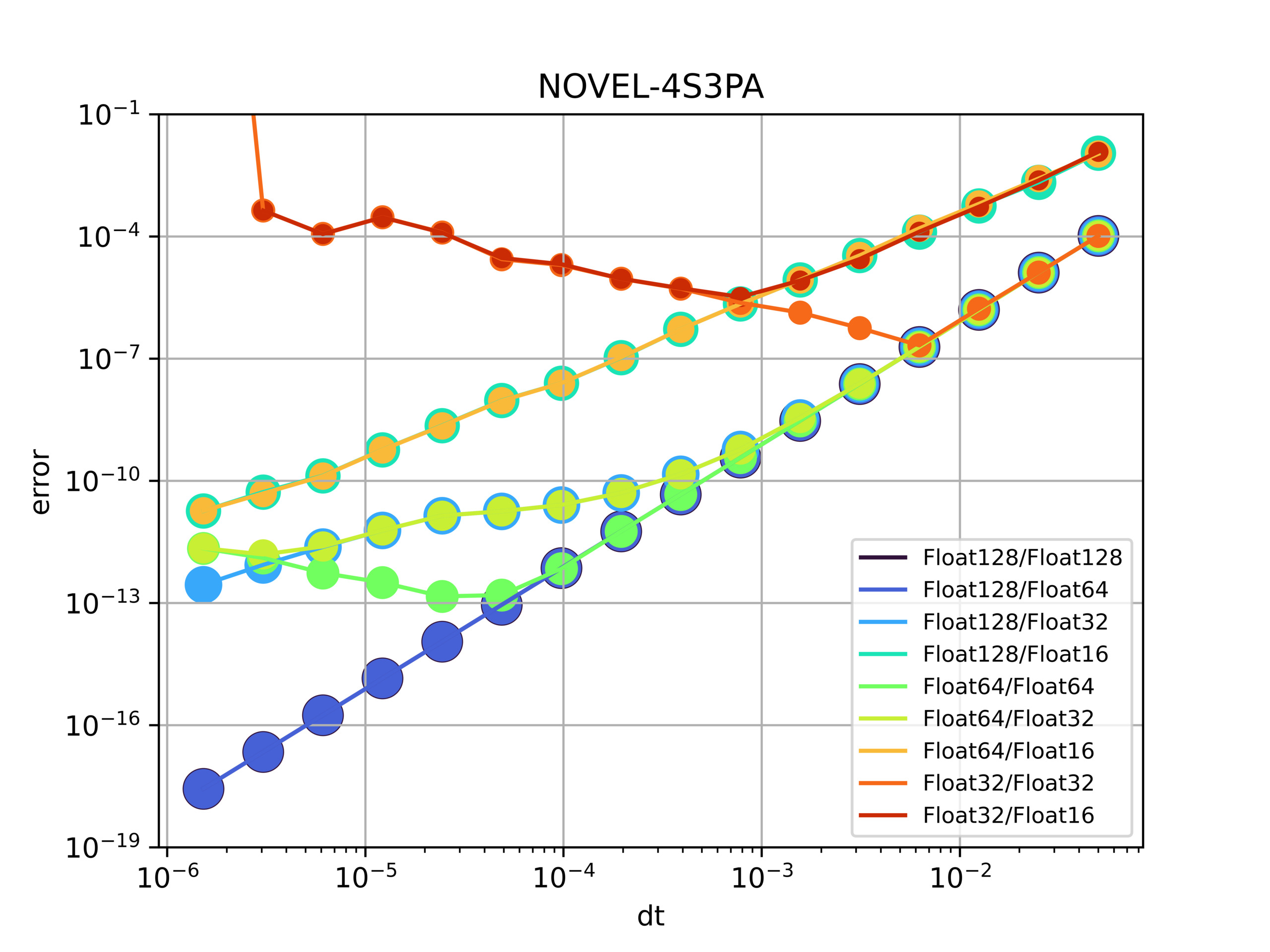}}
{\includegraphics[width=0.32\textwidth]{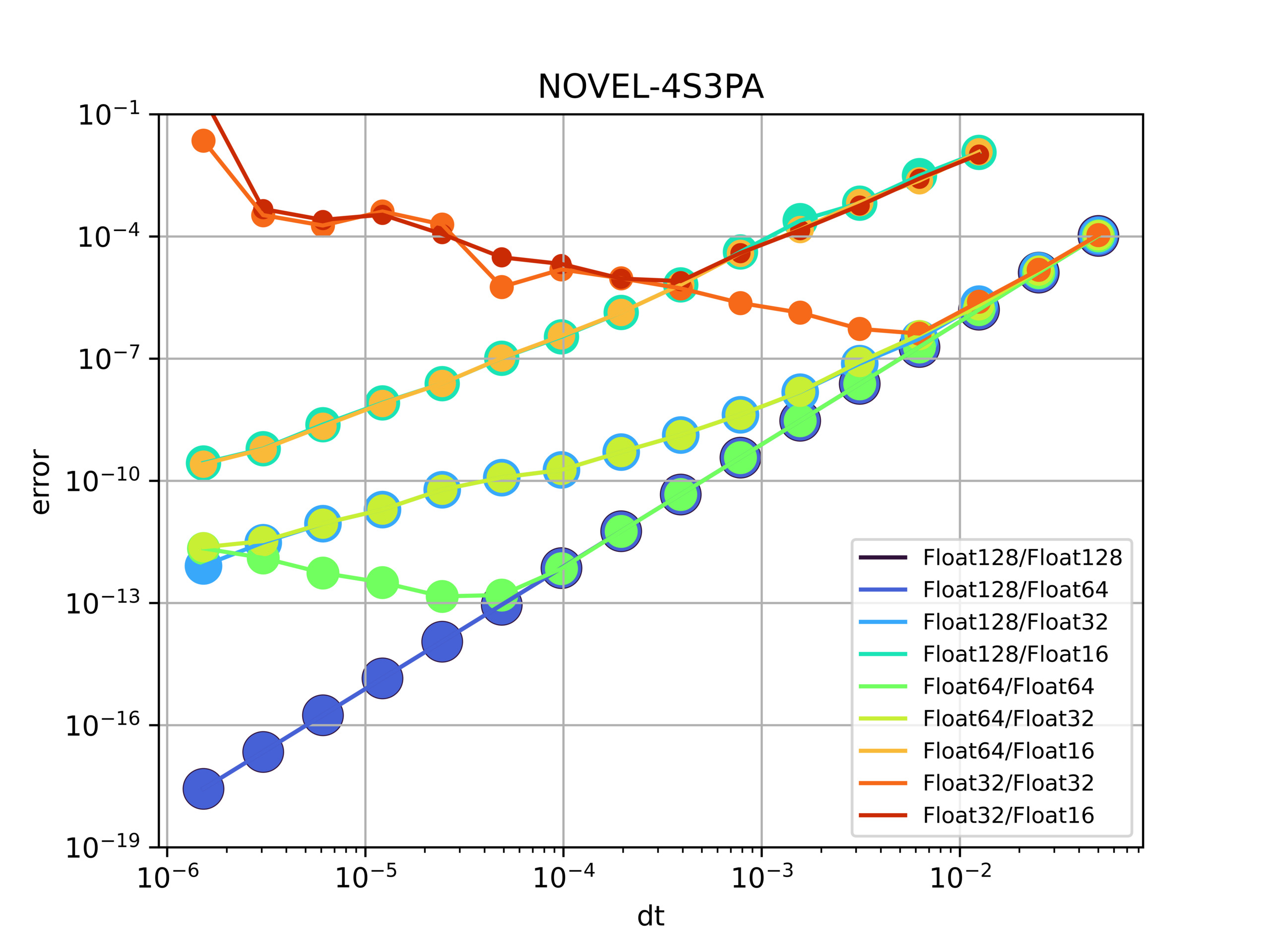}}\\
\caption{The novel method for the viscous Burgers' equation with $N_x=50$ (left), $N_x=100$ (middle), and $N_x=200$ (right).
\label{BurgersNovelConv}}
\end{center}
\end{figure}

Figure \ref{BurgersNovelConv} shows the behavior of the errors of the mixed precision codes using the NovelA
method \eqref{MP-4s3pA} with
$N_x=50$ (left), $N_x=100$ (middle), and $N_x=200$ (right) points in the spatial discretization. 
These results are comparable to those of the  SDIRK method with one correction. 
This figure highlights the differences as the ODE system is  stiffer.
The main difference between  is that the errors from the 
mixed precision quad/single and double/single stop following the quad precision third order convergence line sooner 
as the problem is stiffer: for $N_x =50$, third order convergence is seen down to $  \Delta t \approx 9.7 \times 10^{-5} $
while for $N_x =100$, the third order convergence line is only followed down to $  \Delta t \approx 8 \times 10^{-4} $,
and for  $N_x =200$,  only down to  $\Delta t \approx 6 \times 10^{-3} $.

\section{Linear stability and sensitivity to roundoff of MP-ARK methods}

\subsection{Motivation from mixed model simulation}
To motivate our stability analysis, we begin with a mixed model approach 
rather than with  the mixed-precision approach we have been focusing on.
 In this mixed model approach, we begin with a PDE, and use two different spatial discretizations
to obtain obtain two different ODEs. Assuming we begin with a linear problem, the ODE systems become
\[ u_t = F(u), \; \; \; \mbox{and} \; \; \; u_t = \Fep(u) \]
where $F$ is a costlier, but more accurate, spatial discretization, and $\Fep$ is cheaper, but less accurate.

An example of this was given in \cite{Grant2022}, where we solve the heat equation
\[ u_t = u_{xx}, \; \; \; \mbox{on} \; \; x \in [0, 2\pi ]  \] with periodic boundary conditions and 
initial condition $u(x,0) = \frac{1+\sin(x)}{2}$. We discretize the grid with $N_x$ equidistant points, ans
discretize the spatial derivative in two different ways: using a Fourier spectral method discretization to
create the differentiation matrix $D_s$ and using a simple 3-point centered differencing to create the 
differentiation matrix $D_c$. The spectral approach gives exponential convergence in $N_x$, while the centered difference
gives second order convergence in $N_x$. However, the centered difference matrix is tridiagonal and so the 
implicit steps are inexpensive, while the spectral differentiation matrix is a full matrix and so the implicit steps are
more costly. 

\begin{SCfigure}[][h]
  \includegraphics[width=0.5\textwidth]{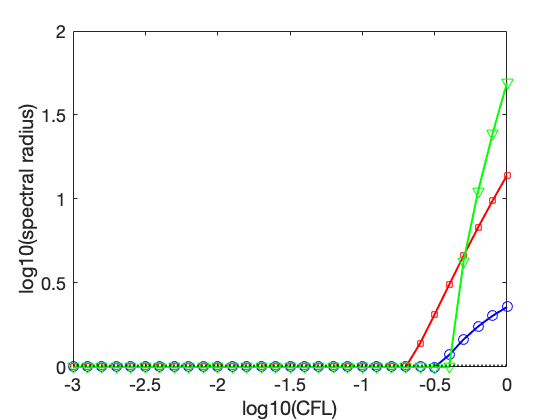}
  \caption{The spectral radius of the stability matrix $P$ \eqref{StabilityMatrix}
   of the mixed model implicit midpoint rule with \\ 
no corrections (blue circles), \\ one correction (red squares), \\
and two corrections (green triangles), \\
  for the heat equation, where the spatial derivative is discretized with a centered difference differentiation 
  matrix for the implicit stage and the spectral method differentiation matrix for the explicit stages.
  }
  \label{IMRmixedmodel}
\end{SCfigure}

We confirmed that for sufficiently small $\dt$ the method converges as designed, and the additional
corrections improve the order of convergence as expected. However, we observed that the 
time-step for which this method is stable is small, and becomes generally worse with more corrections. 
To understand this, we  analyze the stability of this system for the mixed precision implicit midpoint rule, 
 observe that for the mixed model problem the method with no corrections becomes 
\begin{eqnarray*}
y_1 & = & u^n + \frac{1}{2} \Delta t D_c \; y_1 \\
u^{n+1} & = & u^n +  \Delta t D_s \; y_1 \\
\end{eqnarray*}
so that a one-step evolution is 
$ u^{n+1}  =  \left( I + \dt D_s \left(I -  \frac{1}{2} \Delta t D_c \right)^{-1} \right) u^n  $
and to assess the stability we need to look at the size of the maximum absolute value eigenvalue of the stability polynomial.
Figure \ref{IMRmixedmodel} shows the magnitude of the eigenvalues as a function of the CFL coefficient 
$CFL = \frac{\dt}{\Delta x^2}$.
We observe clearly that if $\lambda >0.3$ this mixed model scheme is unstable for this problem. This is interesting, because using
only $D_s$ or only $D_c$ results in an unconditionally stable scheme, 
and even using a mixed model approach but where $D_s$ is used in
the implicit stage and $D_c$ in the explicit stage  is unconditionally stable. 
Additional corrections gives us a one-step evolution stability matrix $P$ with $c$ corrections is given by 
\begin{eqnarray} \label{StabilityMatrix}
 P = \left(I + \Delta t D_s \left( I + \frac{1}{2} \Delta t D_s \right)^c \left(I -  \frac{1}{2} \Delta t D_c \right)^{-1} \right).  
 \end{eqnarray}
 Figure \ref{IMRmixedmodel} shows that the corrections do not significantly improve the stability properties, and 
 in fact when $\dt$ is large additional corrections lead to higher errors.

As we see in Figure \ref{IMRmixedmodel}, use of the mixed model approach is limited by the fact that 
the two discretizations need to work well with each other.
The observed reduction in stability can be explained by the fact that the implicit stages stabilize the method by damping out certain modes.  In the case where we use the centered difference differentiation matrix $D_c$ the 
modes being damped out do not compensate for the growth of the modes introduced by $D_s$ in the explicit stage.

\subsubsection{Linear stability analysis for mixed precision methods} \label{sec:stability}
Given a mixed precision additive Runge--Kutta method \eqref{RK-Butcher}
\begin{eqnarray*}
y^{(i)} & = & u^n + \dt \sum_{j=1}^{s} A_{ij} F(y^{(j)}) +  \dt \sum_{j=1}^{s} \Aep_{ij} \Fep (y^{(j)})\\
u^{n+1} & = & u^n + \dt \sum_{j=1}^{s} b_{j} F(y^{(j)}) +  \dt \sum_{j=1}^{s} \bep_{j} \Fep (y^{(j)}).
\end{eqnarray*}
for linear problems of the form
\[ u_t =  \lambda u, \; \; \; \mbox{and} \; \; \; u_t = \lambda_\epsilon u \]
the stability polynomial is
\begin{eqnarray*}
1 + \left( \dt  \lambda  b+ \dt  \lambda_\epsilon  \bep \right) \left( 1 - \dt A  \lambda - \dt \Aep \lambda_\epsilon \right)^{-1} e .
\end{eqnarray*}

If we define $ \lambda_\epsilon u  =  \left( \lambda  + \epsilon  \tau \right) u $ then this becomes (recall that $\bep=0$ in our cases)
\[ 1 +  \dt  \lambda  b  \left( 1 - \dt ( A  + \Aep) \lambda - \dt  \epsilon \tau \right)^{-1} e .\]
Setting $z= \lambda  \dt$ and $q = \dt \epsilon \tau$  (where we we assume that  $|q| \approx  \dt \epsilon $),
the stability polynomial is 
\begin{eqnarray}
1 +  z  b  \left( 1 - z ( A  + \Aep)  -  z  \frac{\epsilon}{\lambda}  \tau  \right)^{-1} e .
\end{eqnarray}
To highlight the connection between the three variables, we define \[ \tilde{\epsilon} =  \frac{\epsilon}{\lambda} ,\]
and  look at the linear stability regions defined by 
\begin{eqnarray}
\Psi_\epsilon = 1 +  z  b  \left( 1 - z ( A  + \Aep)  - z   \tilde{\epsilon}  \tau  \right)^{-1} e ,
\end{eqnarray}
where we simulate rounding error using stochastic rounding and assign $\tau$ to be a random number between  
$-\frac{1}{2}$ and $\frac{1}{2}$.
This allows us to look at the stability regions of the methods, and the impact of the precision and stiffness on the 
region of stability.

\begin{figure*}[ht!]
  \includegraphics[width=0.325\textwidth]{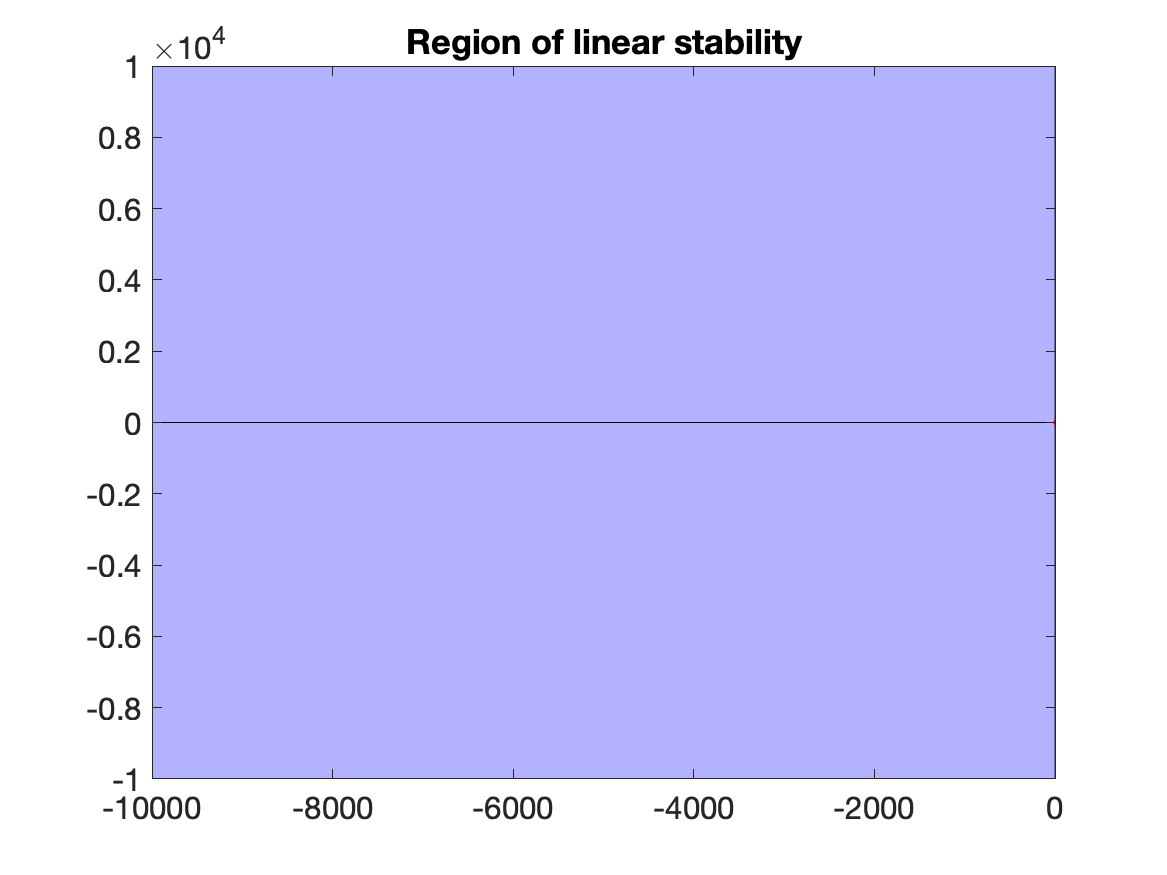}
    \includegraphics[width=0.325\textwidth]{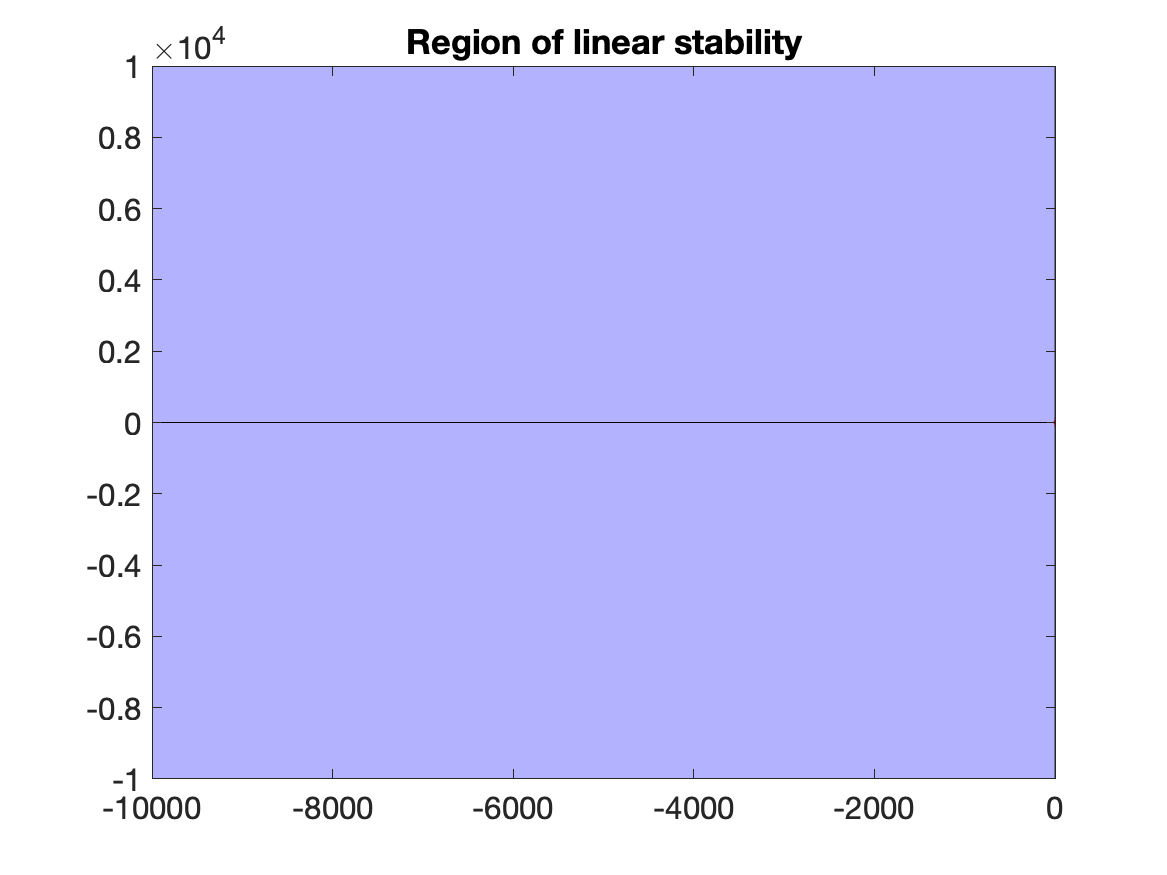}
      \includegraphics[width=0.325\textwidth]{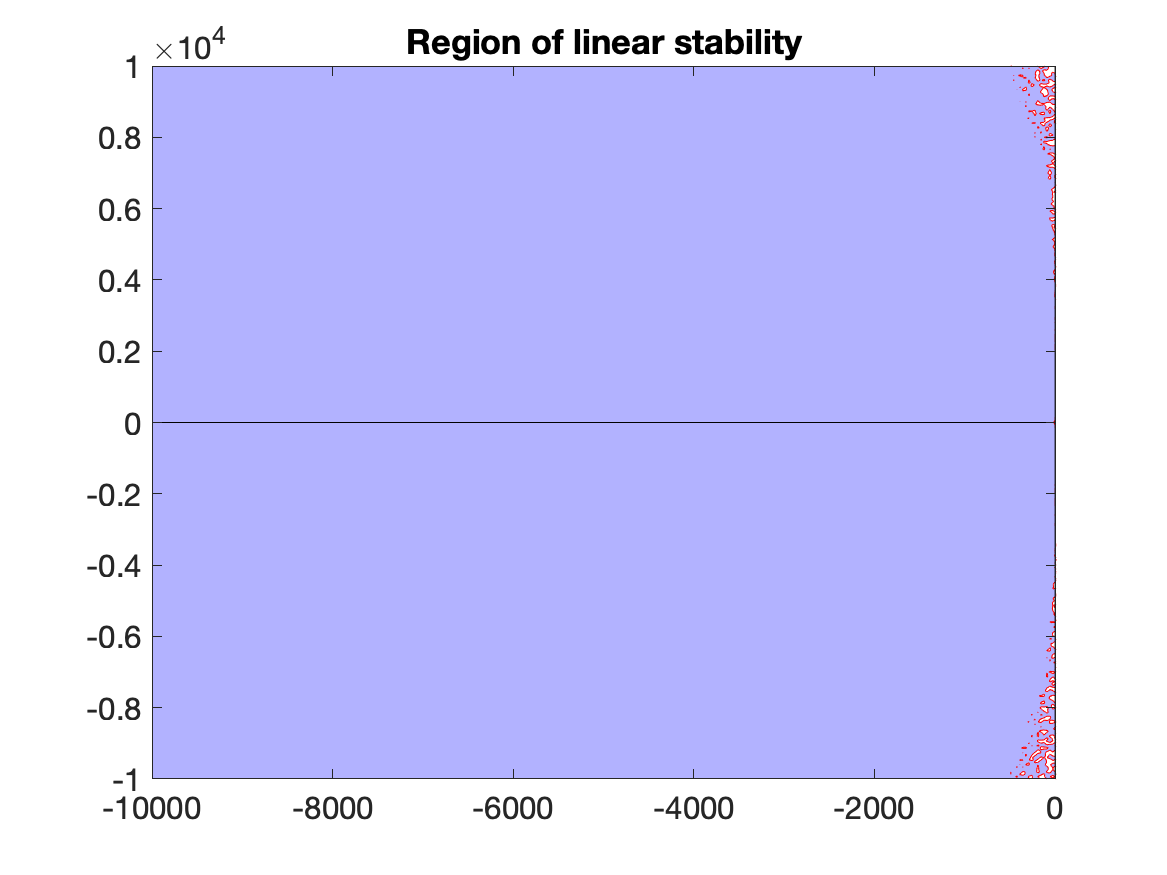} \\
  \includegraphics[width=0.325\textwidth]{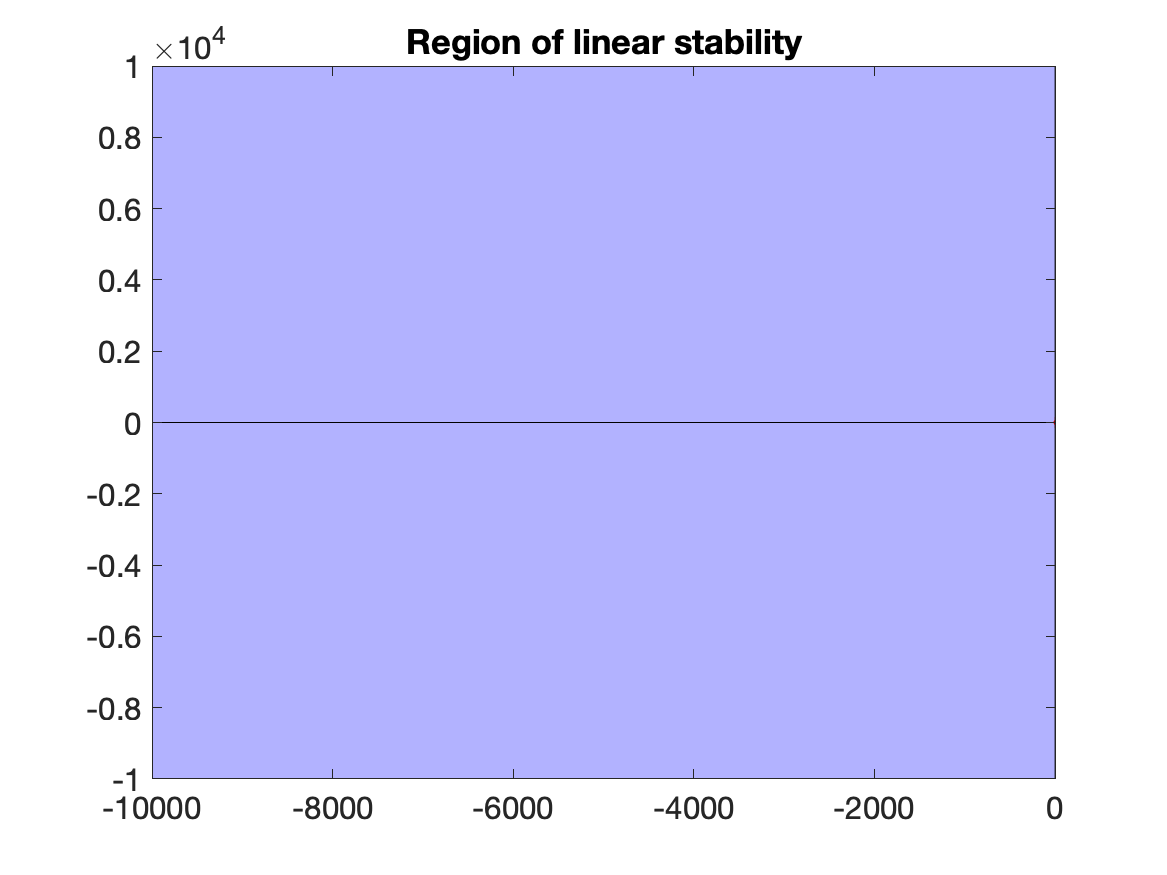}
    \includegraphics[width=0.325\textwidth]{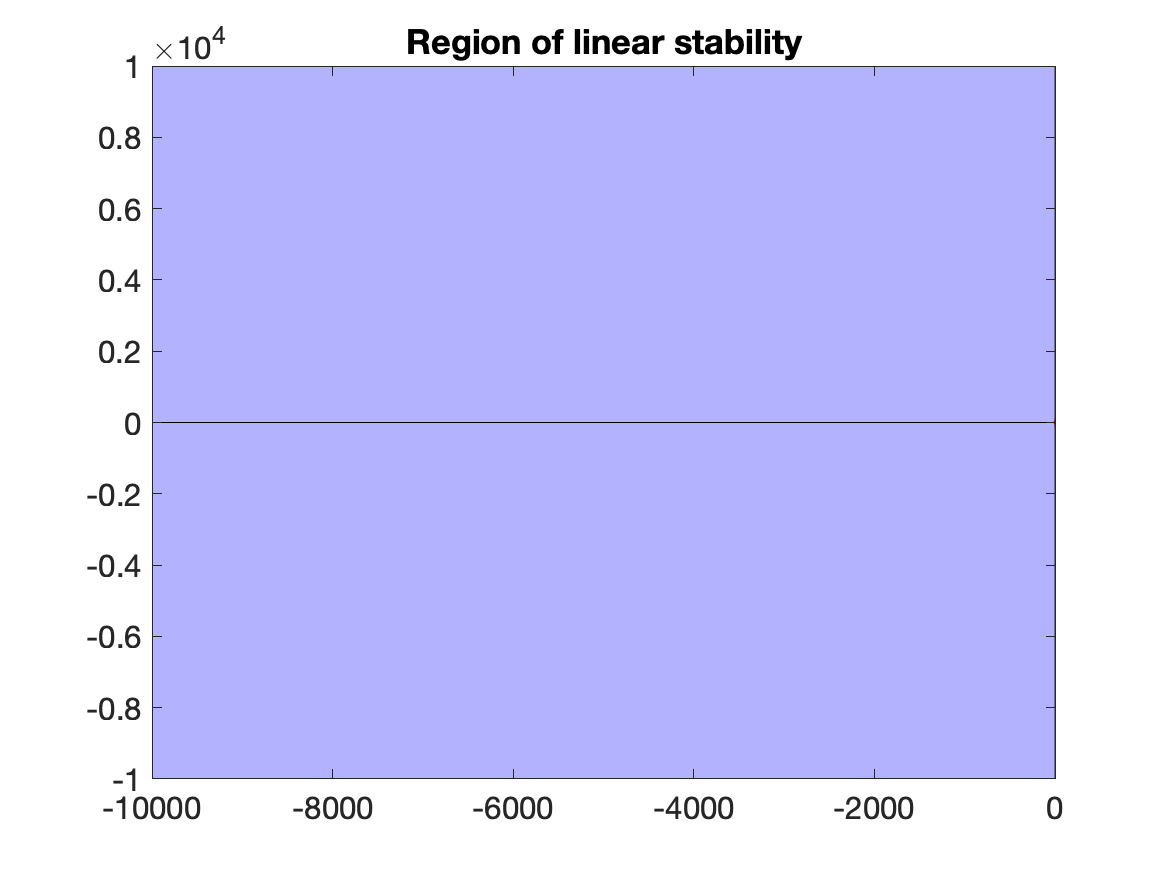}
      \includegraphics[width=0.325\textwidth]{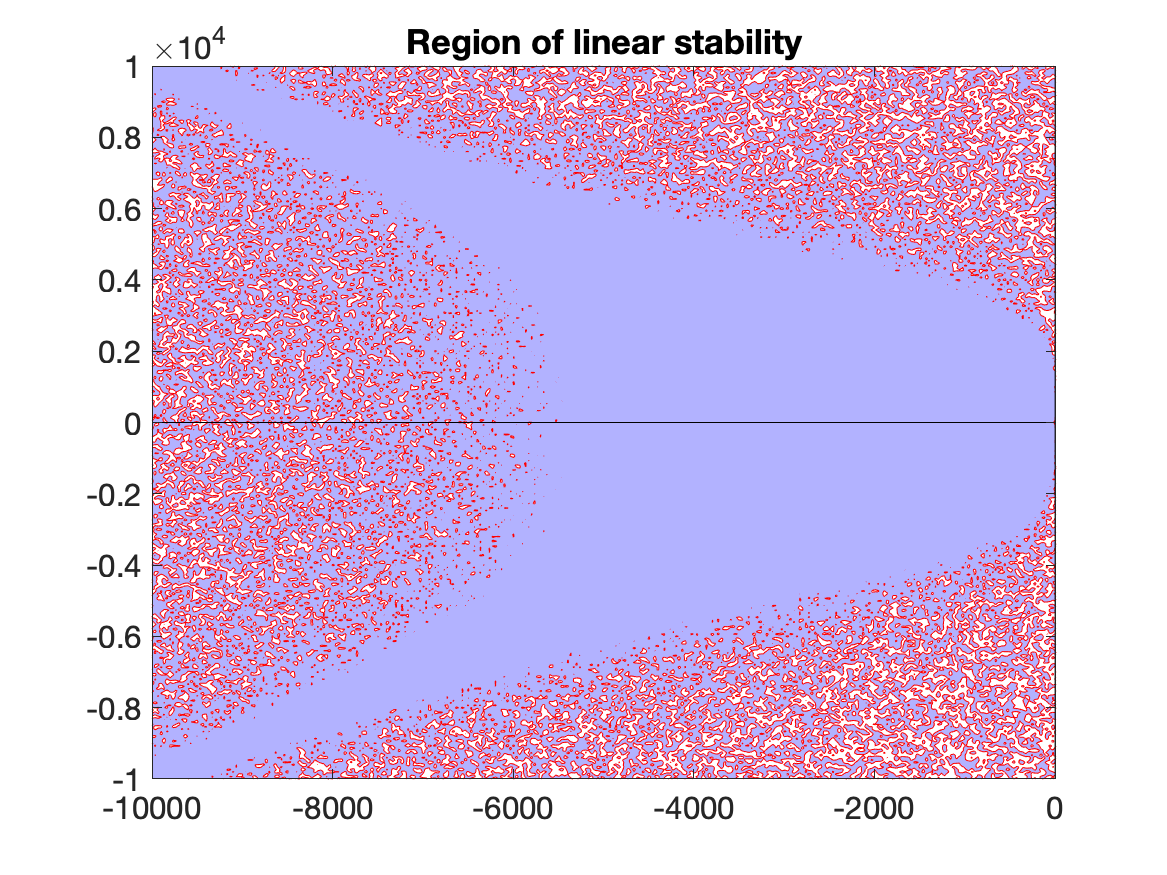} \\
  \includegraphics[width=0.325\textwidth]{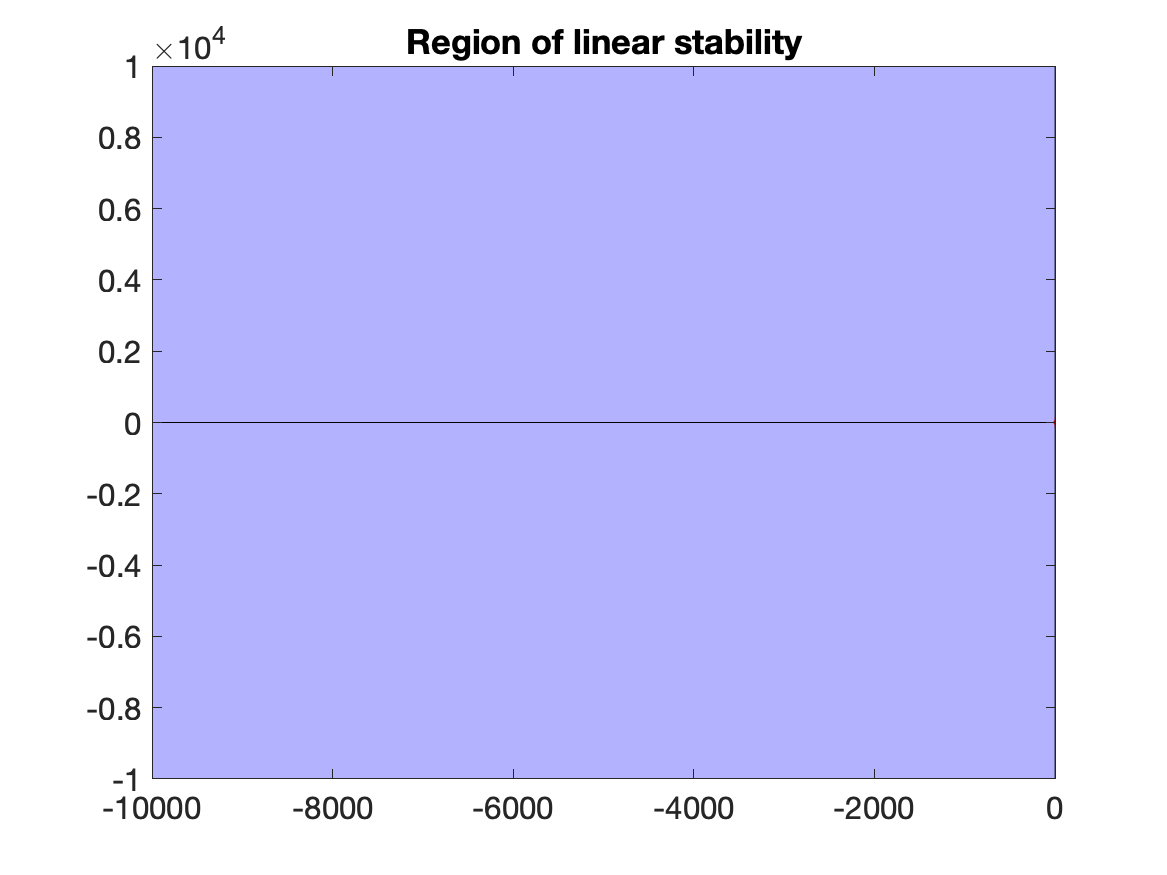}
    \includegraphics[width=0.325\textwidth]{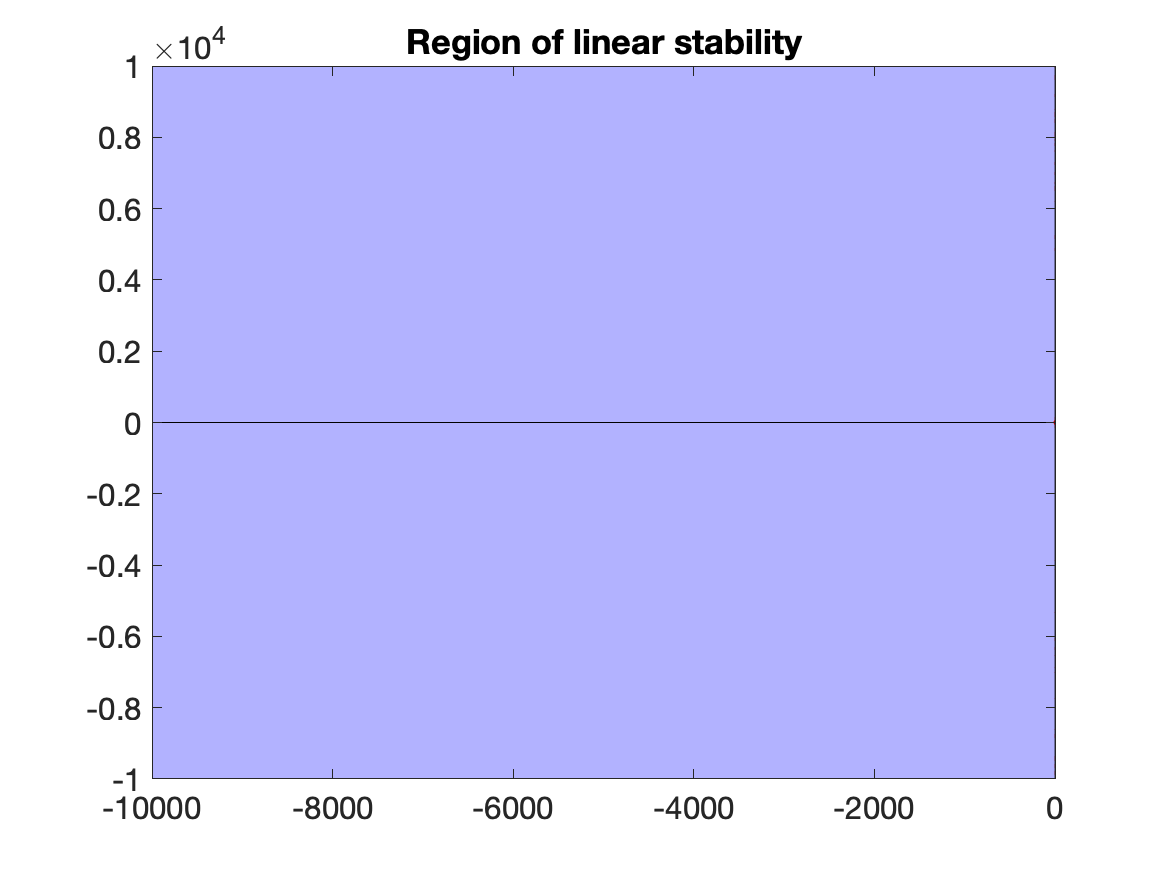}
      \includegraphics[width=0.325\textwidth]{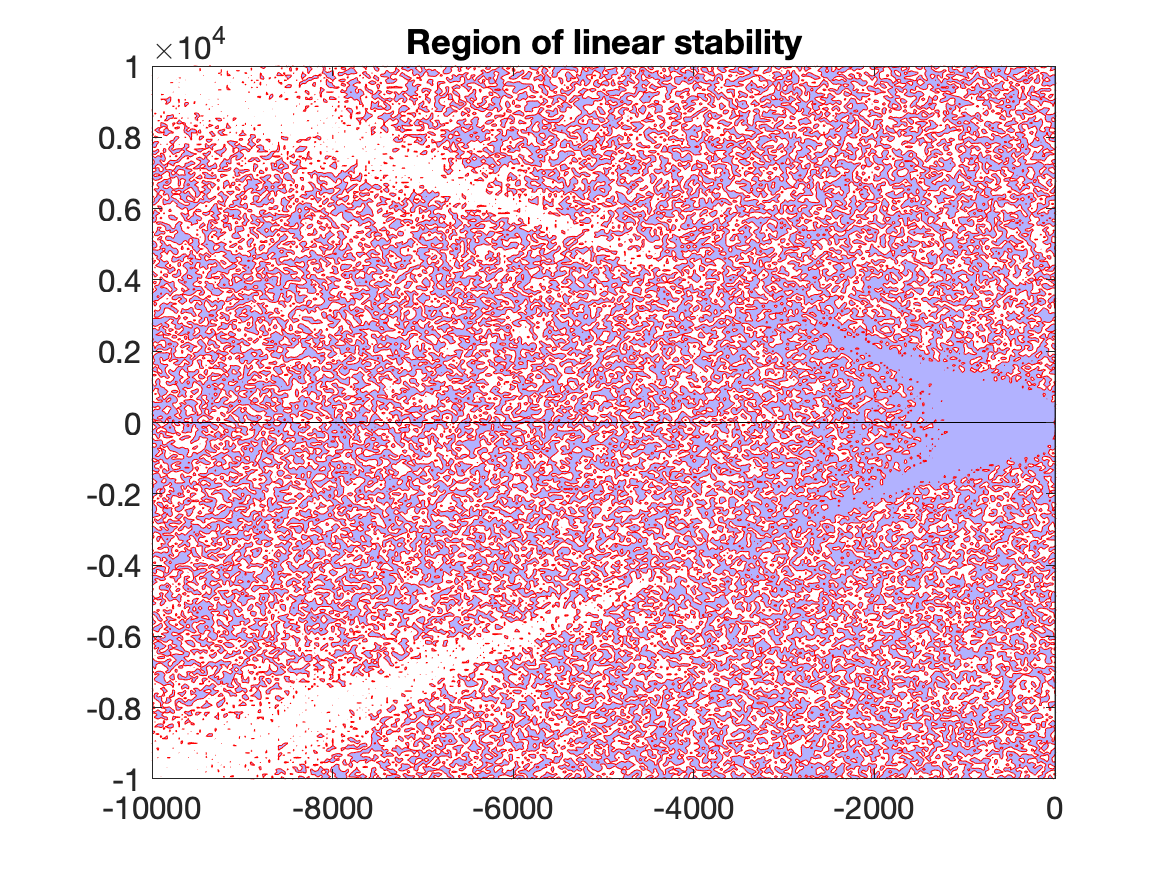} \\
        \includegraphics[width=0.325\textwidth]{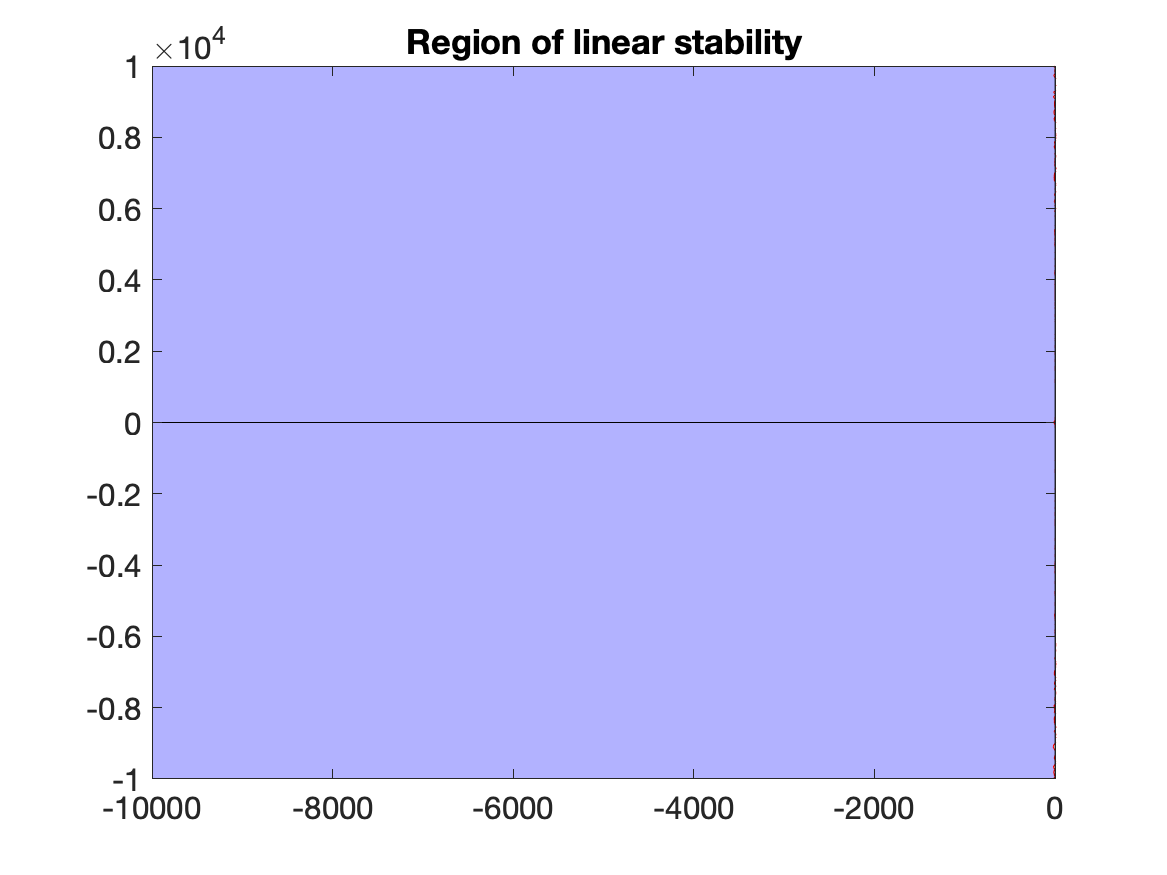}
    \includegraphics[width=0.325\textwidth]{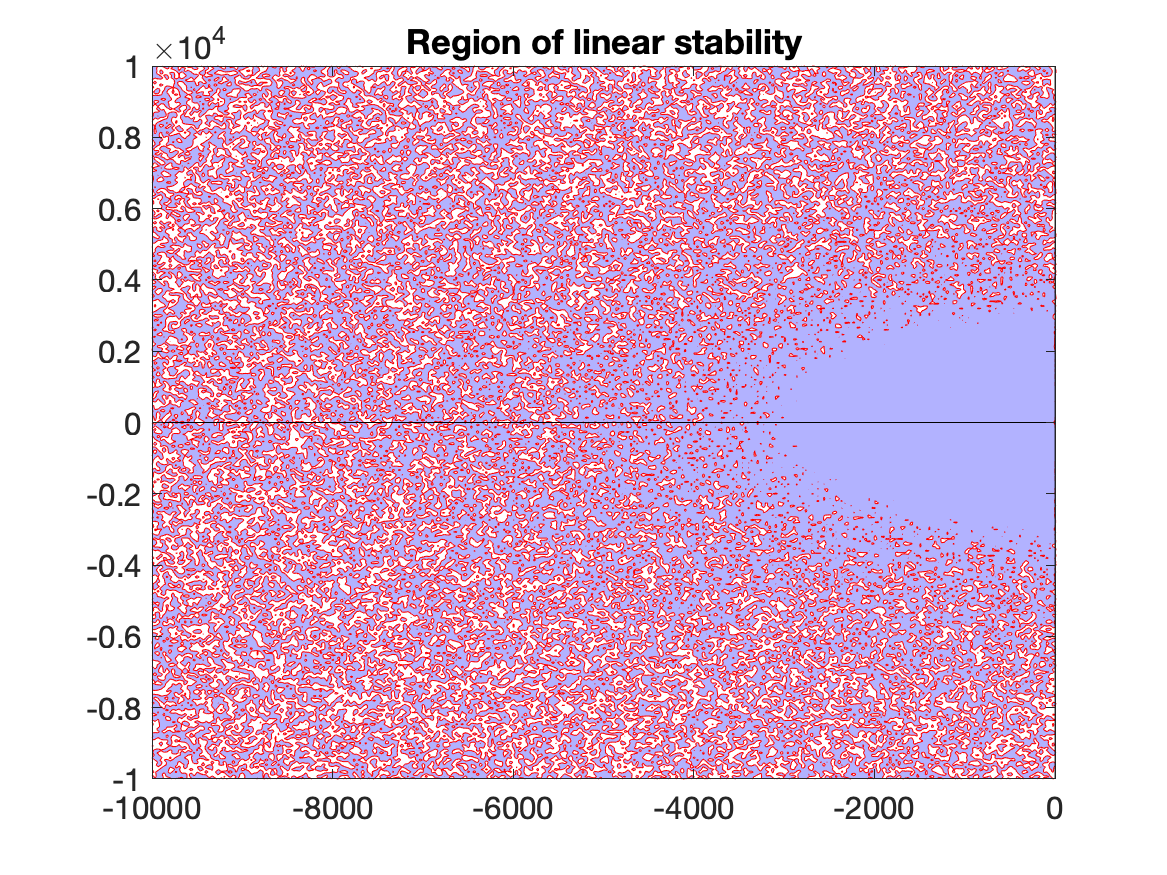}
      \includegraphics[width=0.325\textwidth]{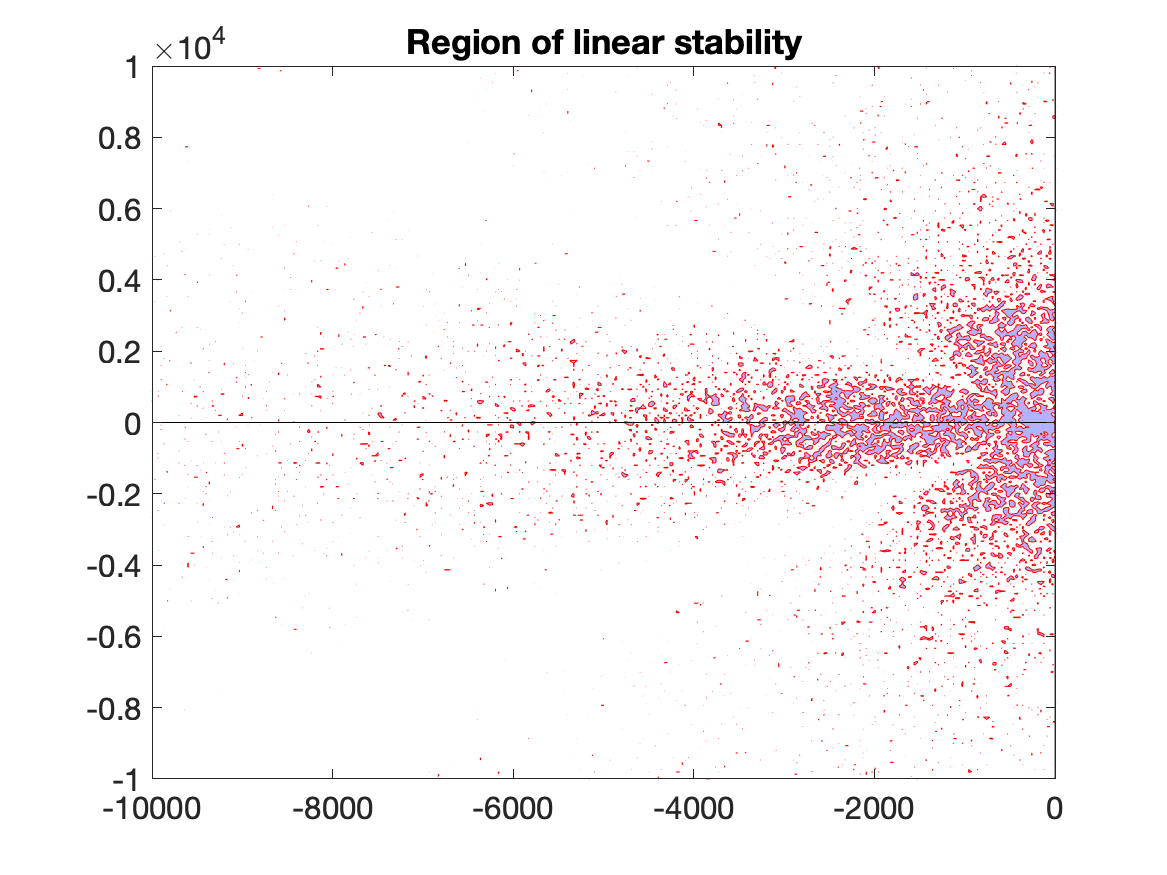} \\
              \includegraphics[width=0.325\textwidth]{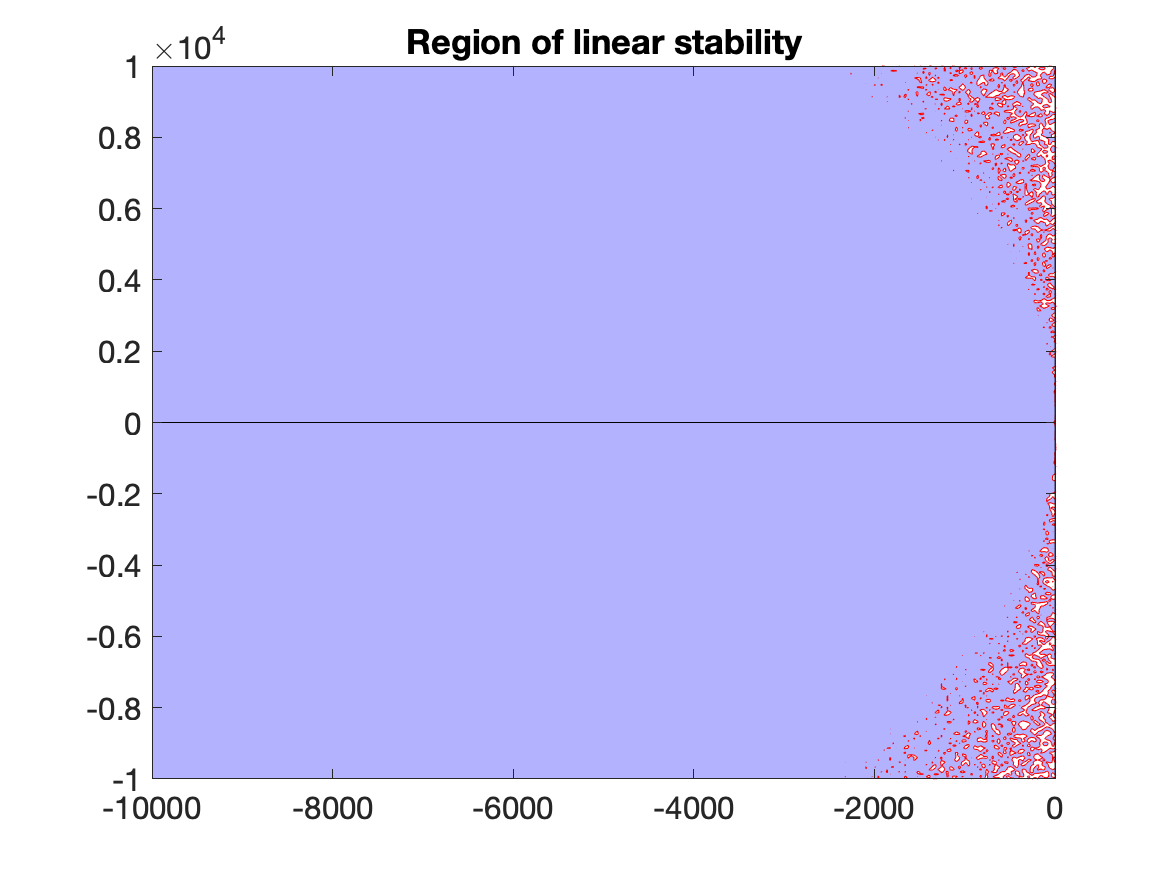}
    \includegraphics[width=0.325\textwidth]{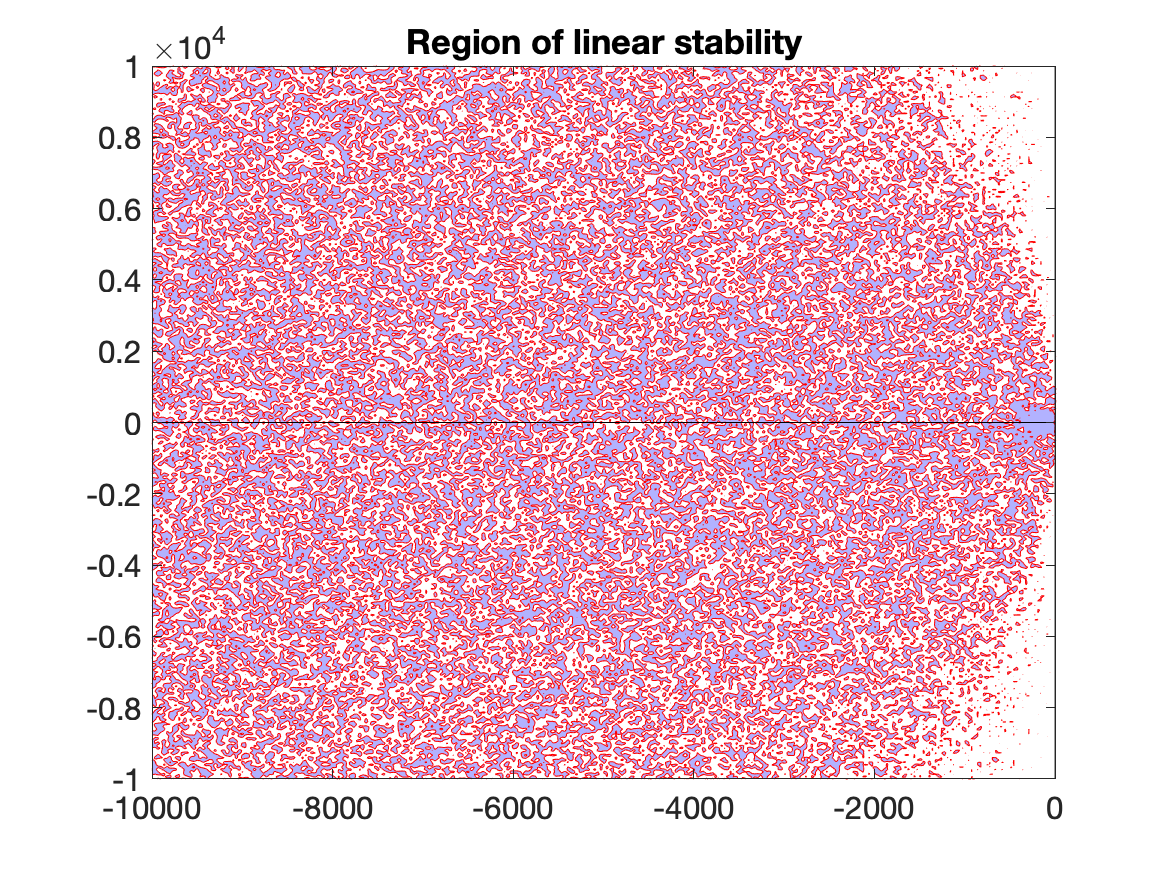}
      \includegraphics[width=0.325\textwidth]{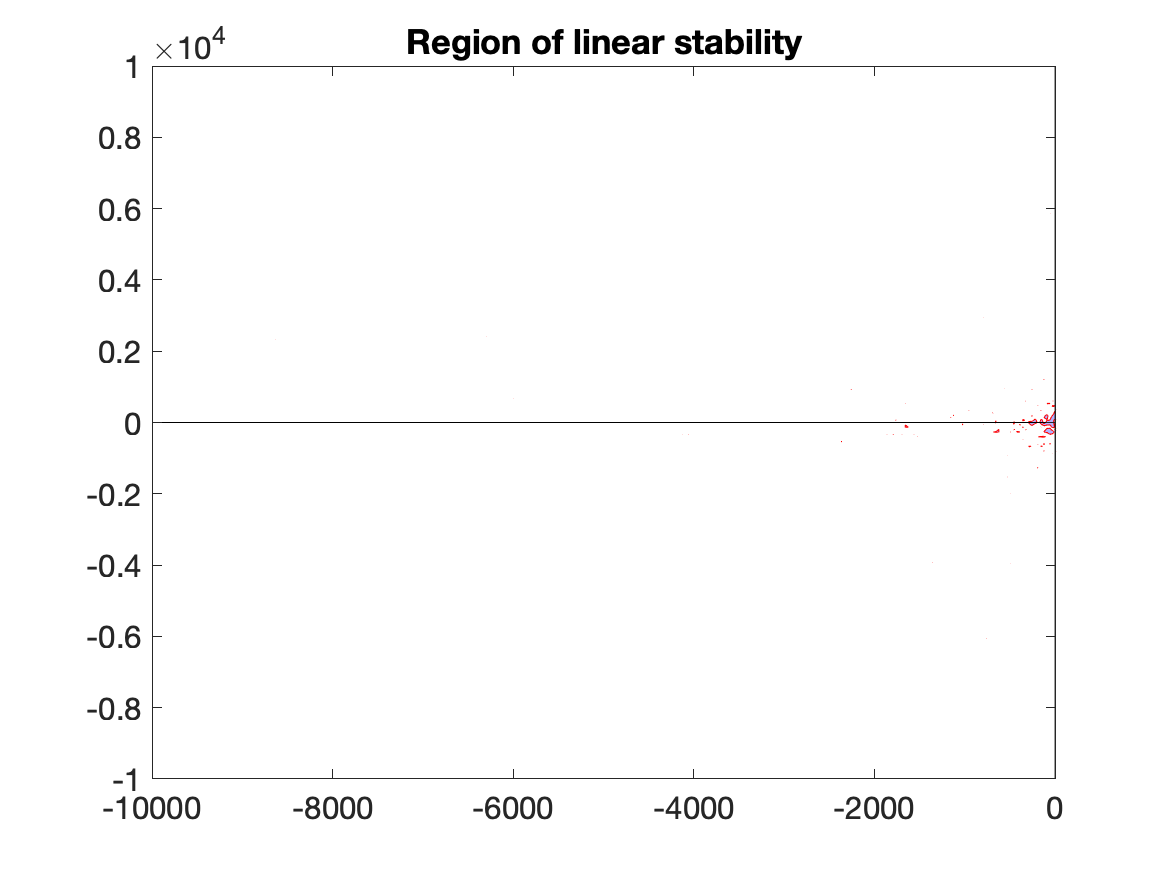} \\
  \caption{Regions of linear stability for the implicit midpoint rule with no corrections
(left), one correction (middle), and two corrections (right), for $\tilde{\epsilon} = 10^{-12}, 10^{-10}, 10^{-8}, 10^{-6}, 10^{-4}$
(top to bottom).}
  \label{IMRstability}
\end{figure*}

Figure \ref{IMRstability} shows the linear stability regions for the implicit midpoint rule with no corrections
(left), one correction (middle), and two corrections (right), for  (from top to bottom) 
$\tilde{\epsilon} = 10^{-12}, 10^{-10}, 10^{-8}, 10^{-6}, 10^{-4}$.
The blue areas are stable, white areas are unstable, and the red marks the contours of these areas.
Looking at the graphs from left to right, we clearly observe that 
for large enough $\tilde{\epsilon}$, corrections greatly reduce the region of stability, i.e.
the values of $z= \dt \lambda$ must be smaller to maintain stability. This means that for
a fixed $\tilde{\epsilon}$ that is large enough, more corrections require a smaller $\dt$ as $\lambda$ gets larger.

Unlike in the usual stability analysis where $z = \lambda \dt$ is the only quantity of interest, in these 
figures $\epsilon$ and $\lambda$ each play a role as well.
In practice, this means that the regions of instability  shown in the figures 
are not as bad as it appears at first glance:  recalling that $\tilde{\epsilon} =\frac{\epsilon}{\lambda} $, 
for a stiff problem  the value of $\tilde{\epsilon} $ is scaled accordingly.
A region of stability for  $\tilde{\epsilon} =10^{-12}$ is relevant for a 
 $\epsilon =10^{-8}$ where the problem has a stiffness of $\lambda =10^{4}$,
but a smaller  $\epsilon =10^{-12}$ where the problem has a lesser stiffness of $\lambda =1$.
So although stiffer problems generally require a smaller time-step for stability, we must look at the 
correct image that corresponds to this, and this is scaled in the opposite direction.

\begin{figure*}[ht!]
  \includegraphics[width=0.325\textwidth]{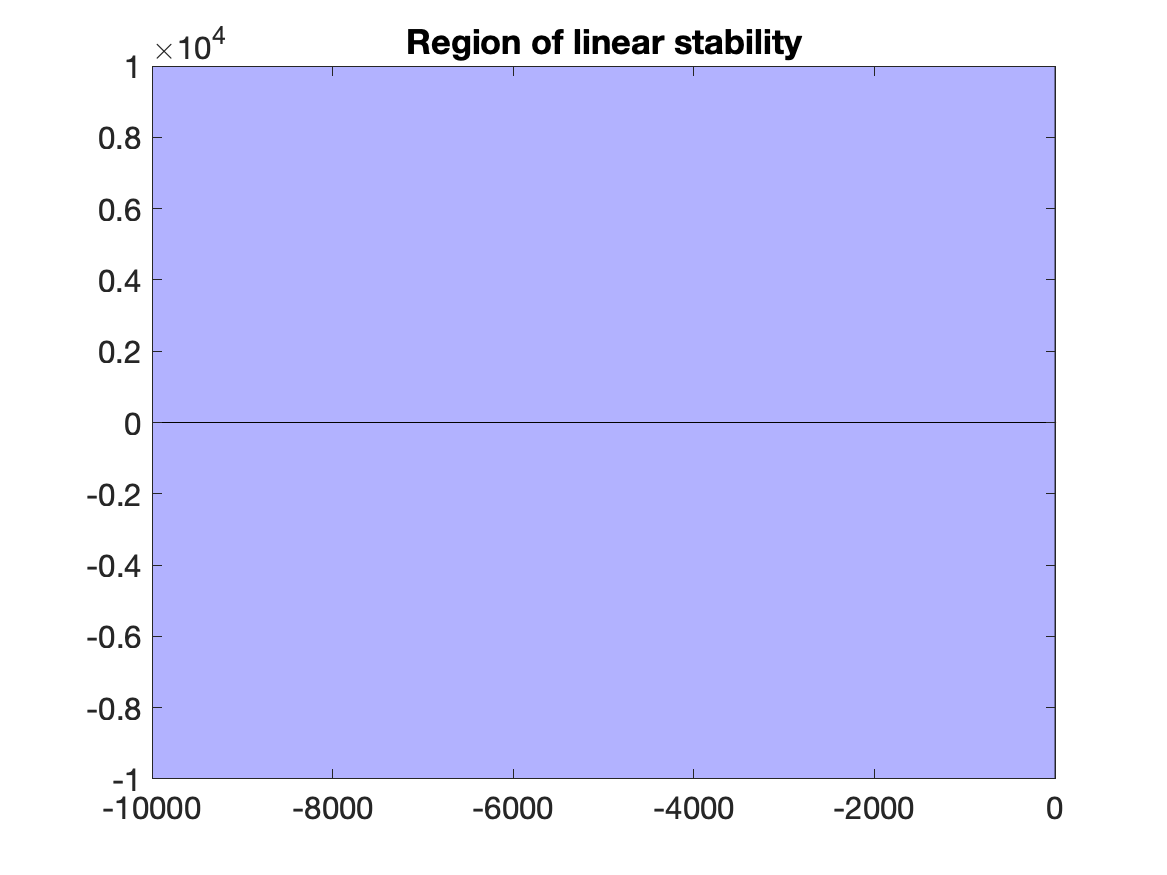}
    \includegraphics[width=0.325\textwidth]{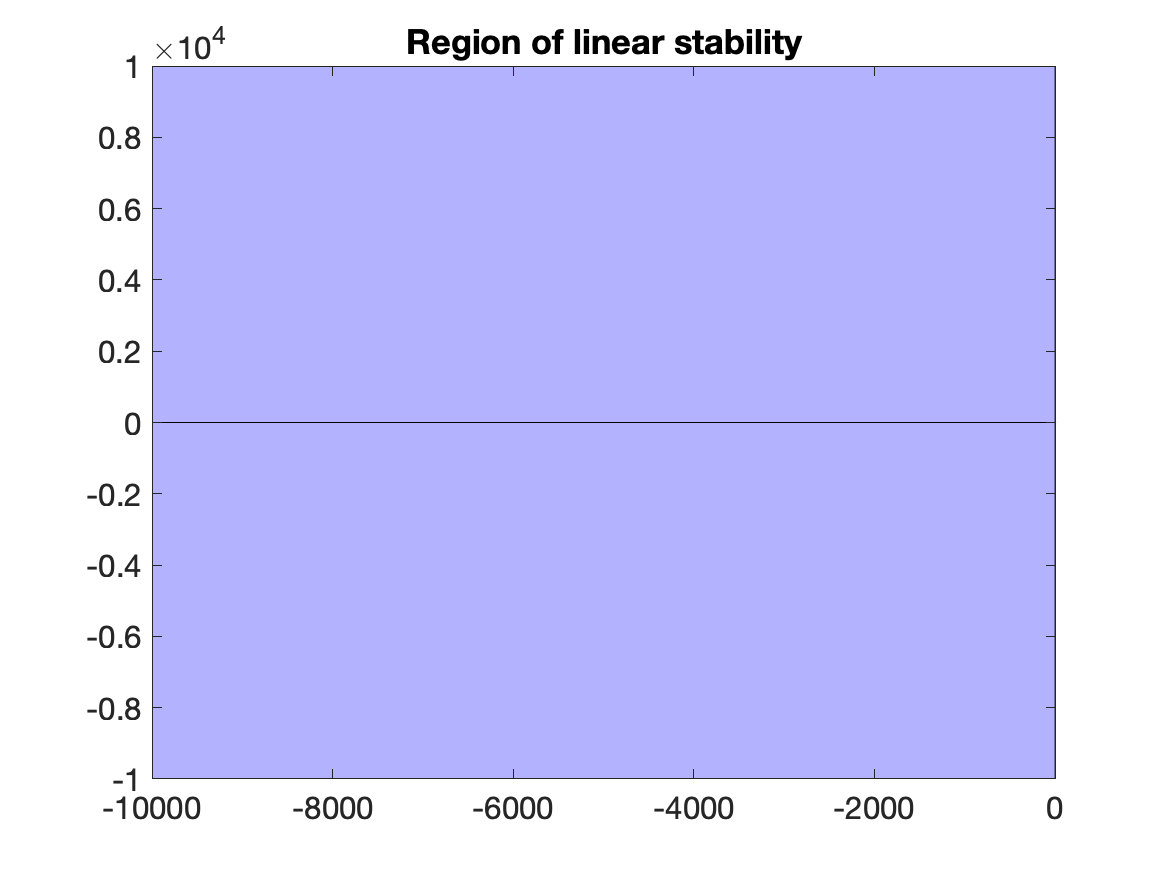}
      \includegraphics[width=0.325\textwidth]{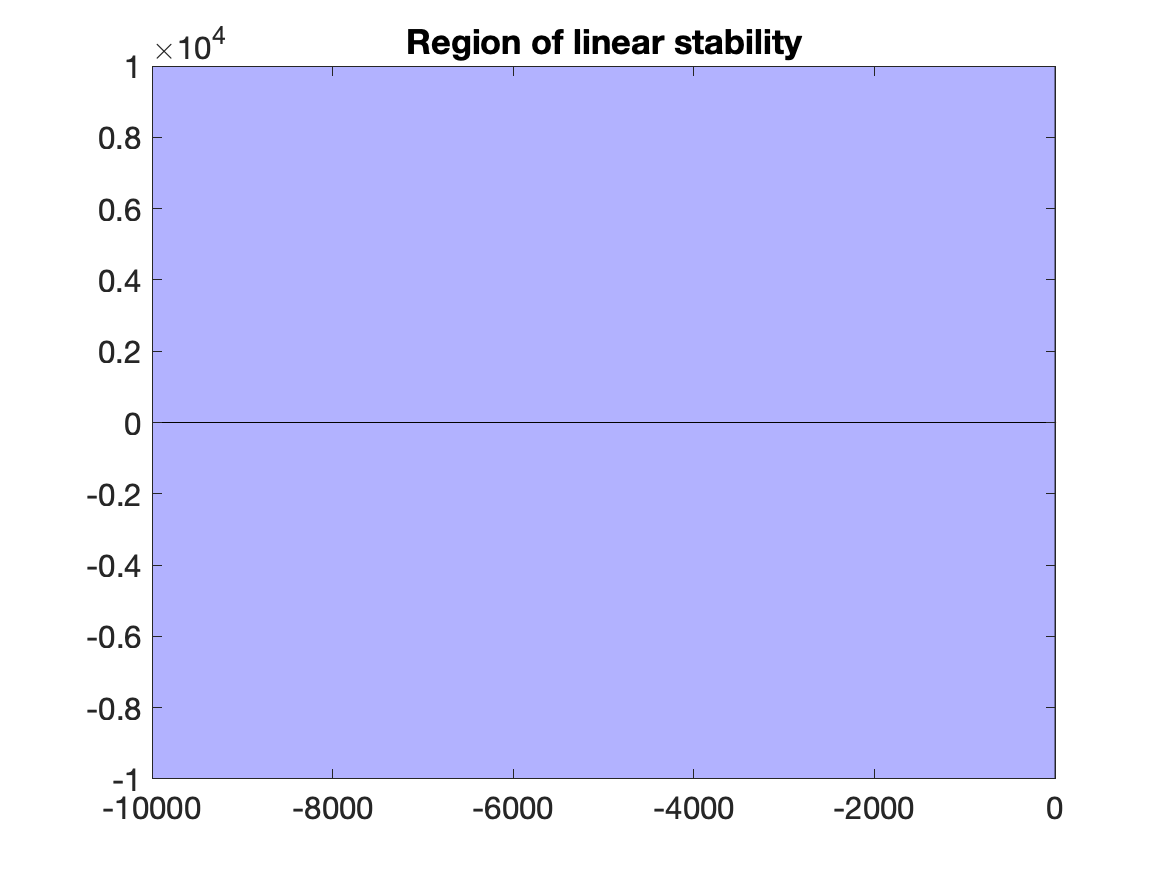} \\
  \includegraphics[width=0.325\textwidth]{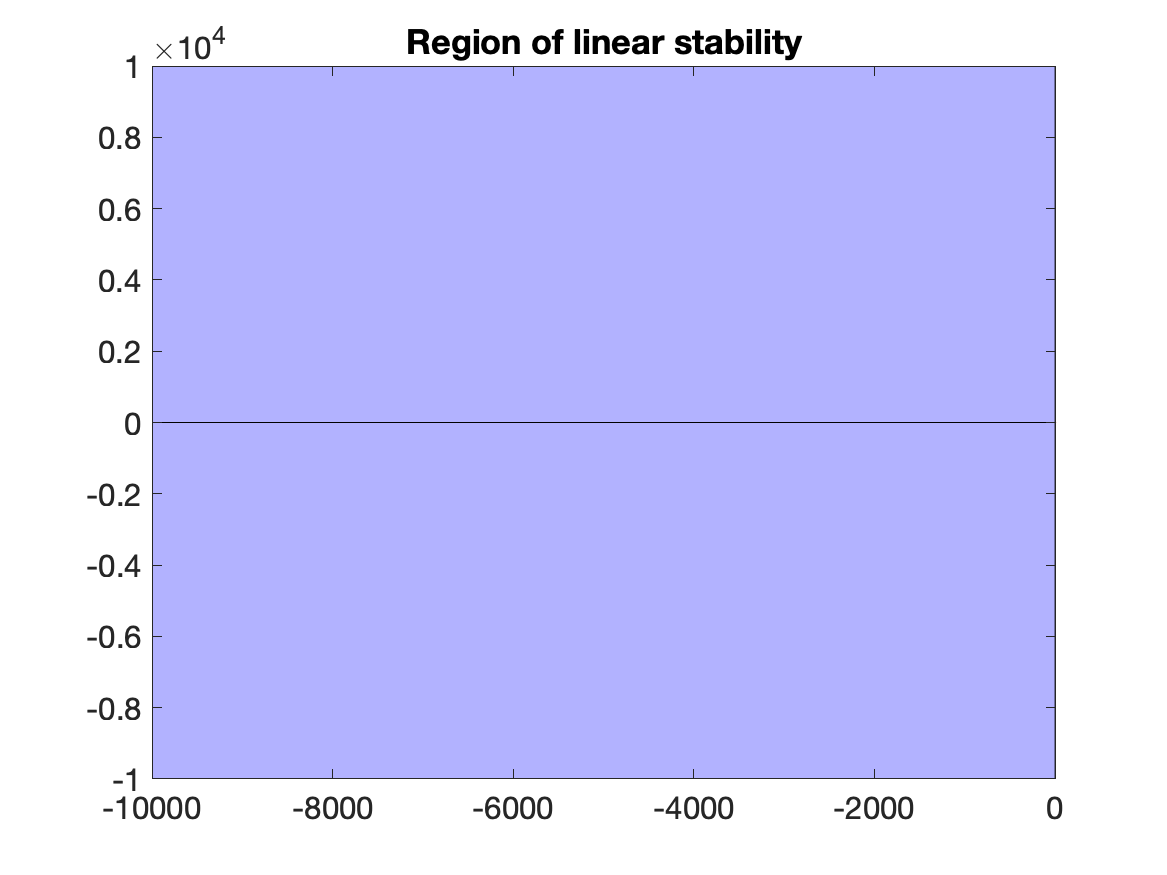}
    \includegraphics[width=0.325\textwidth]{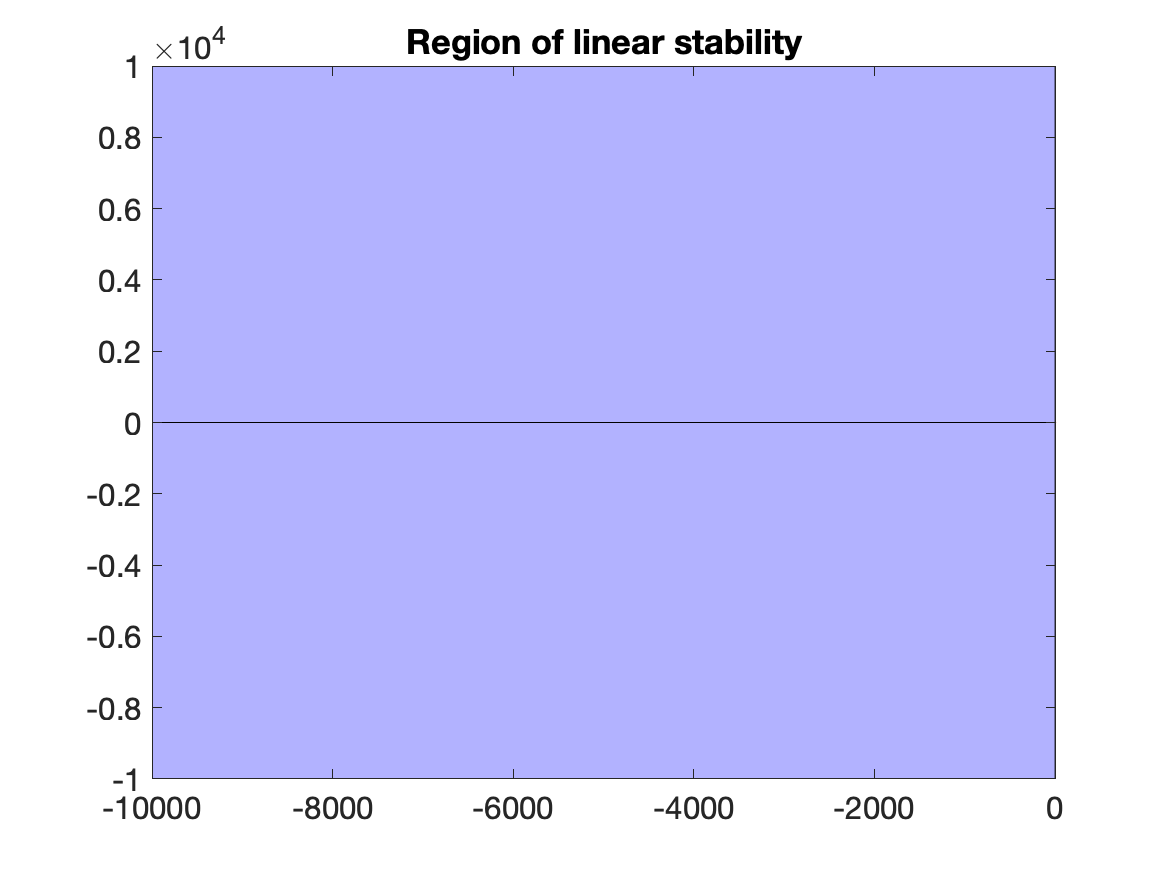}
      \includegraphics[width=0.325\textwidth]{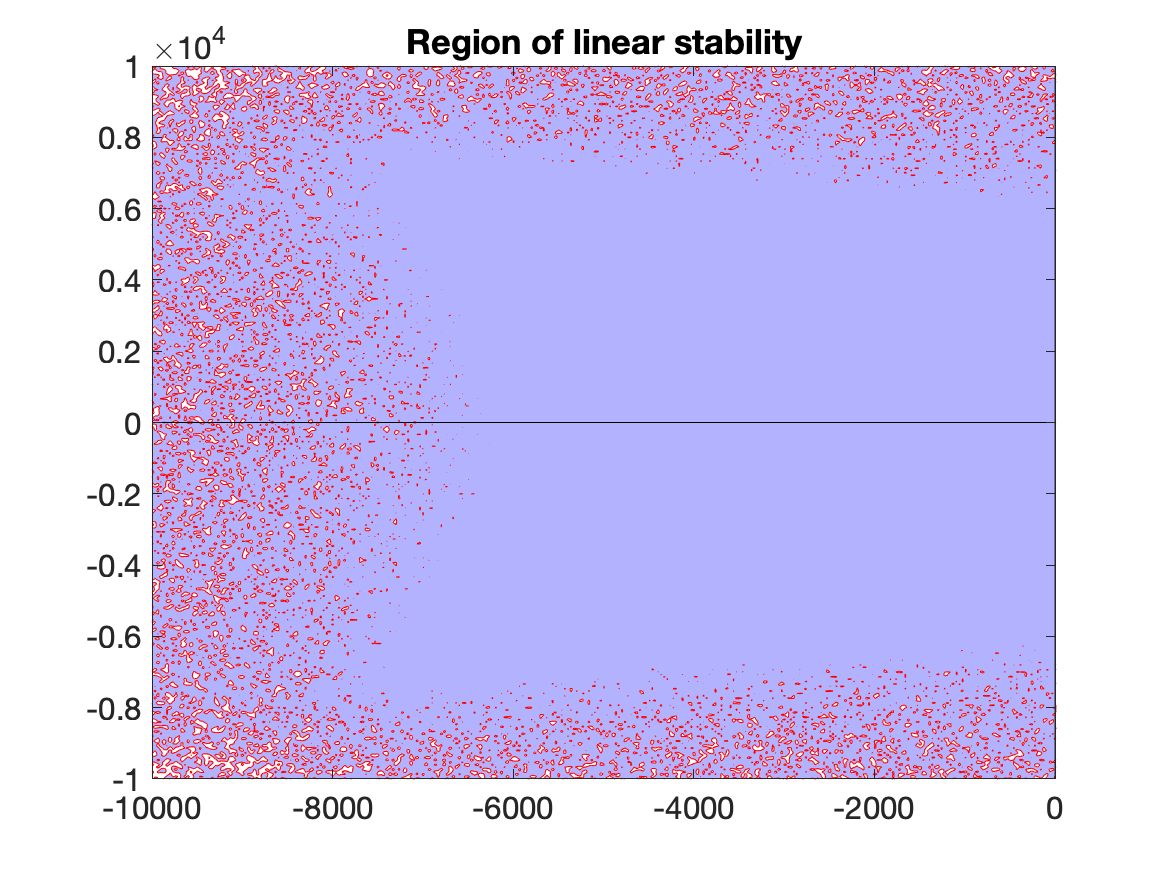} \\
        \includegraphics[width=0.325\textwidth]{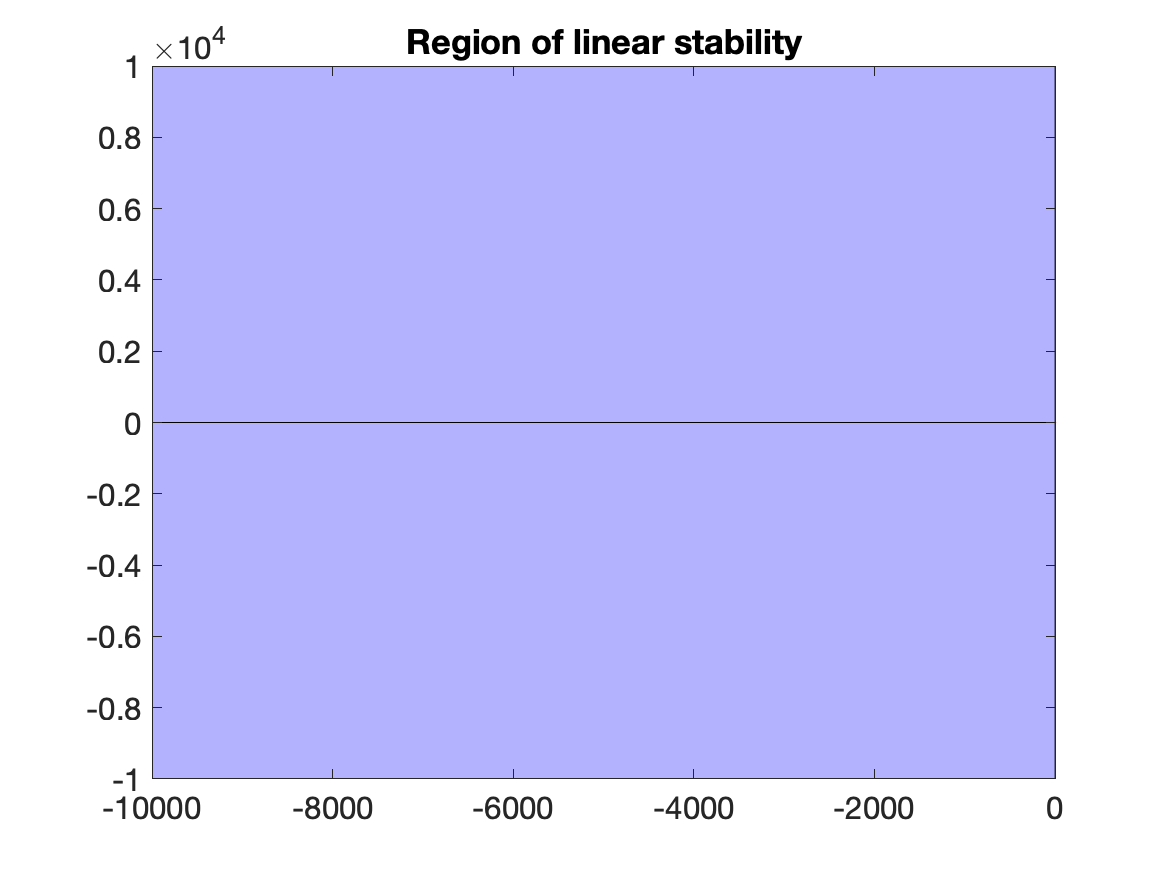}
    \includegraphics[width=0.325\textwidth]{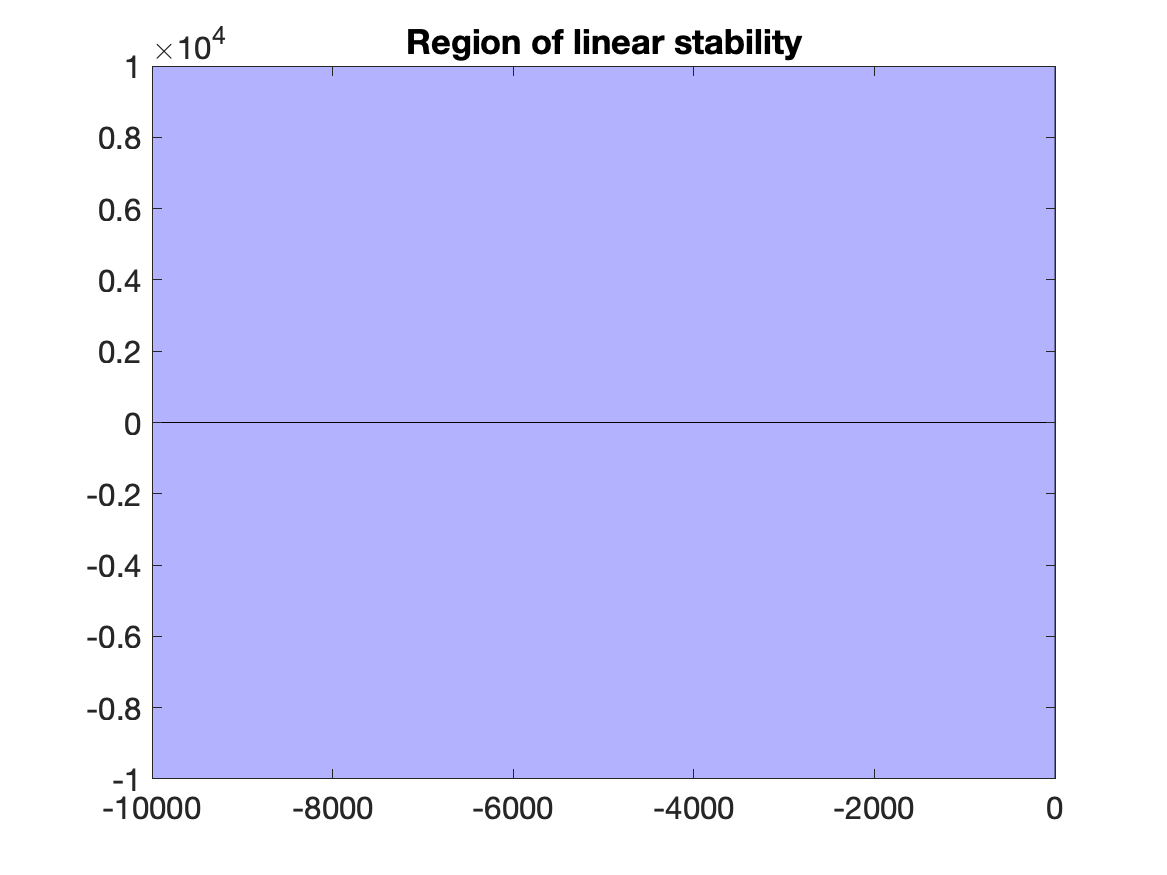}
      \includegraphics[width=0.325\textwidth]{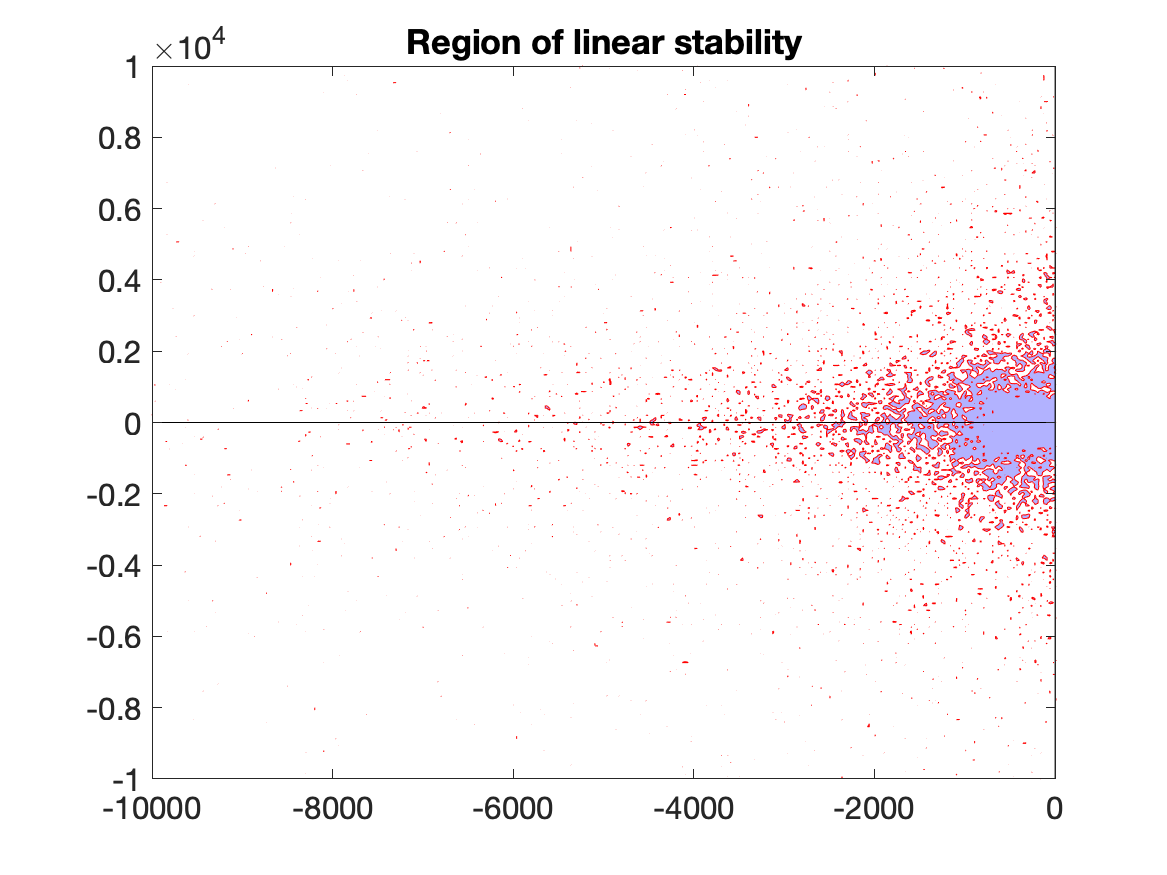} \\
              \includegraphics[width=0.325\textwidth]{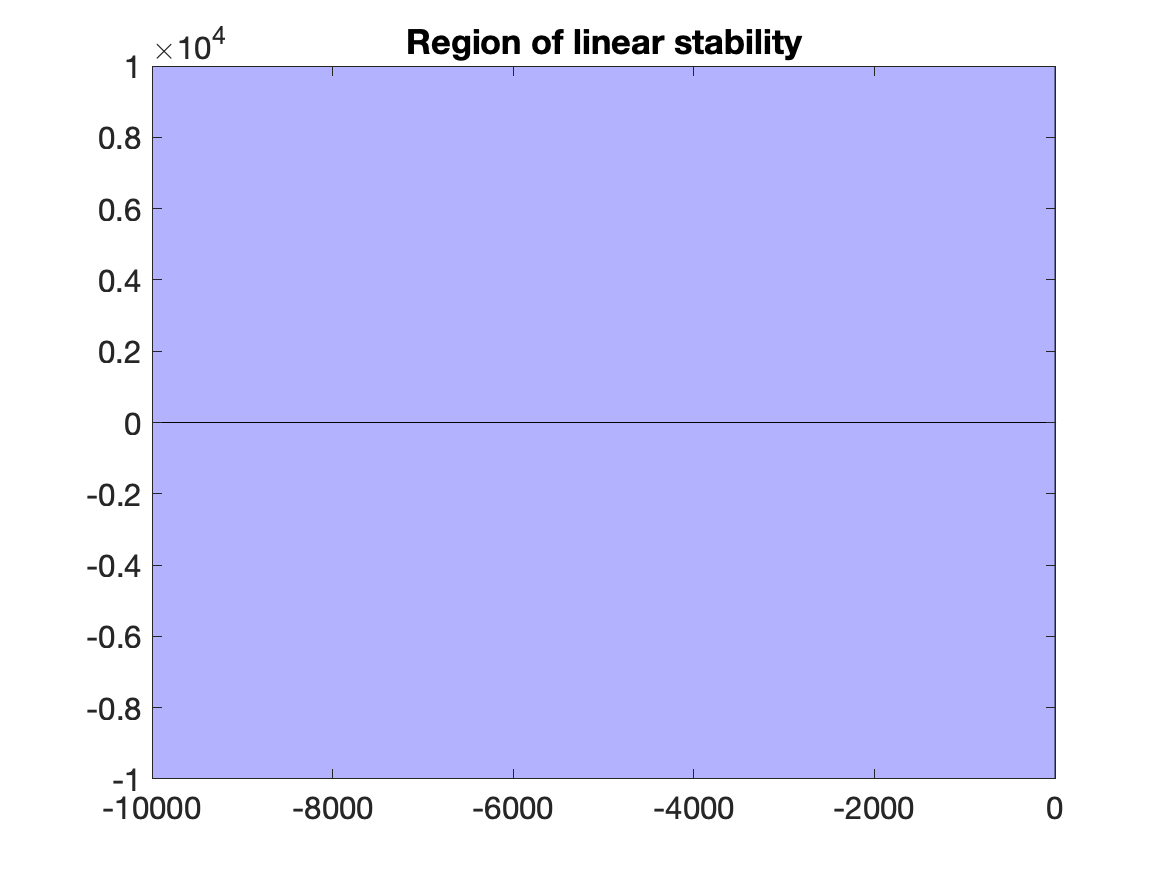}
    \includegraphics[width=0.325\textwidth]{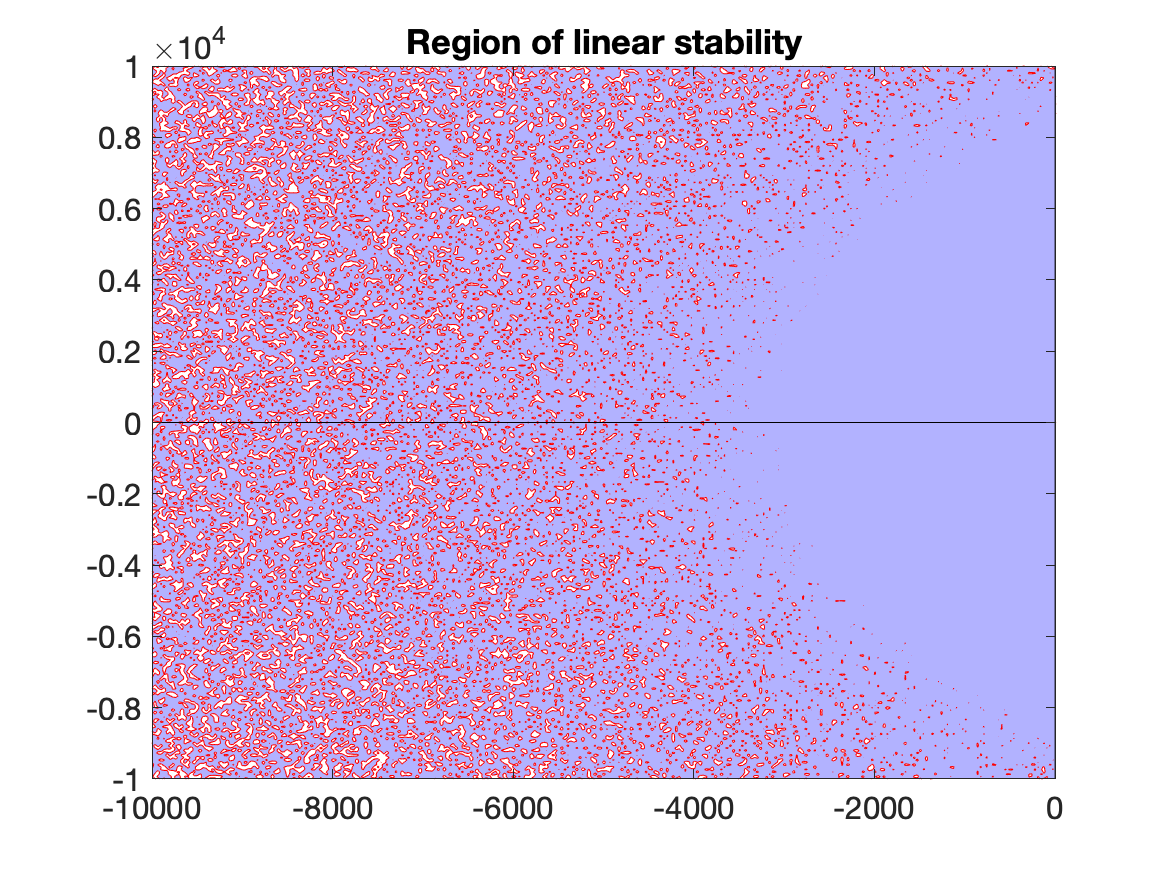}
      \includegraphics[width=0.325\textwidth]{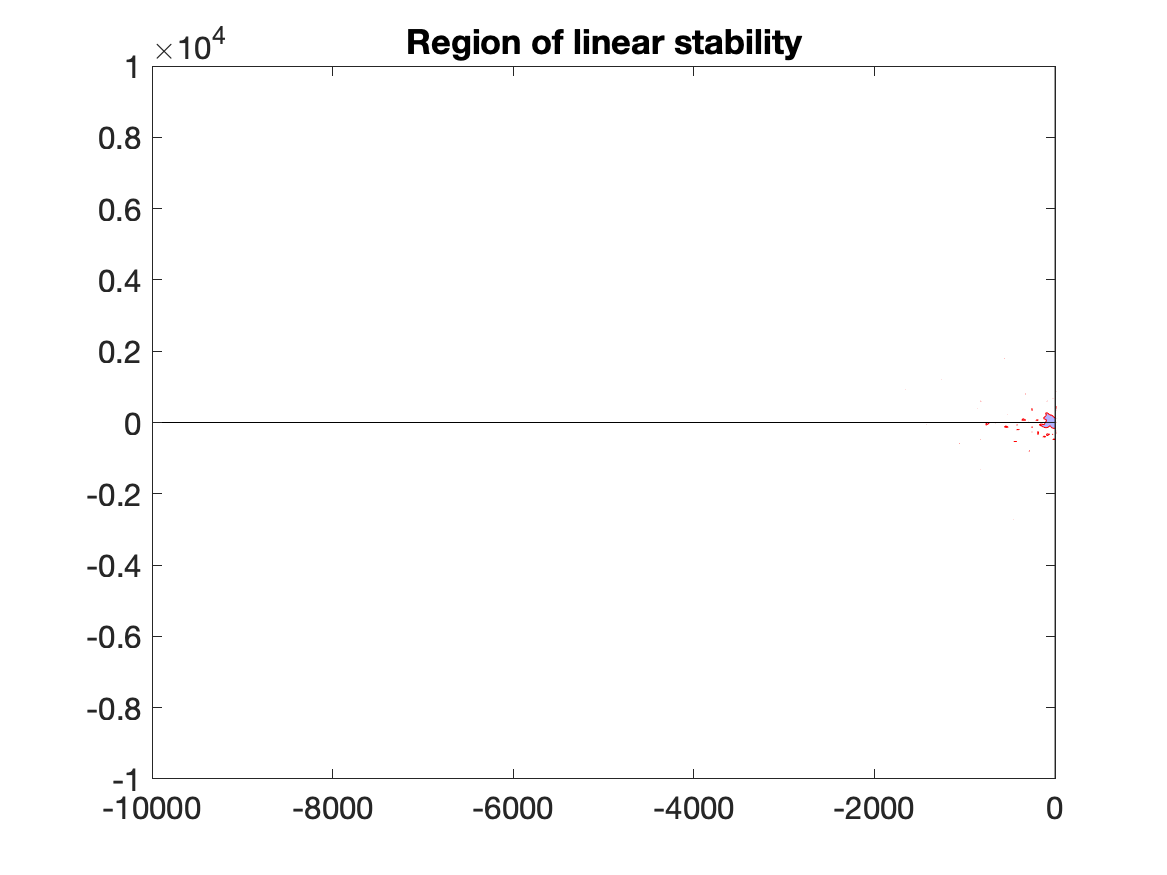} \\
  \caption{Regions of linear stability for the SDIRK with no corrections
(left), one correction (middle), and two corrections (right), for $\tilde{\epsilon} = 10^{-12}, 10^{-8}, 10^{-6}, 10^{-4}$
(top to bottom).}
  \label{SDIRKstability}
\end{figure*}

The linear stability regions for  SDIRK methods have a similar behavior profile, as seen in 
Figure \ref{SDIRKstability}. Here we show the regions of linear stability for the SDIRK with no corrections
(left), one correction (middle), and two corrections (right), for (top to bottom) $\tilde{\epsilon} = 10^{-12}, 10^{-8}, 10^{-6}, 10^{-4}$.
Once again we observe that larger $\tilde{\epsilon}$ and more corrections result in a smaller linear stability region.
The conclusion is that we want to reserve the use of corrections to the case where $\dt$ is small enough
compared to $\epsilon$: this is the regime in which we benefit from corrections in terms of accuracy,
and in which the method is least likely to become unstable. Because $\tilde{\epsilon} = \frac{\epsilon}{\lambda}$, 
the stiffer the problem, the larger is the low precision $\epsilon$ for which stability is ensured. 
This  in turn means that larger $\dt$ will benefit from corrections.

\begin{figure*}[ht!]
      \includegraphics[width=0.325\textwidth]{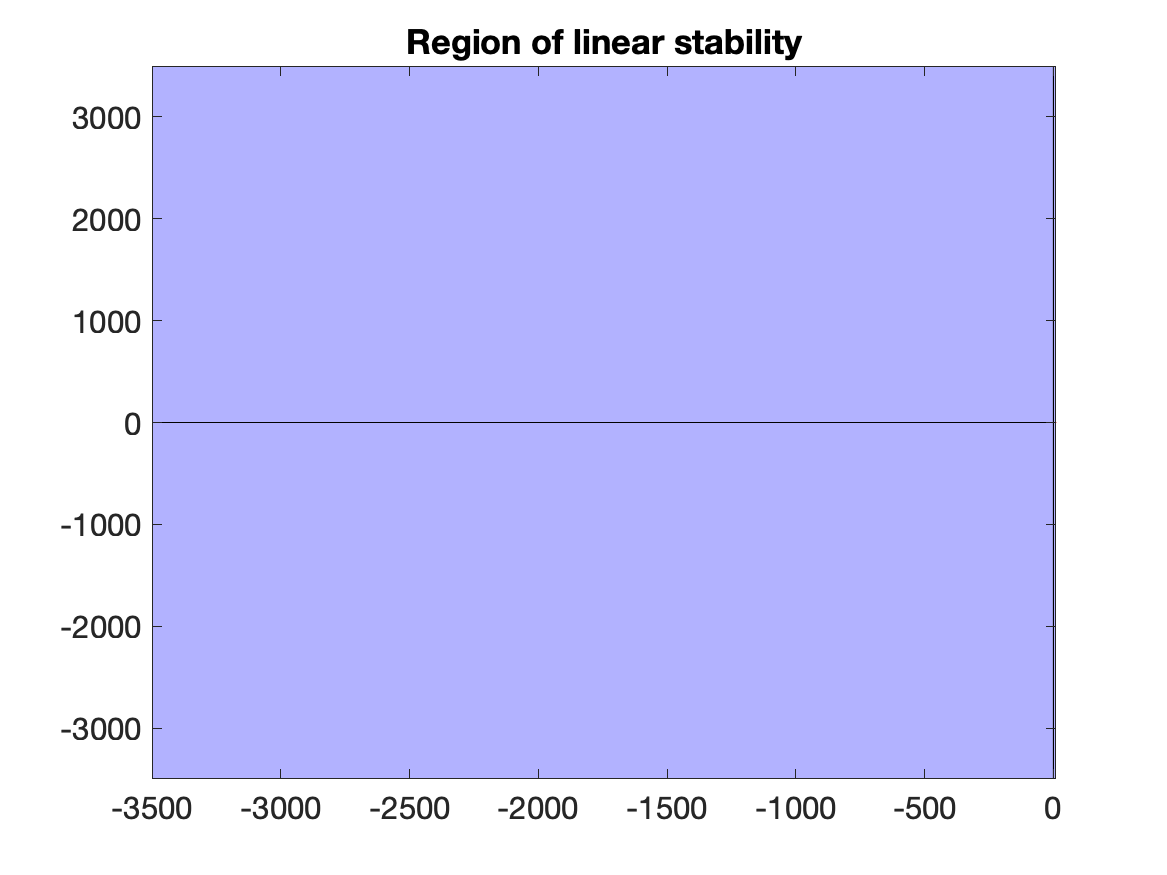}
       \includegraphics[width=0.325\textwidth]{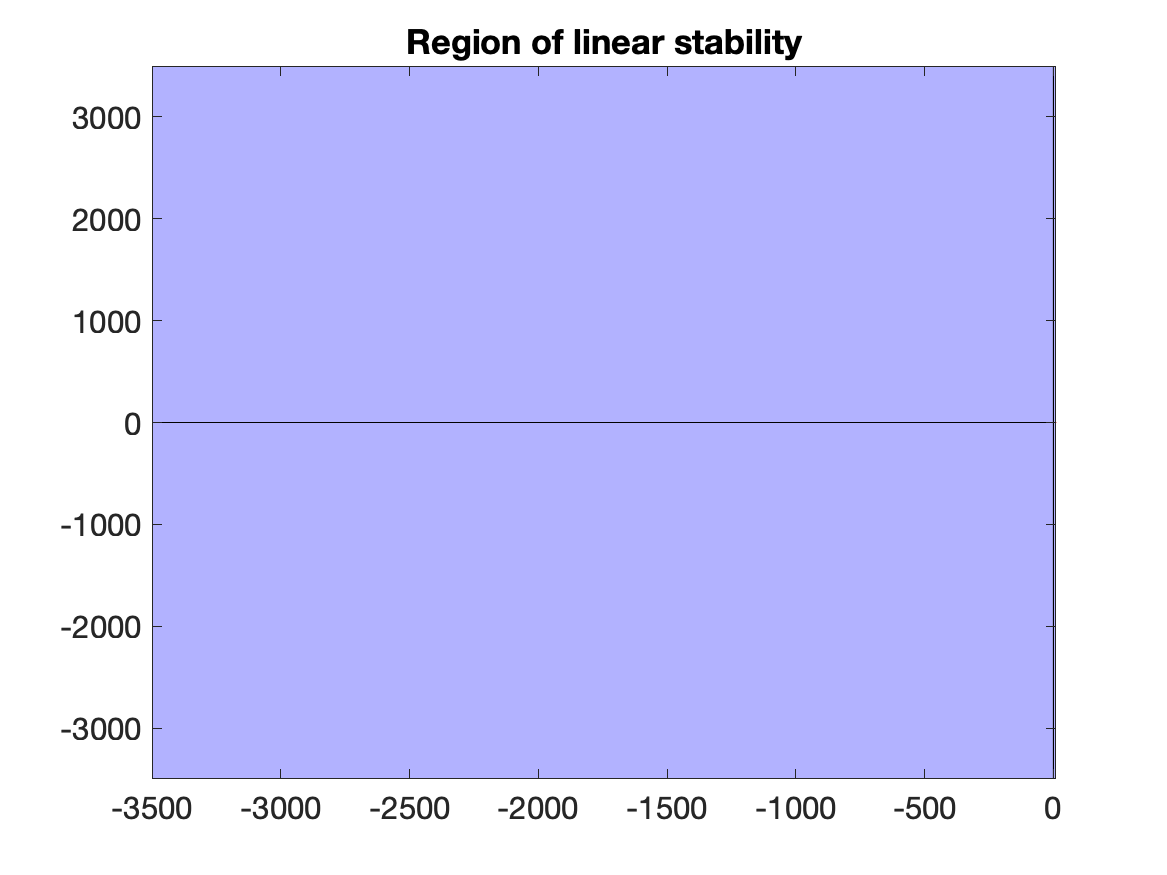}
        \includegraphics[width=0.325\textwidth]{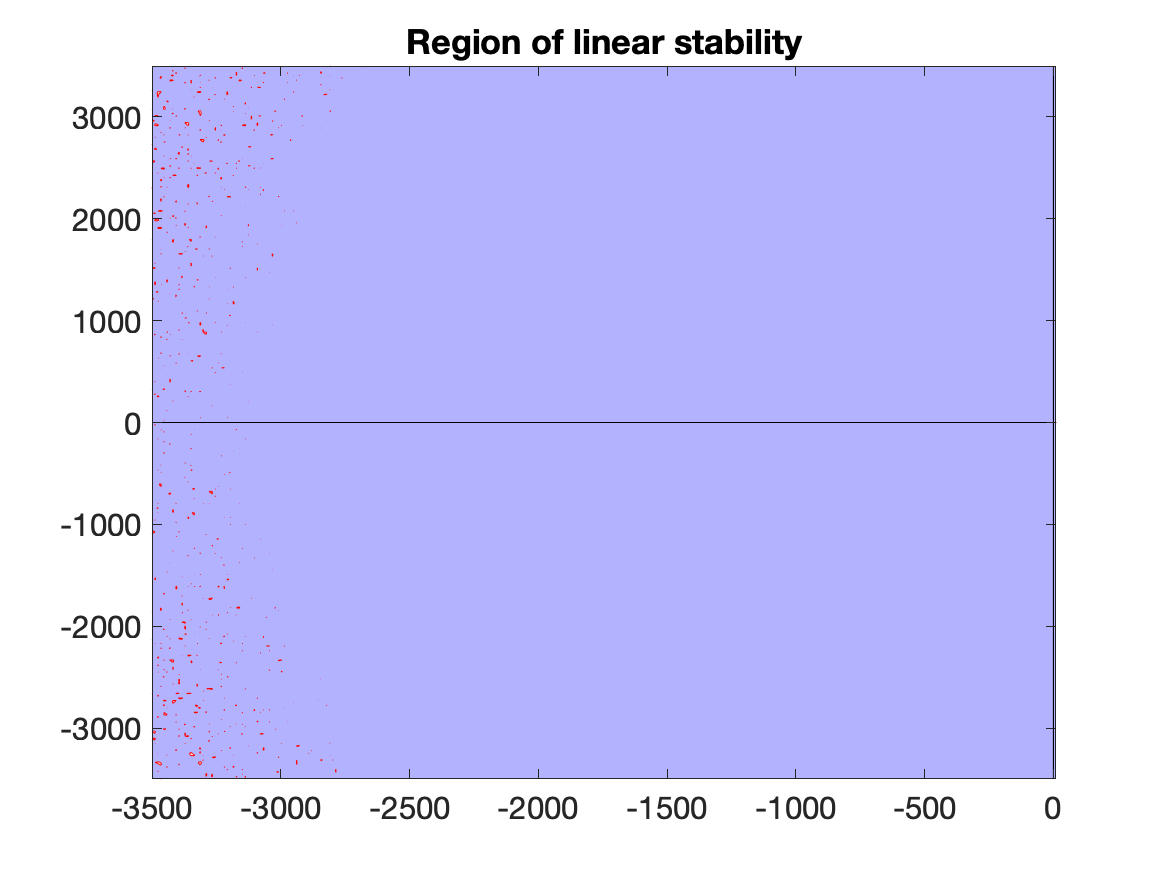} \\
      \includegraphics[width=0.325\textwidth]{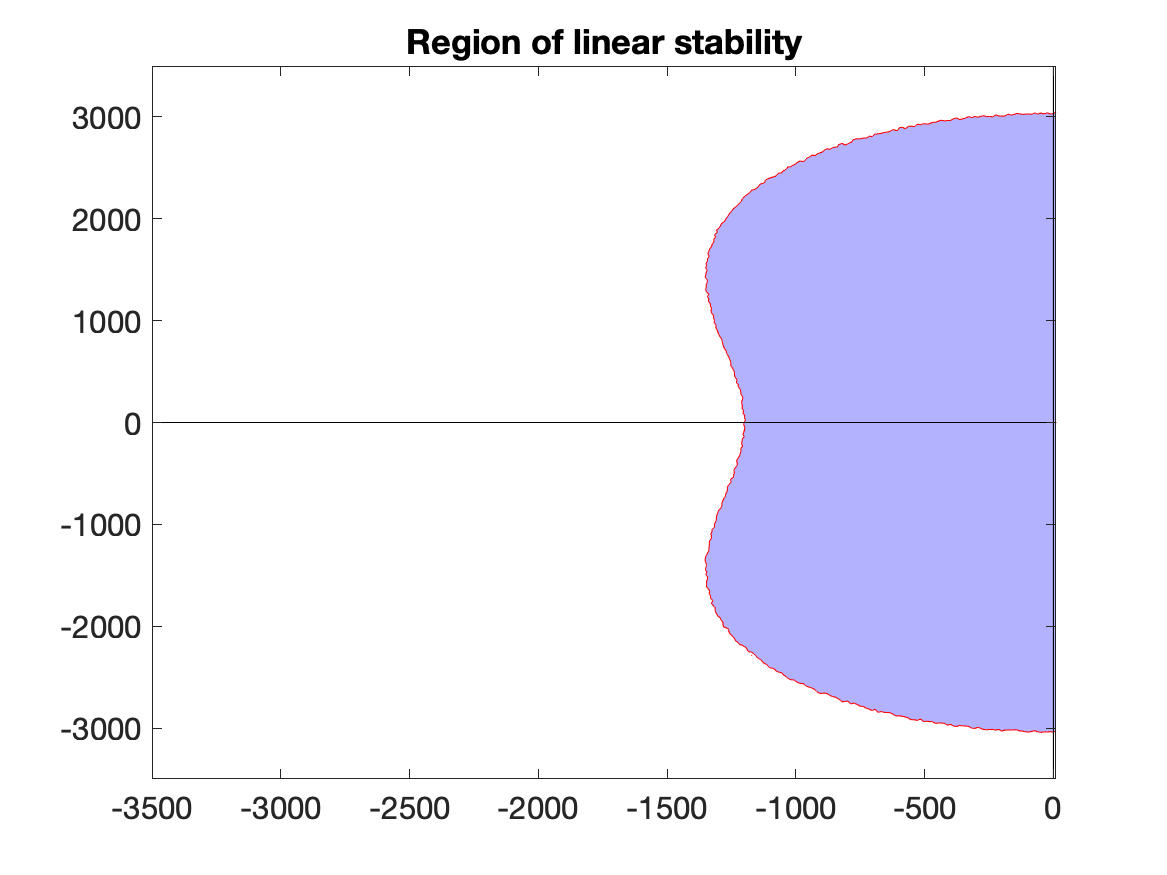}
       \includegraphics[width=0.325\textwidth]{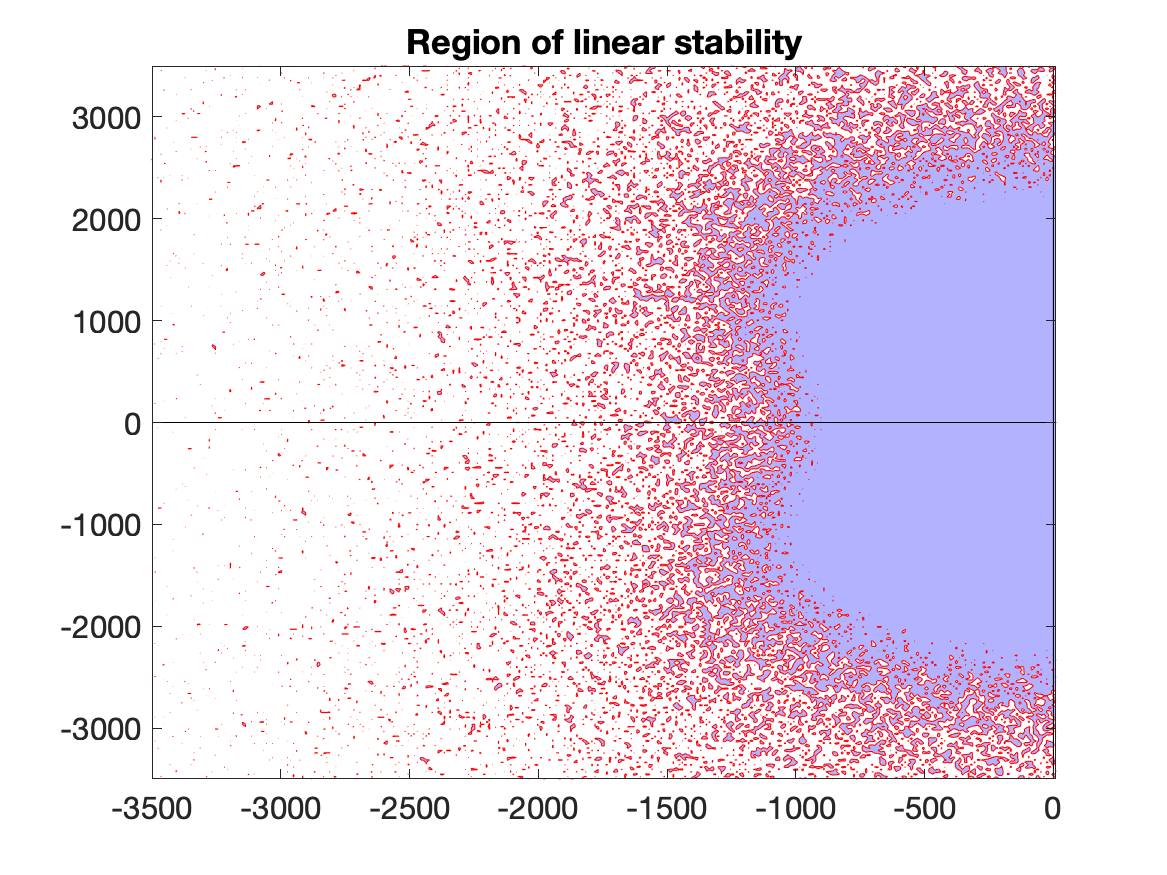}
        \includegraphics[width=0.325\textwidth]{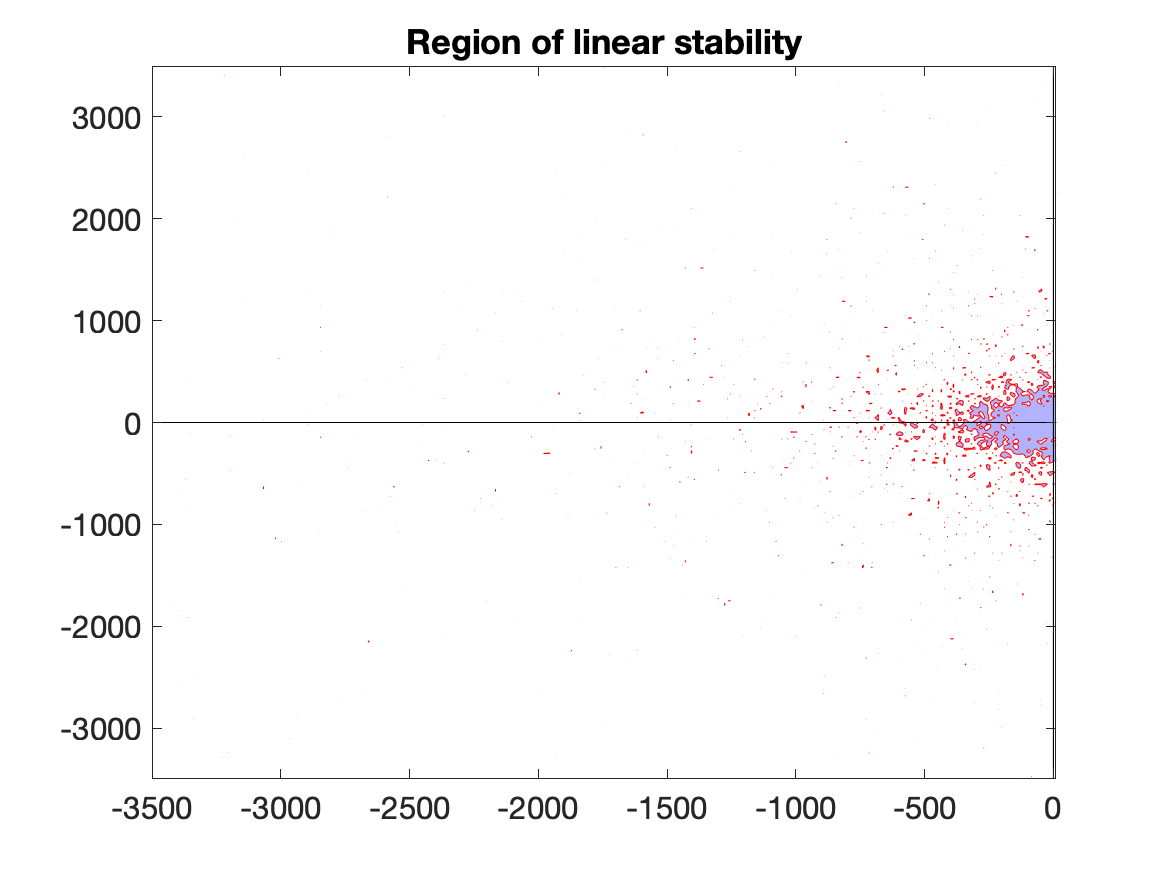}
  \caption{Regions of linear stability for the SDIRK method with one correction (top)
  compared to the Novel A method  \eqref{MP-4s3pA} (bottom),
for $\tilde{\epsilon} = 10^{-8}$ (left) , $\tilde{\epsilon} =  10^{-6}$ (middle) , $\tilde{\epsilon} = 10^{-4}$ (right).}
  \label{NovelAstability}
\end{figure*}

Figure  \ref{NovelAstability} compares the regions of linear stability for the SDIRK method with one correction (top)
 compared to the NovelA method  \eqref{MP-4s3pA}  (bottom), for $\tilde{\epsilon} = 10^{-8}$ (left) , 
 $\tilde{\epsilon} =  10^{-6}$ (middle) , $\tilde{\epsilon} = 10^{-4}$ (right).
We observe that the NovelA method has a significantly smaller region of stability, which  is  more sensitive to 
being adversely impacted by the larger $\tilde{\epsilon}$. These result is certainly significant as we move to 
half and single precision simulations. However, stability is not the only effect that impacts our methods:
in the next subsection we aim to understand the effect of sensitivity to rounding errors
on the behavior of the mixed precision methods.

\subsection{Sensitivity to roundoff errors}  \label{sec:sensitivity}
To understand the sensitivity to rounding errors, we look at the linear stability problem
\[ u_t =  \lambda u, \; \; \; \mbox{and} \; \; \; u_t = \lambda u + \epsilon \tau \]
where $\epsilon$ is the roundoff level, and $\tau$ represents the behavior of the roundoff error
and is $-\frac{1}{2} \leq \tau \leq \frac{1}{2}$. 
We note that this problem is slightly different than the one we considered above, but is more in keeping
with our assumption on the behavior of the truncation errors.
 
Plugging this modified Dahlquist problem into the method \eqref{RK-Butcher} we obtain
\begin{eqnarray}
Y &=& u^n e + \dt A (\lambda Y)  + \dt \Aep ( \lambda Y + \epsilon T)  \\
u^{n+1} &=& u^n  + \dt b [ \lambda Y]
\end{eqnarray}
where $|T| < \frac{1}{2} $ is a vector of roundoff errors.
The first row becomes
\[ Y = [ I - \dt \lambda \At ]^{-1} [ u^n e + \dt \Aep \epsilon T] 
= [ I - \dt \lambda \At ]^{-1} e u^n + \dt \epsilon  [ I - \dt \lambda \At ]^{-1} \Aep T.
 \]
 Plugging this into the second line we get
 \[  u^{n+1} = \left( 1 + \dt   \lambda  b [ I - \dt \lambda \At ]^{-1} e \right) u^n  
 + \dt^2 \epsilon    \lambda  b   [ I - \dt \lambda \At ]^{-1} \Aep T. \]
 Let $\Psi =   z  b [ I - z \At ]^{-1}$ and this becomes
 \[u^{n+1} = \left( 1 +  \Psi e\right) u^n +  \dt  \epsilon  \Psi \Aep T,\] 
 whereas the corresponding high precision evolution can be represented as 
 \[U^{n+1} = \left( 1 +  \Psi e\right) U^n ,\]
 where we ignore rounding errors completely.
 
 The growth of the error between the low and high precision evolution 
 $E^n = \left| u^{n+1} - U^{n+1}  \right| $ is  bounded 
 at each iteration
 \[ E^{n+1} \leq  \left| 1 +  \Psi e \right| E^n +  \dt \left| \Psi \Aep  T \right| \epsilon 
 \leq  \left| 1 +  \Psi e \right| E^n +  \frac{1}{2} \dt \left| \Psi \right| \Aep e \epsilon .\]
 This bound on the growth implies that
 \begin{eqnarray*}
  E^{n} & \leq &  \left| 1 +  \Psi e \right|^n E^0 +  \frac{\epsilon}{2} \dt  
  \left( \left| 1 +  \Psi e \right|^{n-1} + \left| 1 +  \Psi e \right|^{n-2} + \dots \left| 1 +  \Psi e \right| + 1 \right)
\left| \Psi  \right| \Aep e  ,
 \end{eqnarray*}
 where we use that fact that in our methods all the elements of $\Aep $ are non-negative, so $|\Aep  T| \leq \frac{1}{2} \Aep e.$
 We observe that if we set $\epsilon =0$ we get back the usual "exact  precision" method -- in other words, the mixed
precision method is stable in the usual sense that assumes no rounding. The growth of the errors
comes from the buildup of roundoff errors.

 We assume that we are dealing with the case that the method is stable, i.e. $ \left| 1 +  \Psi e \right| \leq 1$, this gives us
  \begin{eqnarray*}
  E^{n}  & =  & \frac{\epsilon \dt }{2}   n \dt
\left| \Psi  \right|  \Aep e  \leq  \frac{\epsilon \dt }{2}    \frac{ 1  }{ 1- \left| 1 +  \Psi e \right| } \left| \Psi \right| \Aep  e .
 \end{eqnarray*}
 Notice that the order conditions may play a role here, in the sense that some of the terms in 
 $\left| \Psi \right| \Aep  e $ zero out and result in higher order terms in $\dt$. This means that, 
 as above, the stability does not depend solely on $z$ but on the value of $\dt$ as well!
 However, if we   ignore that effect for the moment, we see that  the growth of the rounding errors 
 looks like \[ E^{n} \approx  \frac{\epsilon  }{2}  T_f  \left| \Psi \right| \Aep  e.\]
 
\begin{SCfigure}[50][h]
    \vspace{-.2in}    \includegraphics[width=0.5\textwidth]{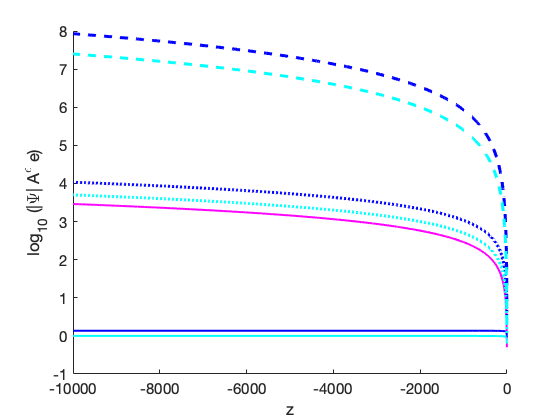} 
  \caption{The values of $log_{10} \left(  \left| \Psi \right| \Aep  e \right) $
 for each of the methods for real valued $z$, with $-10000 \leq z \leq 0$.
 SDIRK methods are in blue, implicit midpoint method in green, and the NovelA method in magenta.\\
 Solid lines: no explicit corrections, \\
 dotted lines: one explicit correction, \\
 dashed lines: two explicit corrections.
  }\vspace{-.1in} 
  \label{MethodSensitivity}
\end{SCfigure}

  Figure \ref{MethodSensitivity} shows the values of $log_{10} \left(  \left| \Psi \right| \Aep  e \right) $
 for each of the methods for real valued $z$, with $-10000 \leq z \leq 0$.
 We observe that the implicit midpoint and SDIRK without corrections have the smallest values.
 Each correction causes the value of $ \left| \Psi \right| \Aep  e$ to grow significantly.
 The NovelA method has values that are close to those of the implicit midpoint and SDIRK methods with
 one correction, but is the smallest of all these, which perhaps explains the behavior of the methods in 
 Table \ref{MPstabilitySDIRKvsNovelA}.

\begin{rmk}
Note that the analysis in this section assumes that each stage that has a low precision computation only introduces a
precision error that looks like $\dt \epsilon$. However,  our implementation does not quite follow the theory exactly. 
All the variables that are computed in an implicit
 step are computed in lower precision and then re-cast to higher precision variables. In this case, the sensitivity analysis 
 above is not quite correct because it is not only the function evaluation that introduces low-precision errors, and these are not
 multiplied by $\dt$. The bound on the error is  better expressed as
  \[ E^{n+1} \leq     \left| 1 +  \Psi e \right| E^n +  \frac{\epsilon}{2}  \left| \Psi \right| \Aep e  ,\]
  and so the growth of the rounding errors  scales with the number  of time-steps
  $ E^{n} \approx   \frac{n \epsilon }{2}  \left| \Psi \right| \Aep  e.$
  This type of roundoff error may grow as $\dt$ is refined.
\end{rmk}

\subsection{Explicit corrections and their impact}
We saw in Figures \ref{SDIRKstability} and \ref{IMRstability} the adverse impact that the explicit corrections 
have on stability, and in   Figure \ref{MethodSensitivity}  we saw that methods with corrections have 
increased sensitivity to roundoff errors. On the other hand, we see in our accuracy analysis that the corrections
confer a significant benefit on our numerical simulations. How do we explain these seemingly contradictory facts?
Consider that accuracy results are asymptotic phenomena, and generally we look at these results as $\dt$ is 
small. On the other hand, the stability and sensitivity results generally focus on the case when $\dt$ is larger.
We see this distinction in practice, while the SDIRK with more corrections had much better accuracy,
it does not always converge for mixed precision where the time-step is large.

Table \ref{MPstabilitySDIRKvsNovelA} shows the maximal stable time-step for the mixed precision SDIRK and 
NovelA method for viscous Burgers' case with $N_x =200$.
While the SDIRK with no corrections   is the most stable for all the mixed precision codes, the NovelA method is
a close second, as it is {\bf stable for all time-steps tested for the mixed/double and mixed/single precision codes}, 
but not for the mixed/half codes. The SDIRK method with one correction is only stable for lower time-steps than the NovelA method,
and increased corrections from one to two require a smaller time-step for stability. 

\begin{table}
\begin{center}
\begin{tabular}{|c|c|ccc|cc|} \hline
Precision & NovelA & \multicolumn{3}{c}{SDIRK} & \multicolumn{2}{c|}{Implicit Midpoint Rule}\\
 &   & $c=0$ &  $c=1$ &  $c=2$ & $c=0$ &  $c=1$\\
128/64 &  all	  &  all & 0.003125   & 0.003125	&  all & 0.003125  \\
128/32 & all	  &  all & 0.003125   & 0.0015625	& all	& 0.003125  \\
128/16 & 0.0125 &  all & 0.0015625 & 0.0015625	& all	& 0.003125  \\
64/32   & all 	  &  all & 0.003125   & 0.0015625	& all	&  0.003125 \\
64/16   & 0.0125 &  all & 0.0015625 & 0.0015625 	& all	& 0.003125 \\ \hline
\end{tabular}
 \caption{Largest stable $\dt$ for each method, for the viscous Burgers' case with $N_x =200$.
 The largest time-step tested was $\dt =0.05$, and going down by powers of two. For comparison,
 the explicit  third order three stage Runge-Kutta is stable for $\dt = 0.0015625$. }
\label{MPstabilitySDIRKvsNovelA}
\end{center}
 \end{table}

What is occurring here is that when the time-step is relatively large compared to $\epsilon$, the 
mixed precision implicit midpoint rule method, which has order $O(\dt^2) + O(\epsilon \dt) $ is accurate enough
because $\dt^2 >> \epsilon \dt$. For such large $\dt$, the correction steps introduce some growth of the 
solution because the local explicit stability limit plays a significant role. To see this, consider the 
implicit midpoint rule method with $c$ corrections:
\begin{eqnarray*}
y^{[0]}_1 & = & u^n + \frac{1}{2} \Delta t F_\epsilon(y^{[0]}_1) \\
y{[c]}_1 & = & u^n + \frac{1}{2} \Delta t F(y^{[c-1]}_1) \\
u^{n+1} & = & u^n +  \Delta t F(y^{[c]}_1) \\
\end{eqnarray*}
With no corrections ($c=0$), we get
\[ y^{[0]}_1 = \frac{1}{1- z/2} u^n + \frac{ \Delta t \epsilon}{2 ( 1- z/2) }  = \frac{1}{1- z/2} u^n + \frac{q/2}{ ( 1- z/2) } ,\]
where $|q| \approx \dt \epsilon $ is a function that represents the rounding error scaled by the time-step, and so varies from line to line.
The one-step impact is 
\[u^{n+1} = u^n +  z \left( \frac{1}{1- z/2} u^n + \frac{q/2}{ ( 1- z/2) } \right)  
\phi u^n +  \frac{z/2}{ ( 1- z/2) } q ,\]
where $ \phi $ is the usual linear stability polynomial.

With one explicit correction
\[ y^{[0]}_1 = \frac{1}{1- z/2} u^n + \frac{ \Delta t \epsilon}{2 ( 1- z/2) }  = \frac{1}{1- z/2} u^n + \frac{q/2}{ ( 1- z/2) } \]
\[y^{[1]}_1 = u^n + \frac{z}{2} \left(  \frac{1}{1- z/2} u^n + \frac{q/2}{ ( 1- z/2) }  \right)
= \left( 1 + \frac{z/2}{1- z/2} \right) u^n + \frac{z/4}{ ( 1- z/2) } q \]
so
\[ u^{n+1} = u^n +    
z \left( 1 + \frac{z/2}{1- z/2} \right) u^n +  \frac{z^2/4}{ ( 1- z/2) } q = \phi u^n +  \frac{(z/2)^2}{ ( 1- z/2) } q .\]
More corrections give
\begin{equation} \label{sensitivity}
u^{n+1} = \phi u^n + \frac{(z/2)^{(c+1)}}{ ( 1- z/2) } q .
\end{equation}
The key observation is that if  $|z| > 2$ then each correction causes the rounding errors at each time-step
to grow. If, on the other hand, we have $|z| < 2$, then the corrections will causes the rounding errors at each time-step
to decay, and have their intended benefit. 

This highlights the fact that the corrections are an explicit process
that are effective within the  explicit stability limit $\dt < \frac{2}{|\lambda|}$. On the other hand, the corrections are needed when
$\epsilon \dt > \dt^2$ (for the implicit midpoint rule), i.e. when $\dt < \epsilon$. We typically have $ \epsilon << \frac{2}{|\lambda|}$,
so that the values of $\dt$ in which the corrections are needed and in which they are beneficial and stable are the same.
For the SDIRK method,  the corrections are needed when $\epsilon \dt > \dt^3$, i.e. when $\dt < \sqrt{\epsilon}$. The explicit corrections 
will be beneficial if  $ \dt \leq  \frac{1}{|\lambda|}$, so these two criteria are consistent except in the 
region where $ \epsilon > \frac{1}{|\lambda|^2} .$ If $\lambda$ is very large, we must use a higher precision method
(so that $ \epsilon < \frac{1}{|\lambda|^2} $) for this approach to work.

\section{Performance analysis of MP-ARK methods} \label{sec:performance} 

For this study we look at the performance of  the three methods we review above in Section 
\ref{sec:methods}. These are the mixed precision implicit midpoint method
with and without correction steps, the mixed-precision SDIRK method with and without correction steps, 
and a novel MP-ARK method from \cite{Grant2022}
designed to suppress the low precision errors without additional correction steps. 

The low precision evaluations of the implicit solvers
require a subroutine that maps the variables low precision, performs a Newton iteration
with tolerances that are consistent with the level of the low precision, and 
returns the low precision value of the stage variable $y^{(i)}$. This value is then cast to
a high precision variable $y^{(i)}$ to be used in the explicit evaluations of $F( y^{(i)})$ in the
explicit correction stages  and explicit final reconstruction stage.

\subsection{Implementation Details}

\noindent{\sc Fortran code}
The mixed-precision Runge-Kutta methods that we explored were implemented in Fortran using the 2008 standard. 
The \verb|iso_fortran_env| intrinsic module was also used to import the \verb|Real32|, \verb|Real64|, and \verb|Real128| derived types.
The implicit solvers were implemented using the Newton-Raphson method for multivariate systems. 

The tolerance of the implicit solvers were set using the Fortran intrinsic machine epsilon. 
The maximum number of iterations for the solver was set to $20$ regardless of the precision, however this fail-safe maximum was not attained. 
These setting do have significant impact on the overall runtime of the method and were determined to be fair for the consideration of mixed precision methods, based on our tests showing that the implicit solvers on average only used $4$ iterations to converge.

Higher precision computations  (specifically quadruple precision but also arbitrary precision) 
are generally evaluated in software on a majority of platforms. 
Notable exceptions for this are the IBM POWER9 
that have hardware level support for quadruple precision. In this work, we compare the runtime performance
on both software based high precision computations and hardware based high precision computations.
The experiments were performed on a system with an Intel Xeon Silver 4116 CPU running 
CentOS 7 and an IBM POWER9 running CentOS 8. 
We will refer to these as x86 and POWER9, respectively.
The experiment performed looks at three levels of precision: single, double, and quadruple. 
The Intel system is limited to hardware evaluation for single and double precision and requires 
quadruple precision to be evaluated using software. The POWER9 system is capable of hardware 
evaluation of all precisions. This was done to show the usefulness of these methods under 
both conditions and that the benefits are platform independent. 

The codes were compiled on both architectures using gfortran version 11.1.0 using no optimization flags. 
On the POWER9 system, the GCC compilers were compiled using the IBM AT14.0 libraries 
in order to enable hardware quadruple precision support while the GCC compilers on the Intel system
were compiled using the default compiler settings. Timing on the method subroutine was done using the 
intrinsic \verb|cpu_time| Fortran subroutine. 

The GitLab repository  \url{https://gitlab.com/bburnett6/mixed-precision-rk}. 
contains all the necessary files to reproduce the study performed here. 
The only dependencies are a working version of the GCC 11.1.0 compilers and an environment 
with the Numpy and MatPlotLib python libraries for the experiment driver and analysis 
files that were implemented in python.

\noindent{\sc Julia code}
Although Julia is typically a dynamically typed language, we developed the codes used 
in this paper to ensure strong typing of all functions and variables. 
Efforts for the Julia code development were focused on maintaining similarity to the Fortran codes used. 
This included using similar selections of parameters to a Newton-Raphson method re-implemented in Julia.

We made use of the built-in language support for half, single, and double precision and used the 
Quadmath.jl Julia wrapper package for the gcc libquadmath quadruple precision support. 
This means that all precisions are IEEE 754 compliant just like the Fortran codes. 
Other packages used were the LinearAlgebra.jl Julia package for linear solves inside 
our Newton-Raphson routine, and TimerOutputs.jl for gathering the run times of our methods. 
The use of a linear solve over a matrix inverse in the Newton-Raphson is likely the source 
of any numerical difference between the Fortran codes and the Julia codes.
The Julia codes can be found at https://gitlab.com/bburnett6/mpark-julia

The codes were run serially, again to maintain comparable results to the Fortran codes, 
on the UMass Dartmouth Carnie Cluster employing Intel Xeon Silver 4116 CPUs using the 
generic Linux Julia binary version 1.7.3. Between the time that the Fortran and the Julia codes were run, 
the cluster had undergone a system upgrade and moved from Centos 7 using a Linux kernel 
version 3.10 to Ubuntu 20.04 using a Linux kernel version 5.15. 
These kernel differences did impact the performance of the codes.

\subsection{Van der Pol system}
We solve the non-stiff van der Pol system  in Eqn. \eqref{vdp}. A reference solution using an explicit fourth order Runge--Kutta at a small time step was used for calculating the errors. The reference solution was computed entirely in quadruple precision and all methods results regardless of computational precision, were cast to quadruple precision in order to compute the error between a method and this reference solution.

\subsubsection{Implicit midpoint rule}

\vspace*{.1in}

\noindent {\sc Fortran code:}
The table below shows the performance of the 
mixed precision IMR method with two corrections, at different levels of errors. 
On the x86  the mixed precision double-single code is the {\em most} efficient. For an error level of 
$\approx 10^{-15}$, the mixed precision quad-single code has a fifteen-fold reduction in runtime.
The mixed precision double-single code gives the same level of error at $\approx 10^{-9}$
as the double precision code with a runtime reduction of 2.6x.
On the POWER9 (the runtime is generally better, but the
runtime   savings are similar. For the  error level of  $\approx 10^{-9}$, we see 3x savings from the mixed precision
double-single over the the double precision. For the  error level of  $\approx 10^{-15}$, 
we see more than 14x savings from the mixed precision quad-single over the the quad precision. 

\smallskip

\begin{center} {\small 
\begin{tabular}{ |c|c|c|c|c|c|c|  }
 \hline
& \multicolumn{6}{|c|}{Mixed precision IMR method with two corrections}\\ \hline
\parbox[t]{2mm}{\multirow{3}{*}{\rotatebox[origin=c]{90}{x86}}}  &  error & 64/64 & 64/32 & 128/128 & 128/64 
 & 128/32 \\
& $\approx 10^{-9}$ & 0.013 & 0.005 & 0.203 & 0.020 & 0.012 \\
& $\approx 10^{-15}$ & N/A & N/A & 108.4 & 13.91 & 7.165 \\ \hline
\parbox[t]{2mm}{\multirow{3}{*}{\rotatebox[origin=c]{90}{{\footnotesize \sc power9}}}}  &  error & 64/64 & 64/32 & 128/128 & 128/64 
 & 128/32 \\
& $\approx 10^{-9}$ & 0.021 & 0.007 & 0.0683 & 0.023 & 0.008 \\
& $\approx 10^{-15}$ & N/A & N/A & 34.64 & 14.00 & 2.395 \\ \hline
 \end{tabular}}
\end{center}

\smallskip

\begin{figure}[t] \vspace{-0.75in}
\begin{center}
  \includegraphics[width=0.495\textwidth]{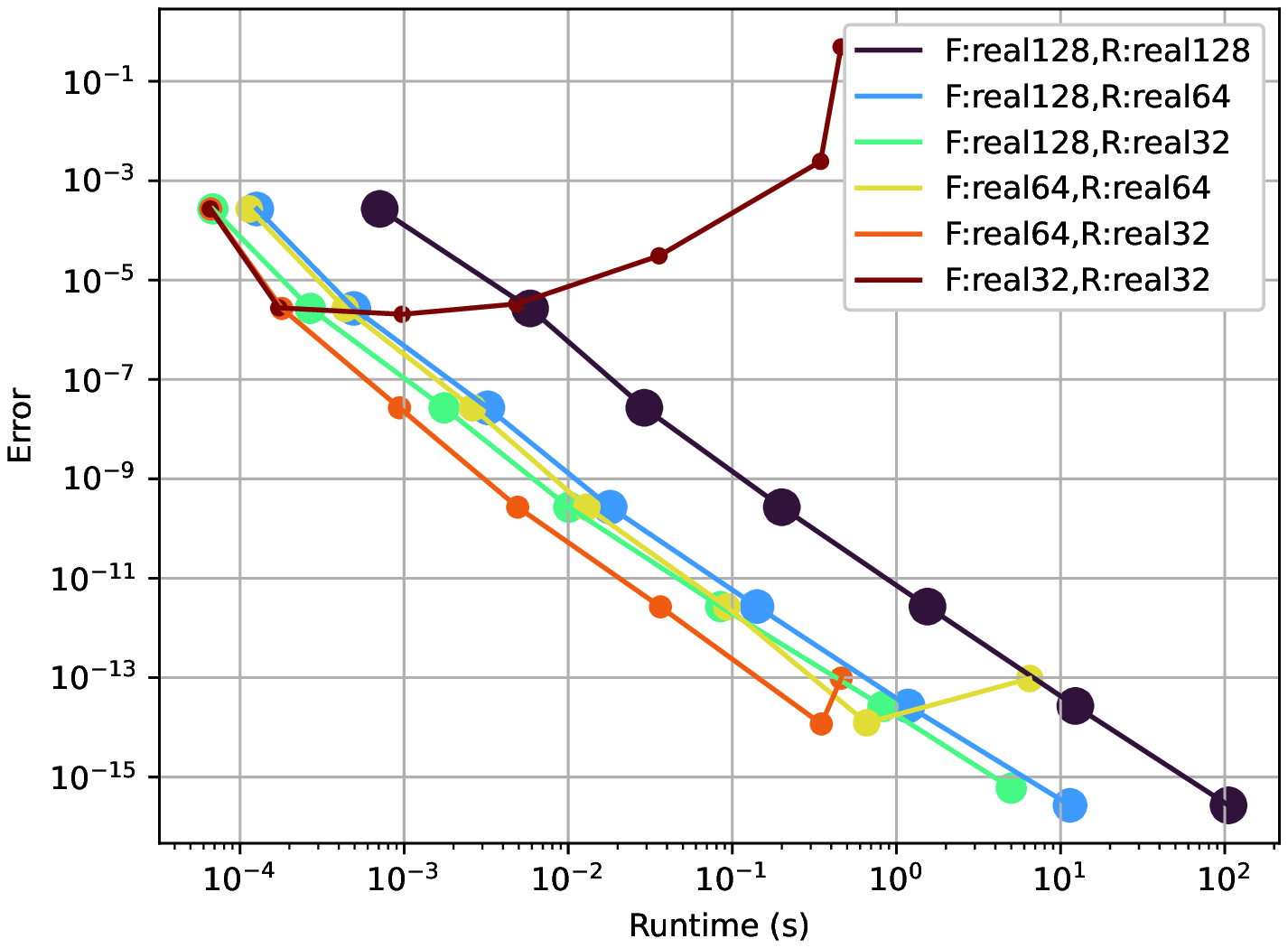}
   \includegraphics[width=0.495\textwidth]{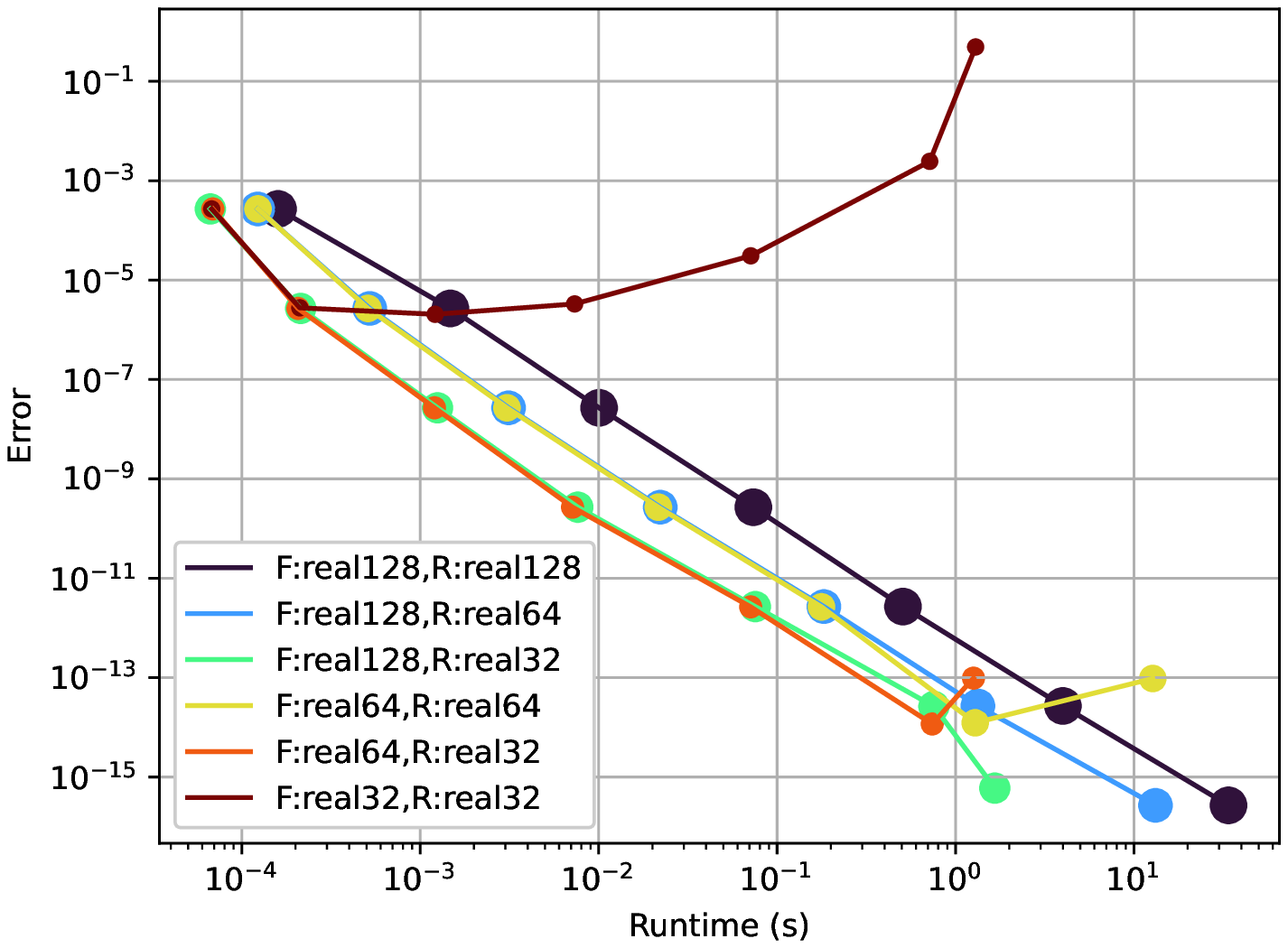}
   \vspace{-.25in}
\caption{Mixed precision implicit midpoint method with one correction step: error vs runtime. Left: x86; right: POWER9.
\label{impmid_1p_time}
}
\end{center} \vspace{-.15in}
\end{figure}

 In Figure \ref{impmid_1p_time} we plot  the numerical error against the runtime  
for the mixed precision implicit midpoint method with one correction. 
Remarkable in this figure is the significant reduction in runtime 
from the quad precision to the mixed quad-double precision for the x86 chip, 
and from the quad precision to the mixed quad-single precision for the POWER9 chip. 
The POWER9 chip has hardware-level support for quad precision, so it is not
surprising to see that the cost of quad precision on the POWER9 chip is significantly less than on the x86 chips,
but the runtime savings is not  impacted, as we see in the table below.

\begin{SCfigure}[][htb]  \vspace{-.1in}
\includegraphics[width=0.5\textwidth]{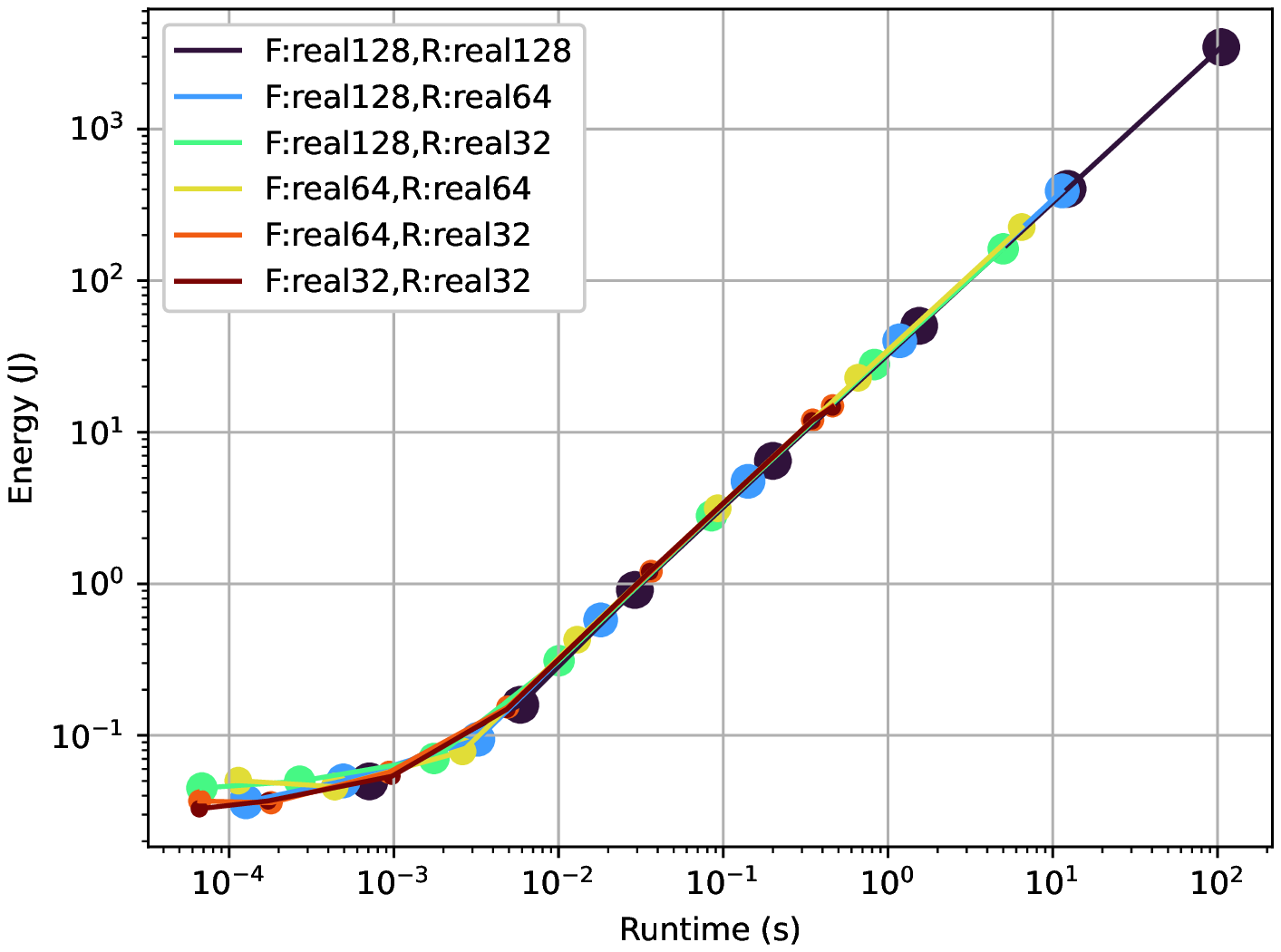}
  \caption{Total energy used by the mixed precision implicit midpoint rule with one correction 
   for the  van der Pol system. \\
   This is for the Fortran code run on the x86 chips. \\
   For larger values of $\dt$, the energy consumptions barely increases
   with the runtime but once the time-step is small enough, the energy
  scales linearly with the runtime regardless of the precision used.
      \label{energy-runtime}}
  \end{SCfigure}

The x86 chips allow us to measure the energy consumption directly. 
The energy consumption results for this simulation follow very closely the runtime results.
  Figure \ref{energy-runtime} shows that for long enough runtimes, the energy consumed scales 
  linearly with the runtime for all precisions. 
  
  \vspace{-.15in}

\noindent {\sc  Julia language  code:} 
Julia language allows us to extend our computation to half precision as well. 
The table below shows the computational time (in seconds) for each of the mixed precision
implicit midpoint rule codes with $c=0,1,2$ corrections, for the van der Pol equation with 
$\alpha=1 $ and $\alpha=6 $. The solutions are different for the two cases, and so we cannot compare the errors
from the $ \alpha=1 $ case to the $\alpha=6 $ case. However, We can observe that the less stiff problem ($\alpha=1 $)
performs better in terms of the mixed precision codes in the sense that without corrections we attain an error of $10^{-10}$
for all the mixed precision codes except for the quad/half (Float128/Float16) code. On the other hand,
for the  stiffer problem ($\alpha=6 $), an error level of $10^{-10}$ is not attained
for all the non-corrected mixed/single and mixed/half codes. Correction terms fix this problem.
The table below shows that more corrections do not always improve efficiency, but tend to do so for most of the mixed/half codes,
and at times for the mixed/single codes. For the quad/double codes and the double/single code, one correction is optimal.

\smallskip

\begin{center}
{\small
\begin{tabular}{ |c|c|ccccccc| }
 \hline
 \multicolumn{9}{|c|}{Implicit midpoint rule $\approx 10^{-10}$ }\\ \hline
& corrections &  128/128 &128/64 & 128/32 & 128/16 &  64/64 & 64/32 & 64/16 \\
\parbox[t]{2mm}{\multirow{3}{*}{\rotatebox[origin=c]{90}{{$\alpha=1 $}}}}  
& $c=0 $ &  0.533 &  0.156 & 0.311 & --       & 0.094 & 0.175 & --  \\
& $c=1 $ &  0.533 & 0.149  & 0.123 & 0.115 & 0.094 & 0.079 & 0.045 \\
& $c=2 $ &  0.533 &  0.207 & 0.136 & 0.102 & 0.094 &  0.106 & 0.044 \\ \hline
\parbox[t]{2mm}{\multirow{3}{*}{\rotatebox[origin=c]{90}{{$\alpha=6 $}}}}  
& $c=0 $ & 0.059 & 0.009 & -- 	     & --        & 0.018 & -- & -- \\
& $c=1 $ & 0.059 & 0.010 & 0.026 & 0.047 & 0.018  & 0.005 & 0.022 \\
& $c=2 $ & 0.059 & 0.011 & 0.018 & 0.004 & 0.018  & 0.003 & 0.002 \\ \hline
  \end{tabular}
  }
  \end{center}
  
 \smallskip
 
\begin{figure}[htb]  \vspace{-.15in}
\begin{center}
{\includegraphics[width=0.45\textwidth]{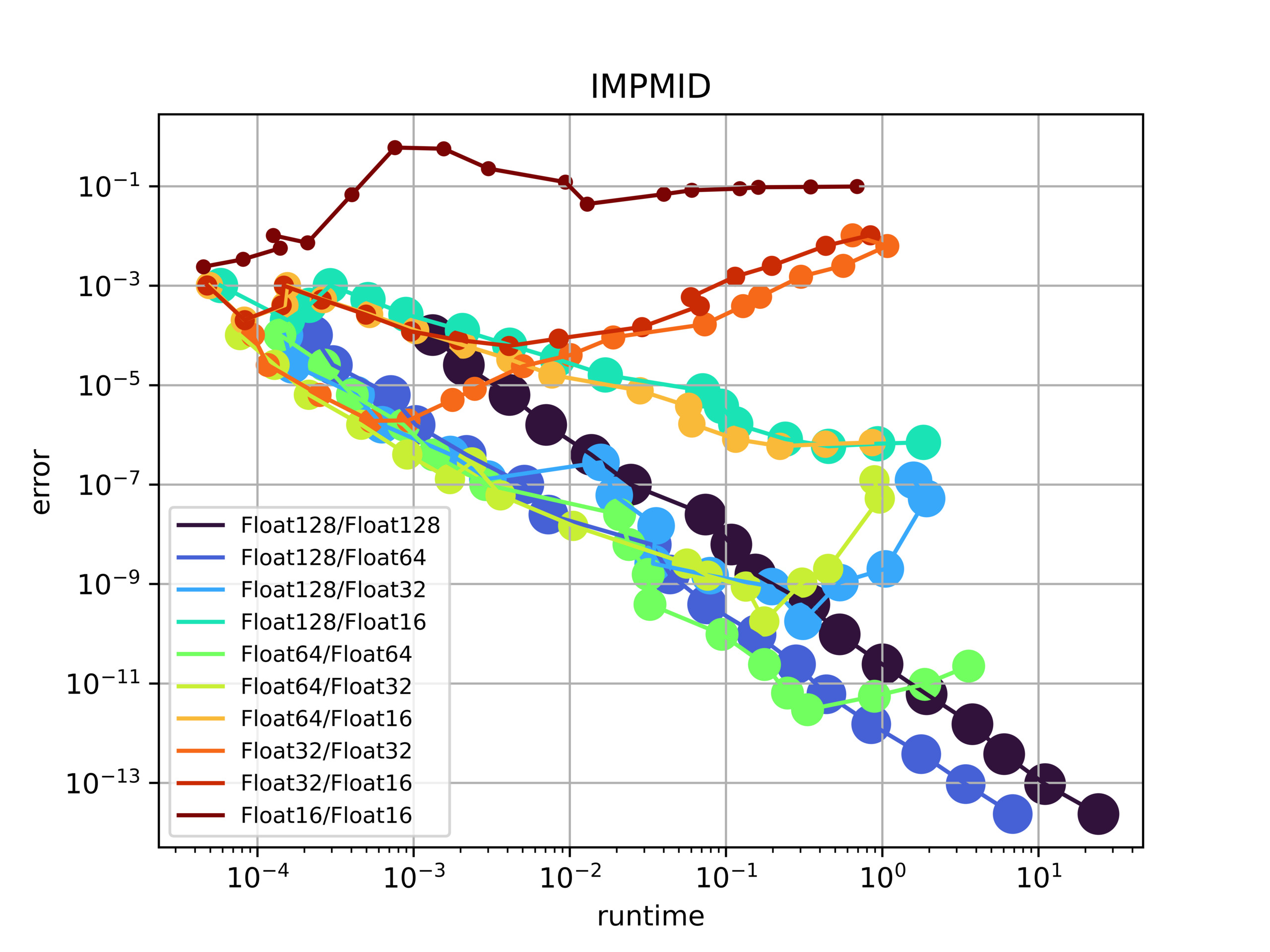}}
{\includegraphics[width=0.45\textwidth]{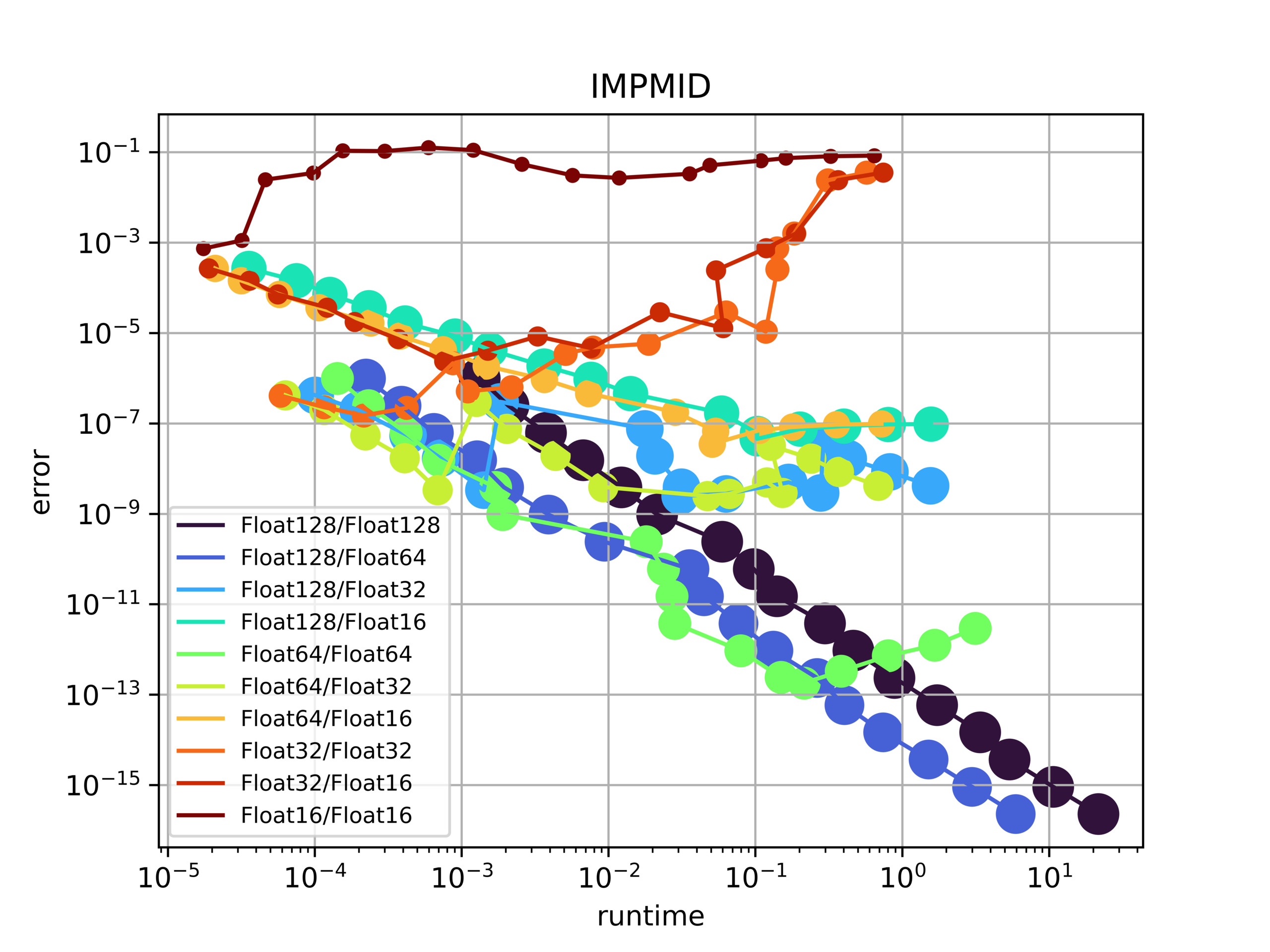}} \\
{\includegraphics[width=0.45\textwidth]{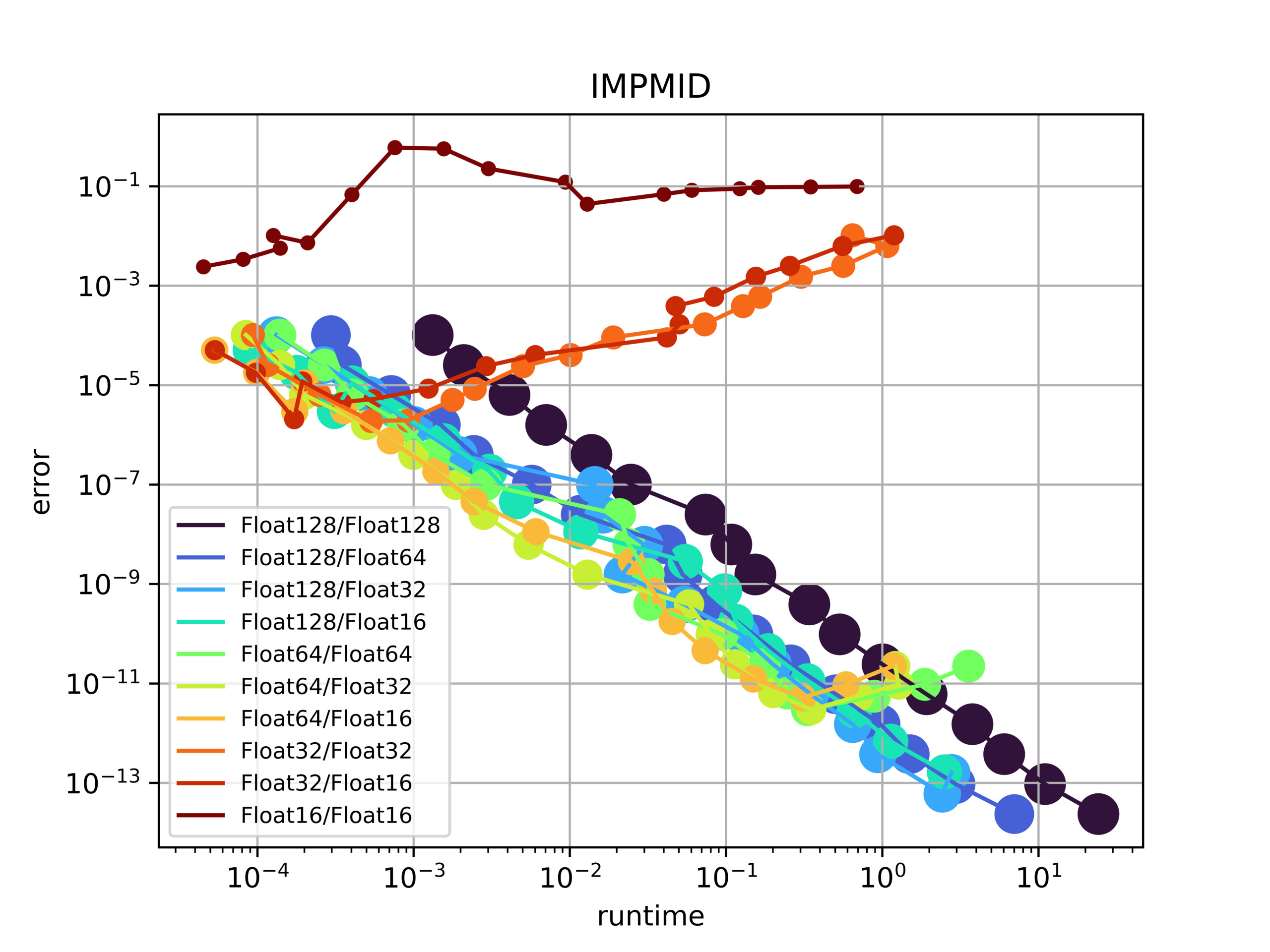}}
{\includegraphics[width=0.45\textwidth]{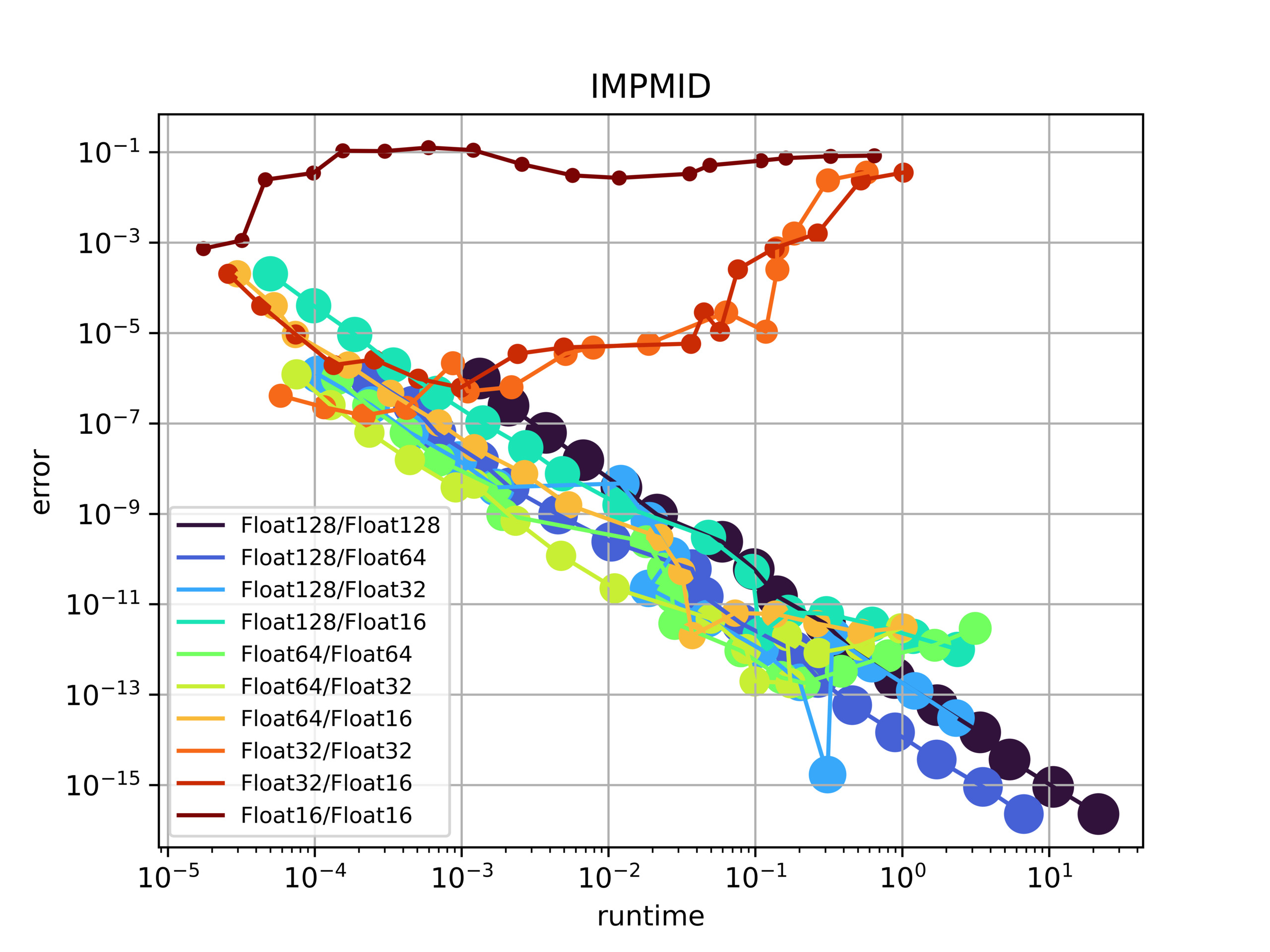}} \\
{\includegraphics[width=0.45\textwidth]{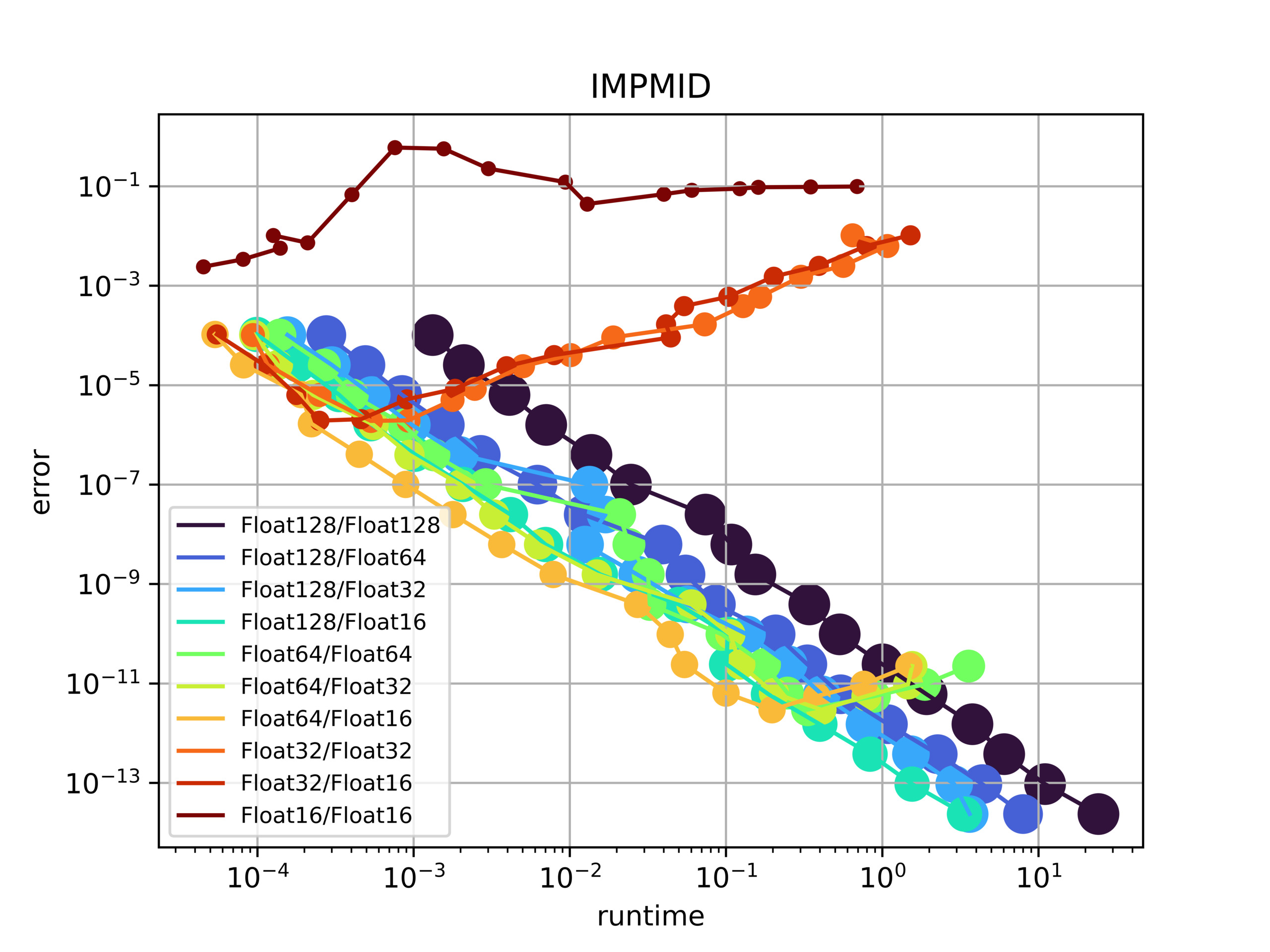}}
{\includegraphics[width=0.45\textwidth]{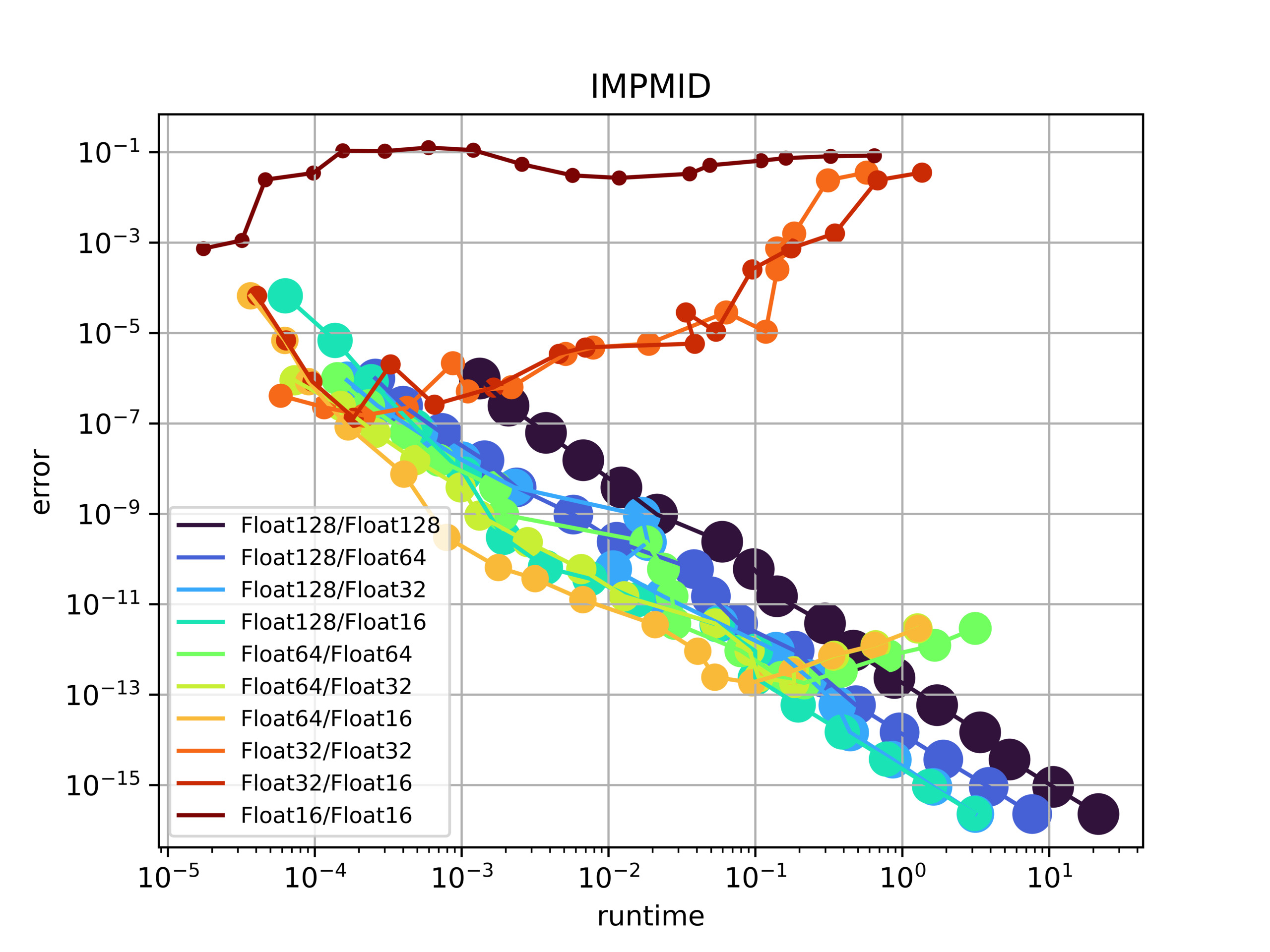}}\\
\caption{Julia code results: the implicit midpoint rule for the van der Pol equation with $\alpha=1$ (left) and $\alpha=6$ (right)
 with no corrections (top), one correction (middle), two corrections (bottom).
\label{vdpIMRefficiency}
}
 \vspace{-.25in}
\end{center}
\end{figure}

Figure  \ref{vdpIMRefficiency}  shows the error vs. runtime for  the implicit midpoint rule 
for the van der Pol equation with $\alpha=1$ (left) and $\alpha=6$ (right),
 with no corrections (top), one correction (middle), two corrections (bottom).
 We observe that successive corrections improve the efficiency of the mixed precision methods.
 From left to right, we see that the light orange line corresponding to the mixed double/half (Float64/Float16)
 has a dramatic improvement with one correction and further improvement with two corrections.
 A similar behavior is seen for the mixed precision quad/half (Float128/Float16) codes in the blue-green line. 
 A similar behavior is seen for the  mixed precision double/single (Float64/Float32 in chartreuse), 
 mixed precision quad/single (Float128/Float32 in light blue). The mixed precision quad/double performs well without
 corrections. This behavior is consistent for both $\alpha=1$ (left) and $\alpha=6$ (right), however we observe
 that for the larger $\alpha$, the mixed precision methods are less consistent initially, but after two corrections
they  end up with smaller errors for equivalent runtimes.

\newpage

\subsubsection{\bf Mixed precision SDIRK Method}

 \begin{figure*}[ht!]
 \begin{center}
  \includegraphics[width=0.45\textwidth]{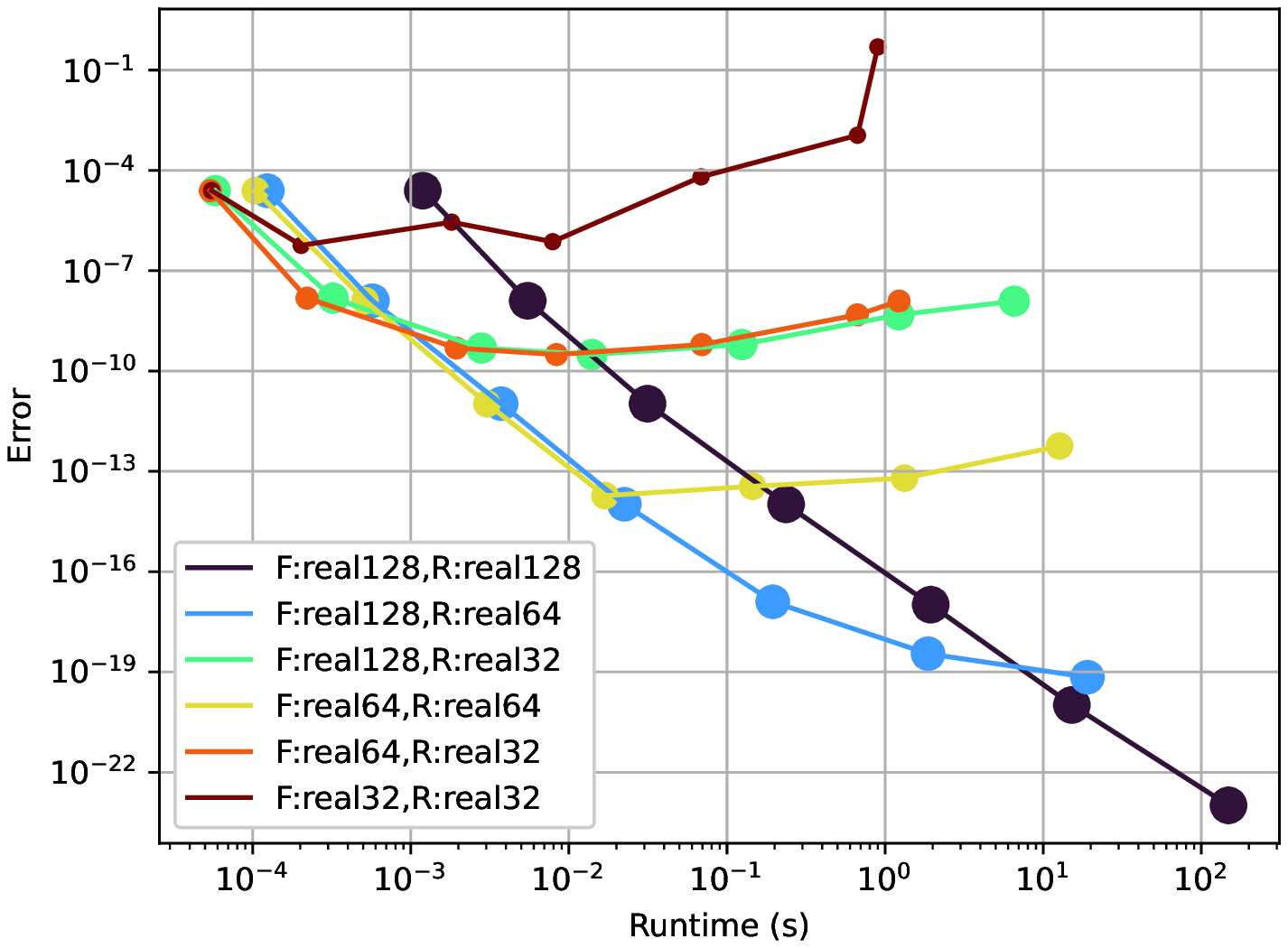}
    \includegraphics[width=0.45\textwidth]{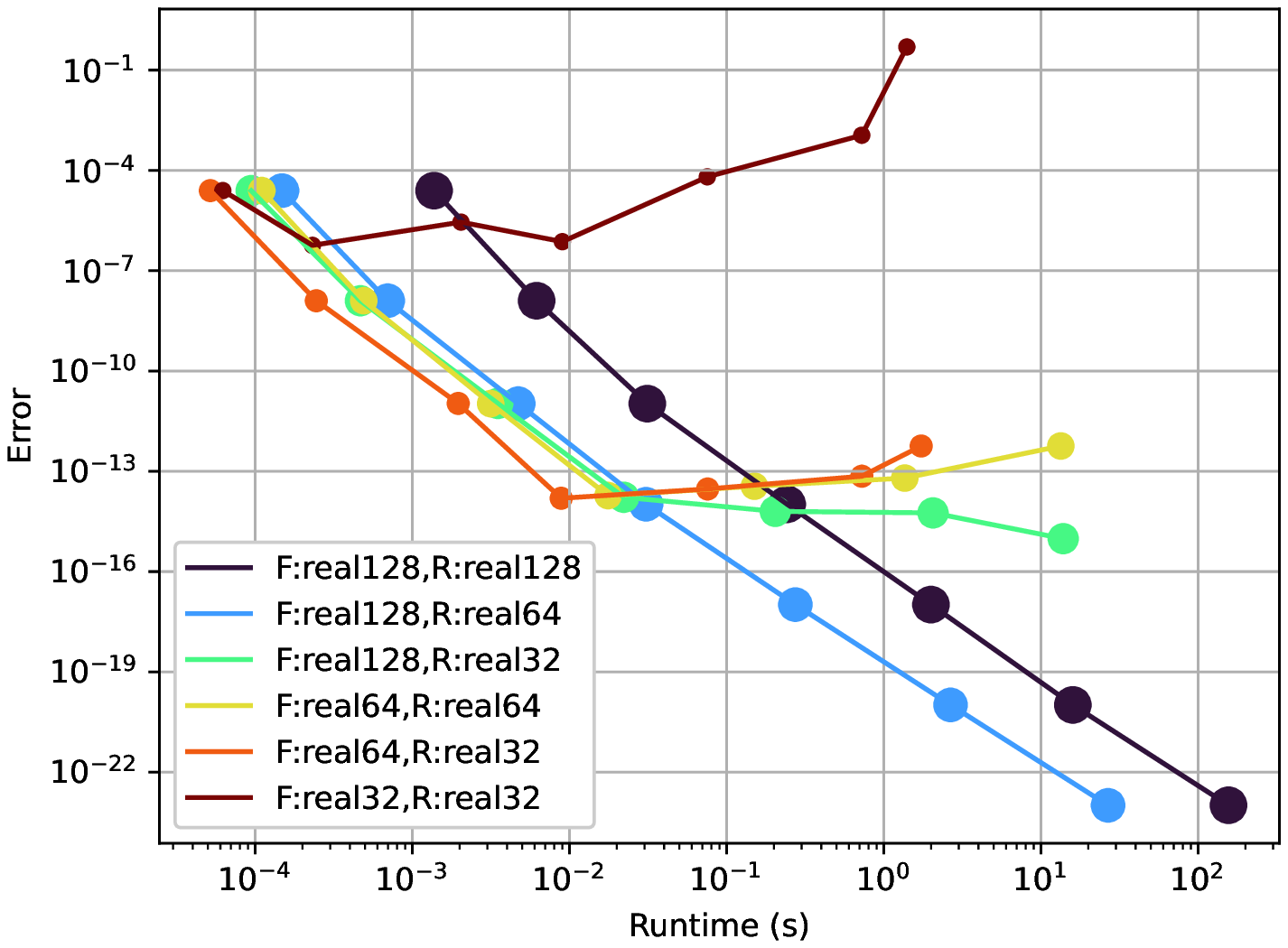}\\
      \includegraphics[width=0.45\textwidth]{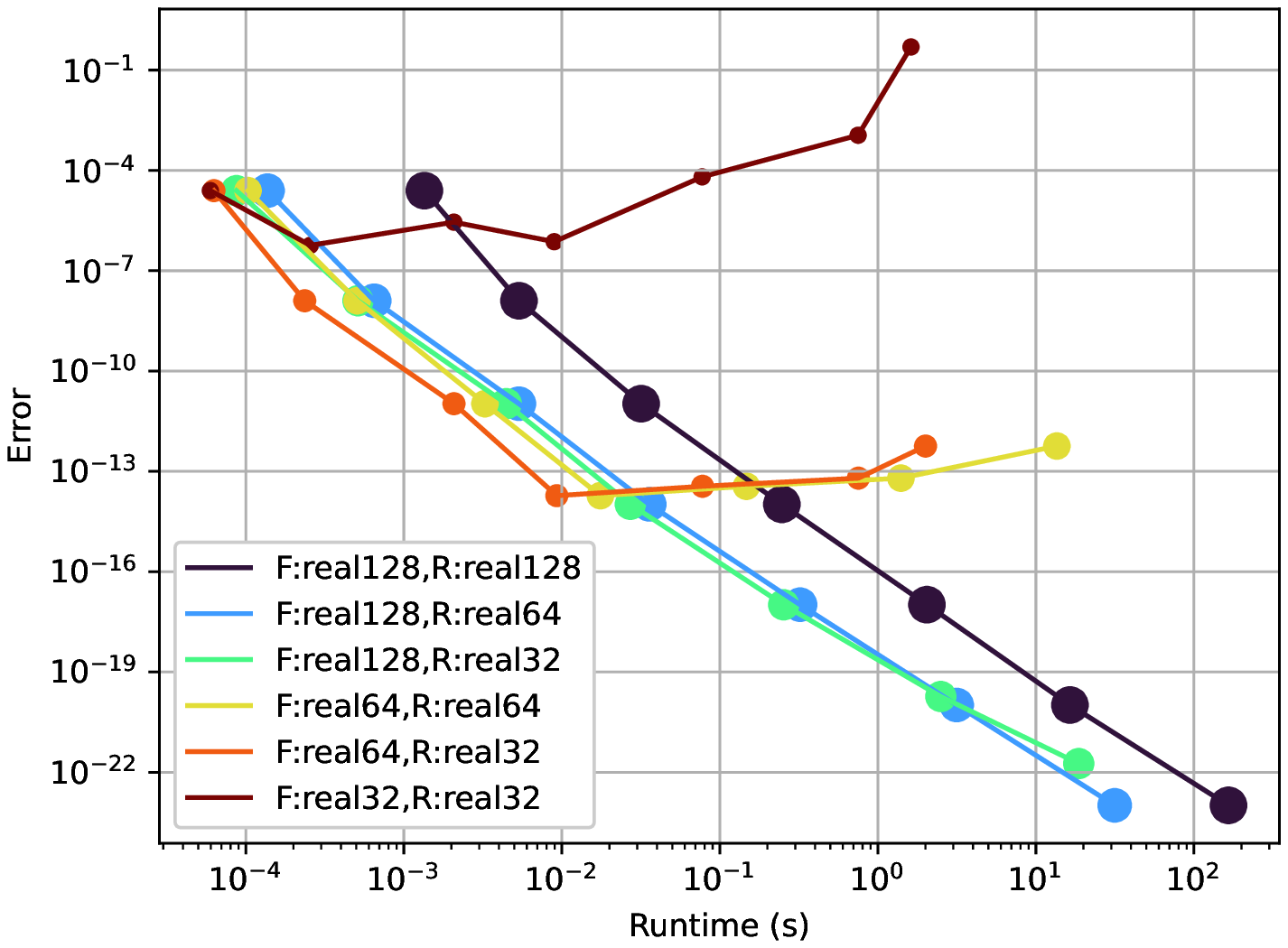}
  \includegraphics[width=0.45\textwidth]{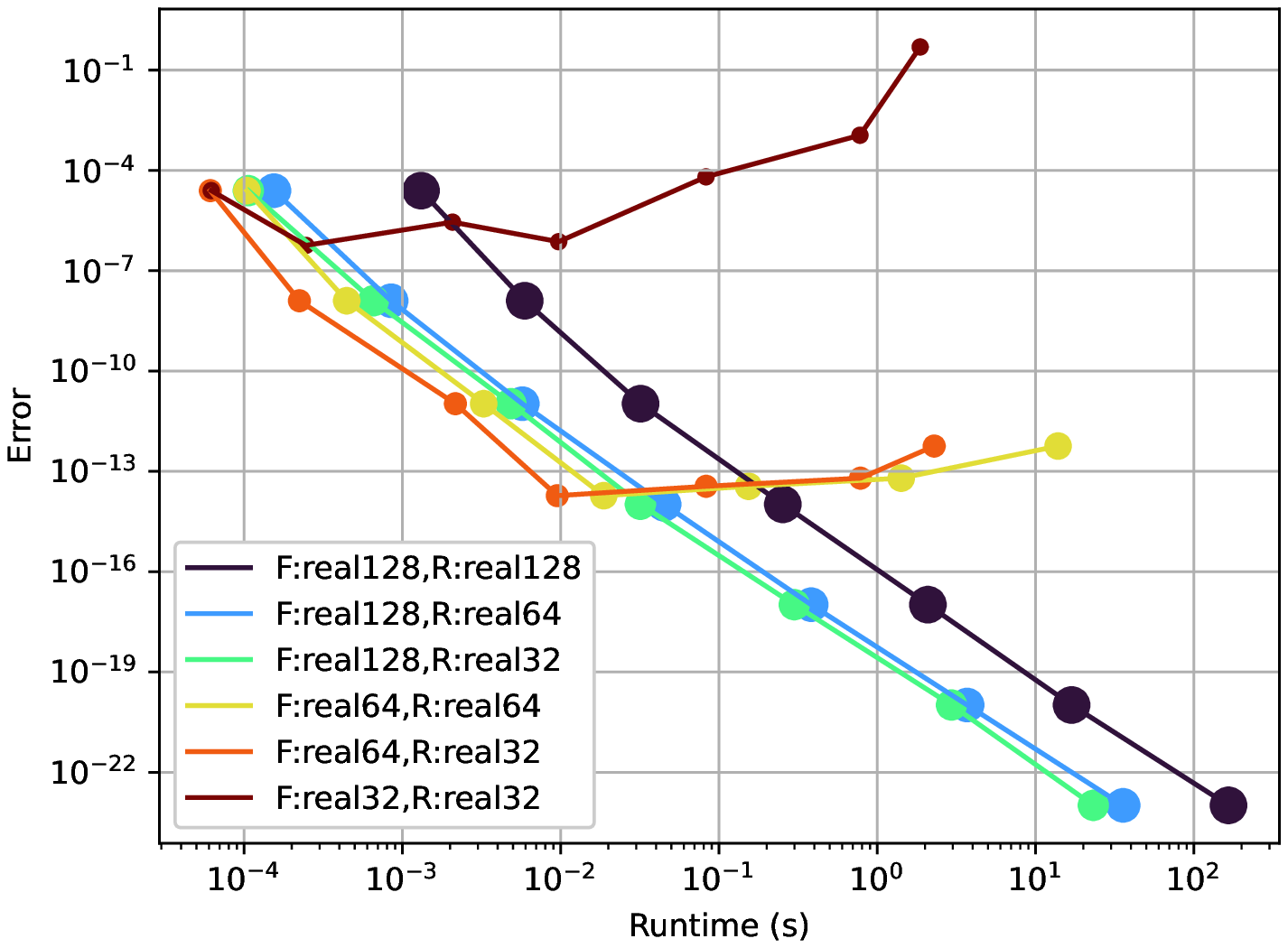}\\
  \caption{Mixed precision SDIRK method errors vs runtime on x86. 
  Top left: no corrections;
  Top right: one correction;
Bottom left: two corrections;
  top right: three corrections.
  \label{sdirk-eff-x86}
  }
  \end{center}
\end{figure*}

\noindent{\sc FORTRAN code}
The results for the mixed precision singly diagonally implicit Runge-Kutta (MP-SDIRK) in 
Section  are similar to those of the implicit midpoint rule.
For an error level of $10^{-13}$ we see that the double/single code runs twice as fast as the 
double precision code for both the x86 and {\sc power9} chips.
The quad precision codes are significantly more efficient -- by a factor of 3 --
on the  {\sc power9} chips as on the  the x86 chips.
For an error level of $10^{-22}$ the  quad/double code is 1.8x  as fast as the full quad precision code, 
and the quad/single code is  more than 6x as fast as the full quad precision code for the {\sc power9}
chip.The x86 chips are less efficient in quad precision than the  {\sc power9}, so we see greater
efficiency in moving to mixed precision.
For the x86 chips the quad/double code is 4.6x as fast as the full quad precision code, 
and the quad/single code is  more than 7x as fast as the full quad precision code. 

\smallskip
 
 \begin{center}
 \begin{tabular}{|c|c|c|c|c|c|c|  }
 \hline
 \multicolumn{7}{|c|}{Mixed precision SDIRK with three corrections}\\ \hline
 \parbox[t]{2mm}{\multirow{3}{*}{\rotatebox[origin=c]{90}{{\footnotesize \sc x86}}}}  &
 error & 64/64 & 64/32 & 128/128 & 128/64 
 & 128/32 \\
& $\approx 10^{-13}$ & 0.018 & 0.009 & 0.253 & 0.045 & 0.031 \\
& $\approx 10^{-22}$ & N/A & N/A & 165.3 & 35.72 & 23.15 \\ \hline 
 \parbox[t]{2mm}{\multirow{3}{*}{\rotatebox[origin=c]{90}{{\footnotesize \sc power9}}}}  &
 error & 64/64 & 64/32 & 128/128 & 128/64 
 & 128/32 \\
& $\approx 10^{-13}$ & 0.031 & 0.016 & 0.092 & 0.034 & 0.019\\
& $\approx 10^{-22}$ & N/A & N/A & 54.33 & 30.06 & 8.151\\ \hline
 \end{tabular}
\end{center}

\smallskip

Figure \ref{sdirk-eff-x86} shows the runtime per error of the mixed precision SDIRK  method with no corrections, 
one correction, two corrections, and three corrections on the x86 chip.
One correction shows a dramatic difference in the accuracy of the method, but more corrections seem to
enhance the efficiency of the method. The figure confirms that mixed precision approach provides significant time-savings.
For errors  above $\approx 10^{-13}$, the mixed precision double-single method performs best. 
For errors below that level, the mixed precision   quad-single method is most efficient, as shown in the table above.

\begin{figure}[htb]
\begin{center}
{\includegraphics[width=0.32\textwidth]{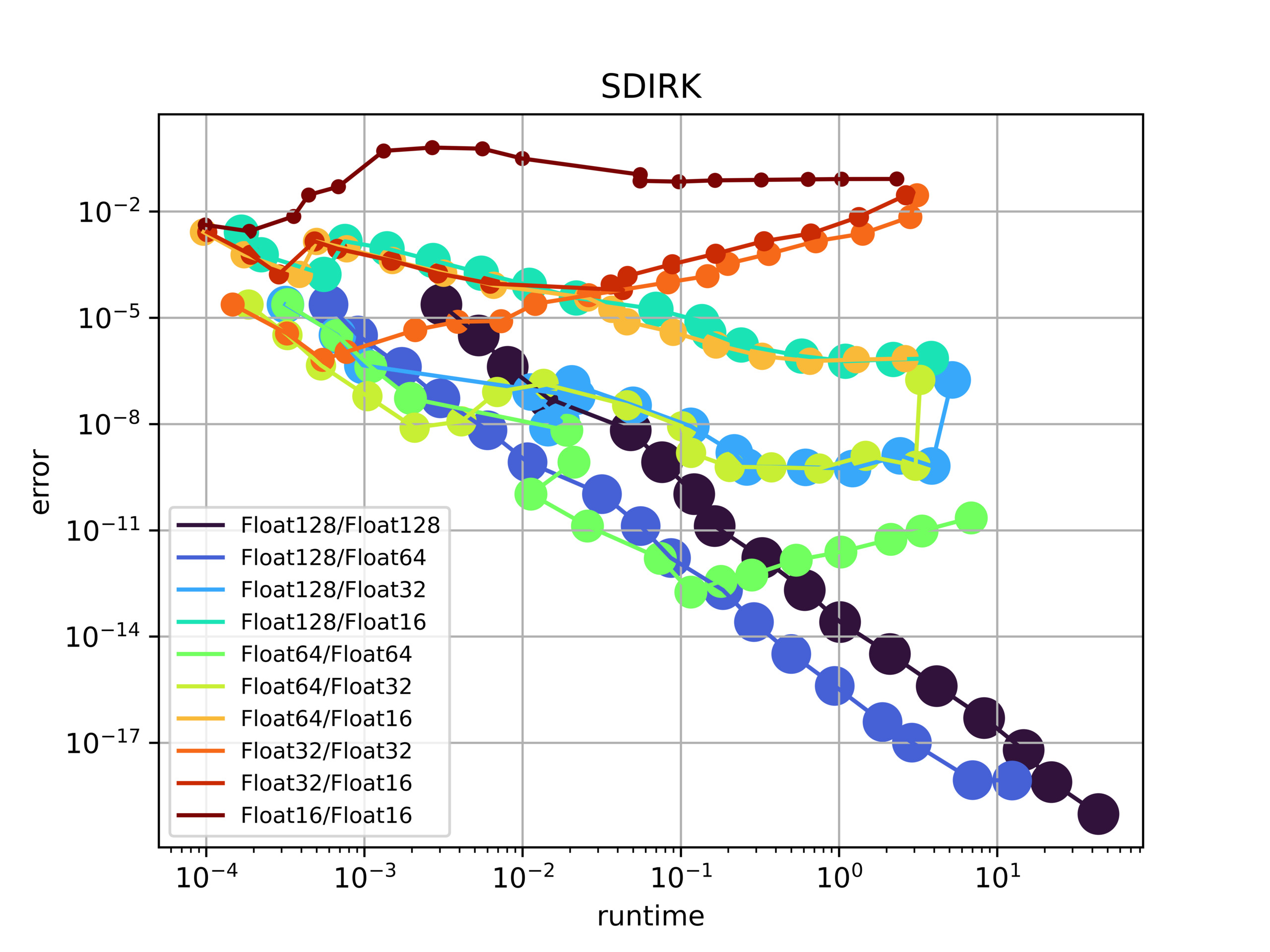}}
{\includegraphics[width=0.32\textwidth]{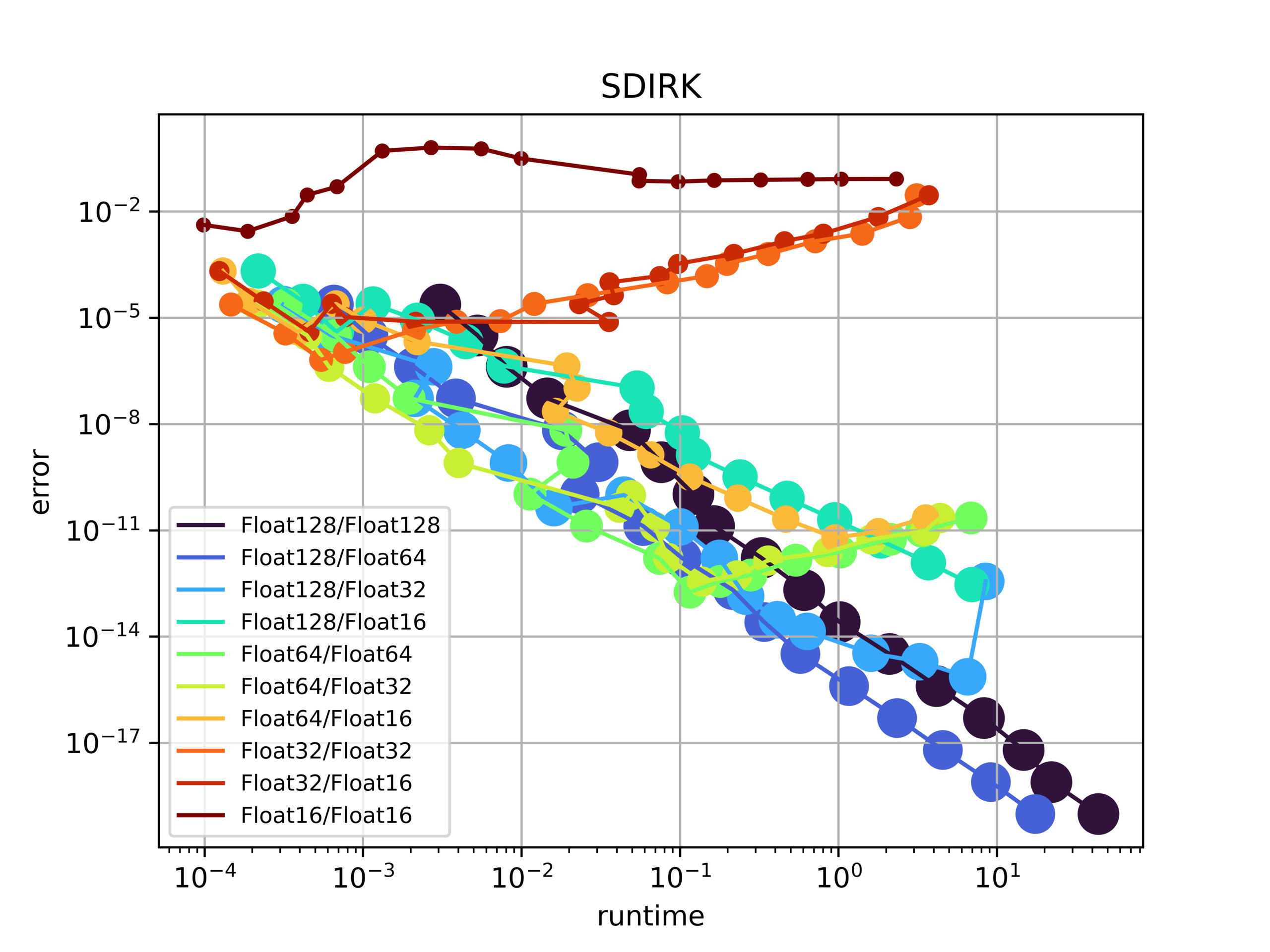}}
{\includegraphics[width=0.32\textwidth]{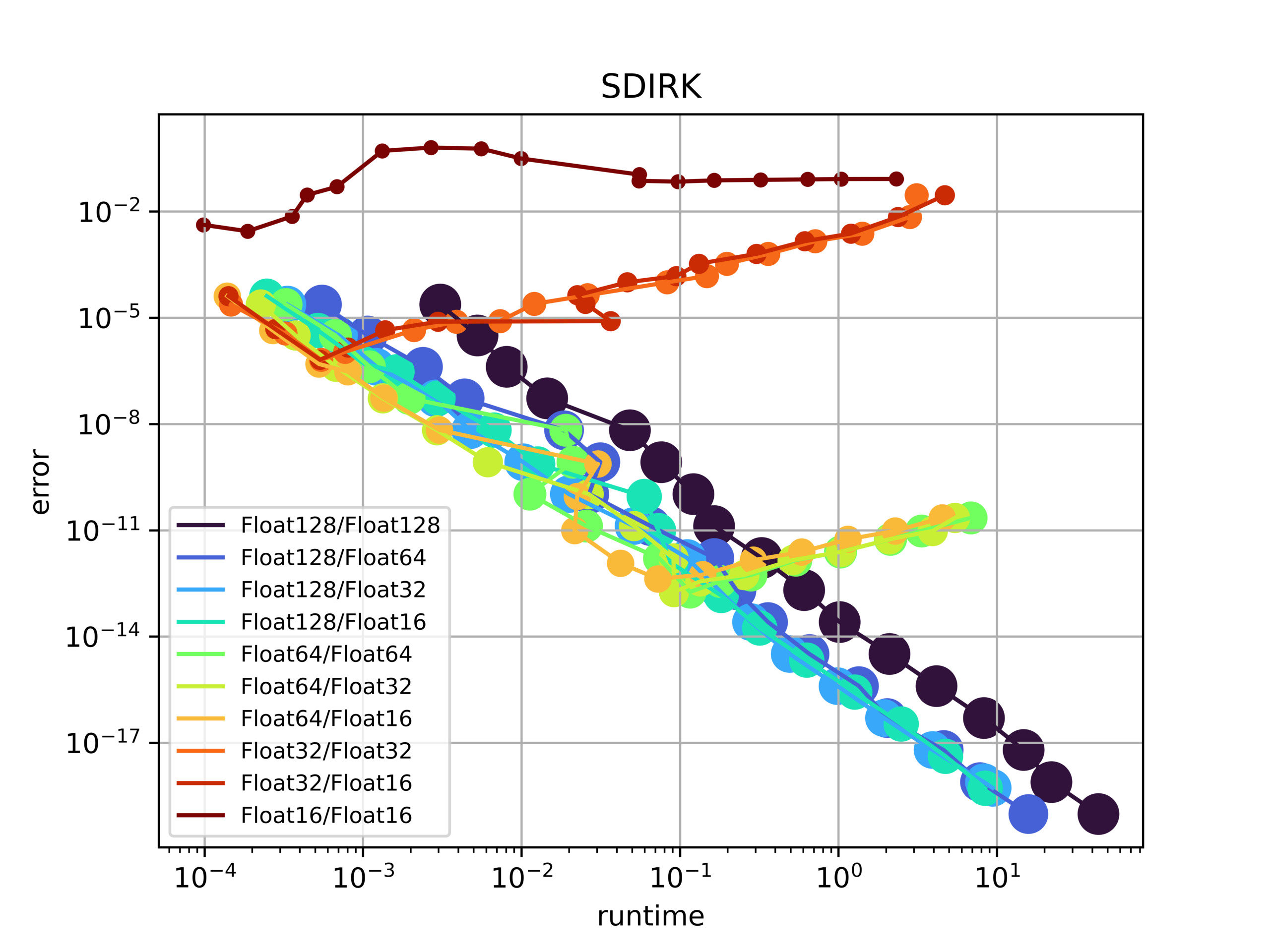}}\\
{\includegraphics[width=0.32\textwidth]{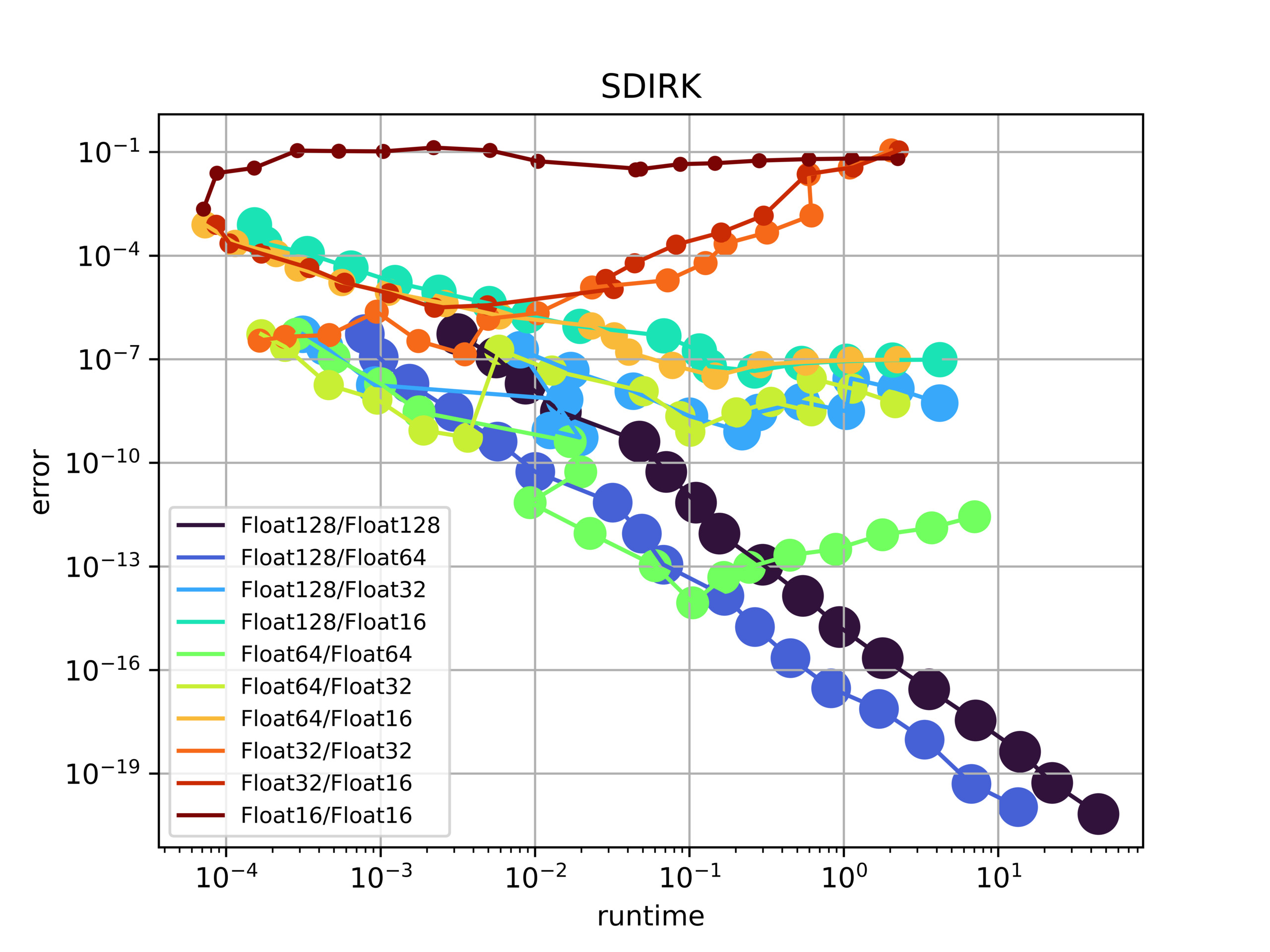}}
{\includegraphics[width=0.32\textwidth]{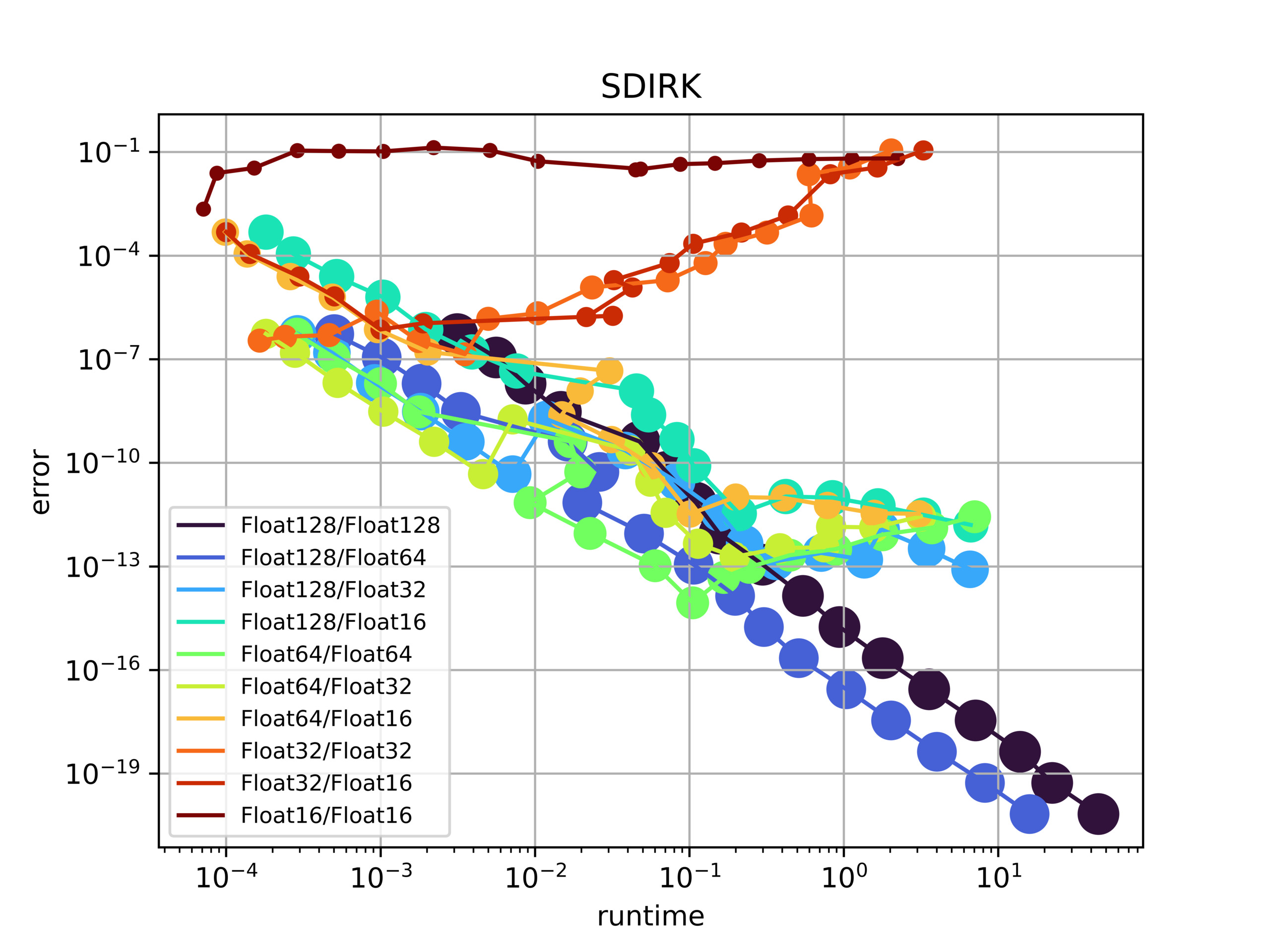}}
{\includegraphics[width=0.32\textwidth]{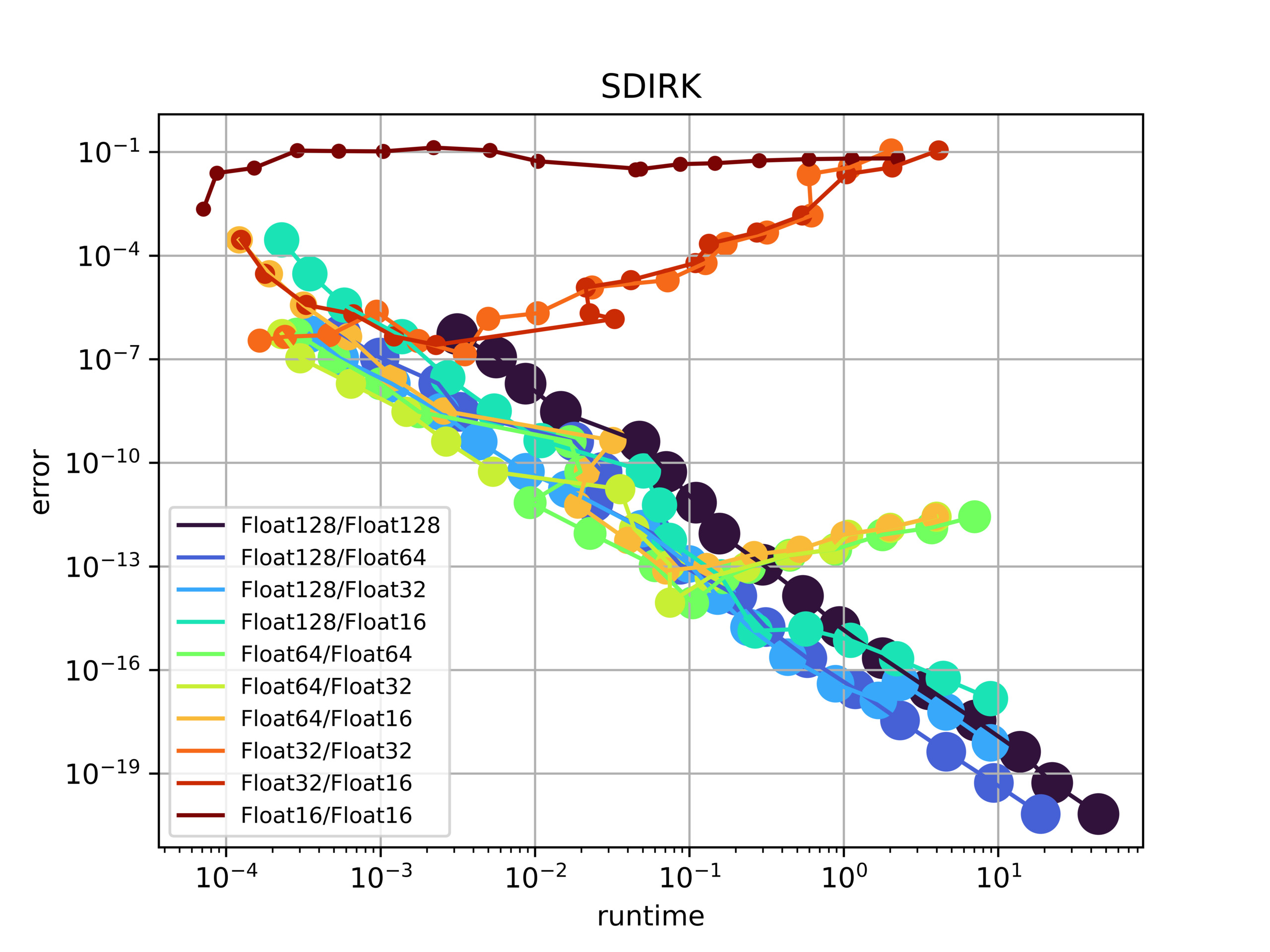}}\\
\caption{Julia code results: Errors vs, runtime for the  SDIRK method for the van der Pol equation with $\alpha=1$ (top) and $\alpha=6$ (bottom)
 with no corrections (left), one correction (middle), two corrections (right).
\label{vdpSDIRKefficiency}
}
\end{center}
\end{figure}

 \noindent {\sc  Julia language  code:} 
 The Julia language codes tell a similar story. Figure \ref{vdpSDIRKefficiency} show the efficiency of the
 SDIRK method for the van der Pol equation with $\alpha=1$ (top) and $\alpha=6$ (bottom)
 with no corrections (left), one correction (middle), two corrections (right).
 The top row (read across from left to right) shows the effect of the successive explicit corrections.
 Notable here, once again are the double/single  (Float64/Float32 in chartreuse),
 double/half  (Float64/Float16 in light orange), and quad/half (Float128/Float16 in blue-green) results, 
 which show the dramatic effect of successive corrections. 
 Notice, too, how one correction fixes the  leveling off for the quad/double code 
 (Float128/Float64 in blue) for $\alpha=1$.
 
 The distinction between the different stiffnesses is expressed most clearly in the 
 quad/half  (Float128/Float16 in blue-green)  and the quad/single (Float128/Float32 in light blue) 
 figures: for the non-stiff case $\alpha=1$, two corrections produce efficiency comparable to the 
  quad/double code (Float128/Float64 in blue) -- note that those three lines are right on top of each other
  for the top right graph. However, for the stiffer problem with  $\alpha=6$, two corrections 
  are not enough to produce comparable efficiency, and the quad/double code remains most efficient
  for very small errors, while the quad/half and quad/single are less so.
  
  The behavior for different stiffnesses is also seen by looking at the 
  double/single  (Float64/Float32 in chartreuse) and double/half  (Float64/Float16 in light orange)
  results. For $\alpha=1$, the  double/half code becomes as efficient as the double/single after two corrections,
  while for $\alpha=6$, the  double/half code does not gain as much in efficiency with two corrections.
  Clearly, the stiffness effect makes it harder for corrections to give the desired efficiency.

\subsubsection{\bf  NovelA Method}

\noindent{\sc Fortran code}
Finally, we present the results from the NovelA  method  \eqref{MP-4s3pA}.
 The table below shows  that for an error level of $10^{-13}$ we see that the double/single code runs $\approx 2.8$ as fast as the 
double precision code for both the x86 and {\sc power9} chips.
\begin{figure}[t] \vspace{-0.75in}
\begin{center}
{\includegraphics[width=0.45\textwidth]{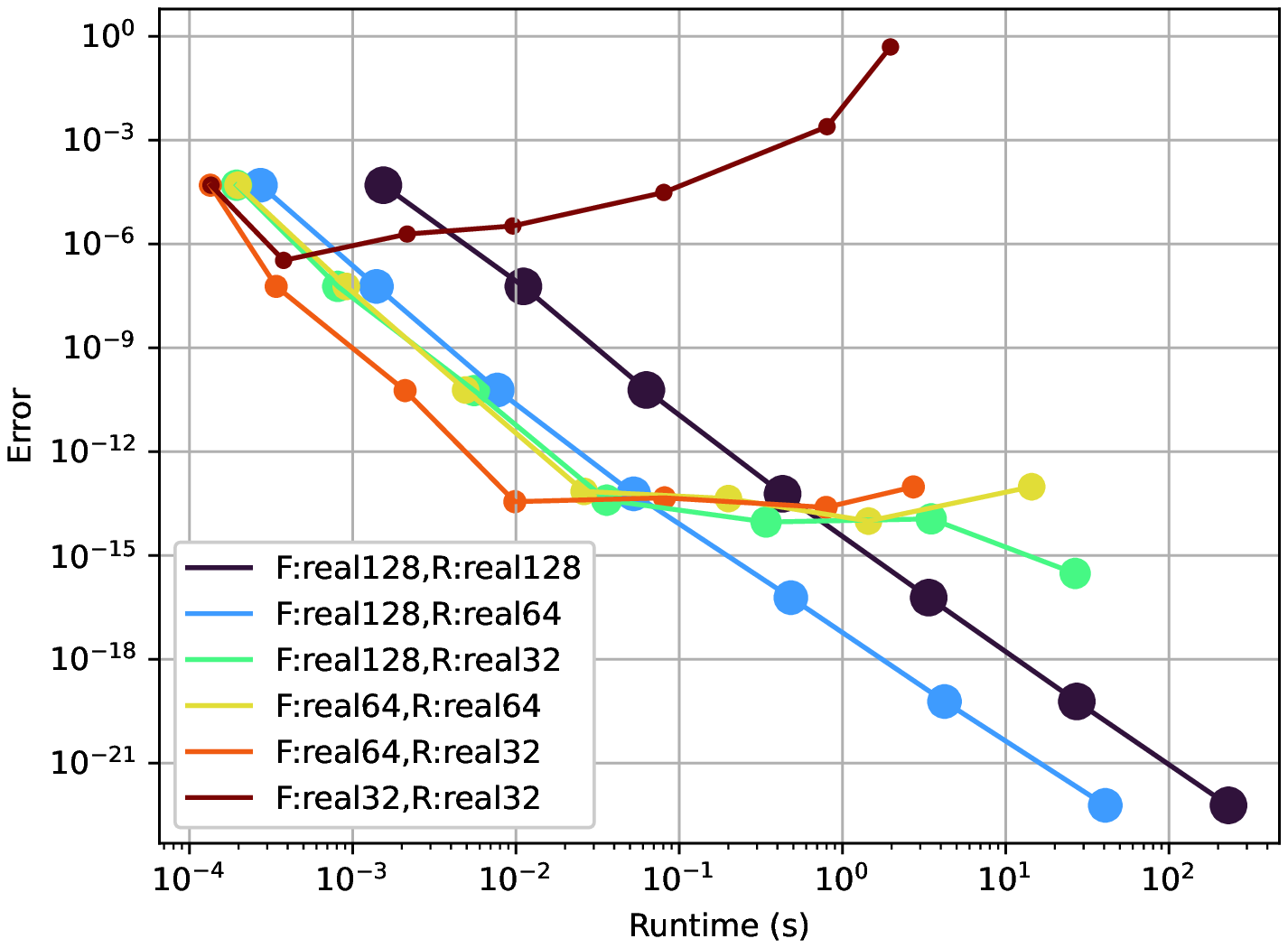}}
{\includegraphics[width=0.45\textwidth]{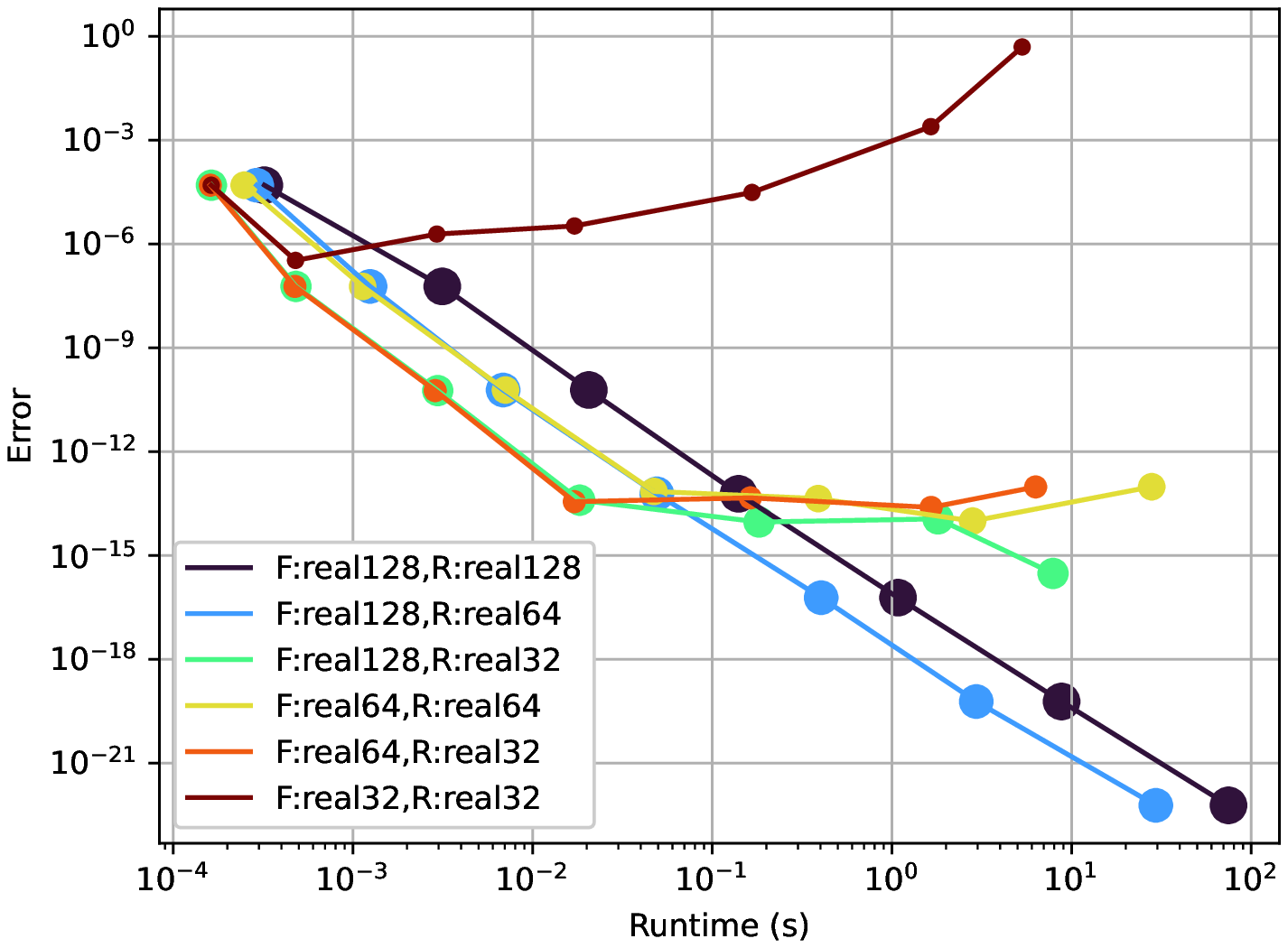}}
\caption{NovelA mixed-precision  method error vs runtime. Left: x86 chip. Right: {\sc power9} chip.
\label{vdpNovelAeffFort}
}
\end{center}
\end{figure}
The quad precision codes are significantly more efficient -- by more than a factor of 3 --
on the  {\sc power9} chips as on the  the x86 chips.
For an error level of $10^{-22}$ the quad/single code does not reach this level of accuracy. 
On the {\sc power9} chip, the  quad/double code is 2.5x as fast as the full quad precision code, 
and on the  x86 chips  the quad/double code is close to 5.7x as fast as the full quad precision code.
In Figure \ref{vdpNovelAeffFort} we see the runtime per level of accuracy for the x86 chip (left)
and power9 chip (right).
We observe that up to $\approx 10^{-13}$, the double/single is the most efficient method,
whereas below that,  the quad/double is the most efficient.
\smallskip

\begin{center}
{\small
\begin{tabular}{ |c|c|c|c|c|c|  }
 \hline
 \multicolumn{6}{|c|}{Mixed precision NovelA method}\\ \hline
  \multicolumn{6}{|c|}{Runtime (s) for x86} \\ \hline
 error & 64/64 & 64/32 & 128/128 & 128/64 
 & 128/32 \\
$\approx 10^{-13}$ & 0.026 & 0.009 & 0.4300 & 0.052 & 0.035\\
$\approx 10^{-22}$ & N/A & N/A & 231.4 & 40.71 & N/A\\ \hline
\multicolumn{6}{|c|}{Runtime (s) for POWER9} \\ \hline
 error & 64/64 & 64/32 & 128/128 & 128/64 
 & 128/32 \\
$\approx 10^{-13}$ & 0.047 & 0.017 & 0.140 & 0.049 & 0.018\\
$\approx 10^{-22}$ & N/A & N/A & 74.68 & 29.42 & N/A\\ \hline
 \end{tabular}
 }
 \end{center}
 
 \smallskip

 \begin{figure}[h]
\begin{center}
{\includegraphics[width=0.45\textwidth]{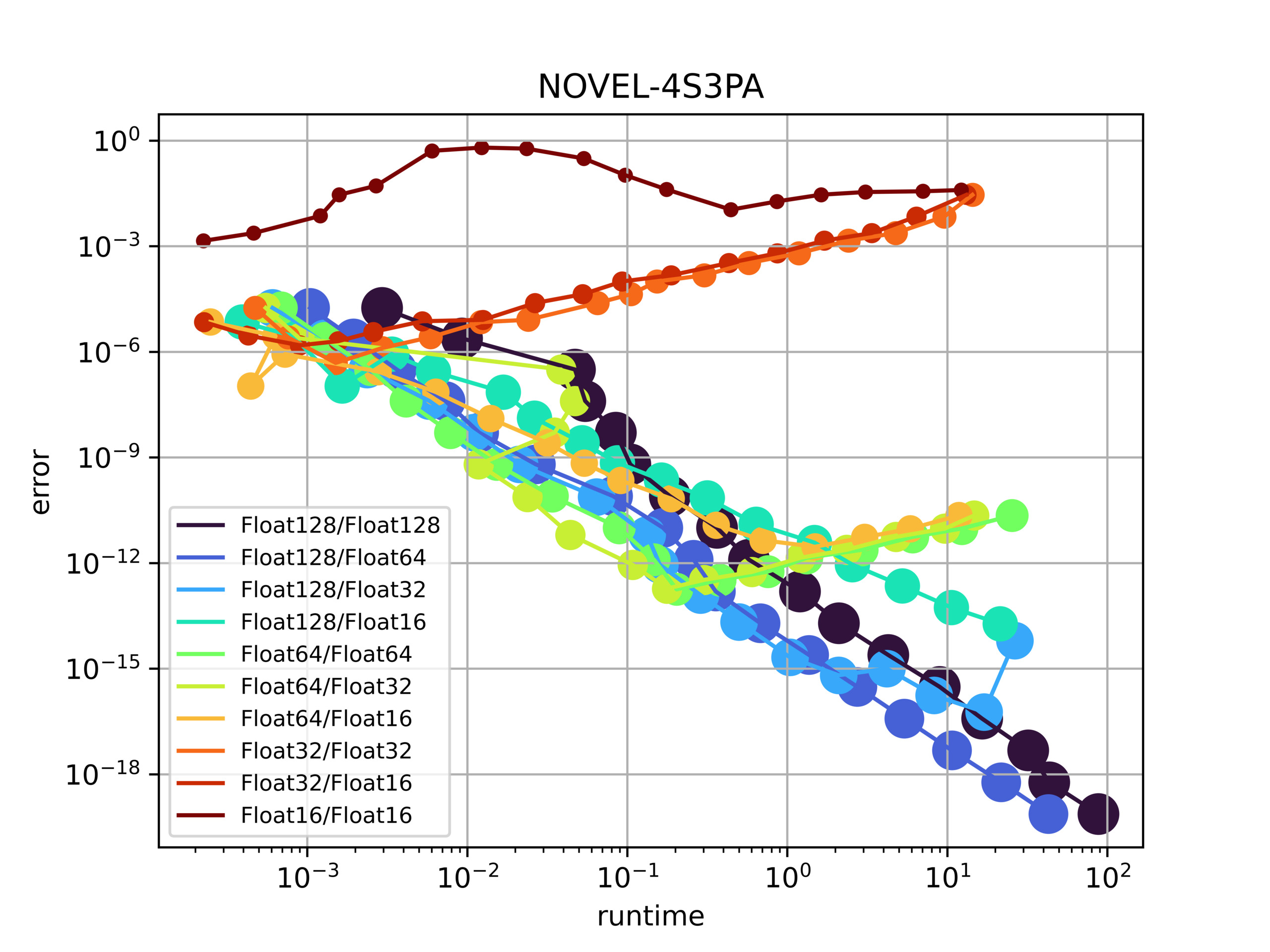}}
{\includegraphics[width=0.45\textwidth]{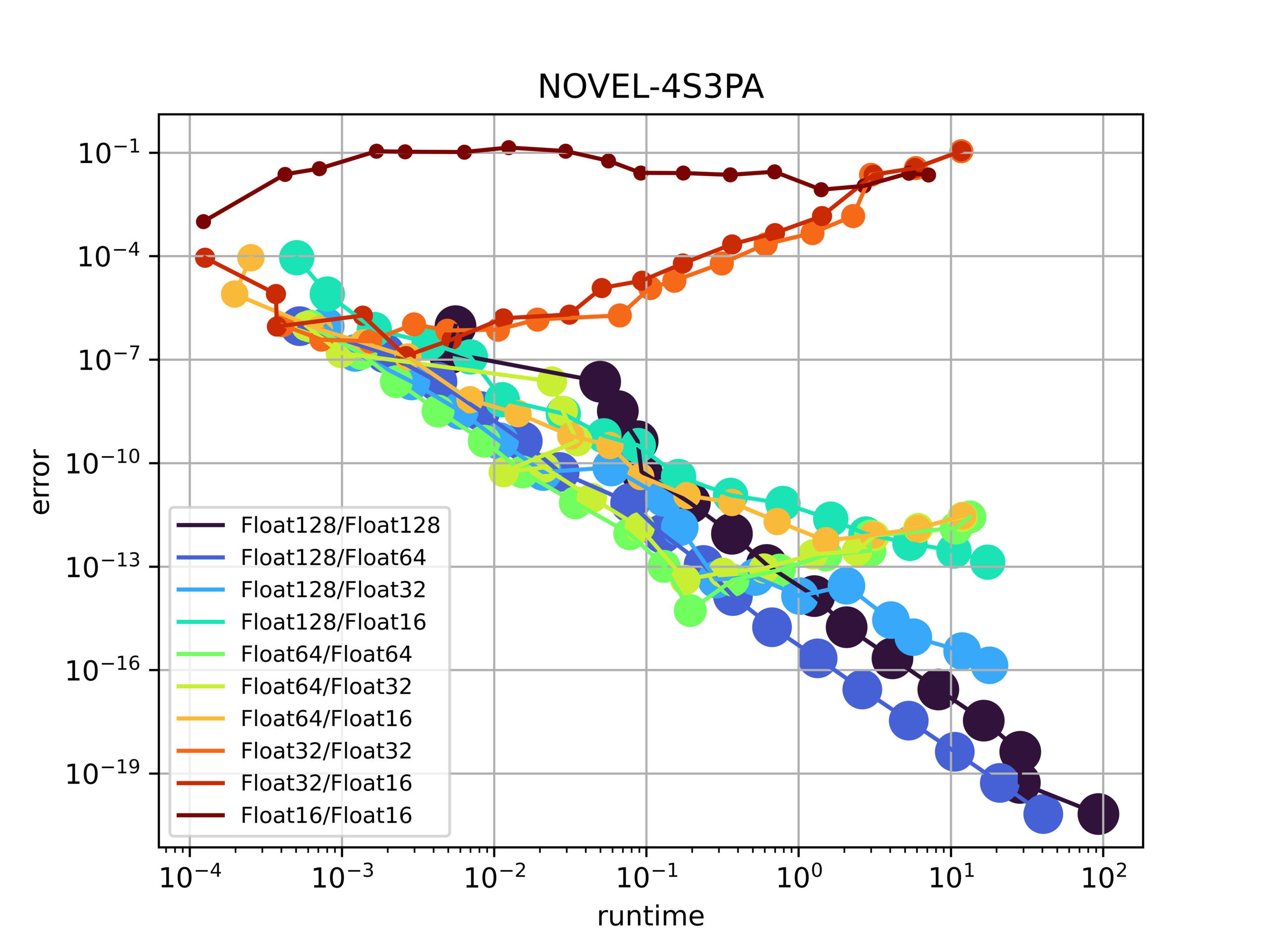}}
\caption{Julia code results: the NovelA method for the van der Pol equation with $\alpha=1$ (left) and $\alpha=6$ (right).
\label{vdpNovelAefficiency}
}
\end{center}
\end{figure}

 \noindent {\sc  Julia language  code:} 
 The results presented in Figure \ref{vdpNovelAefficiency} show that the mixed precision 
 quad/double code provides efficiency over the full precision quad code for both the van der Pol problem 
 with $\alpha =1$ and with $\alpha=6$. The double/single code (Float64/Float32 in chartreuse) 
 is most efficient for error levels between $10^{-9}$ and $10^{-13}$  for the non-stiff problem ($\alpha =1$), 
 but not for the  stiffer problem ($\alpha=6$).  Similarly, the quad/single code (Float128/Float32 in light blue)
 is most efficient for error levels between  $10^{-13}$ and $10^{-15}$  for the non-stiff problem,
 but not for the  stiffer problem. Finally, the mixed/half precision codes do not perform well here: the
 error of $\epsilon \dt^2$ is evident when $\epsilon$ is at half precision level.

\vspace{-.1in}
\subsection{Viscous Burgers' Equation}
The small-scale  example above allows us to see the savings in the mixed precision approach, but to truly
realize these savings, it is best to look at a larger scale nonlinear system resulting from the semi-discretization
of a partial differential equation.  We consider the viscous Burgers' equation presented in Section \ref{sec:convergence}.
In the following, we compare the results of the mixed precision implicit midpoint rule,  SDIRK method, and  NovelA method
for this problem.

\begin{figure}[t!]
\begin{center}
{\includegraphics[width=0.32\textwidth]{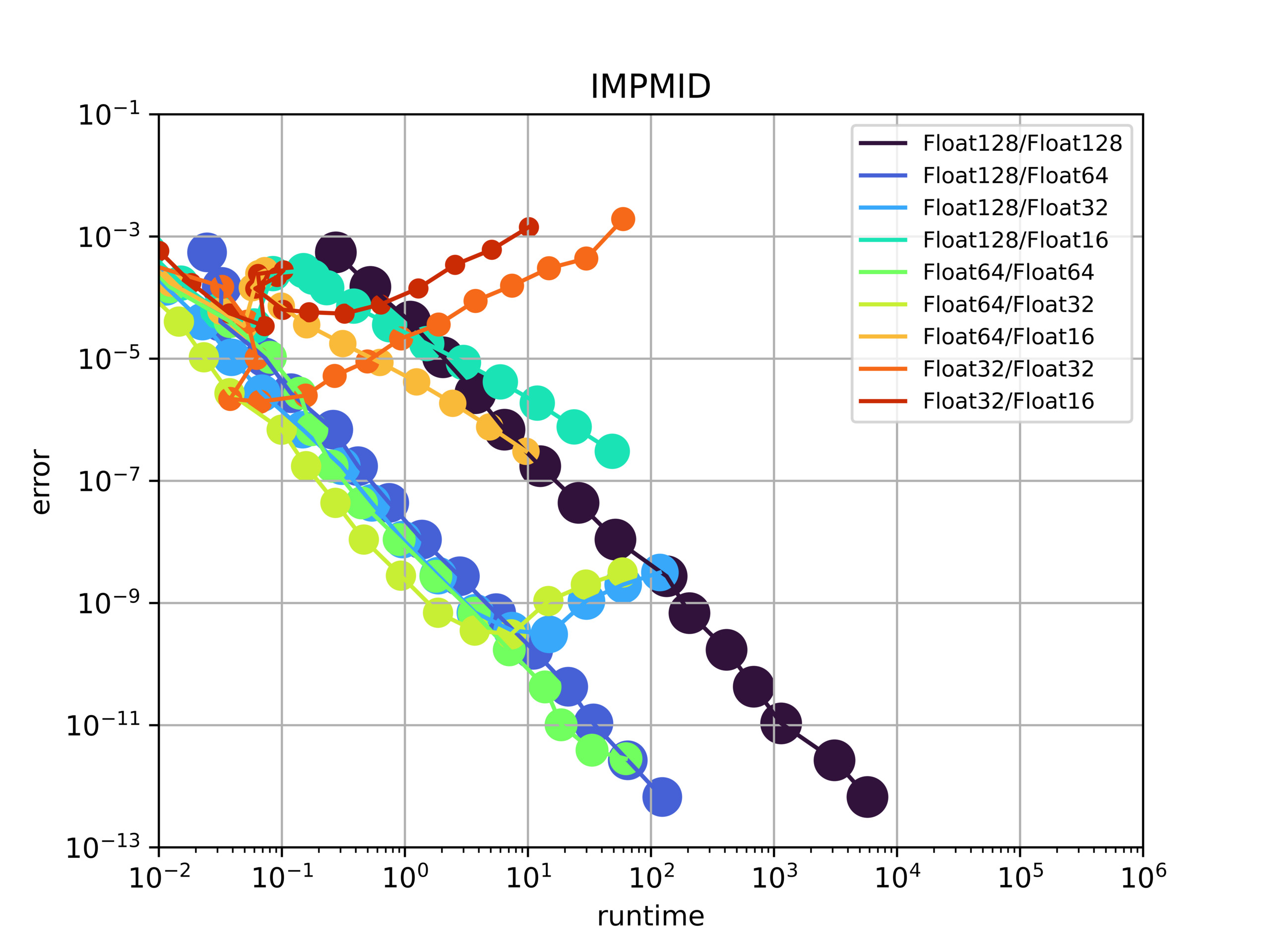}}
{\includegraphics[width=0.32\textwidth]{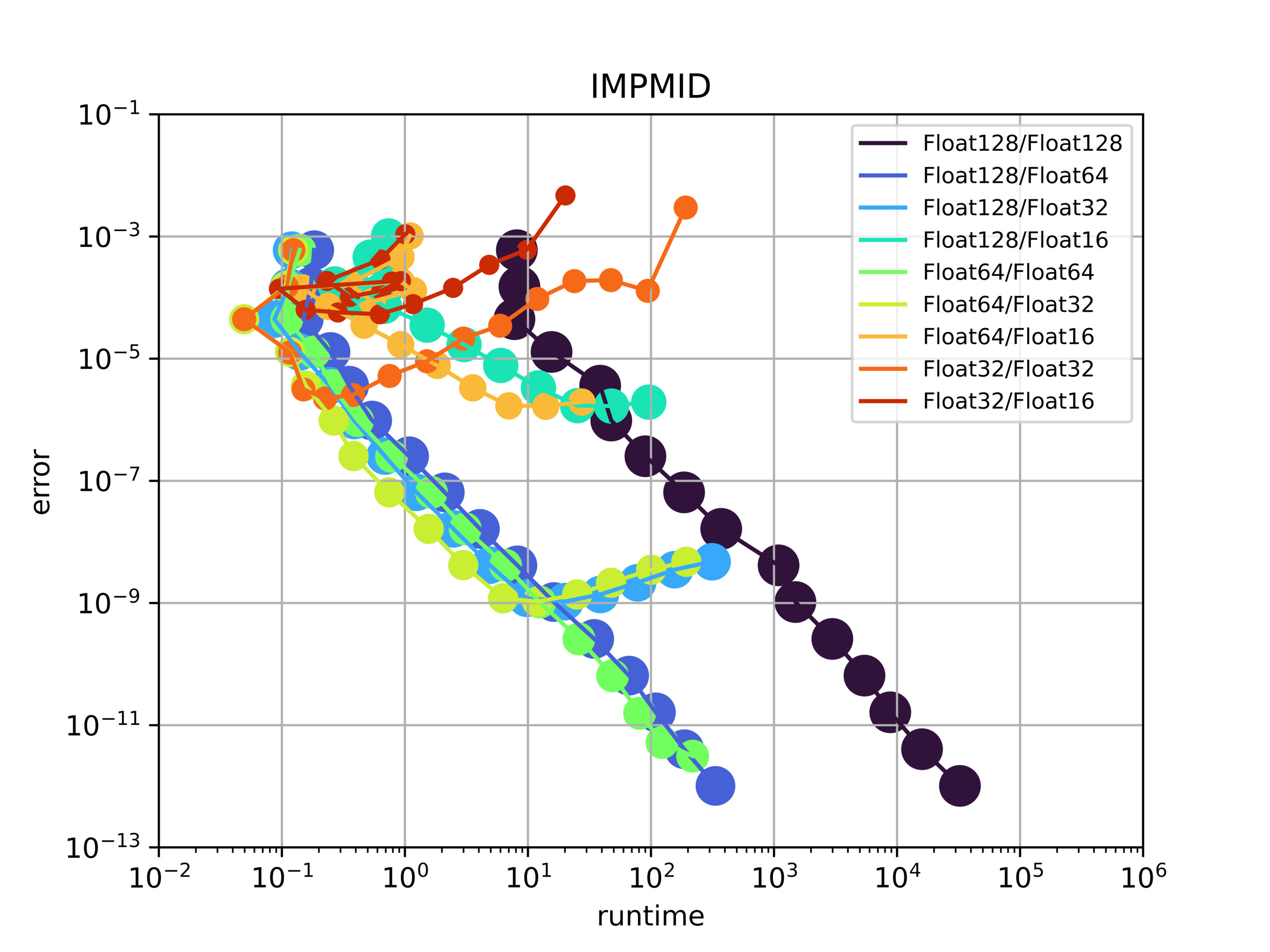}}
{\includegraphics[width=0.32\textwidth]{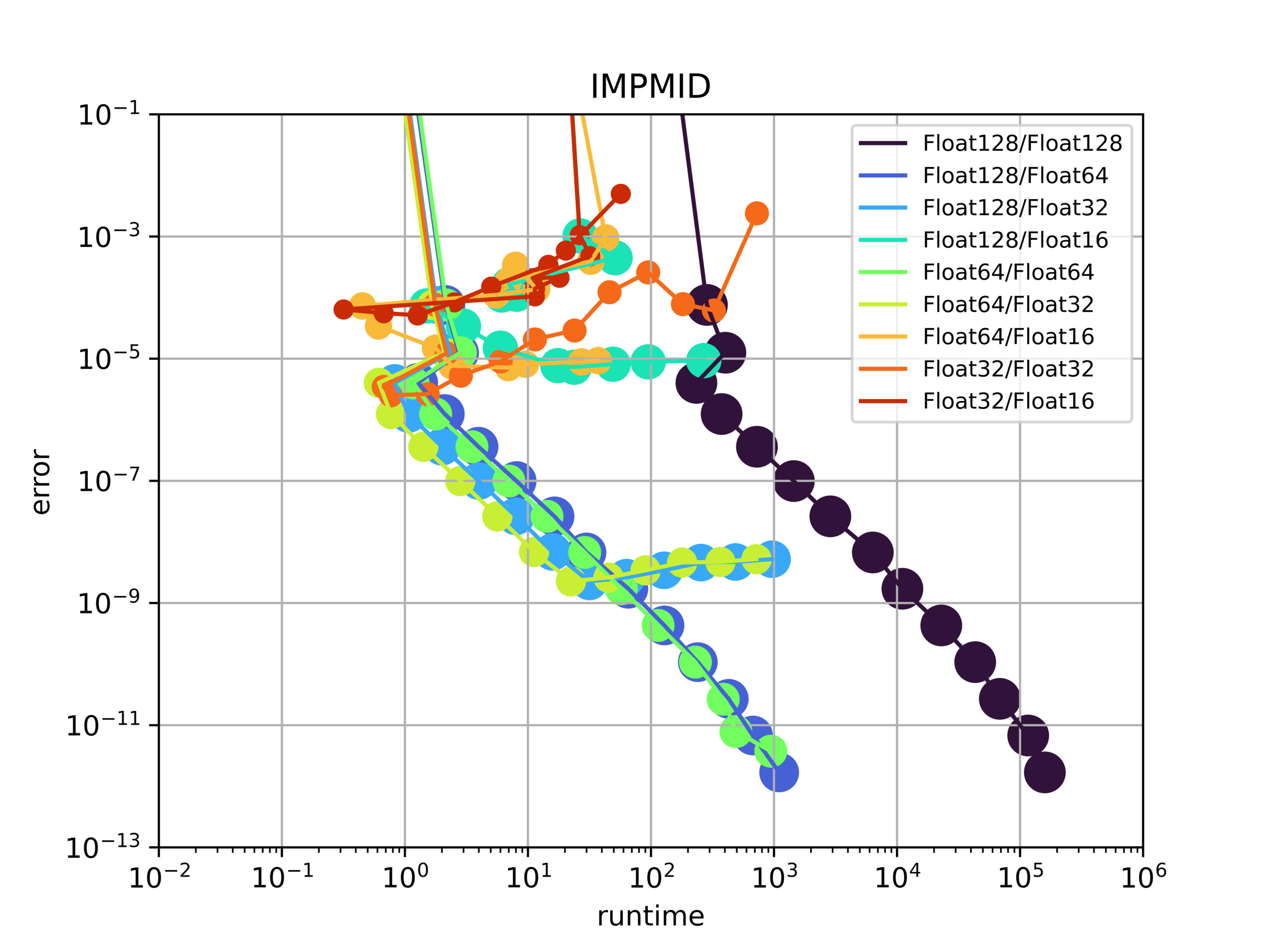}}\\
{\includegraphics[width=0.32\textwidth]{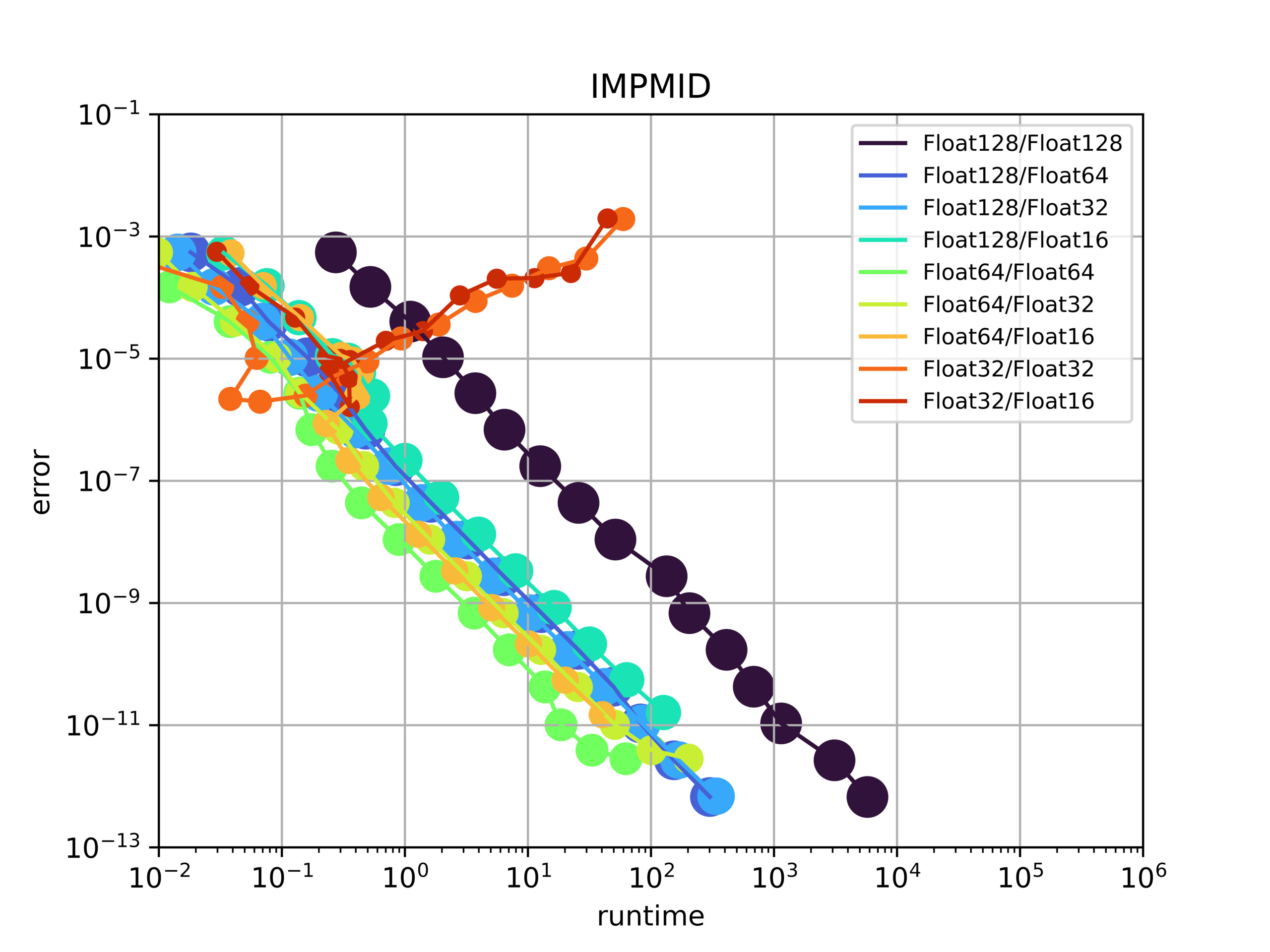}}
{\includegraphics[width=0.32\textwidth]{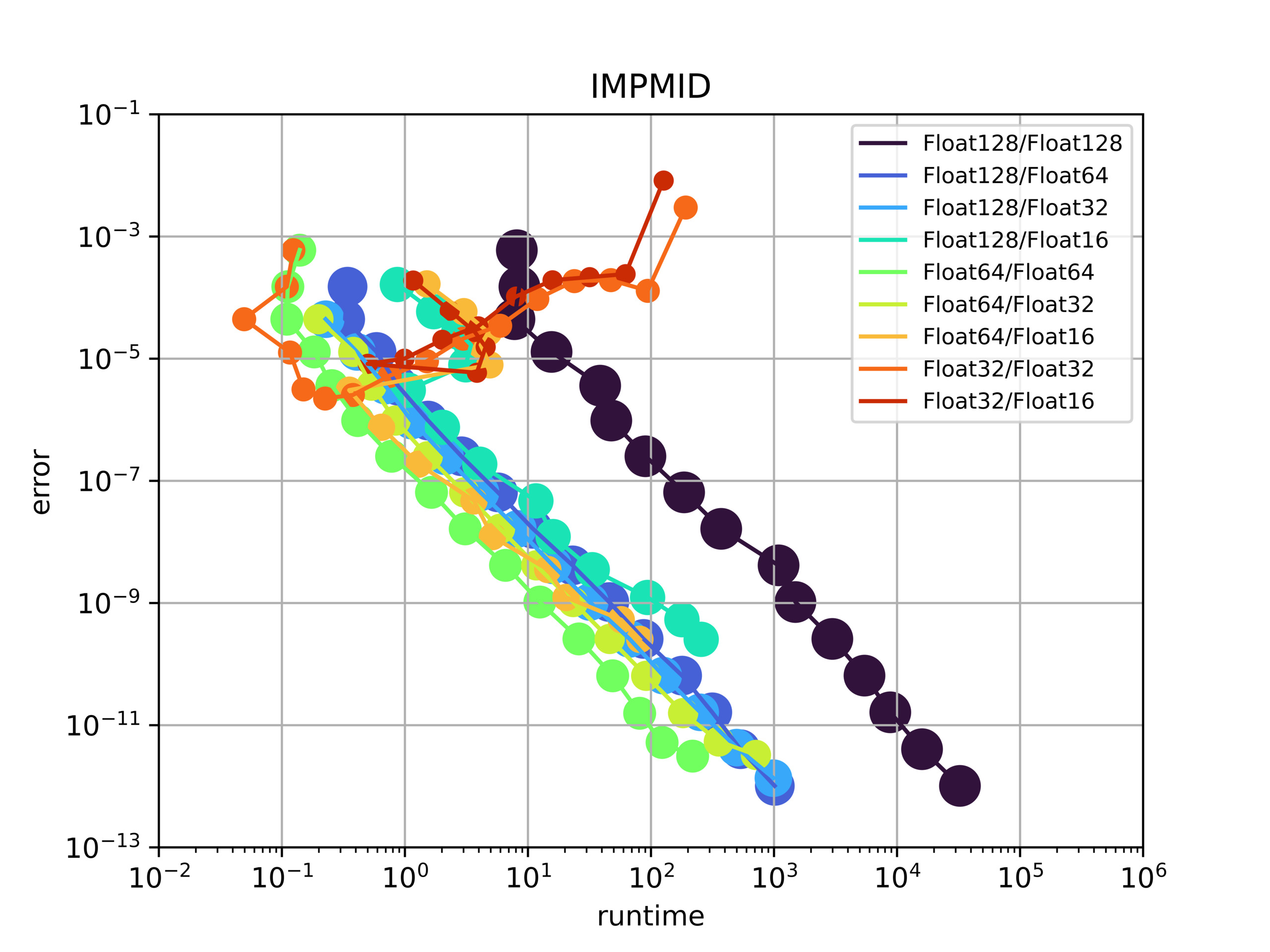}}
{\includegraphics[width=0.32\textwidth]{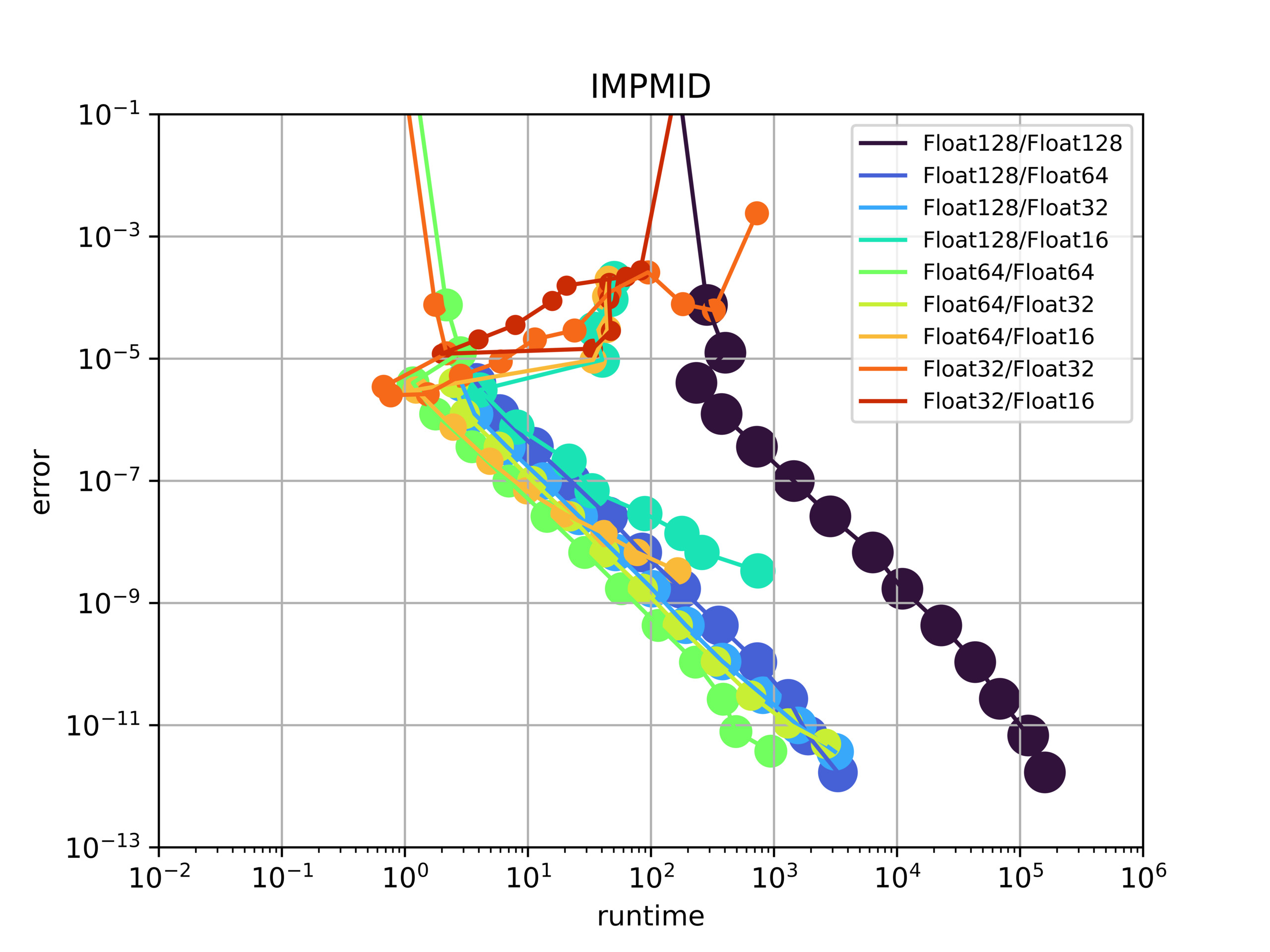}}\\
{\includegraphics[width=0.32\textwidth]{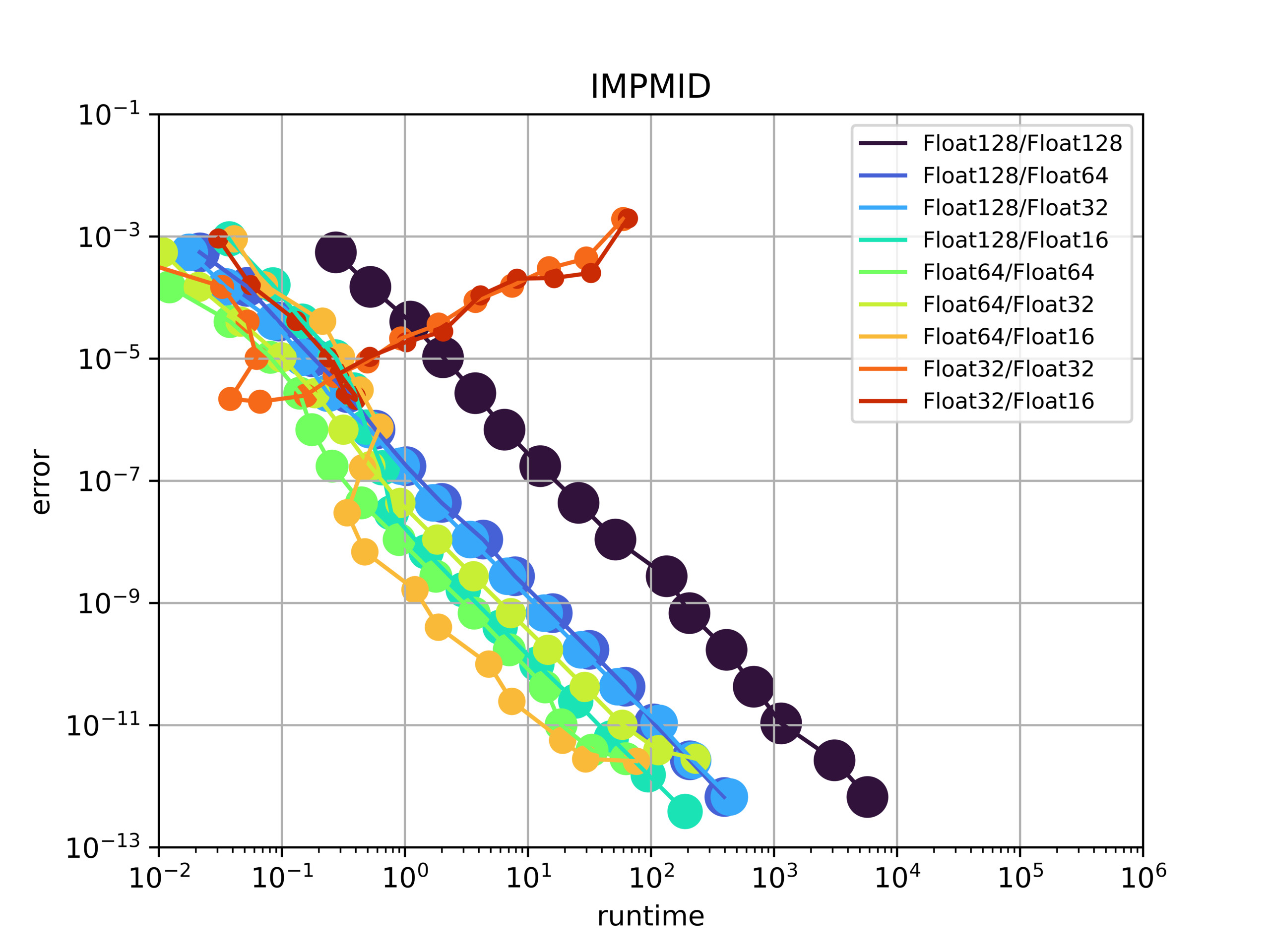}}
{\includegraphics[width=0.32\textwidth]{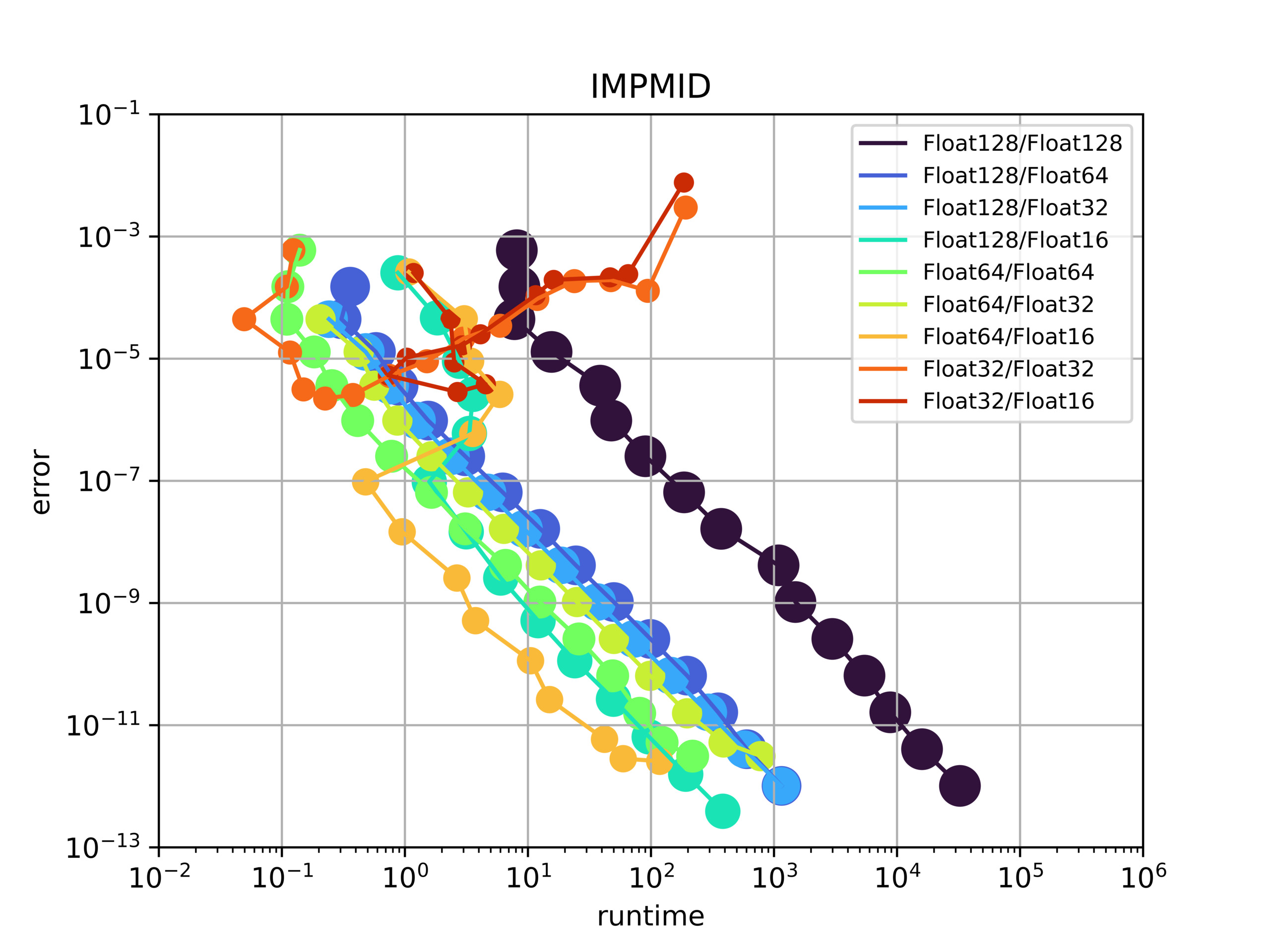}}
{\includegraphics[width=0.32\textwidth]{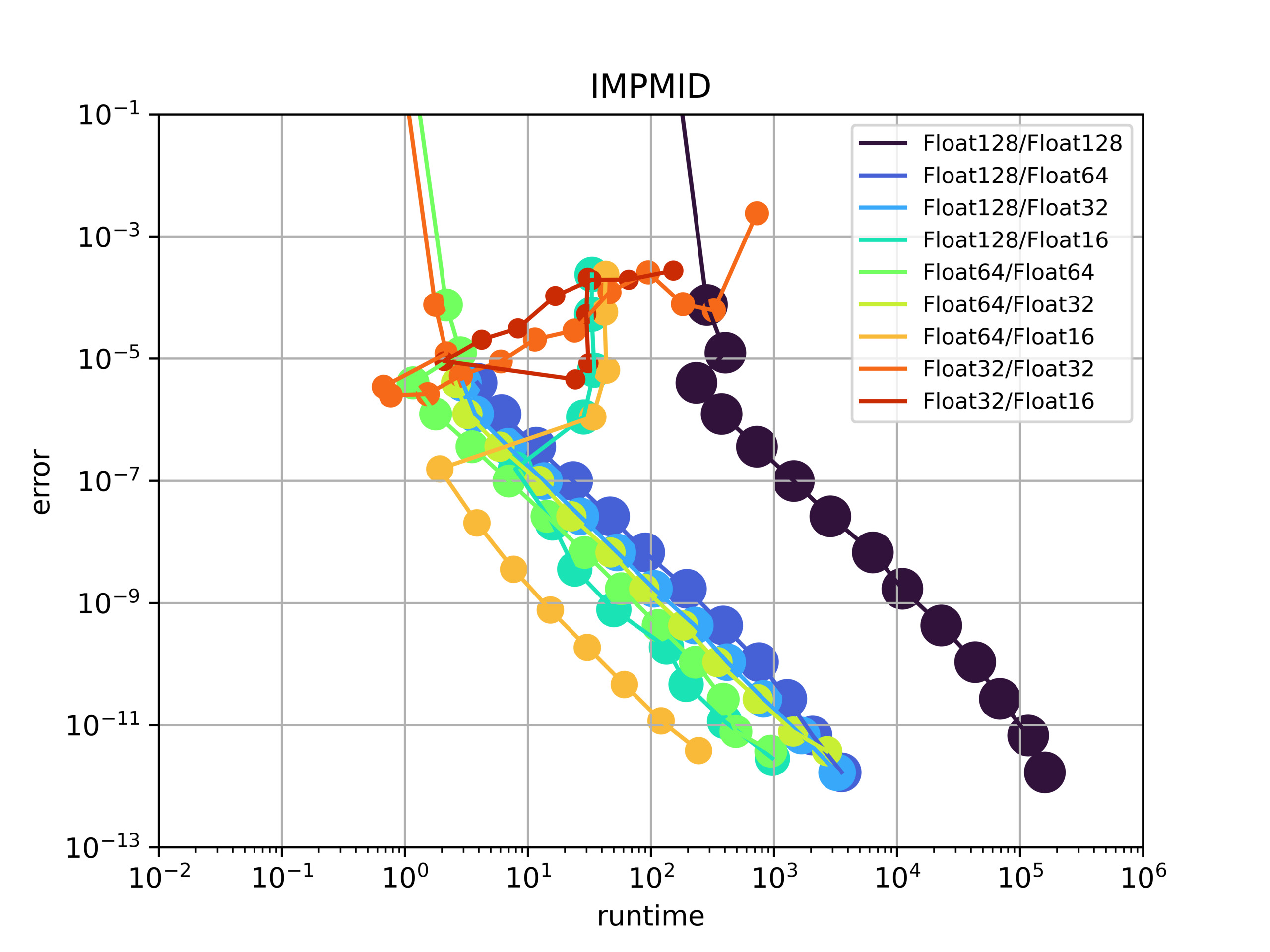}}\\
\caption{Julia code results:  error vs. runtime for the mixed precision 
implicit midpoint rule for the viscous Burgers' equation with $N_x=50$ (left), $N_x=100$ (middle), and $N_x=200$ (right), with no corrections (top row), one correction (middle row), two corrections (bottom row).
\label{BurgersIMR}
}
\end{center}
\end{figure}

\subsubsection{Implicit midpoint rule}
Figure \ref{BurgersIMR} shows the efficiency of the mixed precision 
implicit midpoint rule for the viscous Burgers' equation with $N_x=50$ (left), $N_x=100$ (middle), and $N_x=200$ (right), with no corrections (top row), one correction (middle row), two corrections (bottom row).
We observe that corrections (top to bottom) provide significant gains in efficiency for each problem. Following the behavior of the 
mixed precision double/half code (Float64/Float16 in light orange) shows how each correction improves the efficiency of this code
for all values of $N_x$. For error level between $10^{-7}$ and $10^{-11}$ this method is most efficient.
This code is only second order ($O(\dt^2)$) accurate, so we do not see the accuracy benefits of quad precision. However,
we do clearly observe the relative costs of the mixed methods to the quad  method and can see that 
the mixed precision codes provide close to one order of magnitude savings for $N_x=50$, and significantly more than
one order of magnitude savings for $N_x=100$ and $N_x=200$.

\subsubsection{SDIRK method}
\begin{figure}[htb]
\begin{center}
{\includegraphics[width=0.32\textwidth]{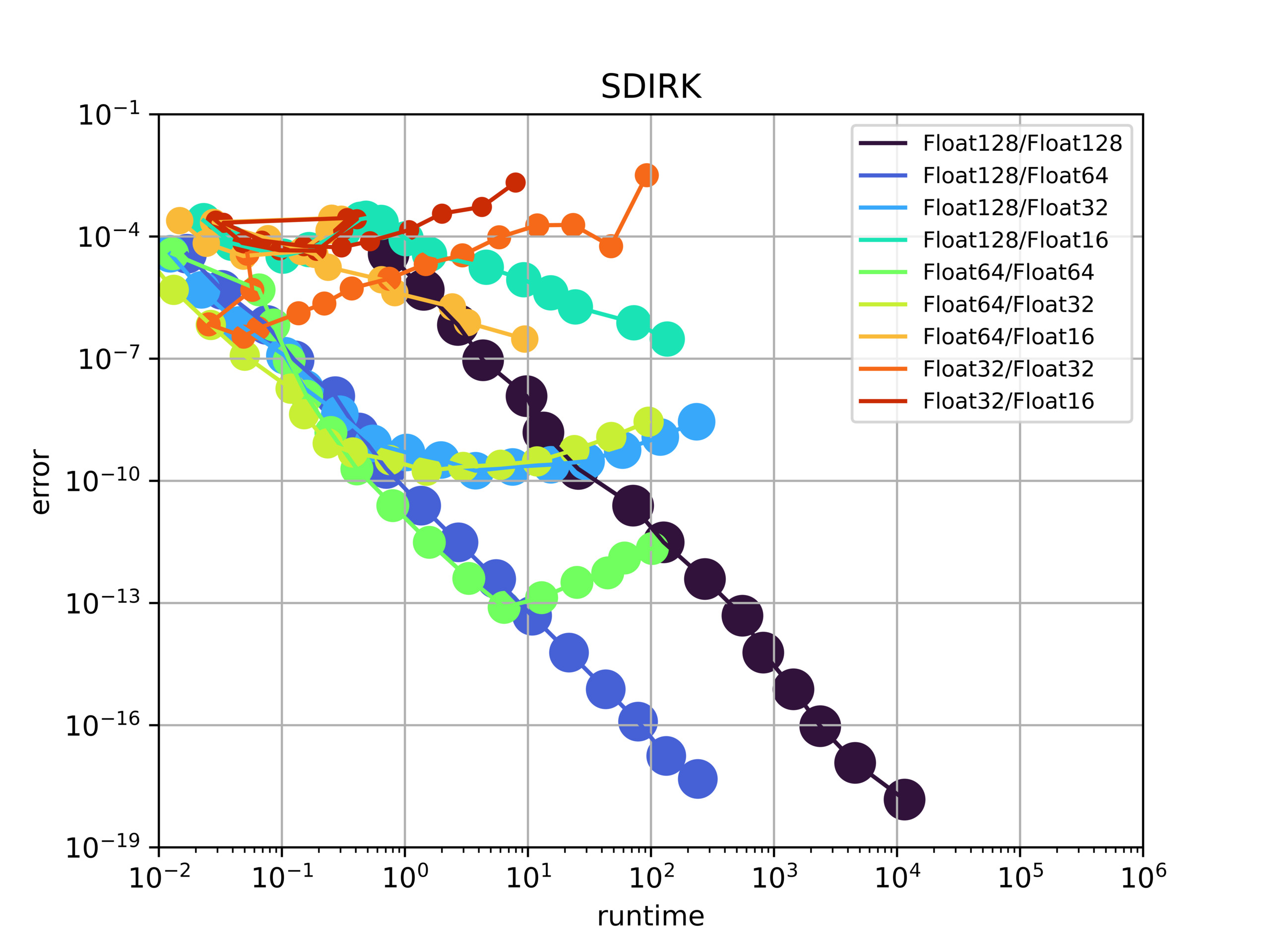}}
{\includegraphics[width=0.32\textwidth]{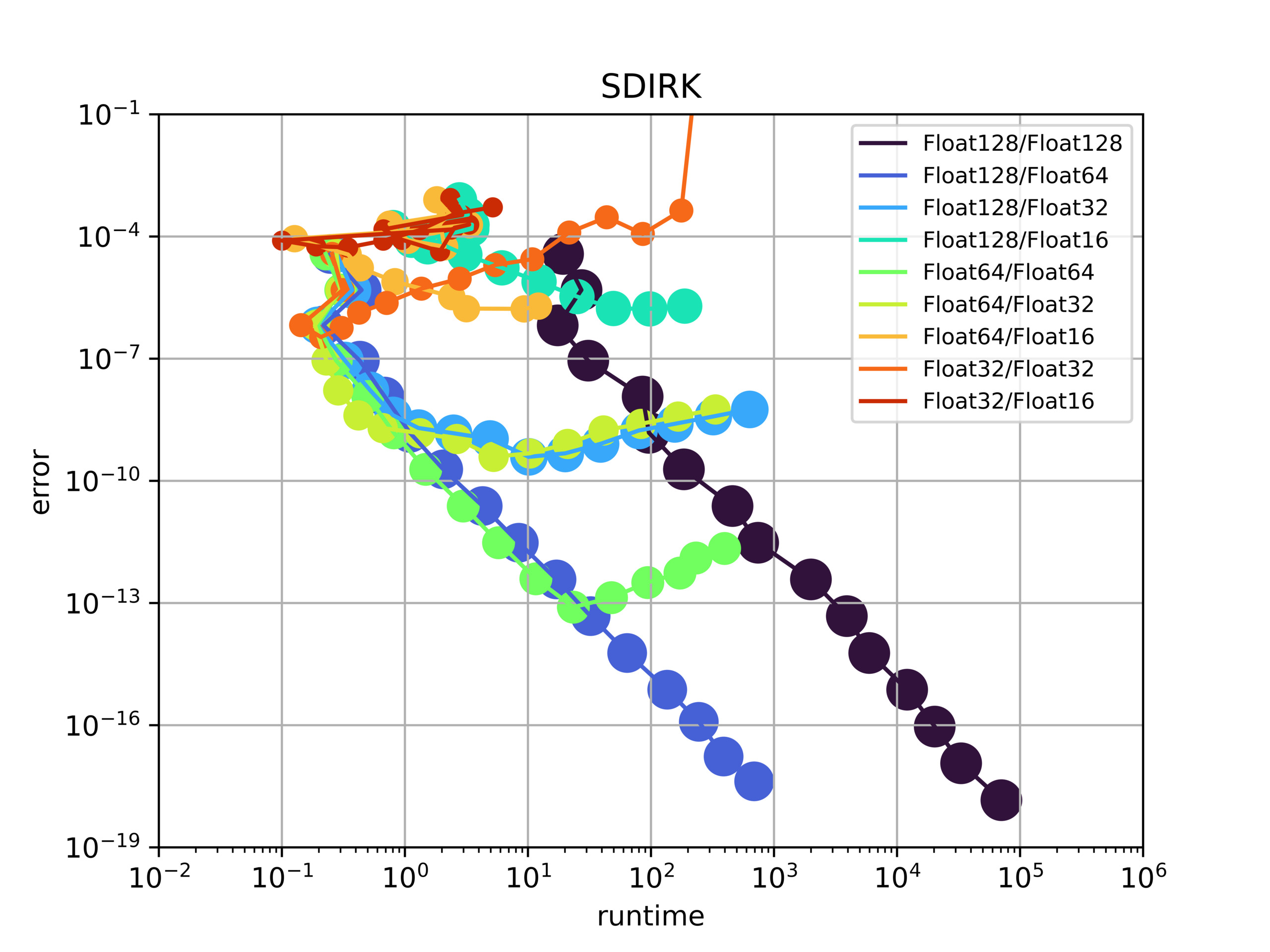}}
{\includegraphics[width=0.32\textwidth]{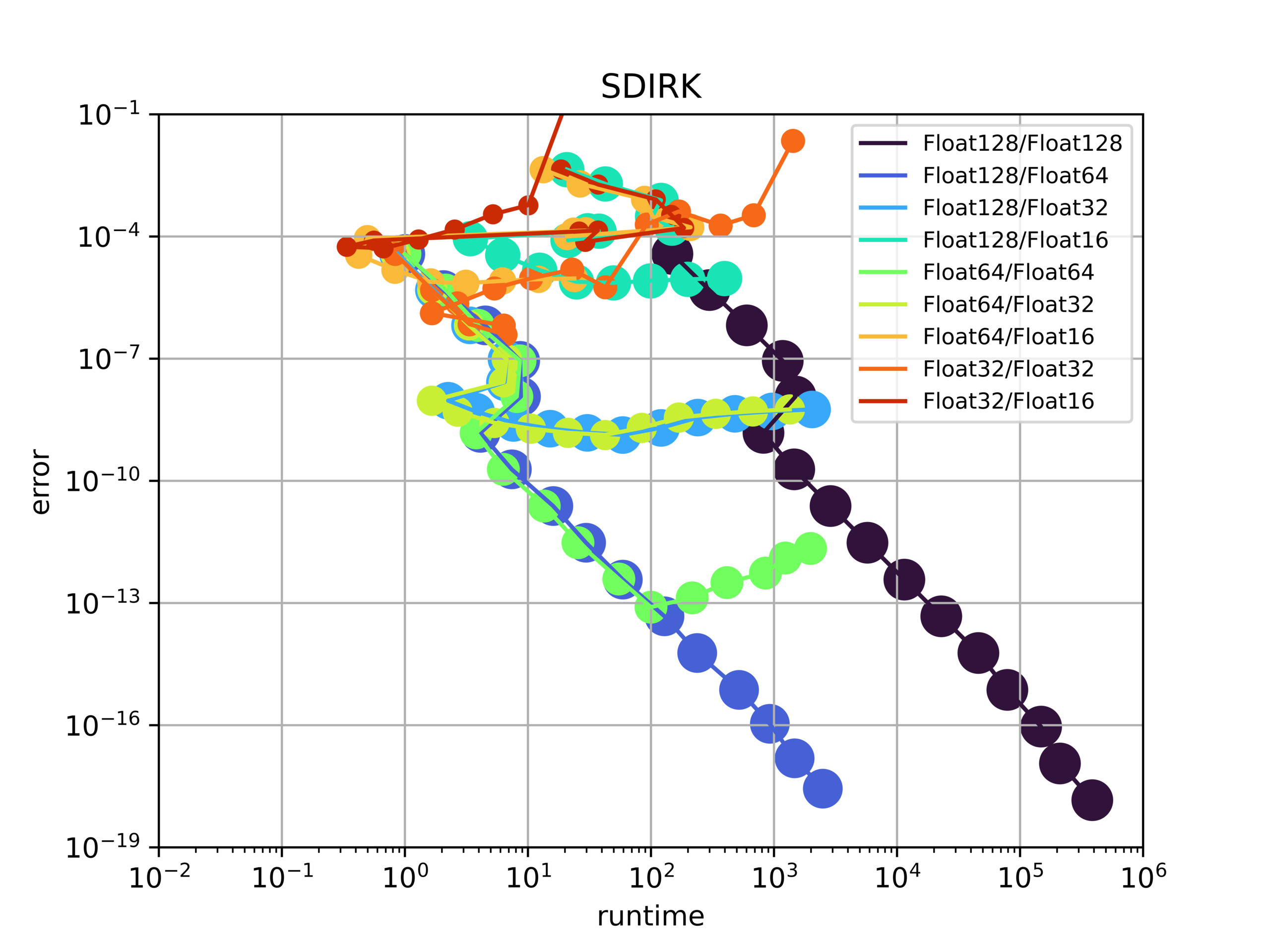}}\\
{\includegraphics[width=0.32\textwidth]{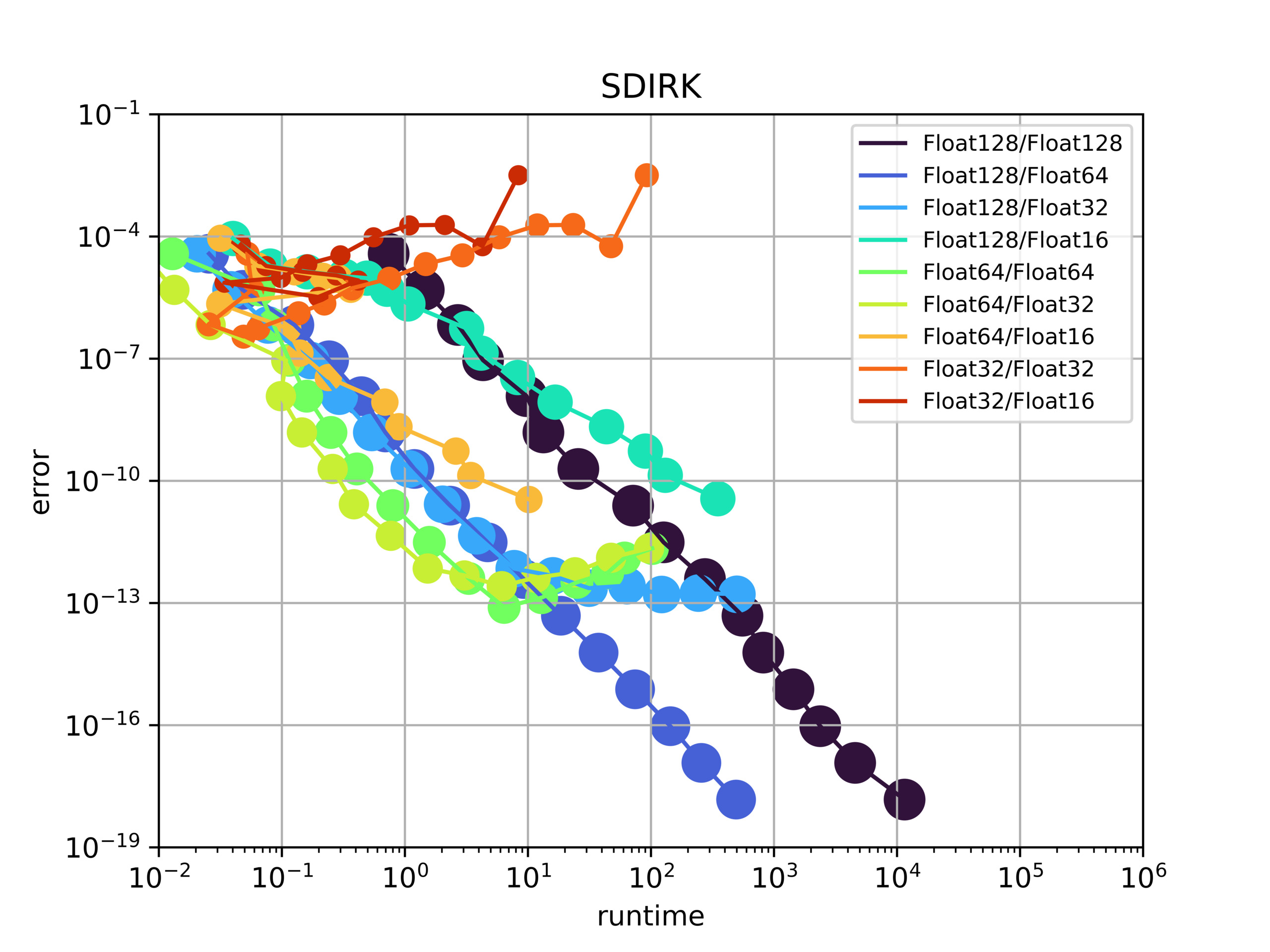}}
{\includegraphics[width=0.32\textwidth]{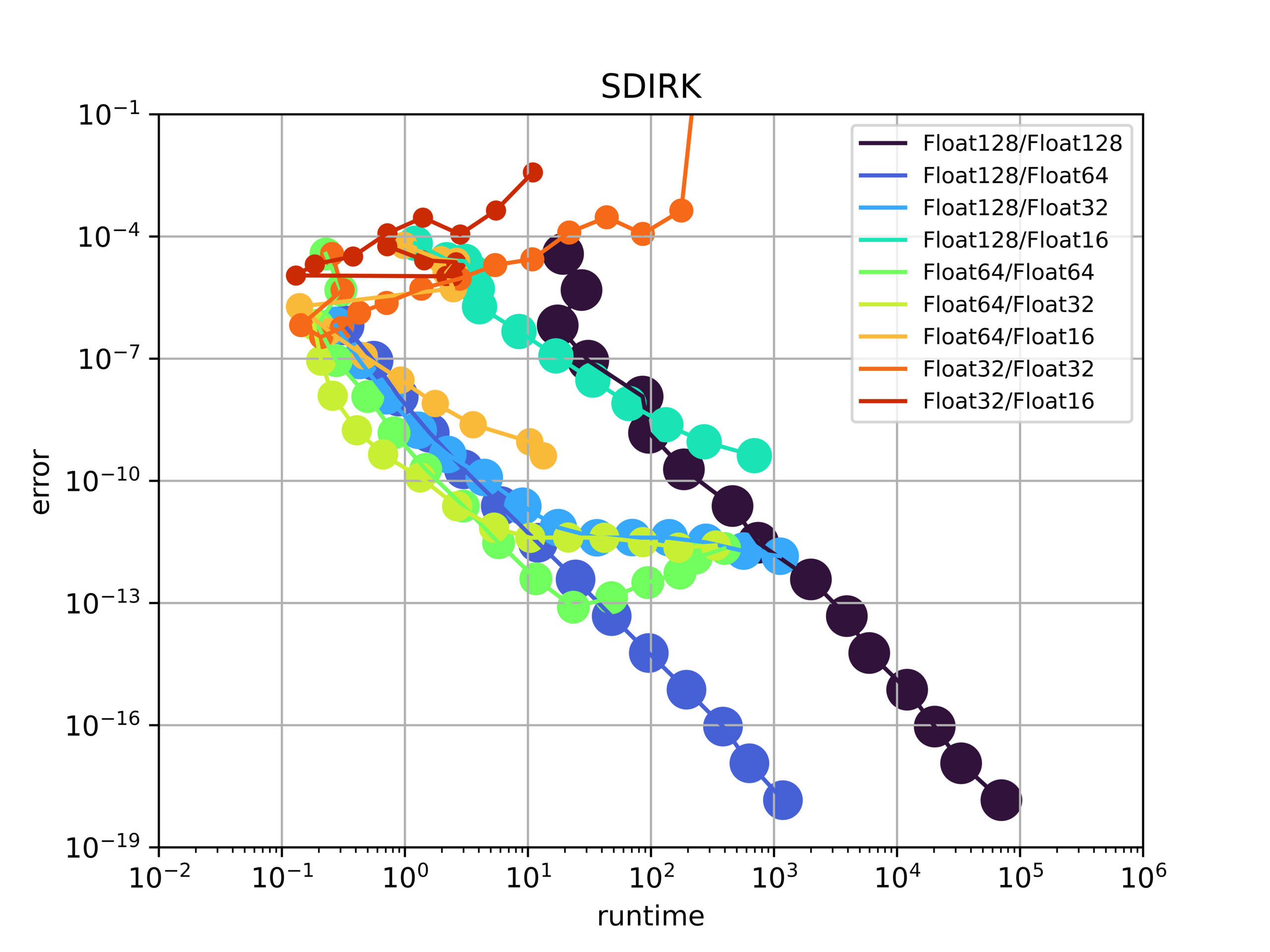}}
{\includegraphics[width=0.32\textwidth]{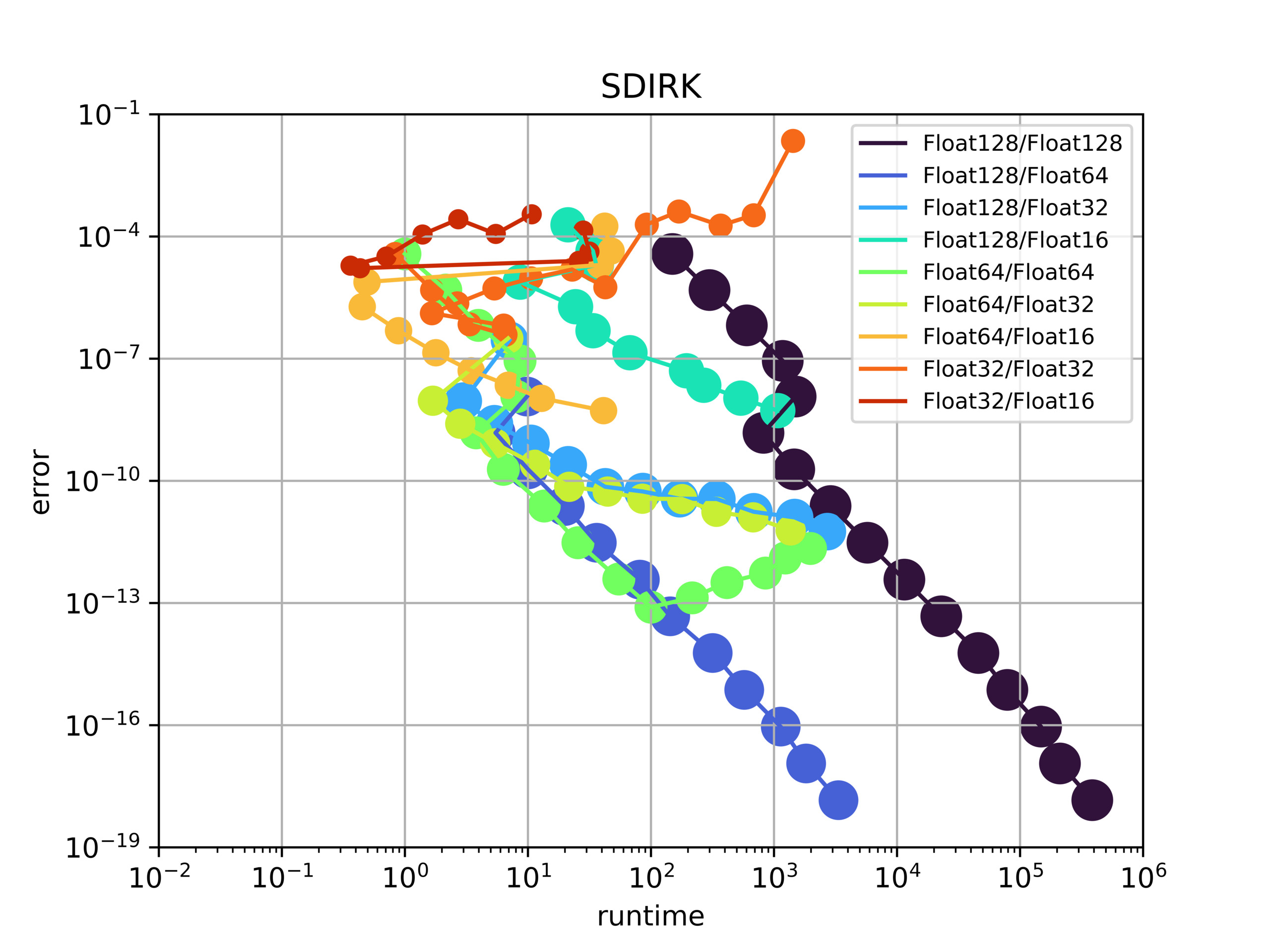}}\\
{\includegraphics[width=0.32\textwidth]{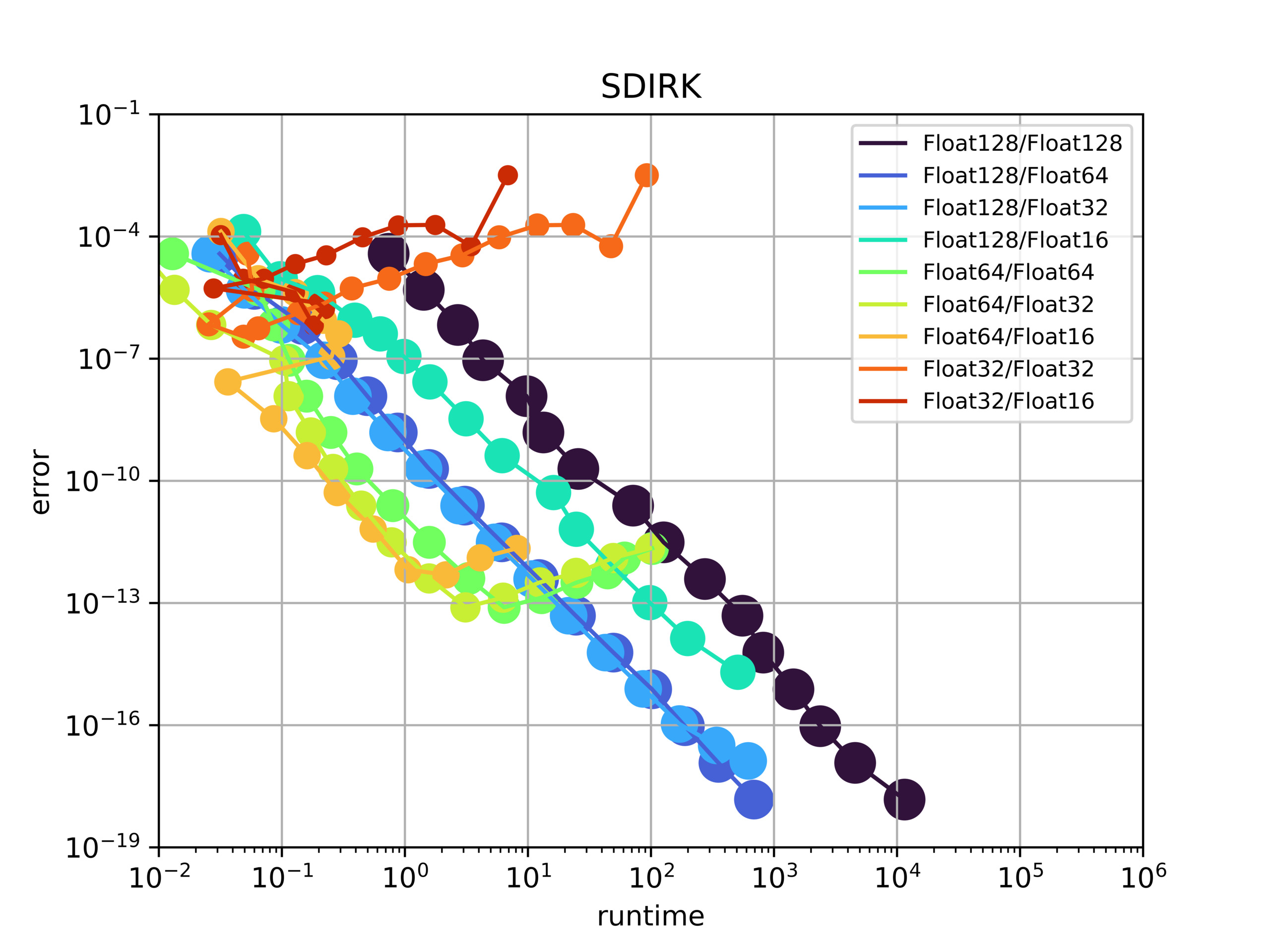}}
{\includegraphics[width=0.32\textwidth]{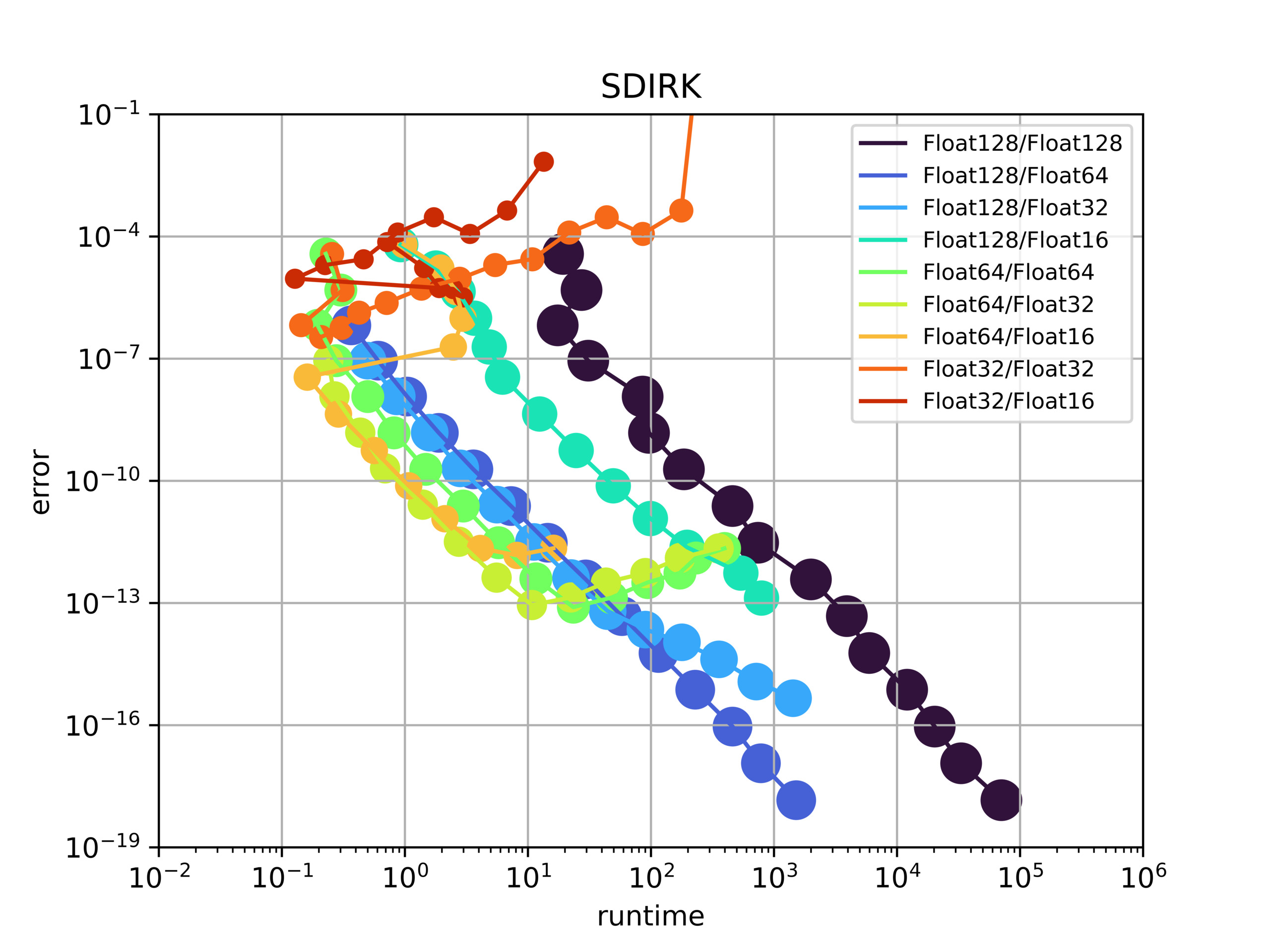}}
{\includegraphics[width=0.32\textwidth]{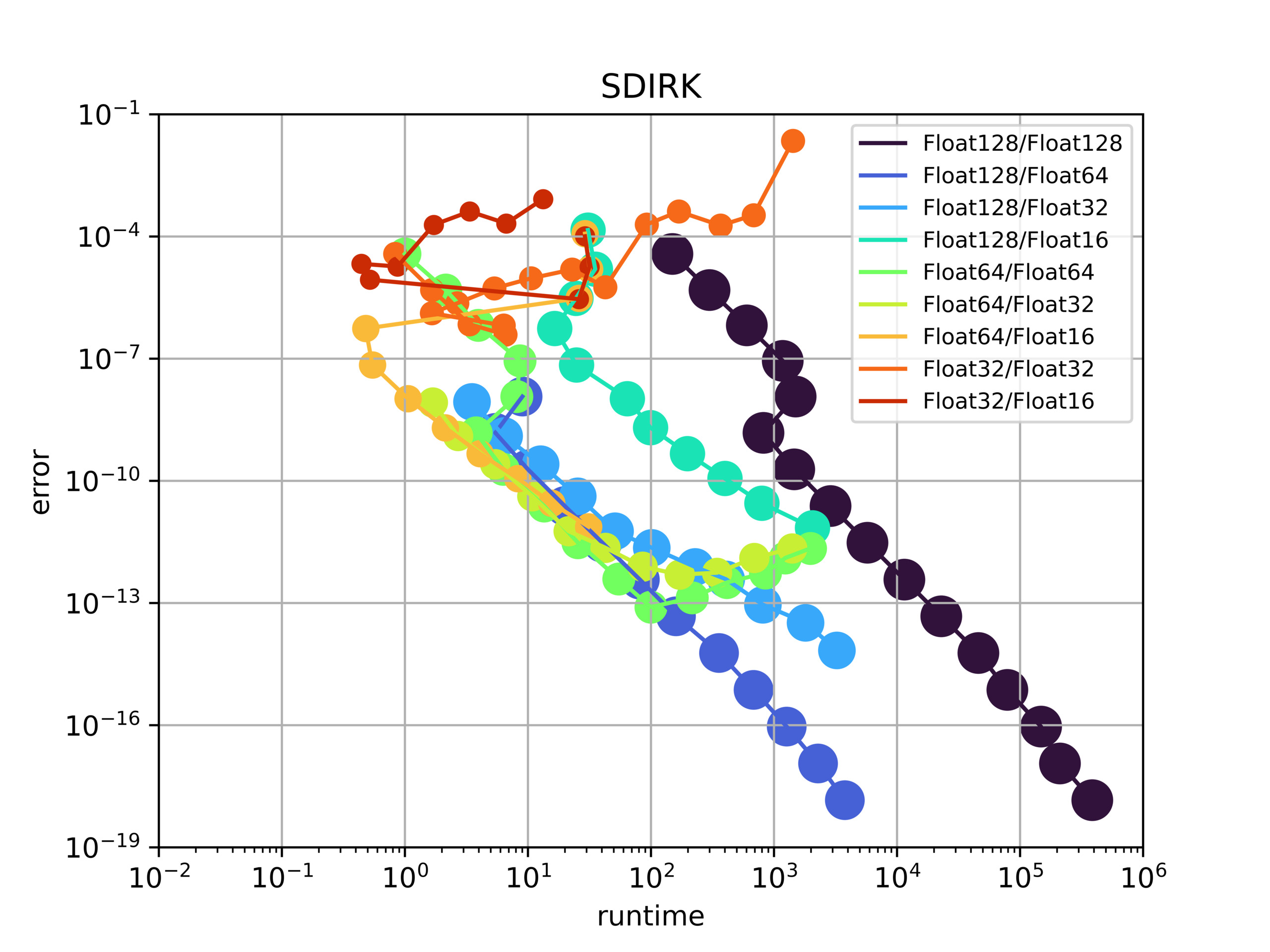}}\\
{\includegraphics[width=0.32\textwidth]{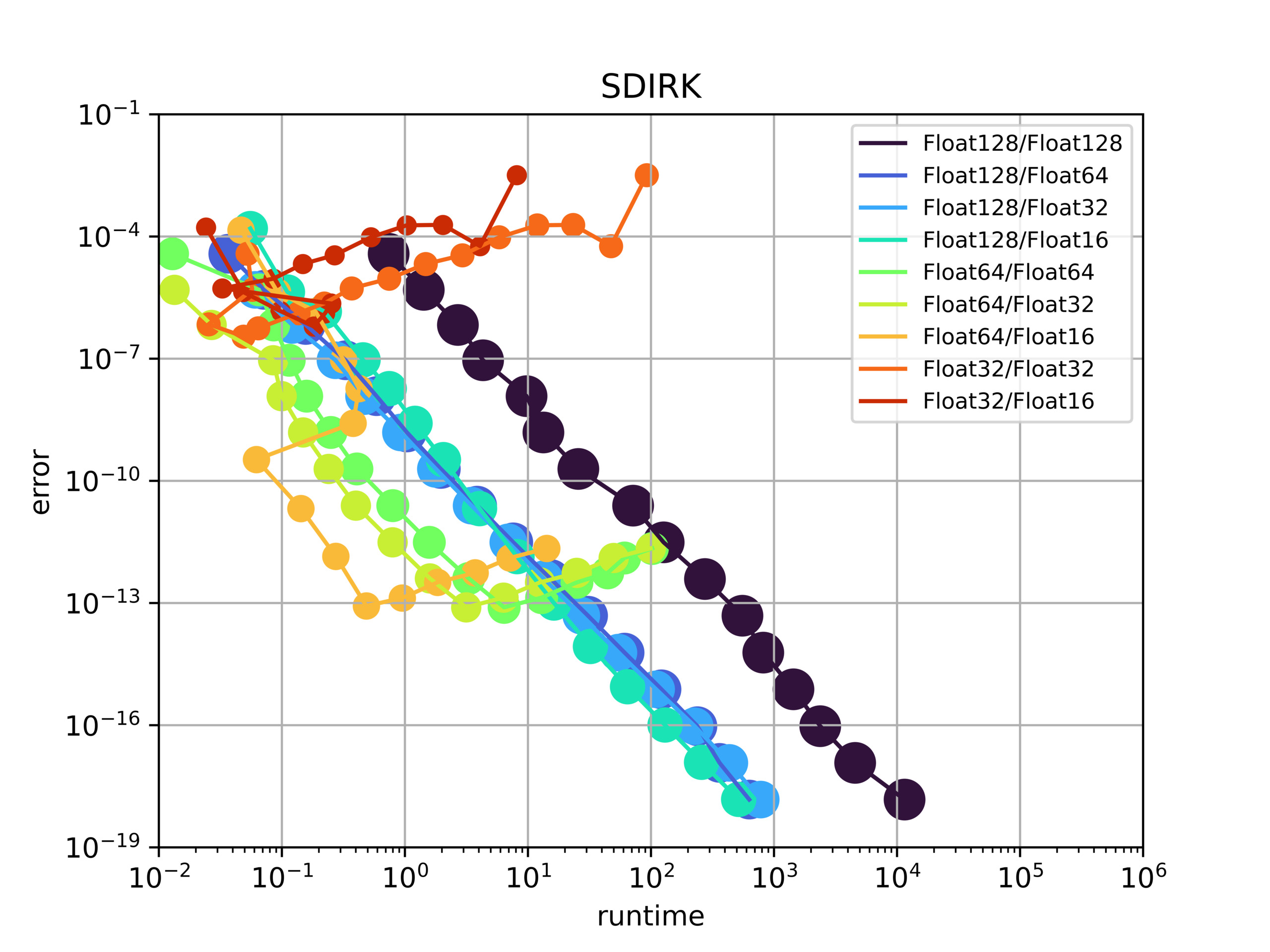}}
{\includegraphics[width=0.32\textwidth]{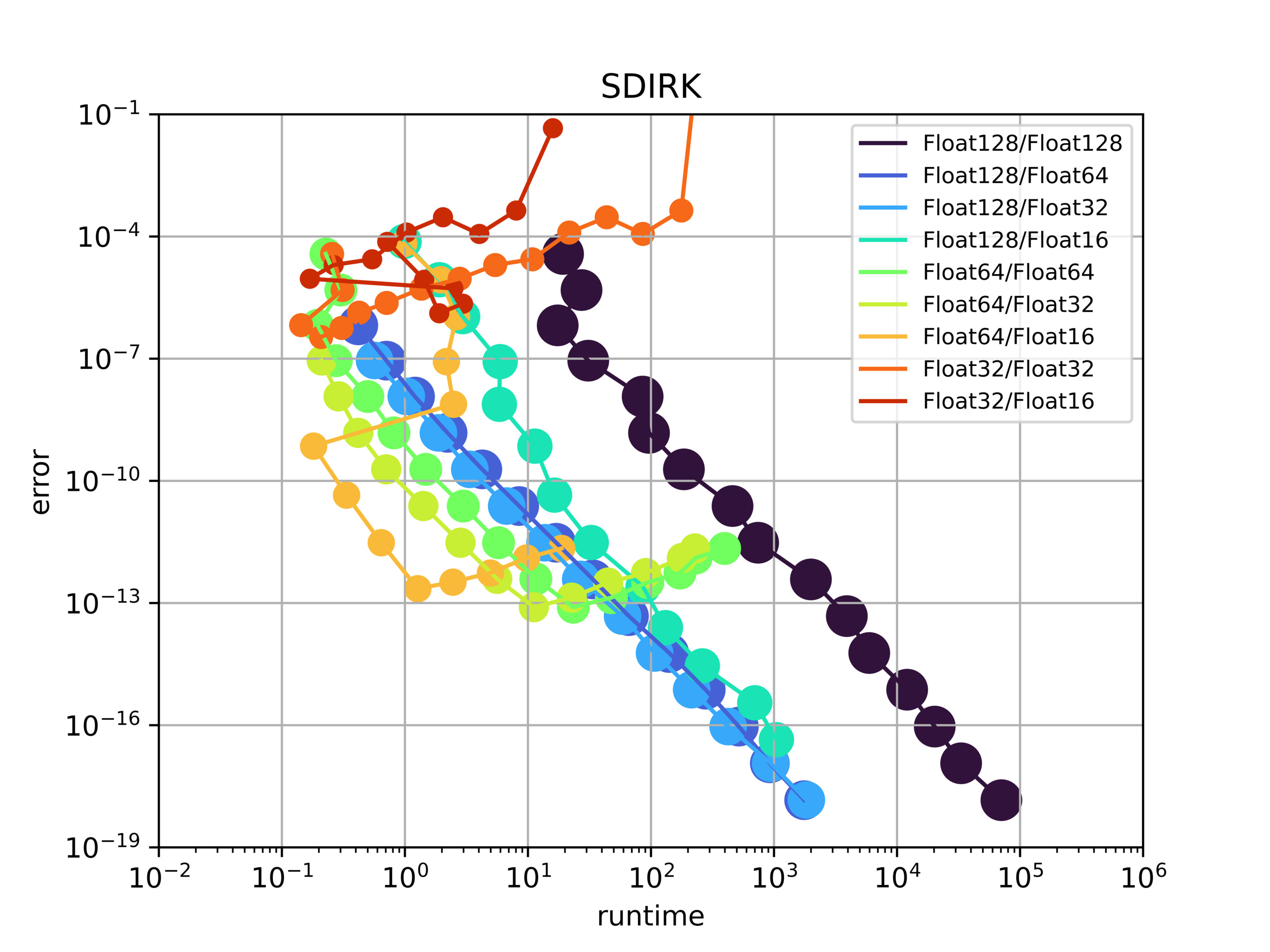}}
{\includegraphics[width=0.32\textwidth]{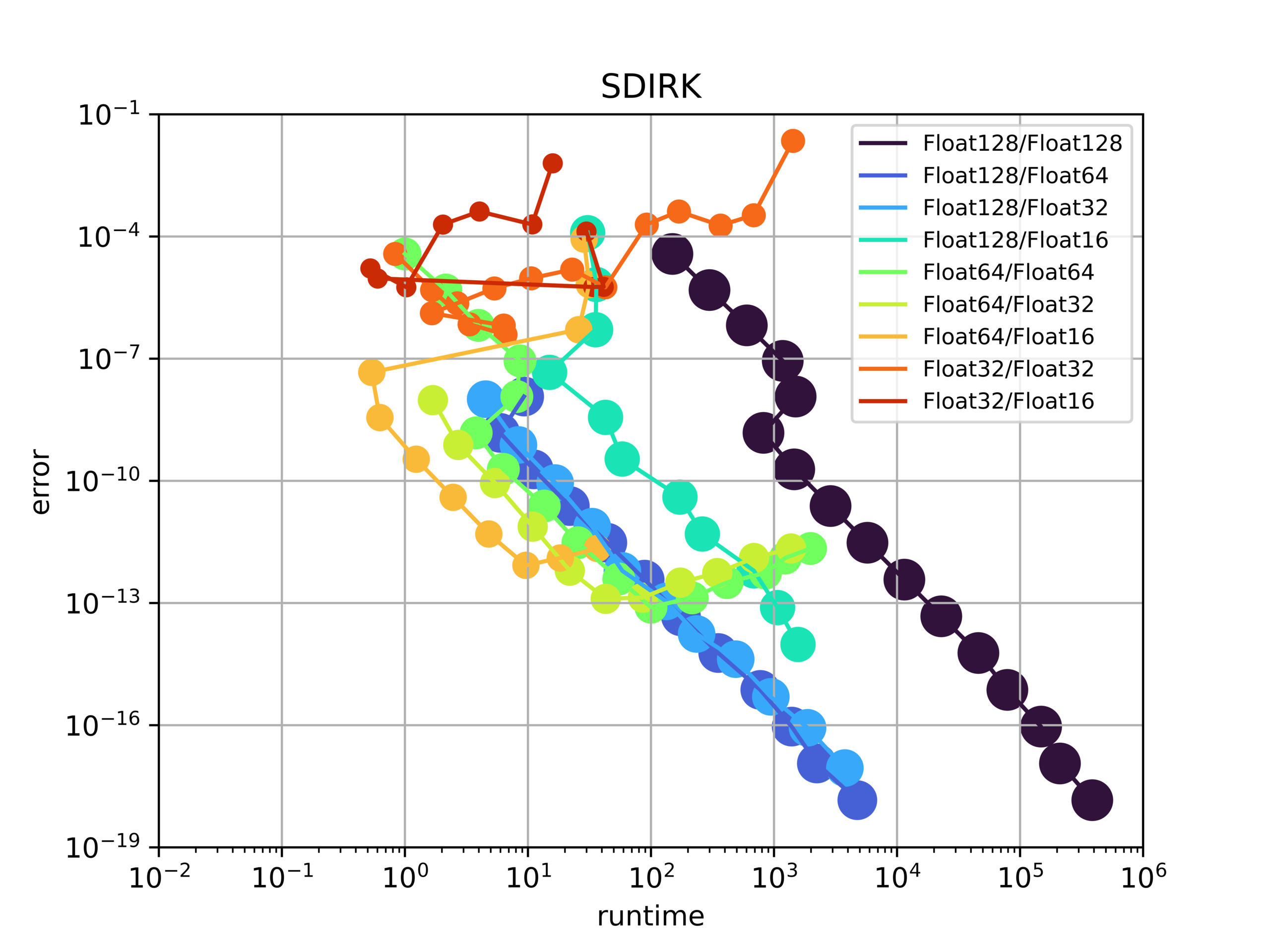}}\\
\caption{Julia code results:  The SDIRK method for the viscous Burgers' equation with $N_x=50$ (left), $N_x=100$ (middle), and $N_x=200$ (right), with no corrections (top row), one correction (second row), two corrections (third row), 
and three corrections (bottom row).
\label{BurgersSDIRK}}
\end{center}
\end{figure}

Figure \ref{BurgersSDIRK} shows the efficiency of the mixed precision 
SDIRK method \eqref{SDIRK-MPfix}  for the viscous Burgers' equation with $N_x=50$ (left), $N_x=100$ (middle), and $N_x=200$ (right), with no corrections (top row), one correction (middle row), two corrections (bottom row).
Notice that the single and half precision (and the mixed precision single/half) codes do not converge.
It is remarkable, then, how well the quad/single, quad/half, double/single, and double/half methods do comparably.
These all converge and with corrections even provide significant efficiency. Once again, reading top to bottom
we see the remarkable effect of the explicit corrections on the efficiency of the methods.

Figure \ref{BurgersSDIRK} shows that the double/half  method (Float64/Float16, light orange)
is the most efficient method for error levels between $10^{-10}$ and $10^{-13}$, beating out 
the double/half by a factor of close to 4, and providing 6.7x speedup over the double precision code (see table below).
For error levels above $10^{-10}$, the double/single  codes are the best performers.
For $N_x = 100$ the results are similar, with the double/half  method with three corrections providing a speedup
factor of  4x over the  double/single method, and a speedup of  over 8x over the double precision.
Finally, for  $N_x = 200$ we see a speedup of 4.4x over the mixed double/single and a speedup of over 5x  
over the double precision code.\\
\begin{center} {\small
\begin{tabular}{ |c|c|c|c| }
 \hline
 \multicolumn{4}{|c|}{SDIRK with 2 corrections for error level $\approx 10^{-9}$ }\\ \hline
  spatial points & 64/64 & 64/32 & 64/16 \\
$N_x = 50$ &   0.25 & 0.17 & 0.086 \\
$N_x = 100$ & 0.81 & 0.434 & 0.29 \\
$N_x = 200$ & 3.79 & 2.71  & 1.06 \\ \hline
\multicolumn{4}{|c|}{SDIRK with 3 corrections for error level $\approx 10^{-10}$ }\\ \hline
$N_x = 50$ &  0.41 & 0.240 & 0.062 \\
$N_x = 100$ & 1.47 & 0.74 & 0.18 \\
$N_x = 200$ &  6.32 &  5.4 & 1.23   \\ \hline
  \end{tabular}}
  \end{center}

\smallskip

Looking at error levels below $10^{-13}$ three corrections result in a quad/double and quad/single 
code that have the same efficiency and are {\em significantly} more efficient than the quad precision code. 
For $N_x =50$, the quad/half code is more efficient than the quad/double and quad/single,
but for $N_x =100$ it is a little less efficient and does not decay down to the final level, 
and for $N_x =200$ much less so. 
The table below shows that for two corrections,  the quad/single and quad/half codes do not get to
an error level of $\approx 10^{-18}$ with two corrections for any of $N_x=50, 100, 200$
but for three corrections all the quad/single codes reach the desired error level, as does the
quad/half code for $N_x=50$.\\
\begin{center} {\small
 \begin{tabular}{ |c|c|c|c|c| }
 \hline
 \multicolumn{5}{|c|}{SDIRK with 2 corrections for error level $\approx 10^{-18}$ }\\ \hline
 spatial points & 128/128 & 128/64 & 128/32 & 128/16 \\
$N_x = 50$ &  11,524 & 686 & -- & --   \\
$N_x = 100$ & 70,449 & 1,517 & -- & --   \\
$N_x = 200$ & 387,140 &3,768 &  -- & --  \\ \hline
\multicolumn{5}{|c|}{SDIRK with 3 corrections for error level $\approx 10^{-18}$ }\\ \hline
$N_x = 50$ &    11,524 & 630 & 778 & 515  \\
$N_x = 100$ &  70,449 & 1,760 & 1,831 & --\\
$N_x = 200$ &  387,140 & 4,760 & 3,765 & --   \\ \hline
  \end{tabular}}
  \end{center}
  
  \smallskip
This table indicates that with three corrections the quad/half results provide a speedup of
over 22x over the quad precision code for  $N_x=50$. For  $N_x=100$ the quad/double 
provides the most efficient result, with a speedup of  40x, and for  $N_x=200$
the quad/single is most efficient with a speedup of two orders of magnitude.

\subsubsection{NovelA method}
\begin{figure}[htb]
\begin{center}
{\includegraphics[width=0.32\textwidth]{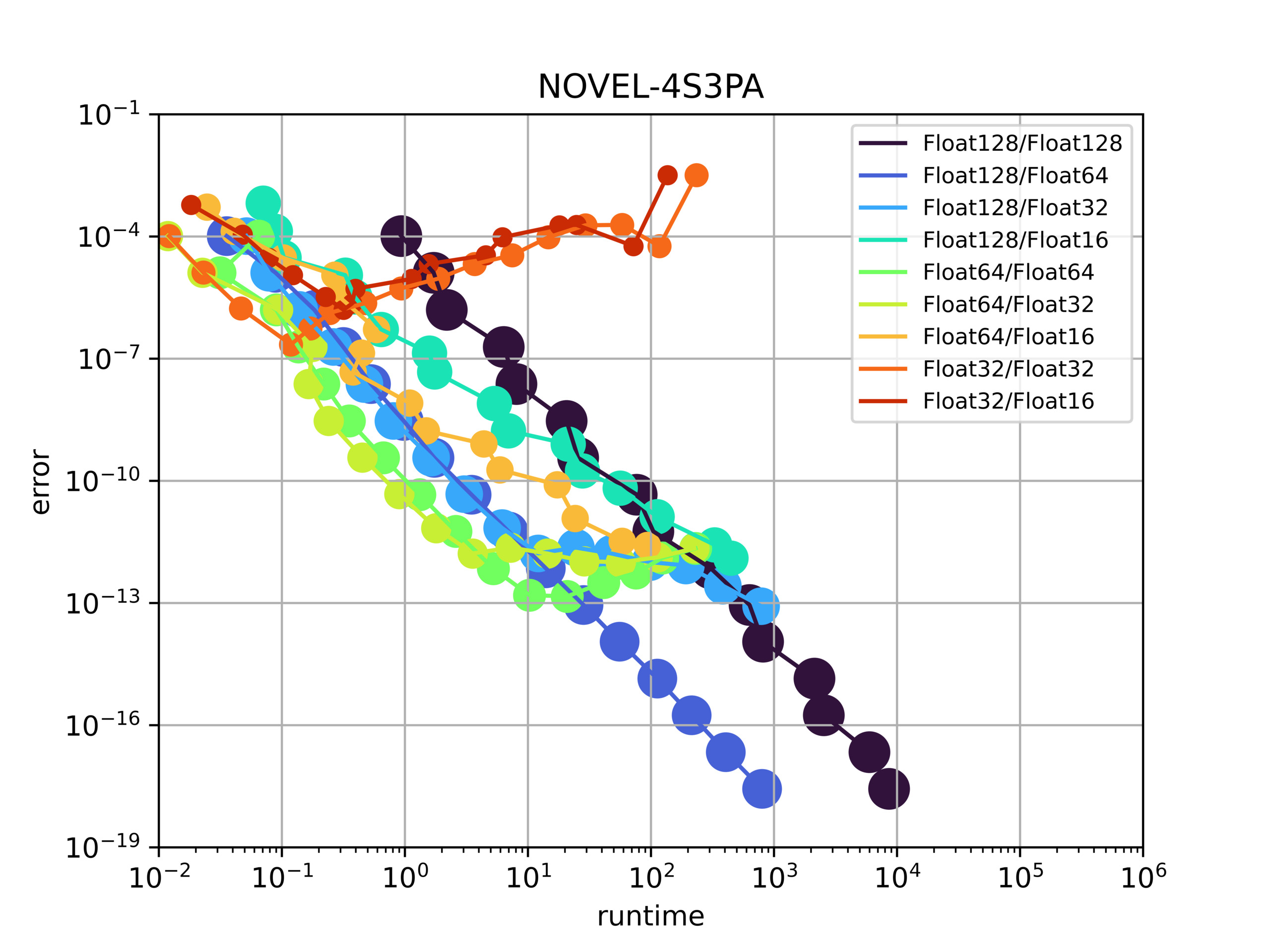}}
{\includegraphics[width=0.32\textwidth]{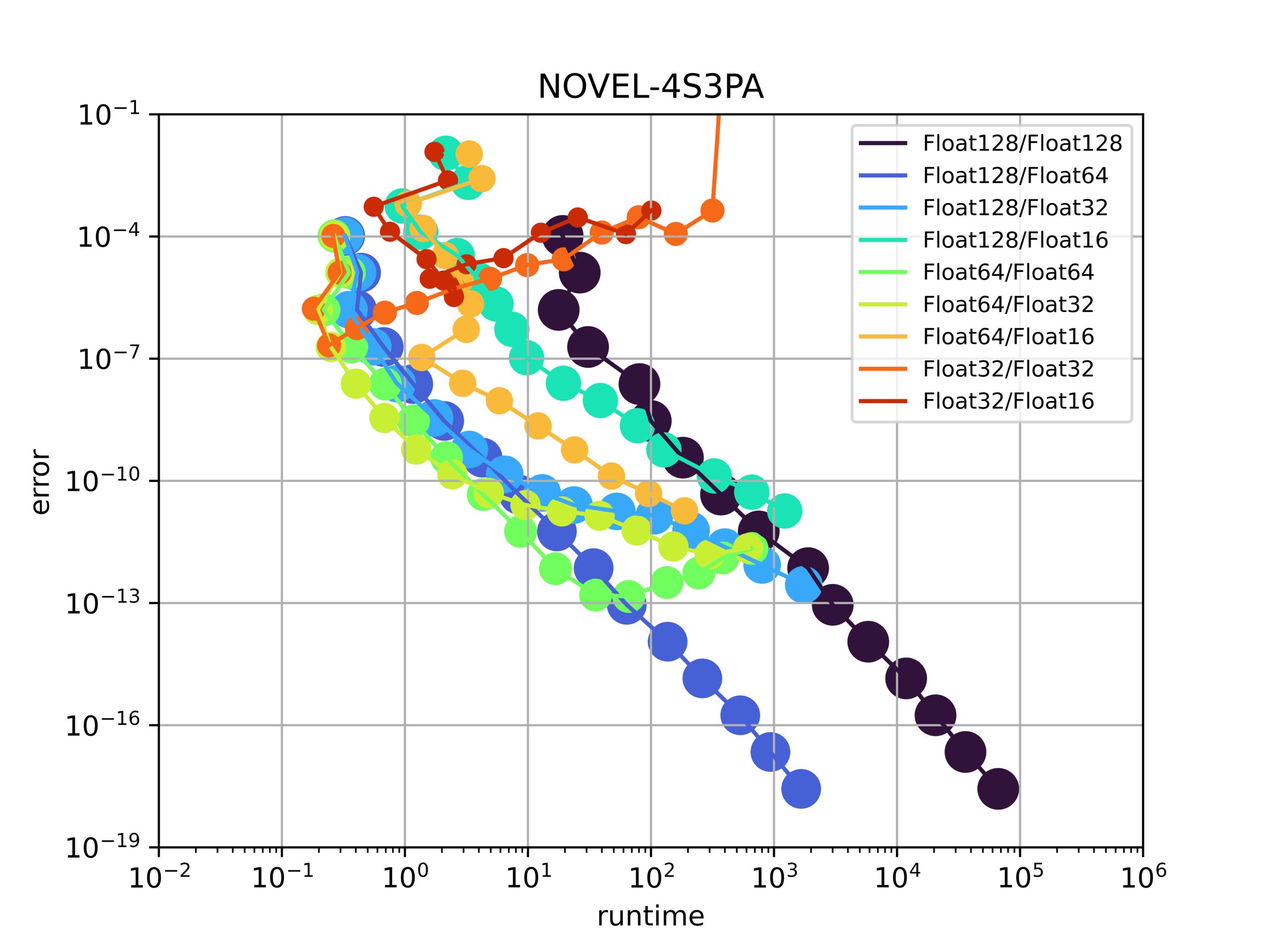}}
{\includegraphics[width=0.32\textwidth]{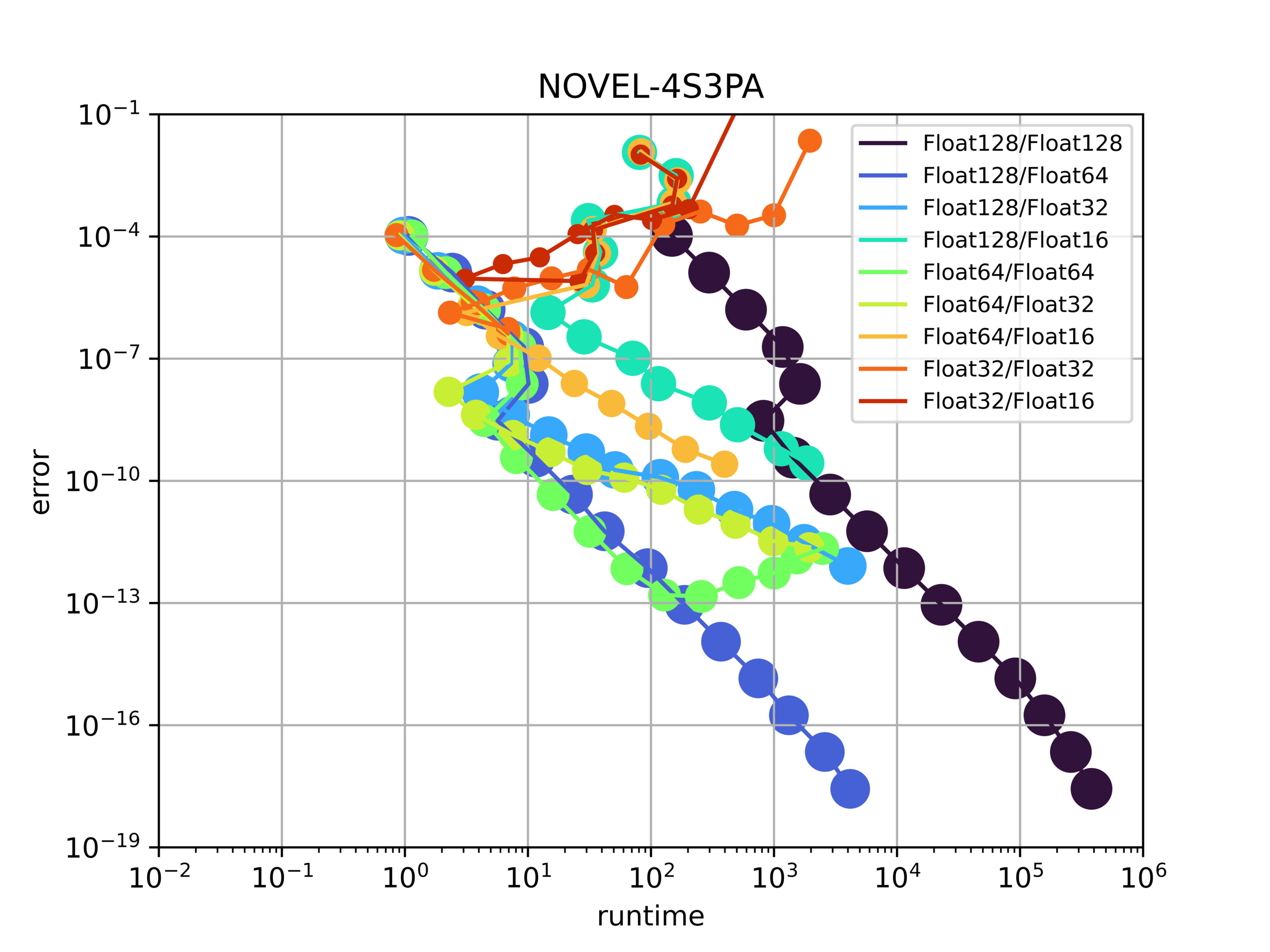}}\\
\caption{Julia code results: The NovelA for the viscous Burgers' equation with  $N_x=50$ (left), $N_x=100$ (middle), 
and $N_x=200$ (right).
\label{BurgersNovel}
}
\end{center}
\end{figure}

Figure \ref{BurgersNovel} shows the efficiency of the mixed precision  NovelA method \eqref{MP-4s3pA}  
for the viscous Burgers' equation with $N_x=50$ (left), $N_x=100$ (middle), and $N_x=200$ (right).
As expected, these methods behave much like SDIRK with one correction. 
A minor efficiency gain is seen using the double/single code for error levels $10^{-7}$ to $10^{-10}$
for $N_x=50$ and  $N_x=100$, but not for $N_x=200$.
Major efficiency gains are seen for the  quad/double code, where  at error levels below $10^{-13}$
we see a speedup factor of 10x for  $N_x=50$, 40x for $N_x=100$, and  90x for $N_x=200$ (see table below). 
We note that it is also possible to add explicit corrections to the NovelA method when needed,
and we expect this to improve the mixed precision errors and efficiency for small $\dt$.

\begin{center} {\small
\begin{tabular}{ |c|c|c|c| }
 \hline
 \multicolumn{4}{|c|}{NovelA for error level $\approx 10^{-18}$ }\\ \hline
 spatial points &   128/128 & 128/64 & savings factor \\ \hline
$N_x = 50$ &  8,622 & 800 & 10.77 \\
$N_x = 100$ & 66,404 & 1662 & 39.95 \\
$N_x = 200$ &  381,190 & 4155 &91.74 \\ \hline
  \end{tabular}}
\end{center}

\section{Conclusions}

In this paper we  presented a stability analysis and numerical study of the 
convergence and efficiency of additive Runge-Kutta Methods for mixed precision
simulations. Following \cite{Grant2022}, the mixed precision implementation
of the Runge-Kutta methods is recast as a pertubation of the function and expressed 
as an additive Runge-Kutta method. This approach was used in  \cite{Grant2022} to derive the 
perturbation order conditions and derive novel methods as well as verify that adding explicit 
corrections to known methods would enhance the perturbation order of accuracy. 
This paper provides confirmation that the approach developed in \cite{Grant2022}
can provide stable, accurate, and efficient simulations in many cases.

The results presented confirm that the  mixed precision Runge-Kutta methods display the 
expected perturbation order, and that the explicit corrections do, as designed, improve the
perturbation order. Recalling that order is an asymptotic property, we verified that these 
explicit corrections are needed and are beneficial for sufficiently small time-step. However,
for large enough time-steps the corrections may in fact harm the performance of the method.
Future work will investigate a systematic approach to including explicit corrections only when 
beneficial, and the impact of these results.

The numerical efficiency studies showed that the mixed precision Runge-Kutta approach
may provide significant to dramatic savings in runtime, depending on the problem.
Our results demonstrate that mixing a single or even half precision function evaluation with the 
double precision corrections allows an accurate and convergent simulation at significant runtime savings, 
whereas the corresponding single or  half precision  codes are unstable.
Furthermore, where quad precision is desirable (for example, for long time simulations) the mixing 
of quad with double precision is a reliable way to significantly improve the runtime of the method  --
by up to two orders of magnitude! -- without adversely impacting the accuracy of the method.

\section*{Author Contribution Statement}

{\bf Ben Burnett} was responsible for all the codes, computations, and preparation of graphs in Sections \ref{sec:convergence} and 
\ref{sec:performance}, was involved in discussions on the presentation of the material, and on the stability and sensitivity 
analysis. He provided the expertise on half precision and quad precision implementation, and on 
running all the codes on the CARNIE cluster. He was the primary writer of Section 5.1, and 
read and commented on the entire manuscript. 

{\bf Sigal Gottlieb} oversaw the preparation of the codes, computations, and preparation of graphs in 
Sections \ref{sec:convergence} and  \ref{sec:performance}, and helped design the numerical experiments and 
debug codes. She was primarily responsible for the sensitivity to roundoff analysis, and for the writing of the 
first draft of the majority of  the manuscript.

{\bf Zachary Grant}  was primarily responsible for the writing of the review section, and for the stability analysis 
and the codes used to produce the stability regions. He was involved in discussions on the sensitivity to roundoff 
analysis, and accuracy studies and their interpretation.
He also helped design the numerical experiments presented in \ref{sec:convergence} and  \ref{sec:performance}
and analyze the results. He collaborated on presentation of the material and read and edited the entire manuscript.

\begin{center} {\em
On behalf of all authors, the corresponding author states that there is no conflict of interest. 
All authors certify that they have no affiliations with or involvement in any organization or entity with any financial interest or non-financial interest in the subject matter or materials discussed in this manuscript. The authors have no financial or proprietary interests in any material discussed in this article.}
\end{center}

\clearpage

\end{document}